\pdfoutput=1
\documentclass[pdftex, noinfoline]{imsart}
 
\usepackage[english]{babel} 
\usepackage{parskip}
\usepackage{amsmath,amsfonts,amsthm, amssymb}
\usepackage{graphicx} 
\usepackage{enumitem} 
\usepackage{xcolor, latexsym} 
\usepackage{subfigure} 

\usepackage{titletoc}   

\usepackage{graphicx}
\usepackage{placeins}
\usepackage{float}

\usepackage[pdftex, 
colorlinks,citecolor=blue,urlcolor=blue,filecolor=blue,linkcolor=blue,backref=page,
bookmarks=true, breaklinks, hypertexnames=false, pdfborder={0 0 0 0}]{hyperref} 
\usepackage[normalem]{ulem} 
\usepackage{float}
\usepackage{placeins}
\usepackage{soul} %for highlighting
\soulregister\cite7 %for allowing cites within highlight!
\soulregister\ref7
\soulregister\pageref7
\DeclareMathOperator*{\KL}{\text{KL}}

\newcommand{\tz}{z} 

\usepackage{xcolor}
\definecolor{forestgreen}{RGB}{34, 139, 34}

\definecolor{blendedblue}{rgb}{0.2,0.2,0.7}

\newcommand{\gambln}{\gamma_{b,L_n}}
\newcommand{\gammln}{\gamma_{M_1, L_n}}
\newcommand{\gamb}{\gamma_b}
\newcommand{\gamm}{\gamma_{M_1}} 

\newcommand{\vpa}{\varphi_a}
\newcommand{\vpb}{\varphi_b}

\newcommand{\rn}{\sqrt{n}}

\newcommand{\tIinv}{\tilde I_{\eta_0}^{-1}}

\newcommand{\Ld}{\Lambda}
\newcommand{\ld}{\lambda}
\newcommand{\Dl}{\Delta}

\newcommand{\given}{\,|\,}

\newcommand{\Ltwo}{{2}}
\newcommand{\Lone}{{1}}
\newcommand{\Linfty}{{\infty}}

\allowdisplaybreaks

\theoremstyle{plain} 
\newtheorem{theorem}{Theorem}
\newtheorem{corollary}{Corollary}
\newtheorem{lemma}{Lemma}
\newtheorem{prop}{Proposition}
\newtheorem{remark}{Remark}

\allowdisplaybreaks

\endlocaldefs
\usepackage{natbib}
\usepackage{multirow}
\usepackage{chngpage}
\usepackage{siunitx}
\usepackage{booktabs}
\usepackage{makecell}
\usepackage[T1]{fontenc}
\usepackage{float}
\usepackage{grffile}
\usepackage{url} 
\usepackage{graphicx,psfrag,epsf}

\usepackage{bbm, dsfont}
\usepackage{mathrsfs} 
\usepackage[ruled,vlined]{algorithm2e}

\usepackage{xr-hyper}

\defcitealias{cast20survival}{CvdP21}

\begin{document}

\begin{frontmatter}

\title{Bayesian multiscale analysis of the Cox model}

\runtitle{Bayesian Cox Model}

\author{\fnms{Bo Y.-C.} \snm{Ning}\thanksref{t1}\ead[label=e2]{bycning@ucdavis.edu}}
\and
\author{\fnms{Isma\"el} \snm{Castillo}\thanksref{t2}\ead[label=e1]{ismael.castillo@upmc.fr}}

\thankstext{t1}{The author gratefully acknowledges the support from the Fondation Sciences Math\'ematiques de Paris postdoctoral fellowship.}
\thankstext{t2}{The author gratefully acknowledges support from the Institut Universitaire de France (IUF) and from the ANR grant ANR-17-CE40-0001 (BASICS).}
 
\affiliation{Sorbonne University \& UC Davis\thanksmark{m1} and Sorbonne University \& IUF\thanksmark{m2}}

\address{University of California, Davis \\
Department of Statistics,\\
1227 Mathematical Science Building,\\
One Shields Avenue, 
Davis, CA 95616 United States\\
\printead{e2}}

\address{Sorbonne Universit\'e \&
Institut Universitaire de France\\
  Laboratoire de Probabilit\'es, 
Statistique et Mod\'elisation\\ 4, Place Jussieu, 75252, Paris cedex 05, France\\
  \printead{e1}}
\runauthor{Ning \& Castillo}

\begin{abstract}
Piecewise constant priors are routinely used in the Bayesian Cox proportional hazards model for survival analysis.
Despite its popularity, large sample properties of this Bayesian method are not yet well understood.
This work provides a unified theory for posterior distributions in this setting, not requiring the priors to be conjugate. We first derive contraction rate results for wide classes of histogram priors on the unknown hazard function and prove asymptotic normality of linear functionals of the posterior hazard in the form of Bernstein--von Mises theorems. Second, using recently developed multiscale techniques, we derive functional limiting results for the cumulative hazard and survival function. Frequentist coverage properties of Bayesian credible sets are investigated: we prove that certain easily computable credible bands for the survival function are optimal frequentist confidence bands. We conduct simulation studies that confirm these predictions, with an excellent behavior particularly in finite samples.
%\sout{showing that even simplest possible Bayesian credible bands for the survival function can outperform state-of-the-art frequentist bands in terms of coverage.} 
{Our results suggest that the Bayesian approach can provide an easy solution to obtain both the coefficients estimate and the credible bands for survival function in practice.}
\end{abstract}



\end{frontmatter}

%\sbl{I feel the title `Bayesian multiscale analysis of the Cox proportional hazards model' would be better. To be discussed} \org{agreed!}

%\org{Multiscale analysis is a bit abstract for applied statisticians. How about change the title to ``\sbl{On the Bernstein--von Mises phenomenon and supremum-norm contraction rate for Bayesian Cox piecewise constant hazards model}''? This is a bit long, but it may attract some applied people to read our paper. }

% INTRODUCTION
\section{Introduction} \label{sec:intro}

The Cox proportional hazards model (hereafter, the Cox model) introduced by \citet{cox72}
is one of the most popular regression models for survival analysis. 
It is a {\em semiparametric} model with two sets of unknown parameters: the regression coefficients, which measure the correlation between the covariates and the
explanatory variables, and the baseline hazard, a nonparametric quantity, which describes the risk of
events happening within given time intervals at baseline levels conditional on the covariates.
A commonly used approach to estimate the two parameters takes two steps: first estimate the regression coefficients from the Cox partial likelihood \citep{cox72} and then derive   
the estimated cumulative hazard function (known as the Breslow estimator \citep{breslow72}) through maximizing the full likelihood via plugging-in the estimated value of the regression coefficients. 

In the past few decades, Bayesian methods for the Cox model have been widely applied for analyzing datasets in, e.g., astronomy \citep{isobe86}, medical and genetics studies \citep{jialiang13}, and engineering \citep{equeter20}.
An advantage of the Bayesian approach is that uncertainty quantification for the parameters of interest is in principle straightforward to obtain once posterior samples are available.  
Contrary to the standard two-step procedure mentioned above, the Bayesian approach provides estimates for the joint distribution of all parameters, which enables to capture dependencies: in particular, as one of the practical applications considered below, one can derive meaningful uncertainty quantification simultaneously for the Cox model parameter and functionals of the hazard rate (e.g. its mean, or the value of the survival function at a point) from corresponding credible sets, in particular automatically capturing the (optimal  `efficient') dependence structure.

The prior for the hazard function needs to be chosen carefully, as it is a nonparametric quantity. 
Two main common approaches to place a prior on hazards have been considered in the literature.
A first approach puts a prior on the cumulative hazard function, modeling this quantity rather than the hazard itself. 
A prominent example is the family of {\it neutral to the right process} priors, which includes the Beta process prior 
\citep{hjort90, damien96},
the Gamma process prior \citep{kalbfleisch78, burridge81}, 
and the Dirichlet process prior \citep{florens99} as special cases.
A second approach, which is the one we follow here, is to put a prior on the baseline hazard function. A commonly used family  is that of piecewise constant priors \citep{ibrahim01}. 
The latter approach is particularly attractive, first because it allows for inference on the hazard rate (assuming it exists), and second because in practice the follow-up period is often split into several intervals, with  
the hazard rate taking a distinct constant value on each sub-interval, making the output of the method easy to interpret for practitioners.

One primary goal of the paper is to validate the practical use of Bayesian credible sets for inference on the Cox model{'s} unknown parameters, for instance credible bands for {the survival function conditional on the covariates.} Indeed, practitioners often treat Bayesian credible sets as confidence sets. Before discussing the possible mathematical validity of this practice, let us conduct a simple illustrative simulation study (see Section \ref{sec:sim} for a detailed description of the simulation setting). 

\begin{figure}
\centering
\includegraphics[width=0.5\textwidth]{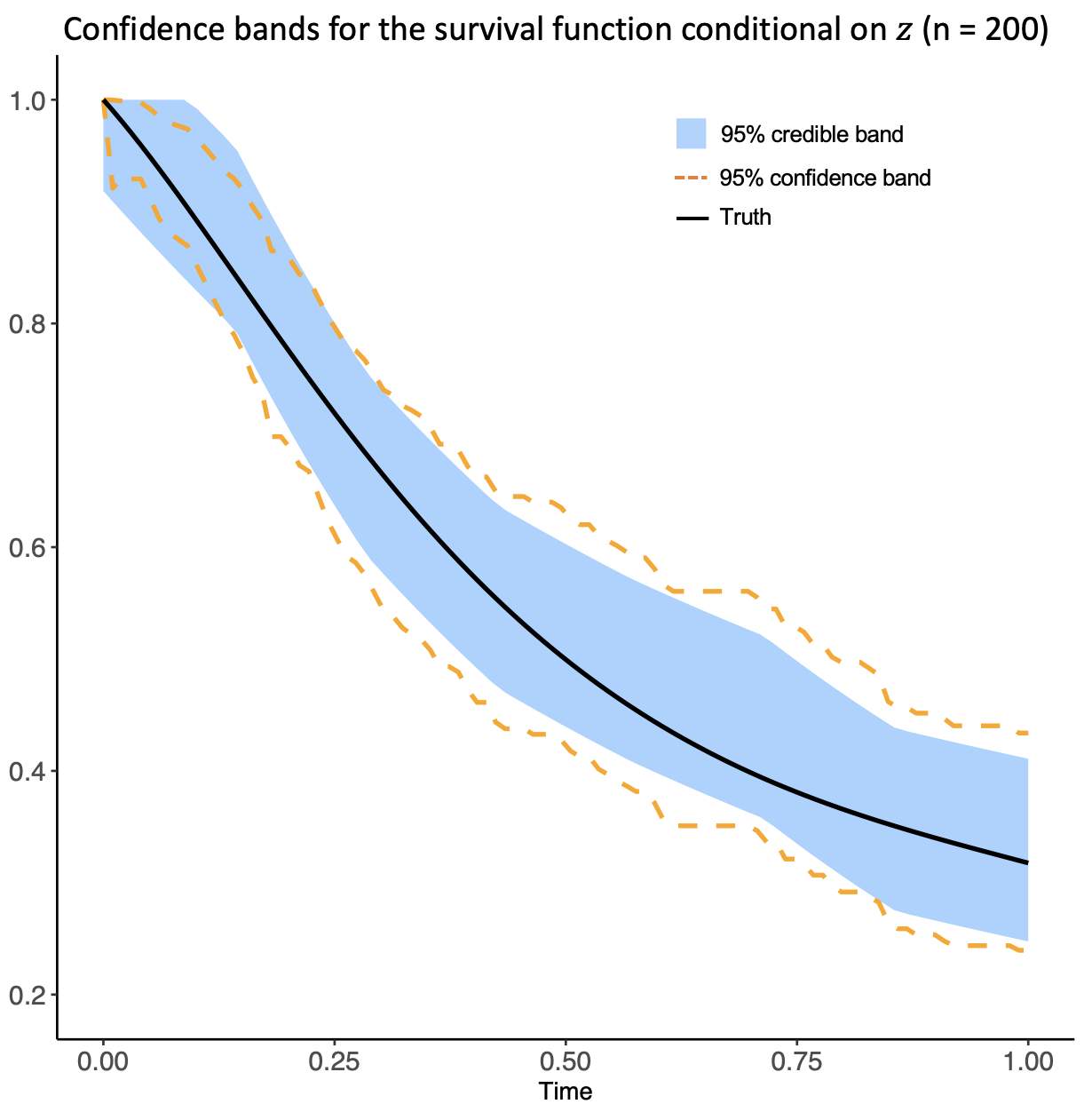}
\vspace{-.2cm}
\caption{
The 95\% Bayesian credible band (blue) using the random histogram prior and the 95\% frequentist confidence band (orange).
The bolded black line is the true {survival function conditional on $z$}, $\exp\{\Ld_0 e^{\theta_0'z}\}$. 
Sample size is $n = 200$. 
}
\label{fig:intro}
\end{figure}

From data simulated from the Cox model, suppose we want to make inference on the {survival function conditional on $z$}, which gives the probability that a patient survives past a certain time $t$ given a covariate $z \in \mathbb{R}^p$, a useful quantity for practitioners. Let us compare a simple $95\%$ credible band of the posterior distribution (with the piecewise constant prior; see Section \ref{sec:sim} for its construction) induced on the {survival function conditional on $z$} and a certain 95\% confidence band---which requires estimation of the covariance structure---of the same function obtained by a commonly used frequentist approach (see Section \ref{sec:sim} for a precise description of how the band is obtained).  
In Figure \ref{fig:intro}, we plot the credible band (blue) and the confidence band (orange).    
The sample size is $n = 200$, and we let $p=1$ and $z = 1$.
One first notes that the true function (black) is contained in both bands, which suggests that both provide a reasonable uncertainty assessment for the {survival function conditional on $z$}. Second, we compared the total area of the two bands; interestingly, we found that the area of the credible band is smaller than that of the second band: the area of the credible band is 0.163 and that of the confidence band is 0.183 ($12\%$ larger). 
A thorough Monte Carlo study, carried out in Section \ref{sec:sim}, confirms that the area of the Bayesian credible band is indeed consistently smaller on average than the size of the frequentist confidence band when the sample size is 200. 
It may be noted that these results hold for the in a sense `simplest possible' Bayesian credible band: as can be seen in Figure \ref{fig:intro}, its width is fixed through the time interval, and even better results are expected for bands that become thinner close to $t=0$; see Section \ref{sec:diss} for more discussion on this.   
This simulation study suggests that, aside from not having to {estimate} covariances, 
using the Bayesian credible set can be particularly advantageous, especially for 
small sample size datasets.

The observations from Figure \ref{fig:intro} raise some interesting questions: can one validate and generalize our findings in the figure, that is, can one provide theory explaining why the credible band is a confidence band, and will the two bands become more similar as sample size increases? 
What can be said in terms of the hazard function: does the Bayesian procedure estimate it in a possibly `optimal' way? We now discuss the existing literature on these questions and the main contributions of the paper. 

In smooth parametric models, taking certain quantile credible sets as confidence sets is justified mathematically by the celebrated Bernstein--von Mises theorem (henceforth BvM, see e.g. \cite{vdvAS}, Chapter 10): a direct consequence thereof is that taking, in dimension $1$ say, the $\alpha/2$ and $1-\alpha/2$ posterior quantiles provides a credible set (by definition) whose frequentist coverage asymptotically goes to $1-\alpha$. Its diameter also asymptotically matches the information bound so is optimal in the frequentist sense from the efficiency perspective. For more complex models, such as the
Cox model, obtaining a semiparametric BvM theorem for the regression parameter is possible, as we see below, but requires non-trivial work. Obtaining analogous results at the level of the {\em survival} function itself is even more challenging. We now review recent advances in the area for such semi- and non-parametric models.
 
Semiparametric BvM theorems where obtained in \citet{cast12} under general conditions on the statistical model using Gaussian process priors on the nuisance parameter. 
\citet{cast15BvM} considered an even more  general framework, allowing for BvMs for linear and non-linear functionals (also generalizing some early results of \cite{RR12} for density estimation).  A multiscale approach was introduced in \citet{cast13, cast14a} in order to derive nonparametric BvM theorems for families of possibly non-conjugate priors, as well as Donsker--type theorems. Yet, the first applicative examples of these works were mostly confined to relatively simple models and/or priors.

%\textcolor{gray}{Theory for convergence of Bayesian posterior distributions in survival models has mostly followed two directions, which we briefly review now (see also \citet{ghosal17}, Section 12.3.3 and Chapter 13). A first series of influential results has been concerned with classes of neutral to the right priors: in the standard nonparametric survival analysis model, \citet{hjort90} showed in particular conjugacy for Beta process priors in the context of cumulative hazard estimation; \citet{kimlee01} obtained posterior consistency, while \citet{Kim2004} derived Donsker theorems for the posterior survival function. In \citet{kim06},  the Cox model was considered and the joint posterior distribution of parameter and survival function was shown to satisfy the Bernstein--von Mises theorem. These results share two common features: they model the cumulative baseline hazard (equivalently, the survival function), not the baseline hazard itself---which can be desirable to model for practitioners---, and they rely on conjugacy of the class of neutral to the right priors, which provides fairly explicit characterisations of the posterior distributions.}
   
{  
Theory for convergence of Bayesian posterior distributions in survival models has mostly followed two directions, which we briefly review now (see also \citet{ghosal17}, Section 12.3.3 and Chapter 13). A first series of influential results has been concerned with classes of neutral to the right priors; e.g., 
\citet{hjort90} and \citet{kimlee01, Kim2004} in the study of the standard nonparametric survival analysis model.  \citet{kim06} studied the Cox model, in which the joint posterior distribution of parameter and survival function was shown to satisfy the Bernstein--von Mises theorem. These results share two common features: they model the cumulative baseline hazard (equivalently, the survival function), not the baseline hazard itself---which can be desirable to model for practitioners---, and they rely on conjugacy of the class of neutral to the right priors, which provides fairly explicit characterisations of the posterior distributions.}

A second series of results, closer in spirit to ours, considers priors on the baseline hazard function. \citet{DeBlasi2009a} used a kernel mixture  with respect to a completely random measure as a prior, and obtained both posterior consistency for the hazard and limit results for linear and nonlinear functionals thereof. The work \cite{DeBlasi2009b} derived  a semiparametric BvM theorem in competing risk models. The present work can be seen as following the footsteps of \cite{cast20survival}, where the simple nonparametric model with right-censoring is treated, and for which  results for the hazard and cumulative hazard are derived. However, the latter model is much simpler than the Cox model, which features both regression coefficients and random covariates: this requires several more delicate bounding of both (semiparametric) bias terms and remainder terms, see Section \ref{sec:support-lemma-thm1} in \citet{ning22supp} for more details.    
 In \citet{cast12}, the Cox model is treated as an application of the general results, which yield the semiparametric BvM theorem for the Cox model parameter for (transformed) Gaussian process priors on the hazard. Although this result has a general flavor, and can be adapted to handle other prior families, it requires a fast enough posterior contraction rate for the baseline hazard, and thus cannot be applied for histogram priors on the hazard (see  Section \ref{sec:hellinger} for more on this). Perhaps more importantly, the later result is confined to the Cox regression parameter, and says nothing about uncertainty quantification for the cumulative hazard, for the survival function, or even simply for linear functionals of the hazard.

The present work obtains the first results for Bayesian uncertainty quantification jointly on regression parameters and survival function in the Cox model using non-conjugate priors. In particular, our results demonstrate that the popular and broadly used histogram priors on the baseline hazard provide not only contraction of the posterior distribution around the true unknown parameters, but  also optimal and efficient uncertainty quantification on those. We adopt the multiscale analysis approach, which is motivated by \citet{cast20survival}'s study of the survival model.
More precisely, we derive 
\begin{itemize}
\item[(a)] a joint Bernstein--von Mises (BvM) theorem for linear functionals of the regression coefficients and of the baseline hazard function;
\item[(b)] a Bayesian Donsker theorem for the conditional cumulative hazard and survival functions;
\item[(c)] a minimax optimal contraction rate for {the hazard function conditional on $z$} in supremum-norm distance.
\end{itemize}

{We would like to highlight that (a) as well as the upper-bound part of (c) are the most important novel contributions of this paper. Given these results are proved, results  (b) are obtained by adopting a similar philosophy as  in \citet{cast20survival}, but still require new arguments, in particular when deriving a joint `nonparametric' BvM result jointly in $(\theta,\lambda)$, see e.g. Proposition \ref{prop-tightness}.}
The Bayesian Donsker theorem, in particular, implies that certain $(1-\alpha)\%$ credible bands for the {survival function conditional on $z$} are asymptotically  $(1-\alpha)\%$ confidence bands.
In addition to these results, a nonparametric BvM theorem for the conditional baseline hazard function 
 is obtained in the Supplemental Material \citep[see][]{ning22supp}. All these results are completely new for non-conjugate (in particular, histogram) priors; also, we derive the first supremum-norm posterior contraction rates for the hazard in the Cox model; we also show the corresponding matching minimax lower bound, which to the best of our knowledge was not yet available in the literature. 
{We also demonstrate that in practice this easy-to-implementable computational algorithm can provide estimates for both the coefficients estimates and confidence bands for the cumulative hazard and survival functions. 
We note, as pointed out to us by a referee, that in current literature there seems to be a lack of practical algorithms in the Cox model model producing such type of bands. 
%This is in contrast with frequentist methods, as confidence bands can not be obtained easily. 
}

Finally, the techniques we introduce have a general flavor. First, they do not rely on conjugacy of the priors considered, so that they can virtually be applied to a wide variety of families, as long as a certain change-of-measure condition is met.  Second, the specific form of the model (here the Cox model) comes in through its local asymptotic normality (LAN) expansion, so similar techniques can be used in more complex settings, as long as a form of local asymptotic normality of the model holds. In particular, the techniques developed herein can serve as a useful tool for future studies of other semiparametric and nonparametric models, in survival analysis and beyond: 
 as further discussed in Section \ref{sec:diss}.

The paper is organized as follows. 
Section \ref{sec:setup} presents the model, the prior families, and key assumptions. 
Main results are presented in Section \ref{sec:main-results}.
Simulation studies are conducted in Section \ref{sec:sim}. 
Section \ref{sec:diss} concludes the paper and discuss a variety extensions for future studies. 
All the relevant proofs for the main results and auxiliary lemmas are left to the Supplementary Material. 

{\em Notation.} For any two real numbers $a$ and $b$, let $a \vee b = \max(a, b)$ and $a \wedge b = \min(a, b)$;
also, let $a \lesssim b$ as $a \leq C b$ for some constant $C$.  
{For a positive semidefinite matrix $A$, we denote $A \succeq 0$.}
Let us denote by $o(1)$ a deterministic sequence going to 
$0$ with $n$ and $o_P(1)$ a sequence of random variables going to $0$ in probability under the distribution $P$. 
For a vector $x$, denote $\|x\|_q$ as the $\ell_q$-norm of $x$ ($q \geq 1$), i.e., $\|x\|_q = (\sum_i |x_i|^q)^{1/q}$. When $q = \infty$, $\|x\|_\infty = \max_i |x_i|$ is the infinity norm of $x$. 
For a matrix $A$, denote $\|A\|_{(\infty, \infty)} = \max_{ij} |a_{ij}|$.
 
For a function $f \in L^{p}[a,b]$ ($p \geq 1$), where $L^p[a,b]$ is the space of functions whose $p$-th power is Lebesgue integrable on $[a,b]$, we define $\|f\|_{p} = (\int_a^b |f|^{p})^{1/p}$ the $L^{p}$-norm of $f$. 
If $p = 2$, $\|f\|_{\Ltwo} = \sqrt{\langle f, f \rangle}$ is the $L^2$-norm
and if $p = \infty$, $\|f \|_{\Linfty}= \sup_t |f(t)|$ is the supremum norm.
The associate inner product between any two functions, $f, g \in L^p[0, 1]$, is denoted as $\langle f, g\rangle = \int_0^1 fg$ and the space of continuous functions on $[a,b]$ is given by $\mathcal{C}[a, b]$ (resp. $L^{\infty}[a, b]$), which is equipped with the supremum norm $\|\cdot\|_\infty$.

For $\beta, D > 0$, let $l = \lfloor \beta \rfloor$ be the largest integer smaller than $\beta$, a standard H\"older-ball on $[a,b]$ can be defined as 
$$
\mathcal{H}(\beta, D) = \{f: |f^{(l)}(x) - f^{(l)}(y)| \leq D|x-y|^{\beta - l}, \  \|f\|_\infty \leq D, \ x, y\in [a, b]\}.
$$
Let $(\mathcal{S}, d)$ be a metric space and $\mu, \nu$ be probability measures of $\mathcal{S}$.
For $F: \mathcal{S} \to \mathbb{R}$, set
$$
\|F\|_{BL} = \sup_{x \in \mathcal{S}} |F(x)| + \sup_{x \neq y} \frac{|F(x) - F(y)|}{d(x,y)},
$$
and denote the bounded Lipschitz metric $\mathcal{B}_{\mathcal{S}}$ as
$$
\mathcal{B}_{\mathcal{S}}(\mu, \nu) = \sup_{F: \|F\|_{BL} \leq 1} 
\left|
\int_{\mathcal{S}} F(x) (d \mu - d \nu ) (x)
\right|.
$$
For two densities $f$ and $g$, denote $h^2(f, g) = \int (\sqrt{f} - \sqrt{g})^2 d\mu$ as their squared Hellinger distance.

\section{Model, prior families and structural assumptions}
\label{sec:setup}

{In this section, we first introduce the Cox model in Section \ref{sec:model}. Information-related quantities of the Cox model and structural assumptions for those quantities are given in Section \ref{sec:infonot} and \ref{sec:assumption} respectively. Last, prior distributions are provided in Section \ref{sec:priors}.}

\subsection{The Cox model with random right censoring}
\label{sec:model}

The observations $X=X^n$ are $n$ independent identically distributed (i.i.d.) triplets given by 
$X = ((Y_1, \delta_1, Z_1), \dots, (Y_n, \delta_n, Z_n))$. The observed $Y_i\in\mathbb{R}^+$ are censored versions of (unobserved) survival times $T_i\in\mathbb{R}^+$,
%\sbl{[indeed, it is $\mathbb{R}^+$, not $\mathbb{R}$!]}, 
with $\delta_i$ indicator variables informing on whether $T_i$ has been observed or not: that is, $Y_i = T_i \wedge C_i$ and $\delta_i = \mathbbm{1}(T_i \leq C_i)$, where $C_i$'s are i.i.d. censoring times.  The  variables $Z_1, \dots, Z_n \in \mathbb{R}^p$, $p$ fixed, are called covariates. 

For a fixed covariate vector $z\in\mathbb{R}^p$ and $t>0$, define the {\it conditional hazard rate}  $\lambda(t\given z) = \lim_{h\to0} h^{-1} P(t \leq T \leq t+h \given T \geq t, Z = z)$. 
The Cox model assumes, for some $\theta\in\mathbb{R}^p$ and denoting by $\theta'z$ the standard inner product in $\mathbb{R}^p$, 
\[ \lambda(t \given z) = e^{\theta'z} \lambda(t),\]
where $\lambda(t)$ is the {\it baseline hazard function}. 
The {\it conditional cumulative hazard function} is defined as $\Lambda(\cdot \given z) = \int_0^\cdot \lambda(u \given z) du = e^{\theta'z} \int_0^\cdot \ld(u) du = e^{\theta'z} \Lambda(\cdot)$ and the {\it survival function conditional on $z$} is denoted by 
$S(\cdot \given z) =  
\exp(-\Ld(\cdot \given z)) = 
\exp(- e^{\theta'z} \Lambda(\cdot) )$. 

Assuming the baseline hazard rate is positive, one can alternatively make inference on the log-hazard. The unknown parameters of the Cox model are then 
\[ \eta = (\theta, r),\qquad \text{where } r = \log \lambda.\] 
The goal is to estimate the pair $\eta=(\theta,r)$. We denote by $\eta_0=(\theta_0,r_0)$ the true values of the parameters (and similarly for the related quantities $\lambda_0, \Lambda_0$).

We now give a set of standard assumptions used in this paper. 
First, we assume both $T$ and $Z$ admit a continuous density function, $f_T(\cdot)$ and $f_Z(\cdot)$ respectively. Given $Z$, the survival time $T$ and the censoring time $C$ are independent. 
At the end of the follow-up, some individuals are still event free and uncensored such that 
$P(T > \varrho \given Z = z) > 0$ and $P(C > \varrho \given Z = z) = P(C = \varrho \given Z = z) = 0$ for some fixed time $\varrho$. 
The censoring $C$ is assumed to follow a distribution $G$ and admit a density such that
    \begin{align*}
        p_C(u) = g_z(u) \mathbbm{1}\{0\leq u < \varrho\} + \bar G_z(\varrho) \mathbbm{1}\{u = \varrho\}
    \end{align*}
    with respect to $\text{Leb}([0, \varrho]) + \delta_\varrho$, where $\text{Leb}(I)$ is the the Lebesgue measure on $I$.
Without loss of generality, we assume $\varrho = 1$ throughout the paper.

Based on this set-up, the joint density function of the triple $(y, \delta, z)$ is given by
\begin{equation}
\begin{split}
\label{eqn:density-function}
    f_\eta (y, \delta, z) = 
   & \left(
        g_z(y) e^{- \Lambda(y) e^{\theta'z}}
    \right)^{1-\delta} 
    \left(
        \bar{G}_z (y) \lambda(y) e^{\theta' z - \Lambda(y) e^{\theta'z}}
    \right)^\delta 
    f_Z(z)
    \mathbbm{1}(y < t)\\
    & + 
    \left(
        \bar{G}_z(t) e^{-\Lambda(t) e^{\theta'z}}
    \right) 
    f_Z(z)
    \mathbbm{1}(\delta = 0, y = t),
\end{split}
\end{equation}
where $g_z(\cdot)$ is a continuous density of $C$ given $Z = z$,
$\bar{G}_z(\cdot) = 1 - G_z(\cdot -)$, where $G_z(\cdot)$ is the cumulative distribution function of $g_z(\cdot)$.

Let $\ell_n(\eta) = \sum_{i=1}^n \log f_\eta (X_i)$ be the log-likelihood function,
the likelihood ratio is given by 
\begin{align}
\label{log-likelihood}
\ell_n(\eta) - \ell_n(\eta_0)
= 
\sum_{i=1}^n \left(\delta_i \{(\theta-\theta_0)'Z_i + (r-r_0)(Y_i) \} - \Ld(Y_i) e^{\theta'Z_i} + \Ld_0(Y_i) e^{\theta_0'Z_i}\right).
\end{align}
From (\ref{log-likelihood}), one sees that the log-likelihood ratio does not depend on $g_z(y)$ and $\bar G_z(\cdot)$, thus one does not need to model $g_z(y)$ in order to make inference on $\eta$. 

\subsection{Information--related quantities} \label{sec:infonot}

We now introduce some of the key quantities arising in the study of the Cox model: these are all related to the information operator (extending the usual Fisher information in parametric models) arising from the LAN--expansion in the model (see also Section \ref{sec:bg}).  
For a bounded function $b(\cdot)$ on $[0, 1]$ and a cumulative hazard function $\Lambda(\cdot)$, 
we denote 
$\Lambda\{b\}(\cdot) = \int_0^\cdot b(u) d\Ld(u)$ and, in slight abuse of notation, we set $\Lambda \{b\} = \Lambda \{b\}(1)$. 
The following notations are commonly used in the literature for the Cox model (see e.g., Section VIII 4.3 of \citet{andersen93} and Section 12.3.3 of \citet{ghosal17}):
\begin{align}
\label{M0}
    M_0(u) &= \mathbb{E}_{\eta_0} \left(e^{\theta_0'Z} \mathbbm{1}_{u \leq T}\right) = 
    \int \bar{G}_z(u) e^{\theta_0'z - \Lambda_0(u) e^{\theta_0'z}} f_Z(z) dz,\\
    \label{M1}
    M_1(u) & = \mathbb{E}_{\eta_0} \left(Z e^{\theta_0'Z} \mathbbm{1}_{u \leq T}\right) 
    = \int z \bar{G}_z(u) e^{\theta_0'z - \Lambda_0(u) e^{\theta_0'z}} f_Z(z) dz,\\
    \label{M2}
	M_2(u) & = \mathbb{E}_{\eta_0} \left(ZZ' e^{\theta_0'Z} \mathbbm{1}_{u \leq T}\right) 
    = \int zz' \bar{G}_z(u) e^{\theta_0'z - \Lambda_0(u) e^{\theta_0'z}} f_Z(z) dz.
\end{align}
The least favorable direction is defined as $\gamma_{M_1}(\cdot) = (M_1/M_0)(\cdot)$. 
For a bounded function $b(\cdot)$ on $[0,1]$, we let $\gamma_b(\cdot) = (b/M_0)(\cdot)$.
The efficient information matrix of the Cox model is denoted by $\tilde I_{\eta_0}$ and is given by
\begin{align}
\tilde I_{\eta_0} = & \Lambda_0\{M_2(\cdot) - \gamm(\cdot) \gamm'(\cdot) M_0(\cdot)\}.
\label{eff-info-mtx}
\end{align}
For any $\vartheta \in \mathbb{R}^p$ and $g \in L^2\{\Ld_0\}$ (i.e., $\int g^2d\Ld_0<\infty$), we define
\begin{align}
\label{Wn}
W_n(\vartheta, g)
    =
    \frac{1}{\sqrt{n}}
    \sum_{i=1}^n \left\{
    \delta_i \left(
    \vartheta'Z_i + g(Y_i)
    \right)
    - e^{\theta_0'Z_i}
    \left(
    \vartheta'Z_i \Lambda_0(Y_i) + (\Lambda_0 g)(Y_i)
    \right)
    \right\}.
\end{align}
This quantity corresponds to the empirical process part of the LAN expansion of the Cox model (see (\ref{lan}) in Section \ref{sec:bg} of \citet{ning22supp}).

\subsection{Structural assumptions}
\label{sec:assumption}

The following fairly mild conditions are assumed on the unknown quantities of the model.  For some positive constants $c_1, \dots, c_8$,
\begin{enumerate}[label=\textbf{(i)}]
    \item \label{asp:i} \text{the random variable }$Z\in\mathbb{R}^p$ is bounded (i.e. $\|Z\|_\infty\le c_1$ a.s.);

\end{enumerate}    
\begin{enumerate}[label=\textbf{(ii)}]    
	\item \label{asp:ii} 
         $\|\theta_0\|_\infty \le c_2$;
\end{enumerate}	
Note that from \ref{asp:i} and \ref{asp:ii}, one can bound $e^{\theta_0'z} \le e^{pc_1c_2}$. Also, suppose
\begin{enumerate}[label=\textbf{(iii)}]    
    \item \label{asp:iii} $c_3 \leq \inf_{t \in [0, \varrho]} \lambda_0(t) \leq \sup_{t \in [0, \varrho]} \lambda_0(t) \leq c_4$ and $r_0 = \log \ld_0 \in \mathcal{H}(\beta, D)$, where $\beta, D > 0$;
\end{enumerate}    
\begin{enumerate}[label=\textbf{(iv)}]
    \item \label{asp:iv} $g_z$ is a continuous density and
    $c_5 \leq \inf_{t \in [0, \varrho]} g_z(t) \leq \sup_{t \in [0, \varrho]} g_z(t) \leq c_6$;
\end{enumerate}

Assumptions \ref{asp:i}--\ref{asp:iv} are common in the related literature; e.g., see \citet{cast12} (p.17). 
As $\theta_0, \ld_0, g_z(u)$ are all assumed to be bounded from above and below, \ref{asp:i}--\ref{asp:iv} imply
 \begin{enumerate}[label=\textbf{(v)}]    
	\item \label{asp:v} For $M_2(u)$ and $\tIinv$ in (\ref{M2}) and (\ref{eff-info-mtx}) respectively,
	{$\Ld_0\{M_2(\cdot)\} - c_7 I_p \succeq 0$} and $\|\tIinv\|_{(\infty,\infty)} \leq c_8$.
\end{enumerate} 

The above conditions are mostly assumed for technical simplicity: boundedness of the true vector $\theta$ enables one to take a prior with bounded support, which is particularly helpful in order to carry out likelihood expansions. Attempting to remove this condition would lead to delicate questions on likelihood remainder terms and is beyond the scope of this work. The smoothness condition on the (log--) hazard is quite mild for smooth hazards: it includes for instance the case of Lipschitz hazards. Another interesting setting would be the one of piecewise constant hazards. It could be treated with the techniques developed of this paper (a given histogram can always be approximated arbitrarily well in the $L^2$--sense by a Haar-histogram, and then the problem becomes --nearly-- parametric) although it would require a somewhat separate treatment: we refrain from providing theory here for this case, although we consider it as one of the examples of the simulations study in Section \ref{sec:sim}, where simulations show very good behavior in this setting as well.  Henceforth we assume that \ref{asp:i}--\ref{asp:iv}  hold without explicit reference. 

\subsection{Prior distributions}
\label{sec:priors}

Priors for $\theta$ and $\lambda$ are chosen independently. 
The prior for $\theta$ is chosen as follows: 
\begin{enumerate}[label=\textbf{(T)}]
\item \label{prior:T} 
Let $\theta_j$ be the $j$-th coordinate of $\theta$, 
$\pi(\theta) = \bigotimes_{j=1}^p \pi(\theta_j) = \bigotimes_{j=1}^p f(\theta_j)\mathbbm{1}_{[-C, C]}$ for some constant $C > {c_2}$.
Examples for $f(\theta_j)$ include the uniform density, i.e., $f(\theta_j) = 1$, the {\em truncated} (to $[-C,C]$) $\tau$--Subbotin density \citep{subb23}, which includes the {\em truncated} Laplace (when $\tau = 1$) density and {\em truncated} normal (when $\tau = 2$) density as special cases. 
\end{enumerate}

Imposing the truncation does not seem necessary in practice, as we found in the simulation study. However, for deriving our theoretical results, and similar to \citet{cast12} and \citet{ghosal17}, we assume that at least some upper-bound on $\|\theta\|_\infty$
 is known, as discussed in the previous subsection (so that one can take e.g. $C=c_2$ in \ref{asp:ii}). 

For the prior on $\ld$, two classes of piecewise constant priors are considered throughout the paper:
\begin{enumerate}[label=\textbf{(H)}]
\item \label{prior:H} {\it Random histogram prior.} 
For $k \geq 1, L' = L+1$, and $L$ an $n$-dependent deterministic value to be specified below,
let
\begin{align}
\label{histogram-prior}
\lambda_H = \sum_{k=0}^{2^{L+1} - 1} \lambda_k \mathbbm{1}_{I_k^{L+1}}, 
\end{align}
where 
$I_0^{L'} = [0, 2^{-L'}], I_k^{L'} = (k 2^{-L'}, (k+1)2^{-L'}]$,
and ($\lambda_k$) are independent random variables. We consider putting the following two types of priors on the $k$-th histogram height, $\lambda_k$:
\end{enumerate}

\begin{itemize}
\item[(i)] An independent Gamma prior: $\lambda_k \sim \text{Gamma}(\alpha_0, \beta_0)$  are i.i.d. variables, for some fixed positive $\alpha_0,\beta_0$.
\item[(ii)] A dependent Gamma prior: let $\lambda_0 \sim \text{Gamma}(\alpha_0, \beta_0)$ and $\lambda_k| \lambda_{k-1} \sim \text{Gamma}(\alpha, \alpha/\lambda_{k-1})$ for some positive constant $\alpha$. 
Then for $k \geq 1$, $E(\lambda_k \given \ld_{k-1}) = \ld_{k-1}$ and $\text{Var}(\ld_k \given \ld_{k-1}) = \ld_{k-1}^2/\alpha$. 
\end{itemize}

\begin{enumerate}[label=\textbf{(W)}]
\item \label{prior:W} {\it Haar wavelet prior.} 
Let $r_S = \log \lambda_S$ and again for $L$ to be specified below, let us set
\begin{align} 
\label{wavelet-prior}
r_S = \sum_{l = -1}^{L} \sum_{k=0}^{2^l - 1} \sigma_l Z_{lk} \psi_{lk},
\end{align}
where $Z_{lk}$ are random variables, $\sigma_l = 1$, $0 \leq l \leq L$, and $(\psi_{lk})$ are Haar wavelet basis. 
\end{enumerate}
Although a variety of densities can be considered for $Z_{lk}$, we specifically consider for simplicity the standard normal density and the standard Laplace density (other choices of $\sigma_l$ e.g. $2^{-l/2}$ are possible). Note that both densities give a non-conjugate prior for $r_S$. 
 
Also note that the random histogram prior can be viewed as a special case of the Haar wavelet prior if one allows for possibly dependent variables $Z_{lk}$ (and possibly different values for $\sigma_l$). 
{The above priors can mostly be chosen free of dependence in the constants, except for $c_2$ for which we need to know at least an upper-bound for $\theta_0$. Such assumption is unavoidable for using bounded priors, as they put no mass outside their support.}

{\em Choice of the parameter $L$}. For results on the specific priors as above, we consider the choice of cut-off $L=L_n$ defined as, for $\beta>0$ the assumed regularity of $r_0=\log\lambda_0$ (see \ref{asp:iii}),
\begin{align} \label{choicel}
2^{L_n}=2^{L_n(\beta)} & \circeq \left(\frac{n}{\log{n}}\right)^{\frac{1}{2\beta+1}},
\end{align}
where $\circeq$ means that one picks a closest integer solution in $L_n$ of the equation. If the regularity of the true $r_0$ is not known in advance, as is usually the case in practice, then all the limiting shape (Bernstein--von Mises) results below go through if one replaces $\beta$ by $1/2$ in \eqref{choicel} (note also that all the main results, except the Hellinger rate which requires no minimal smoothness, require a regularity $\beta>1/2$). In other words, for the semiparametric results, it is enough to `undersmooth'. The strict knowledge of $\beta$ is only required if one wishes to obtain an optimal minimax supremum-norm contraction rate (see Section \ref{sec:sup-norm} for more on this).

%%%%%%%%%%%%%%%%%%%%%%%%%%%%%%%%%%%%%%%%%%%%%%%%%%%%
%%%%%%%%%%%%%%%%%%%%%%%%%%%%%%%%%%%%%%%%%%%%%%%%%%%%
%%%%%%%%%%%%%%%%%%%%%%%%%%%%%%%%%%%%%%%%%%%%%%%%%%%%
\section{Main results}
\label{sec:main-results}

Let us give a brief outline of our results.
Section \ref{sec:hellinger} provides a preliminary contraction result in Hellinger distance.
Section \ref{sec:jointBvM} presents a joint Bernstein--von Mises (BvM) theorem for linear functionals of $\theta$ and $\lambda$. 
Section \ref{sec:npBvM} derives a Donsker theorem for the joint posterior distribution of $\theta$ and the cumulative hazard function $\Ld$. 
This result leads to the Donsker theorem for the posterior of the conditional cumulative hazard and survival functions.
A supremum-norm convergence rate for {the hazard function conditional on $z$} is obtained in Section \ref{sec:sup-norm}.
In Sections \ref{sec:jointBvM}-\ref{sec:sup-norm}, we provide generic conditions that are suitable for a wide range of priors. In Section \ref{sec:application}, we verify those conditions for those specific choices of priors listed in Section \ref{sec:priors}.

\subsection{A key preliminary contraction result}
\label{sec:hellinger}

We start by obtaining a preliminary Hellinger contraction rate for the posterior distribution for the priors considered above.  

Define the rate, for $\beta\in(0,1]$,
\begin{equation} \label{rateveps}
{{\nu}}_n={{\nu}}_{n,\beta}=\left(\frac{\log{n}}{n}\right)^{\frac{\beta}{2\beta+1}}.
\end{equation}

Let $\epsilon_n=o(1)$ be a sequence such that $n\epsilon_n^2\to\infty$ as $n \to \infty$. Define $\zeta_n=\zeta_n(\epsilon_n)= 2^{L_n/2} \epsilon_n + 2^{-\beta L_n}$ and 
\begin{equation}\label{defan}
A_n = \{\eta:\  \|\theta - \theta_0\| \leq \epsilon_n, \|\lambda - \lambda_0\|_\Lone \leq \epsilon_n, \|\ld - \ld_0\|_\infty \leq \zeta_n\}
\end{equation}
Let us consider the following condition:
\begin{enumerate}[label=\textbf{(P)}]
\item \label{cond:P} 
The sequences $L_n,\epsilon_n$ verify $L_n = o(\rn \epsilon_n)$, $L_n^2 = o(1/\epsilon_n)$, and $\rn \epsilon_n^2 L_n = o(1)$ and, for $A_n$ as in \eqref{defan},
\[ \Pi(A_n^c \given X) = o_{P_{\eta_0}}(1).\]
\end{enumerate}
Condition \ref{cond:P} requires the posterior distribution to contract in a certain sense around the true pair $(\theta_0,r_0=\log{\lambda_0})$. In order to derive such a result, one may first apply the general contraction rate theorem of \cite{ghosal00}. This, however, entails a rate for the overall density $f_\eta$ in the Cox model only, not the parameters themselves. The main difficulty here is then to derive results on $\theta$ and $\lambda$ separately. The rate $\epsilon_n$ can be thought of as a typical (possibly optimal) nonparametric rate. We call condition  \ref{cond:P} a preliminary contraction result, because faster  rates both for $\theta$ and for $\lambda$ in the supremum norm can be derived, as will be seen below. In fact, for $\theta$, it is expected that the posterior contracts at parametric, near $1/\sqrt{n}$ rate; a much more precise result is obtained in Section \ref{sec:jointBvM} below in the form of a BvM theorem. 

The next lemma shows that condition \ref{cond:P} is indeed satisfied for the examples of priors introduced in Section \ref{sec:priors}. 
\begin{lemma}
\label{lemma1}
Consider the Cox model 
with priors as specified in {\rm \ref{prior:T}}, and {\rm \ref{prior:H}} or {\rm \ref{prior:W}} with $L=L_n$ as in \eqref{choicel}. 
Then for any $\beta\in(0,1]$ and ${{\nu}}_n$ as in \eqref{rateveps},
\begin{align*}
& \Pi[\{\eta:\ h^2(f_{\eta_0}, f_\eta) \gtrsim {{\nu}}_{n}^2\} \given X]=o_{P_{\eta_0}}(1).
\end{align*}
Further, for any $\beta\in(1/2,1]$, condition {\rm \ref{cond:P}} is satisfied for these priors for $\epsilon_n={{\nu}}_n$ as in \eqref{rateveps}.
\end{lemma}

Let us now briefly comment on the preliminary supremum-norm rate $\zeta_n$ for $\ld$ entailed by \ref{cond:P}. For some cut--offs $L_n$, the rate can be slow. However, it is a $o(1)$ as soon as $\epsilon_n=o(2^{-L_n/2})$. For the typical choice of $L_n$ in \eqref{choicel}, this only requires that $\beta>1/2$, which corresponds to a preliminary rate faster than $n^{-1/4}$. This is much less than what is required for Theorem 5 in \cite{cast12} or Theorem 12.12 in \cite{ghosal17}, where a preliminary rate faster than $n^{-3/8}$ is needed (note that the latter rate rules out the use of regular histograms as priors, since these can get only a rate $n^{-1/3}$ at best). In Section \ref{sec:sup-norm}, we show the rate $\zeta_n$ can be improved by adopting a multiscale analysis approach. For this, a BvM theorem for linear functionals of $\ld$ will be needed: it is a consequence of the joint BvM derived in the next section. 

From Section \ref{sec:jointBvM}  to Section \ref{sec:sup-norm}, we work with a generic histogram prior of the form \eqref{wavelet-prior} with cut--off $L=L_n$ (which includes both {\rm \ref{prior:H}} and {\rm \ref{prior:W}} as special cases), under the above condition \ref{cond:P}. This way, the reader can directly see what generic conditions underpin our results, and adapt these to other relevant families of priors not considered here for the sake of brevity. For instance, smoother wavelet bases $(\psi_{lk})$ can be used in the prior definition and require only minor adaptations of the proofs (in a similar way as in  \cite{cast20survival} for the simple nonparametric survival model); although we do not prove this here for brevity, using these priors would enable one to derive optimal contraction rates in the supremum norm for arbitrary regularities $\beta>1/2$. We come back to the specific examples of priors {\rm \ref{prior:H}} and {\rm \ref{prior:W}} in Section \ref{sec:application}. 

\subsection{The joint BvM theorem for the linear functionals of $\theta$ and $\lambda$} 
\label{sec:jointBvM}

Let us consider the joint estimation of the two linear functionals defined by
\begin{align}
\label{linear-fns}
\varphi_a(\theta) = \theta'a, \quad 
\varphi_b(\lambda) = \langle b, \lambda \rangle = \int_0^1 b(u) \lambda(u) du = \Lambda\{b\},
\end{align}
for fixed $a \in \mathbb{R}^p$ and $b \in L^2(\Lambda)$.
 Let us recall that we work under the generic form of prior  \eqref{wavelet-prior} with cut--off $L=L_n$ and generic condition \ref{cond:P}. Let us also recall the notation from Section \ref{sec:infonot}: $M_0, M_1, M_2$, $\gamma_b(\cdot) = (b/M_0)(\cdot)$ for a bounded function $b$, and $\tilde I_{\eta_0}$. 

Consider the following conditions:
\begin{enumerate}[label=\textbf{(B)}]
\item \label{cond:B} 
Let $P_{L_n}(\cdot)$ be the orthogonal projection onto $\mathcal{V}_{L_n} := \text{Vect}\{\psi_{lk}, \ l \leq L_n, \ 0 \leq k < 2^l\}$, the subspace of $L^2[0,1]$ spanned by the first $L_n$ wavelet levels, and denote $\gambln = P_{L_n}(\gamb)$ and $\gammln = P_{L_n}(\gamm)$.
For any fixed $b \in L^\infty[0,1]$ and $\epsilon_n$ in {\ref{cond:P}}, 
$$
\rn \epsilon_n \|\gamb - \gambln\|_\Linfty = o(1).
$$
\end{enumerate}

Condition \ref{cond:B} is sometimes called {\it no-bias} condition, and holds true if $b$ is sufficiently smooth and (or) the preliminary contraction rate $\epsilon_n$ is fast enough. 
Next, let $h = (t, s) \in \mathbb{R}^2$,
for fixed $a \in \mathbb{R}^p$ and $b \in L^2(\Ld)$, consider the two local paths:
\begin{align}
&\theta_h = \theta - \frac{t \tIinv a}{\rn} + \frac{s \tIinv \Lambda_0\{b \gamm\}}{\rn},
\label{path-theta}\\
r_h = r +& \frac{t \gammln' \tIinv a}{\rn} - \frac{s\gambln}{\rn} - \frac{s\gammln' \tIinv \Lambda_0\{b\gamm\}}{\rn}.
\label{path-r}
\end{align}
\begin{enumerate}[label=\textbf{(C1)}]
\item \label{cond:C1} ({\it Change of variables condition}) with $r = \log \lambda$ and $r_0 = \log \lambda_0$, let $\eta_0 = (\theta_0, r_0)$ and $\eta_h = (\theta_h, r_h)$, 
where $\theta_h$ in (\ref{path-theta}) and $r_h$ in (\ref{path-r}) with $a, b$ to be specified below
and $A_n$ as in \ref{cond:P}, suppose 
$$
\frac{
\int_{A_n} e^{\ell_n(\eta_h) - \ell_n(\eta_0)} d\Pi(\eta)
}{
\int e^{\ell_n(\eta) - \ell_n(\eta_0)} d\Pi(\eta)
} 
= 1 + o_{P_{\eta_0}}(1).
$$
\end{enumerate}

Condition \ref{cond:C1} is often called {\em change of variables condition}: indeed, one natural way to check it is via controlling the change in distribution from $\eta\sim \Pi$, to the distribution induced on $\eta_h$. For priors such as \ref{prior:H} and \ref{prior:W}, this can be checked by posing a change of variables with respect to the Lebesgue measure on $\mathbb{R}^{L_n}$. 
The verification of these conditions for these priors is given in  Section \ref{verify-cov}. 

For $\eta \sim \Pi(\cdot \given X)$, $\mu = (\mu_1, \mu_2) \in \mathbb{R}^2$, $a \in \mathbb{R}^p$, and $b \in L^2(\Lambda)$,
let us define the map
$$
\tau_\mu: \eta \to \rn 
(\theta'a - \mu_1,
\langle \lambda, b \rangle - \mu_2),
$$
and let $\Pi(\cdot \given X) \circ \tau_\mu^{-1}$ be the distribution induced on $\rn (\theta'a - \mu_1,
\langle \lambda, b \rangle - \mu_2)$.
We are ready to present the joint BvM theorem for the bivariate functions $\varphi_a(\theta)$ and $\varphi_b(\ld)$. 

\begin{theorem}
\label{thm:bvm-linear-fns}
Let $a$ and $b$ be fixed elements of $\mathbb{R}^p$ and $L^2(\Lambda_0)$ respectively, for any $b$ that satisfies {\rm \ref{cond:B}},
suppose the prior for $\eta = (\theta, \ld)$ is chosen such that {\rm \ref{cond:P}} and {\rm \ref{cond:C1}} hold. Then
\begin{align}
\label{thm1-1}
\mathcal{B}_{\mathbb{R}^2} \left( \Pi(\cdot \given X) \circ \tau_{\hat \varphi}^{-1}, 
\mathcal{L}(a'\mathbb{V}, \Upsilon_b - \mathbb{V}\Ld_0\{b\gamm\} )
\right) 
\stackrel{P_{\eta_0}}{\to} 0,
\end{align}
where $\mathbb{V}$ and $\Upsilon_b$ are independent, 
$\mathbb{V} \sim N\left(0, \tIinv\right)$ and $\Upsilon_b \sim N(0, \Lambda_0\{b \gamb\})$, 
$\mathcal{B}_{\mathbb{R}^2}$ is the bounded Lipschitz metric between distributions on $\mathbb{R}^2$, and
$\hat\varphi = (\hat\varphi_a(\theta), \hat\varphi_b(\ld))$ is given by 
\begin{align*}
\hat \varphi_a(\theta) & := \varphi_a(\hat \theta) = \varphi_a(\theta_0) 
+ \frac{1}{\rn}W_n\left(\tIinv a, \ - \gamm' \tIinv a\right),\\
\hat \varphi_b(\ld) : = \varphi_b(\hat \lambda) & = \varphi_b(\lambda_0) + 
\frac{1}{\rn} W_n \left(- \tIinv \Ld_0\{b \gamm\}, \gamb + \gamm' \tIinv \Ld_0 \{b\gamm\}\right).
\end{align*}
\end{theorem}

The centering sequences in the last display of the statement can be seen to be `efficient' ones from the semiparametric perspective (see e.g. \cite{vdvAS}, Chapter 25). An important added value to the {\em joint} BvM (in contrast to individual limiting statement for marginal coordinates) is that it captures the dependence between $\theta$ and $\lambda$: a practical application is given in Figure \ref{fig:joint-density}.    
The result enables to consider many combinations of functionals by choosing specific $a \in \mathbb{R}^p$ and $b \in L^2(\Ld_0)$. For example, let $a = (1, 0, \dots, 0)$ and $b = \mathbbm{1}_{[0,1]}$: Theorem \ref{thm:bvm-linear-fns} implies a joint joint BvM for $(\theta_1, \Ld(1))$, where $\theta_1$ is the first coordinate of $\theta$ and $\Ld(1) = \int_0^1 \ld$ is the cumulative hazard function at time one.
The limiting distribution is given in the next corollary, where, in addition, we center the joint posterior at efficient frequentist estimators for $\theta$ and $\Ld(1)$.

\begin{corollary}[Joint BvM for $\theta_1$ and $\Ld(1)$]
\label{cor-1}
Consider the Cox model with the density function in (\ref{eqn:density-function}), 
let $a = (1, 0, \dots, 0)$ and $b = 1$, and $\tau_{(\hat \theta_1, \hat \Ld(1))}$ be the map such that $$
\tau_{(\hat \theta_1, \hat \Ld(1))}: \eta \rightarrow \rn \left(\theta_1 - \hat \theta_1, \Ld(1) - \hat \Ld(1)\right),
$$
where $\hat \theta$ is the maximum Cox partial likelihood estimator and $\hat \Ld(1)$ is the Breslow estimator. 
Denote $\Pi(\cdot \given X) \circ \tau_{(\hat\theta_1, \hat \Ld(1))}^{-1}$ as the distribution induced on 
$\rn (\theta_1 - \hat \theta_1, \Ld(1) - \hat \Ld(1))$, then 
under the same conditions as in Theorem \ref{thm:bvm-linear-fns},  
\begin{align}
\label{BvM-2}
\mathcal{B}_{\mathbb{R}^2}
\left(
\Pi(\cdot \given X) \circ \tau_{(\hat \theta_1, \hat \Ld(1))}^{-1}, \
\mathcal{L}(a' \mathbb{V}, \Upsilon_1 - \mathbb{V}\Ld_0\{\gamm\})
\right) \stackrel{P_{\eta_0}}{\to} 0,
\end{align}
where 
$\Upsilon_1 \sim N(0, \Ld_0\{M_0^{-1}\})$ and $\mathbb{V} \sim N(0, \tIinv)$ are independent. 
\end{corollary}

An immediate practical implication of the BvM theorem in Corollary \ref{cor-1} is that two-sided quantile credible sets for $\hat\theta_1$ (or more generally for any given coordinate $\theta_j$, $j = 1, \dots, p$) are asymptotically optimal confidence sets from the perspective. Results in this vein can also be derived for the survival function in the functional sense: this is the object of the next section.

%\begin{corollary}[BvM for $\theta \in \mathbb{R}^p$]
%\label{cor-2}
%{
%Consider the Cox model with the density function in (\ref{eqn:density-function}), 
%let $\hat \theta$ be the maximum Cox partial likelihood estimator
%and $\mathbb{V} \sim N(0, \tIinv)$, then 
%under the same conditions as in Theorem \ref{thm:bvm-linear-fns},  
%\begin{align*}
%%\label{BvM-theta}
%\mathcal{B}_{\mathbb{R}^p}
%\left(
%\Pi(\rn (\theta - \hat \theta) \given X), \ \mathbb{V}
%\right) \stackrel{P_{\eta_0}}{\to} 0.
%\end{align*}
%}
%\end{corollary}

%%%%%%%%%%%%%%%%%%%%%%%%%%%%%%%%%%%%%%%%%%%
\subsection{Joint Bayesian Donsker theorems}
\label{sec:npBvM}

We now present the second main result in this paper, the Bayesian Donsker theorem for the joint posterior distribution of $\theta$ and the cumulative hazard function $\Ld(\cdot)$. 

Let us denote
\begin{align} 
W_n^{(1)} = W_n^\star \left(\tIinv, \ -\gamm' \tIinv\right),
\label{Wn-1}
\end{align}
where $W_n^{(1)}$ is a $p$-dimensional vector and
\begin{equation*}
\begin{split}
W_n^\star \left(\tIinv,\ -\gamm' \tIinv\right) 
= 
\frac{1}{\rn} \sum_{i=1}^n
\left\{
\delta_i \tIinv \left(Z_i - \gamm \right)
- e^{\theta_0'Z_i} \tIinv \left(Z_i \Ld_0(Y_i) - (\Ld_0{\gamm})(Y_i) 
\right)
\right\},
\end{split}
\end{equation*}
and given $b\in L^2(\Lambda_0)$,
\begin{align} 
W_n^{(2)}(b) = W_n \left(-\tIinv \Ld_0\{b \gamm\}, \gamb + \gamm' \tIinv \Ld_0 \{b \gamm\} \right).
\label{Wn-2}
\end{align}

Define the centering sequences for $\theta$ and $\ld$ as follows:
$$
T_n^\theta = \theta_0 + W_n^{(1)}/\rn,
$$
and, for a given sequence $L_n$, 
\begin{align*}
\langle T_{n}^\ld, \psi_{lk} \rangle 
= 
\begin{cases}
\langle \lambda_0, \psi_{lk} \rangle 
+
W_n^{(2)}(\psi_{lk})/\rn & \quad \text{if} \ l\leq L_n,\\
\ 0 & \quad \text{if} \ l > L_n.
\end{cases}
\end{align*}

The Donsker theorem requires, in addition to \ref{cond:C1}, a similar condition,  where $t$ and $s$ are allowed to increase with $n$. This condition is stated as follows:

\begin{enumerate}[label=\textbf{(C2)}]
\item \label{cond:C2} (Change of variables condition, version 2) with the same notation as in \ref{cond:C1}, let $\eta_h = (\theta_h, r_h)$, $\theta_h$ and $r_h$ be the local paths in \eqref{path-theta} and \eqref{path-r} respectively with $a, b$ to be specified below,
for $n$ large enough and any $|t|, |s| \leq \log n$, one assumes, for $A_n$ as in \ref{cond:P} and some constant $C_1 > 0$,

$$
\frac{
\int_{A_n} e^{\ell_n(\eta_h) - \ell_n(\eta_0)} d\Pi(\eta)
}
{
\int e^{\ell_n(\eta) - \ell_n(\eta_0)} d\Pi(\eta)
}
\leq e^{C_1(1 +t^2 + s^2)}.
$$
\end{enumerate}

Condition \ref{cond:C2} is similar to \ref{cond:C1}. A major difference between the two conditions is that in \ref{cond:C2}, $t$ and $s$ are allowed to increase with $n$; however, in \ref{cond:C1}, $t$ and $s$ are fixed. 

We further require the rates $\epsilon_n$ and $\zeta_n$ and the cut-off $L_n$ in \ref{cond:P} to satisfy
\begin{align}
\label{ez-rate-condition}
\rn \epsilon_n 2^{-L_n} = o(L_n^{-5/2}), \quad
\zeta_n L_n^2 = o(1),
\end{align}

\begin{theorem}[Joint Bayesian Donsker theorem]
\label{thm:donsker-joint}
 Suppose the prior for $\eta = (\theta, \eta)$ is chosen such that both
{\rm  \ref{cond:P}} and (\ref{ez-rate-condition}) hold. Suppose  {\rm \ref{cond:C1}} holds for $a = z$, any fixed $z \in \mathbb{R}^p$, and any $b \in \mathcal{V}_{\mathcal{L}} = \text{Vect}\{\psi_{lk}, \ l \leq L_n, \ 0 \leq k < 2^l \}$ 
for a fixed $\mathcal{L}$ and {\rm \ref{cond:C2}} holds uniformly for $a = z$, any fixed $z \in \mathbb{R}^p$, and any $b = \psi_{LK}$ with $0 \leq L \leq L_n$ and $0 \leq K < 2^L$.

Let $\mathcal{L}((\theta, \Ld(\cdot)) \in \cdot \given X)$ be the distribution induced on $\theta$ and $\Ld(\cdot) = \int_0^\cdot \ld$ and $\mathbb{T}_n^\ld(\cdot) = \int_0^\cdot T_n^\ld$ be the centering for $\Ld(\cdot)$. 
Denote $\mathbb{B}(\cdot)$ as standard Brownian motion and set $U_0(\cdot) = \int_0^\cdot (\ld_0/M_0)(u) du$.
Let $\mathbb{V} \sim N(0, \tIinv)$ that is independent of $\mathbb{B}(\cdot)$.
Then, as $n \to \infty$,
\begin{align}
\label{eqn:donsker-joint}
\mathcal{B}_{\mathbb{R}^p \times \mathcal{C}([0,1])} \left(
\mathcal{L} \left(\rn (\theta - T_n^\theta, \Ld(\cdot) - \mathbb{T}_n^\ld(\cdot) ) \given X \right), \
\mathcal{L} \left( \mathbb{V}, \mathbb{B}(U_0(\cdot)) - \mathbb{V}' \Ld_0\{\gamm\}(\cdot) \right)
\right) \stackrel{P_{\eta_0}}{\to} 0,
\end{align}
where $\mathcal{B}_{\mathbb{R}^p \times \mathcal{C}([0,1])}$ is the bounded-Lipschitz metric on $\mathbb{R}^p \times \mathcal{C}([0,1])$.
\end{theorem}

\begin{remark}
While the proof is left to the supplemental material \citep{ning22supp}, a key step for obtaining the joint Bayesian Donsker theorem, following ideas from \cite{cast14a}, is to establish first a BvM for $(\theta,\lambda)$ in an appropriate space. However, unlike in \cite{cast14a} where one can work directly on the nonparametric quantity of interest, here due to the split semiparametric model at hand, one needs to prove a {\em joint} nonparametric BvM for the pair $(\theta,\lambda)$, see Proposition \ref{prop-tightness}.
This result is new in this context and is of independent interest for proving similar results in other semiparametric models.  
\end{remark}

The centerings $T_n^\theta$ and $\mathbb{T}_n^\ld$ in Theorem \ref{thm:donsker-joint} can be replaced with any efficient estimators for $\theta$ and $\Ld$: the next corollary formalizes this with centering at standard frequentist estimators.

\begin{corollary}
\label{cor-donsker}
Let $\hat \theta$ be the maximum Cox partial likelihood estimator and $\hat \Ld(\cdot)$ be the Breslow estimator,
then,
under the same conditions as in Theorem \ref{thm:donsker-joint}, 
as $n \to \infty$,
\begin{align}
\label{cor-donsker-1}
& \mathcal{B}_{\mathbb{R}^p \times \mathcal{D}([0,1])} \left(
\mathcal{L}\left(\rn (\theta - \hat \theta, \Ld(\cdot) - \hat \Ld(\cdot)) \given X \right), \
\mathcal{L} \left( \mathbb{V}, \mathbb{B}(U_0(\cdot)) - \mathbb{V}' \Ld_0\{\gamm\}(\cdot) \right)
\right) \stackrel{P_{\eta_0}}{\to} 0,
\end{align}
where $\mathcal{D}([0,1])$ is the Skorokhod space on $[0,1]$. 
\end{corollary}

{Corollary \ref{cor-donsker} immediately implies the Bernstein-von Mises theorem for the marginal posterior distribution of $\theta$: $\mathcal{B}_{\mathbb{R}^p} \left( \mathcal{L} (\rn (\theta - \hat \theta) \given X ), 
 \mathbb{V} 
\right) \stackrel{P_{\eta_0}}{\to} 0$.}
As an application of Corollary \ref{cor-donsker}, one obtains the Bayesian Donsker theorem for the conditional hazard and survival functions by simply applying the functional delta method \citep[Chapter 20 in][]{vdvAS}. 
Let $\tz$ be a fixed element in $\mathbb{R}^p$, and recall we define $S(\cdot \given \tz) = \exp(-\Ld(\cdot) e^{\theta' \tz})$, the survival function conditional on $\tz$. 
Denote $\hat S(\cdot \given \tz) = \exp(-\hat \Ld(\cdot) e^{\hat \theta' \tz})$ with $\hat\theta$ and $\hat \Ld(\cdot)$  the frequentist estimators as above, then, as $n \to \infty$,
\begin{align}
& \mathcal{B}_{\mathcal{D}([0,1])} 
\left(
	 \mathcal{L}\left(\rn (\Ld(\cdot) e^{\theta'\tz} - \hat \Ld(\cdot) e^{\hat \theta'\tz}) \given X \right), 
	\ \mathcal{L} (\mathbb{H}_1)
\right) 
\stackrel{P_{\eta_0}}{\to} 0,
\label{donsker-cond-hazard}
\\
& \mathcal{B}_{\mathcal{D}([0,1])} 
\left(
	 \mathcal{L}\left(\rn (S(\cdot \given \tz) - \hat S(\cdot \given \tz)) \given X \right), \
	\mathcal{L}(\mathbb{H}_2)
\right) 
\stackrel{P_{\eta_0}}{\to} 0,
\label{donsker-cond-survival}
\end{align}
where $\mathbb{H}_1$ and $\mathbb{H}_2$ are the transformed processes obtained after applying the functional delta method from (\ref{cor-donsker-1}).
Moreover, by applying the continuous mapping theorem and noting that the map for any function $f \to \|f\|_\infty$ is continuous from $\mathcal{D}([0,1])$, equipped with the supremum norm, to $\mathbb{R}^+$, (\ref{donsker-cond-hazard}) and (\ref{donsker-cond-survival}) imply 
\begin{align*}
&	\mathcal{B}_{\mathbb{R}} 
\left(
	 \mathcal{L} \left(\rn \|\Ld(\cdot) e^{\theta'\tz} - \hat \Ld(\cdot) e^{\hat \theta'\tz}\|_\Linfty \given X \right), 
	\ \mathcal{L} \left(\|\mathbb{H}_1\|_\Linfty\right)
\right) 
\stackrel{P_{\eta_0}}{\to} 0, \\
&	\mathcal{B}_{\mathbb{R}} 
\left(
	 \mathcal{L} \left(\rn \|S(\cdot \given \tz) - \hat S(\cdot \given \tz) \|_\Linfty \given X \right), \
	\mathcal{L} \left(\|\mathbb{H}_2\|_\Linfty\right)
\right) 
\stackrel{P_{\eta_0}}{\to} 0.
\end{align*}
A simple consequence of the last display is that the two-sided $(1-\alpha)\%$ quantile credible band for the conditional hazard (resp. the {survival function conditional on $z$}) function is asymptotically a two-sided $(1-\alpha)\%$ confidence band (see \cite{cast14a}, Corollary 2).

\subsection{The supremum-norm convergence rate for {the hazard function conditional on $z$}}
\label{sec:sup-norm}

In this section, we present the third main result: a faster supremum-norm posterior contraction rate for {the hazard function conditional on $z$} than the rate $\zeta_n$ in \ref{cond:P}. We denote this rate as 
$\xi_n$, which depends on $L_n$, a diverging sequence, such that $L_n 2^{L_n} \lesssim \sqrt{n}$,  where
$$
\xi_n: = \xi_n(\beta, L_n, \epsilon_n) = \sqrt{\frac{L_n 2^{L_n}}{n}} + 2^{-\beta L_n} + \epsilon_n.
$$

\begin{theorem}
\label{thm:supnorm-rate}
Suppose $r_0 \in \mathcal{H}(\beta, D)$ with $1/2 < \beta \leq 1$. 
Let $L_n$ be a diverging sequence such that $L_n 2^{L_n} \lesssim \sqrt{n}$ and let $\tz \in \mathbb{R}^p$ be fixed.
For the prior of $\eta = (\theta, \ld)$ chosen such that both {\rm \ref{cond:P}} and (\ref{ez-rate-condition}) hold, and {\rm \ref{cond:C2}} also holds uniformly for $a = z$ and any $b = \psi_{LK}$, with $(\psi_{LK})$ the Haar wavelet basis, $0 \leq L \leq L_n$ and $0 \leq K < 2^L$,  
then for $\xi_n = o(1)$ and $n \xi_n^2 \to \infty$, and an arbitrary sequence $M_n \to \infty$, 
\begin{align}
\label{supnorm-rate}
\Pi \left(\eta: \|\ld e^{\theta'\tz} - \ld_0 e^{\theta_0'\tz} \|_\Linfty > M_n \xi_n \given X \right) = o_{P_{\eta_0}}(1).
\end{align}
\end{theorem}

We will show that in the next section, with a specific choice of the value of $L_n$, the rate $\xi_n$
is within the same order of the Hellinger rate ${{\nu}}_n$ in \eqref{rateveps}. 
 
\subsection{Results for specific priors}
\label{sec:application}

In this section, we apply the generic results in Sections \ref{sec:jointBvM}, \ref{sec:npBvM}, and \ref{sec:sup-norm} to study the specific priors considered in Section \ref{sec:priors}.
The result is stated in the following theorem.
\begin{theorem}
\label{thm:specific-sup-rate}
Consider the Cox model with the priors as specified in {\rm \ref{prior:T}, and \ref{prior:H}} or {\rm \ref{prior:W}} with $L = L_n$ in (\ref{choicel}) and ${{\nu}}_n$ given in \eqref{rateveps}. For any $\beta \in (1/2, 1]$,
\begin{enumerate}
\item conditions {\rm \ref{cond:P}} and {\rm \ref{cond:C1}} hold, {\rm \ref{cond:B}} holds for any $b \in \mathcal{H}(\mu, D)$ with $\mu > 1/2$ and $D > 0$, then, (\ref{thm1-1}) in Theorem \ref{thm:bvm-linear-fns} holds;
\item condition {\rm \ref{cond:C2}} also holds, thus (\ref{eqn:donsker-joint}) in Theorem \ref{thm:donsker-joint} holds;
\item the supremum-norm rate $\xi_n$ in \eqref{supnorm-rate} can be taken to be $\xi_n = {{\nu}}_n$ as in \eqref{rateveps}. 
\end{enumerate}
\end{theorem}

The proof of this result, implying that conditions \ref{cond:P}, \ref{cond:B}, \ref{cond:C1}, and \ref{cond:C2} hold for priors given in Section \ref{sec:priors}, can be found in \citet{ning22supp}.

\begin{remark}
If $\beta > 1$, the first and second points in Theorem \ref{thm:specific-sup-rate} still hold. The third point also holds but the supremum-norm rate becomes $(\log n/n)^{1/3}$.
\end{remark}

Let us compare the supremum rate ${{\nu}}_n$ in the third point of the theorem and the rate $\zeta_n$ obtained in Lemma \ref{lemma1}. Obviously, ${{\nu}}_n < \zeta_n$, as $\zeta_n \geq 2^{L_n/2} {{\nu}}_n$, and $L_n \to \infty$ as $n\to \infty$. In fact, by plugging-in the value of $L_n$ in (\ref{choicel}), 
one obtains $\zeta_n = (\log n/n)^{\frac{2\beta - 1}{2(2\beta+1)}}$ which can become extremely slow when $\beta$ is close to $1/2$. In Lemma \ref{lower-bound} in \citet{ning22supp}, we derive a lower bound for the minimax rate in the supremum norm for the hazard which shows that the rate ${{\nu}}_n$ is sharp. 
To our best knowledge, this is the first sharp supremum-norm result for the hazard obtained for the Cox model. 

The cut-off $L_n$ in our theorems is chosen to be a deterministic sequence depending on $n$ and the smoothness level $\beta$. As noted below \eqref{choicel}, for semiparametric-type results, including Donsker theorems, it is enough to `undersmooth', and all such results hold for a smoothness parameter taken to be $1/2$ in \eqref{choicel} whenever the true smoothness $\beta$ is larger than 1/2, and this choice  already provides a contraction rate of $n^{-1/4}$ for the posterior of the conditional hazard. It is natural to ask whether the cut-off parameter $L$ can itself be taken random in a hierarchical Bayes approach. Although often used in practice too, we underline that particular caution must be taken with such an `adaptive' prior: indeed, as demonstrated in \citet{cast15BvM} (Section 4.3) in the density estimation model, BvM results may fail to hold for such a prior. This phenomenon would appear in the Cox model too if the regularities of the hazard and of the least favorable direction are too far apart. Regarding adaptive supremum norm rate (or nonparametric BvM) results, it is conceivable that spike-and-slab type priors would work, as in \citet{ray17}, although unlike in the white noise setting considered in \citet{ray17}, one could not use conjugacy here, so this is beyond of the scope of this paper and left for future investigation.

\section{Simulation studies}
\label{sec:sim}

Two simulation studies are conducted in this section. 
The first study, described in Section \ref{study-I}, compares the limiting distribution given in Corollary \ref{cor-1} with the empirical distributions obtained from the MCMC algorithm, which is given in Section \ref{generate-data}.
The second study compares the coverage and the area of the 95\% credible bands for the MCMC algorithm to the 95\% confidence bands for a commonly used frequentist method by varying the sample size and changing the censoring distribution.
{We choose the two random histogram priors for $\ld$ as given in Section \ref{sec:priors}. The prior for $\theta$ is chosen as the standard normal distribution. If $\theta$ is multivariate, we use the standard multivariate normal density instead.} 
In Section \ref{generate-data}, we describe how we generate the simulated data and the MCMC sampler. Section \ref{study-I} presents results for the first study, and Section \ref{study-II} summarizes results for the second study.

\subsection{Generating the data and the MCMC sampler}
\label{generate-data}

The data are generated from the ``true'' conditional hazard function $\ld_0(t) e^{\theta_0'z}$, {where $\ld_0$ and $\theta_0$ will be specified below}.
%\sout{which is chosen as follows: let $\theta_0 = -0.5$, a scalar, and generate the covariate $z$ randomly from the standard normal distribution. Two different baseline hazard functions are used to generate the data:}
%\begin{enumerate}[label=\textbf{(1)}]
%\item \label{chzf} \sout{$\lambda_0(t) = 6 \times ((t + 0.05)^3 - 2(t + 0.05)^2 + t + 0.05) + 0.7, \ t \in [0,1]$,}
%\end{enumerate}
%\begin{enumerate}[label=\textbf{(2)}]
%\item \label{pwhzf} \sout{$\lambda_0(t) = 3 \times \mathbbm{1}_{[0,0.4)}(t) + 1.5 \times \mathbbm{1}_{[0.4, 0.6)}(t) + 2 \times \mathbbm{1}_{[0.6, 1]}(t)$,}
%\end{enumerate}
%\sout{The first one is a smooth function and the second one is piecewise constant. 
%Plots of the two functions are given in Figures \ref{fig:pwhzf} and \ref{fig:chzf} respectively. These numerical choices are for illustration purposes and otherwise fairly arbitrary. Similar simulation results would hold if choosing other either smoothly varying or piecewise constant functions (we note once again that, although our theoretical results assume H\"older smoothness of the true log-hazard, the techniques go through for histogram true hazards as well).}
The ``observations'' $X^n = (Y^n, \delta^n)$ are generated using the ``{\it simsurv}'' function in $\mathsf{R}$ \citep{brilleman20}. 
We consider the following two types of censoring: 
\begin{enumerate}
\item {\it Administrative censoring only}. Time points are censored at a fixed time point $t = 1$;
\item {\it Administrative censoring $+$ uniform censoring}. The censoring time is generated from the uniform distribution on $[0,1]$. Any time point beyond $t = 1$ is also censored. 
\end{enumerate}

Although the first type of censoring violates our assumption in Section \ref{sec:assumption}, as we assumed the censoring follows a random distribution, it is interesting to find out in the next two sections that the empirical results still match with our theoretical results quite well.

Posterior draws are obtained using the MCMC algorithm given as follows:
\begin{enumerate}
\item For the {\it independent gamma prior,} since it is conjugate with the posterior distribution given $\theta$, we sample each $\ld_k \sim \text{Gamma}(d_k + \alpha, T_k(\theta) + \beta)$,
where $d_k = \sum_{i=1}^n \delta_{ki}$ is the number of events in $k$-th interval and 
$T_k(\theta) = \sum_{i=1}^n Y_{ik} e^{\theta Z_i}$, $\alpha$ and $\beta$ are the hyperparameters, and we chose them to 1.
After obtaining samples for $(\ld_{lk})$, we draw $\theta$. Since $\pi(\theta)$ is not conjugate, we first draw a candidate from the proposal density, i.e., $\theta^{\text{prop}} \sim N(\theta^{\text{prev}}, 1)$, where $\theta^{\text{prev}}$ stands for the draw from the previous iteration, and then use the Metropolis algorithm to accept or reject this candidate.

\item For the {\it dependent gamma prior}, as it is non-conjugate, we thus draw each $\lambda_k$ from the proposal density as follows: $\ld_1^{\text{prop}} \sim \text{Gamma}(d_1 + \alpha_0 - \alpha, \ T_1(\theta) + \beta_0)$ and
$\ld_k^{\text{prop}} \sim \text{Gamma}(d_k + \varepsilon, \ \alpha/\ld_{k-1} + T_k(\theta))$ for $k = 2, \dots, L_n-1$. The last interval $\ld_K \sim \text{Gamma}(d_K + \alpha,\  \alpha\ld_{K-1} + T_K(\theta))$ for $K=L_n$. In practice, we choose $\varepsilon = 10^{-6}$, $\alpha_0 = 1.5$ and $\alpha = \beta_0 = 1$. 
The proposal density for $\theta$ is the same.
\end{enumerate}

To initialize the MCMC algorithm, we choose the initial values for $\theta$ and $\ld_k$ as their frequentist estimators (the same as in Corollary \ref{cor-donsker}). We choose $L_n$ as in (\ref{choicel}) and $\beta = 1/2$.
For each simulation, we run 10,000 iterations and discard the first 2,000 draws as burn-in. 

Let us now discuss in more detail the simulation in Figure \ref{fig:intro}. 
%\sout{{The dataset is} generated by choosing \ref{chzf}.}
{The dataset is generated by choosing $\theta_0 = -0.5$ and $\lambda_0(t) = 6 \times ((t + 0.05)^3 - 2(t + 0.05)^2 + t + 0.05) + 0.7, \ t \in [0,1]$. We generate the covariate $z$ randomly from the standard normal distribution.
Here, the true function $\ld_0$ is chosen the same as it in the simulation of \citet{cast20survival}; see also in their $\mathsf{R}$ package `BayesSurvival' \citep{vdp21R}. However, 
this choice is for illustration purposes and otherwise fairly arbitrary. Similar simulation results would hold if choosing other either smoothly varying or piecewise constant functions (see Section \ref{study-II}, which we chose different $\theta_0$ and $\ld_0$). We note once again that, although our theoretical results assume H\"older smoothness of the true log-hazard, the techniques go through for histogram true hazards as well).
In Figure \ref{fig:intro},}
only administrative censoring is considered. 
The prior is chosen to be the independent gamma prior (choosing the dependent gamma prior won't change the result dramatically, as can be seen in Table \ref{tab:coverage-1} below).
The 95\% credible band is a fixed width band whose width is constant with the time. 
The width is determined such that the posterior probability is 95\%.
The 95\% confidence band, on the other hand, is obtained using the $\mathsf{predictCox}$ function of the `riskRegression' package in $\mathsf{R}$ \citep{gerds21}. Its width varies with time. 

Here we briefly describe the approach used in their package. 
We refer the interested readers to read \citet{lin94} and \citet{scheike08} for more details. 
Using the fact that the frequentist estimator $(\hat \theta, \hat \Ld)$ converges to the same limiting process as the joint distribution $(W_n^{(1)}, \int_0^t W_n^{(2)}(\mathbbm{1}_{u \leq t}) du)$, where 
$W_n^{(1)}$ and $W_n^{(2)}$ are given in \eqref{Wn-1} and \eqref{Wn-2}, their approach first defines another process, depending on $W_n^{(1)}$ and $W_n^{(2)}$, that is asymptotically equivalent to $\rn \left(\hat \Ld(t)e^{\hat\theta'z} - \Ld(t) e^{\theta' z}\right)$; see equation 2.1 in \citet{lin94} for the exact expression of that process.
Their approach then further approximates that process by a summation of independent normal variables whose distribution can be easily generated through Monte Carlo simulation and replaces other unknown quantities with their sample estimators.
After large enough samples are generated, the last step is to obtain the size of the 95\% confidence band for the conditional cumulative hazard function by choosing the 95th quantile from those samples. 
The confidence band for the {survival function conditional on $z$} can be obtained similarly, except that one first needs to apply the functional delta method to obtain the limiting process for the {survival function conditional on $z$}. The remaining parts are the same. 
 
\subsection{Study I: Comparing the empirical posterior distributions and the limiting distributions of $\theta$ and $\Ld(1)$}
\label{study-I}

We compare the limiting distribution given in Corollary \ref{cor-1} with the empirical distribution obtained using the MCMC sampler in Section \ref{generate-data}.
{We will study the joint posterior distribution of $\theta$ and $\ld$. For simplicity, we let $\theta_0 = -0.5$ (for now, we simply choose it to be univariate; simulation results for using a multivariate $\theta$ are given in the next section) and generate the covariate $z$ randomly from the standard normal distribution. We also choose $\lambda_0(t) = 6 \times ((t + 0.05)^3 - 2(t + 0.05)^2 + t + 0.05) + 0.7, \ t \in [0,1]$ and generate the hazard rate with a sample of 1,000.}
{We choose the prior for $\theta$ as the standard normal distribution and the prior for $\ld$ as the independent gamma prior.} Results for choosing the dependent gamma prior are similar. To obtain draws, we run the MCMC algorithms in parallel for 1,000 times. 
Each time we only record the last pair draw for $\theta$ and $\Ld(1)$, where $\Ld(1) = \sum_{k=1}^{L_n} \ld_k$.  Therefore, the 1,000 draws we obtained are independent.

\begin{figure}
\centering
\begin{minipage}{0.45\textwidth}
\centering
\includegraphics[width = 1\textwidth]{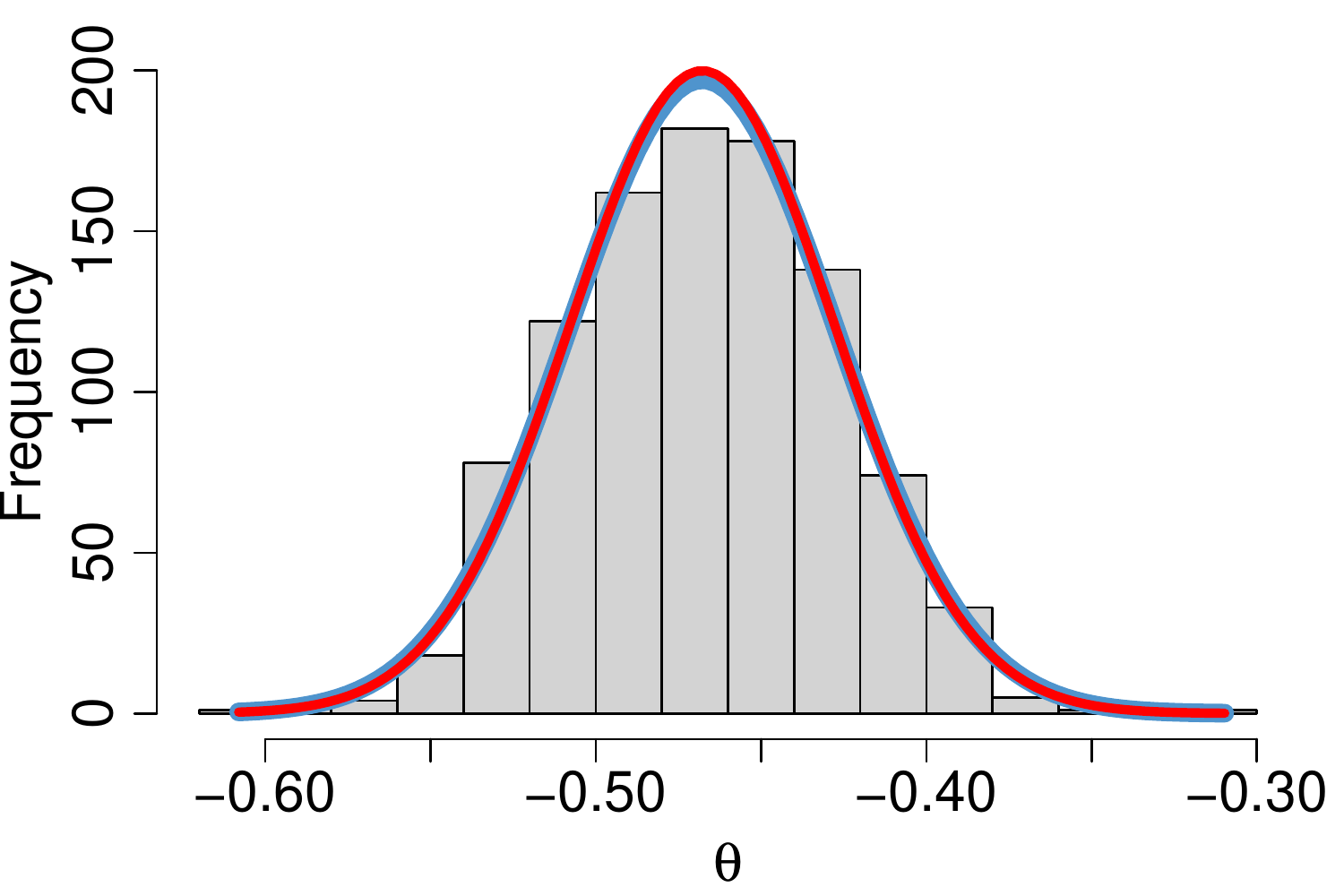}
\caption{Plot of the empirical histogram, the empirical distribution (blue), and the limiting distribution (red) for the marginal posterior distribution of $\theta$. True value of $\theta$ is $-0.5$.}
\label{fig:marginal-theta}
\end{minipage}%
\hspace{0.04\textwidth}
\begin{minipage}{0.45\textwidth}
\includegraphics[width=1\textwidth]{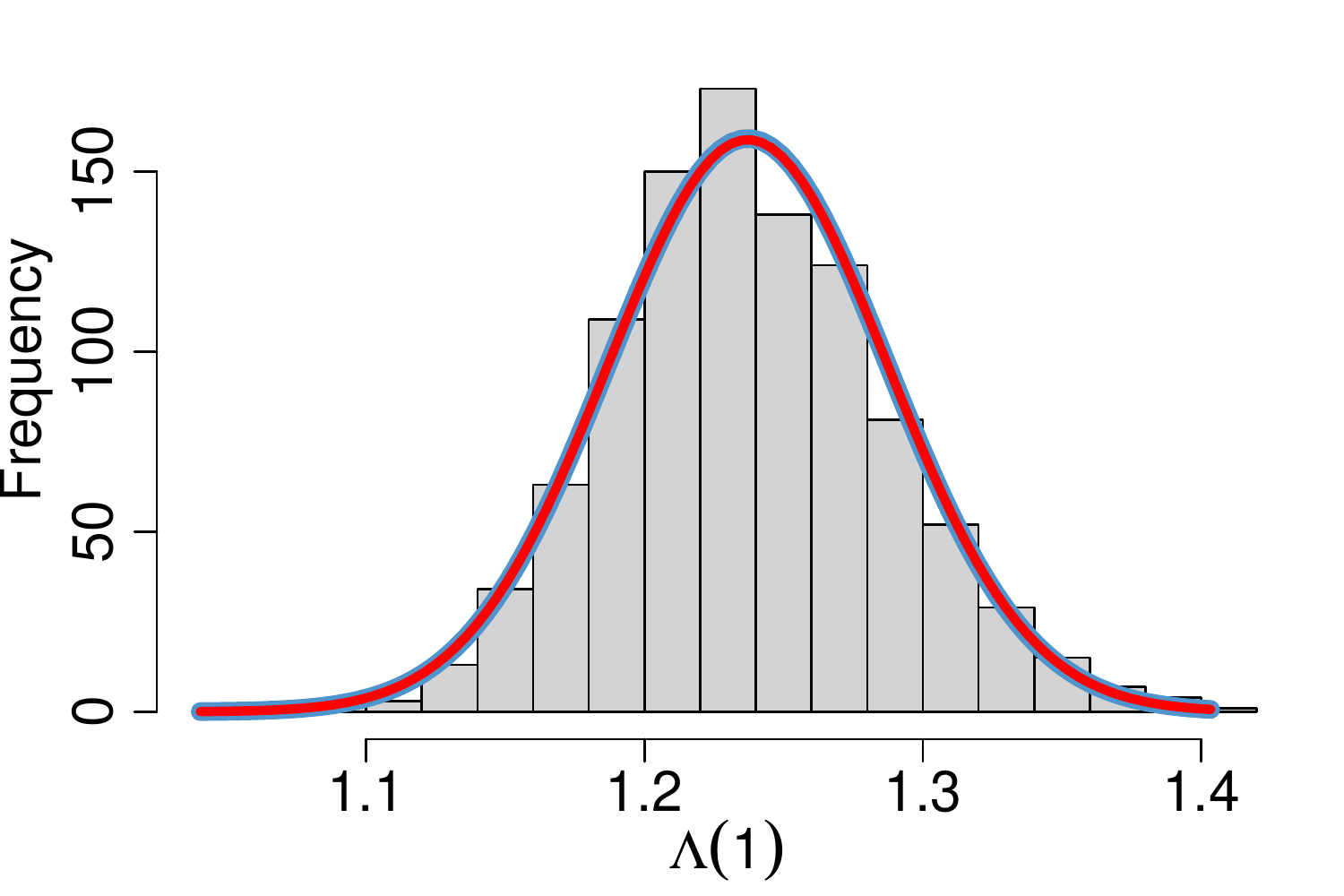}
\caption{Plot of the empirical histogram, the empirical distribution (blue), and the limiting distribution (red) for the marginal posterior distribution of $\Ld(1)$. True value of $\Ld(1)$ is $1.2$.}
\label{fig:marginal-Lambda}
\end{minipage}
\end{figure}

We first study the marginal posterior distributions for $\theta$ and $\Ld(1)$. In 
Figures \ref{fig:marginal-theta} and \ref{fig:marginal-Lambda}, 
we first draw their empirical histogram from the 1,000 independent draws. 
We then draw a normal density with blue color centered at the posterior mean, and its variance is estimated from those draws. 
Last, we draw another normal density with red color, which has the same centering as the blue one, but its variance is chosen as the theoretical value from the limiting distribution in Corollary \ref{cor-1}. 
We observe that in either the left plot (i.e., for $\theta$) or the right plot (i.e., for $\Ld(1)$), 
the density with blue color is well aligned with the one with red color. This finding suggests that empirical variances are close to their theoretical variances obtained from the corollary.
We also found that both empirical histograms show similar shapes as their corresponding normal density, which verifies their limiting distributions should be normal. 
Last, both the true values of $\theta_0$ and $\Ld_0(1)$, $-0.5$ and $1.2$ respectively, are contained inside the corresponding 95\% credible intervals.
For $\theta$, the interval is $[-0.54, -0.39]$, and for $\Ld(1)$, it is $[1.14, 1.34]$.

\begin{figure}[!]
\centering
\subfigure[ ]{\includegraphics[width=0.34\textwidth]{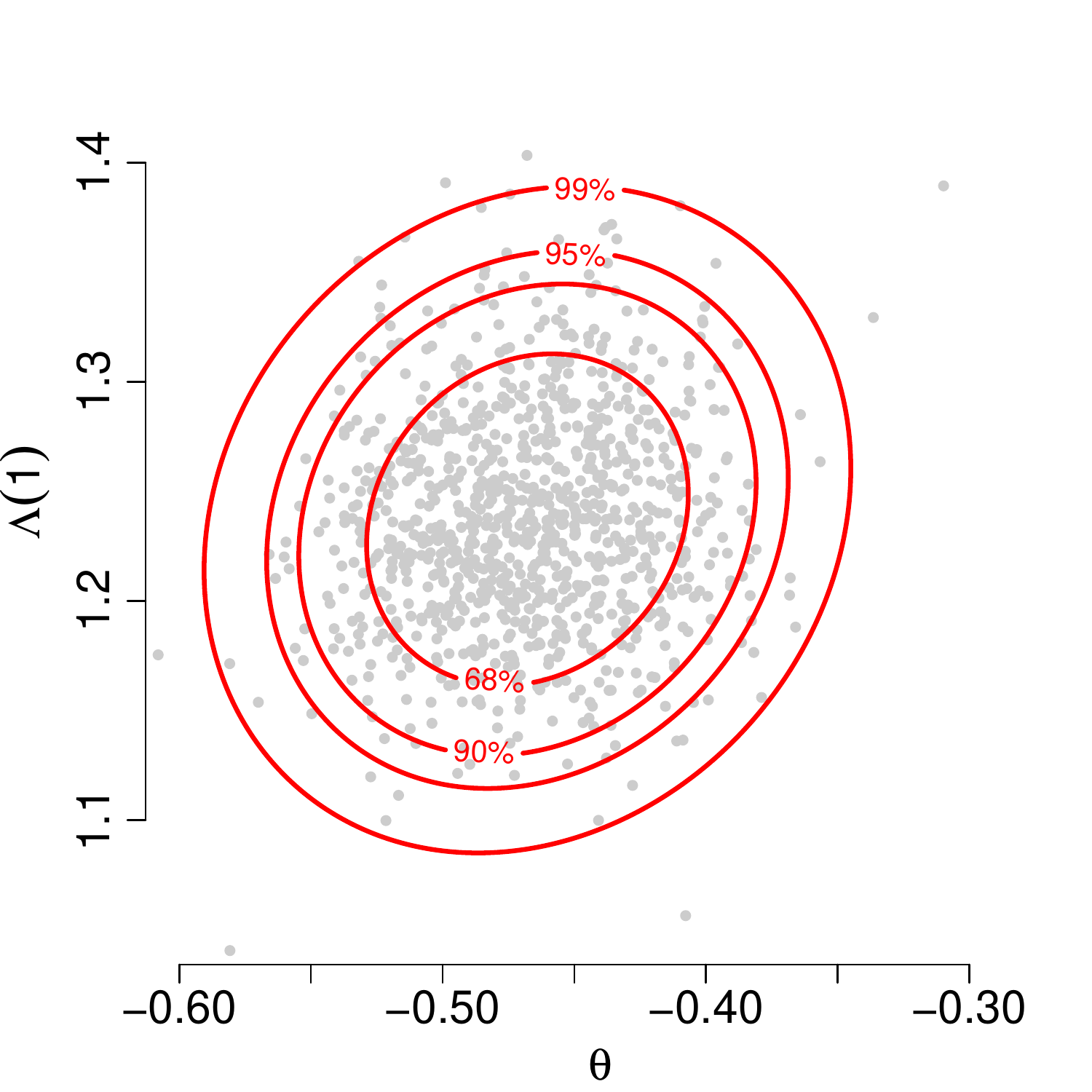}}%
\subfigure[ ]{\includegraphics[width=0.34\textwidth]{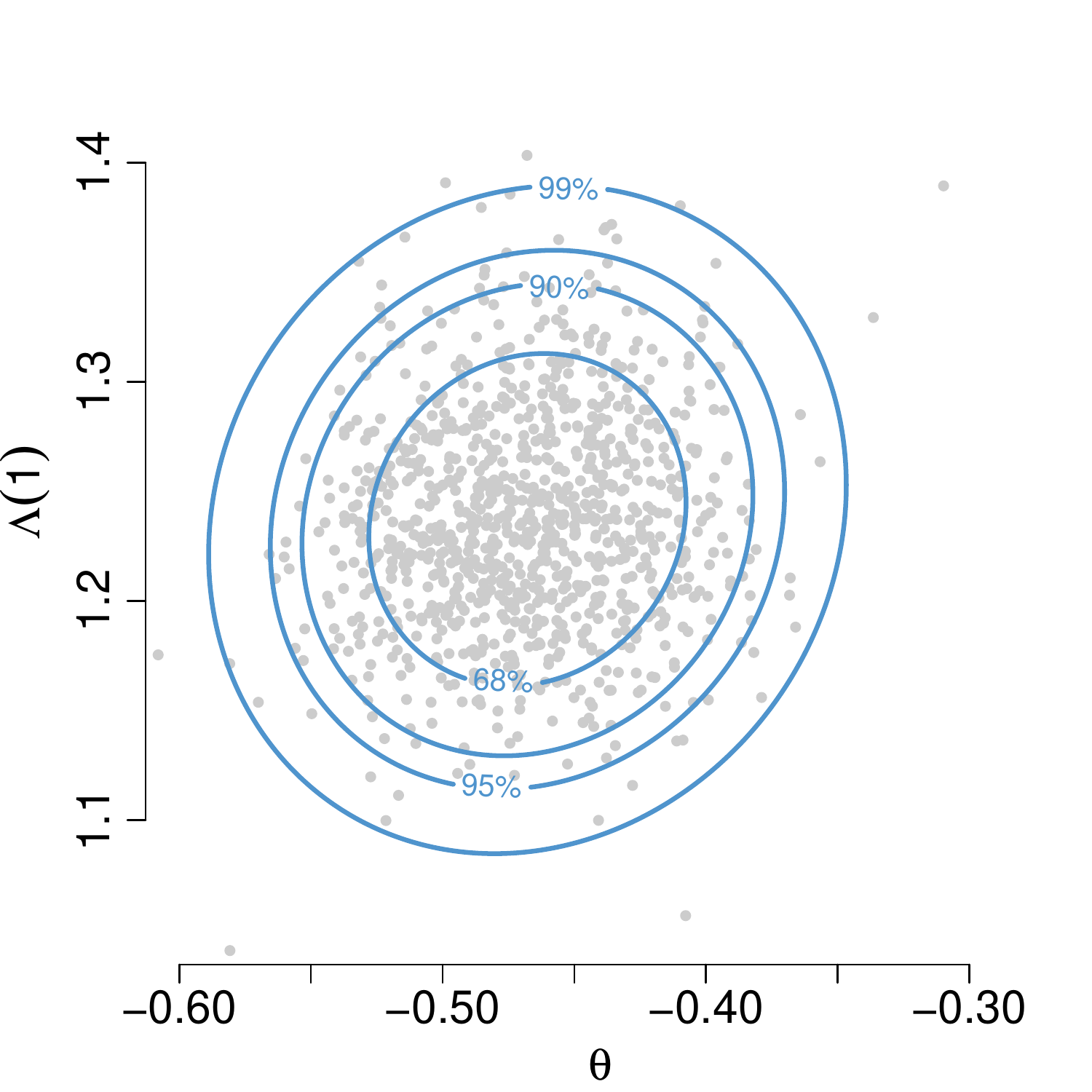}}%
\subfigure[ ]{\includegraphics[width=0.34\textwidth]{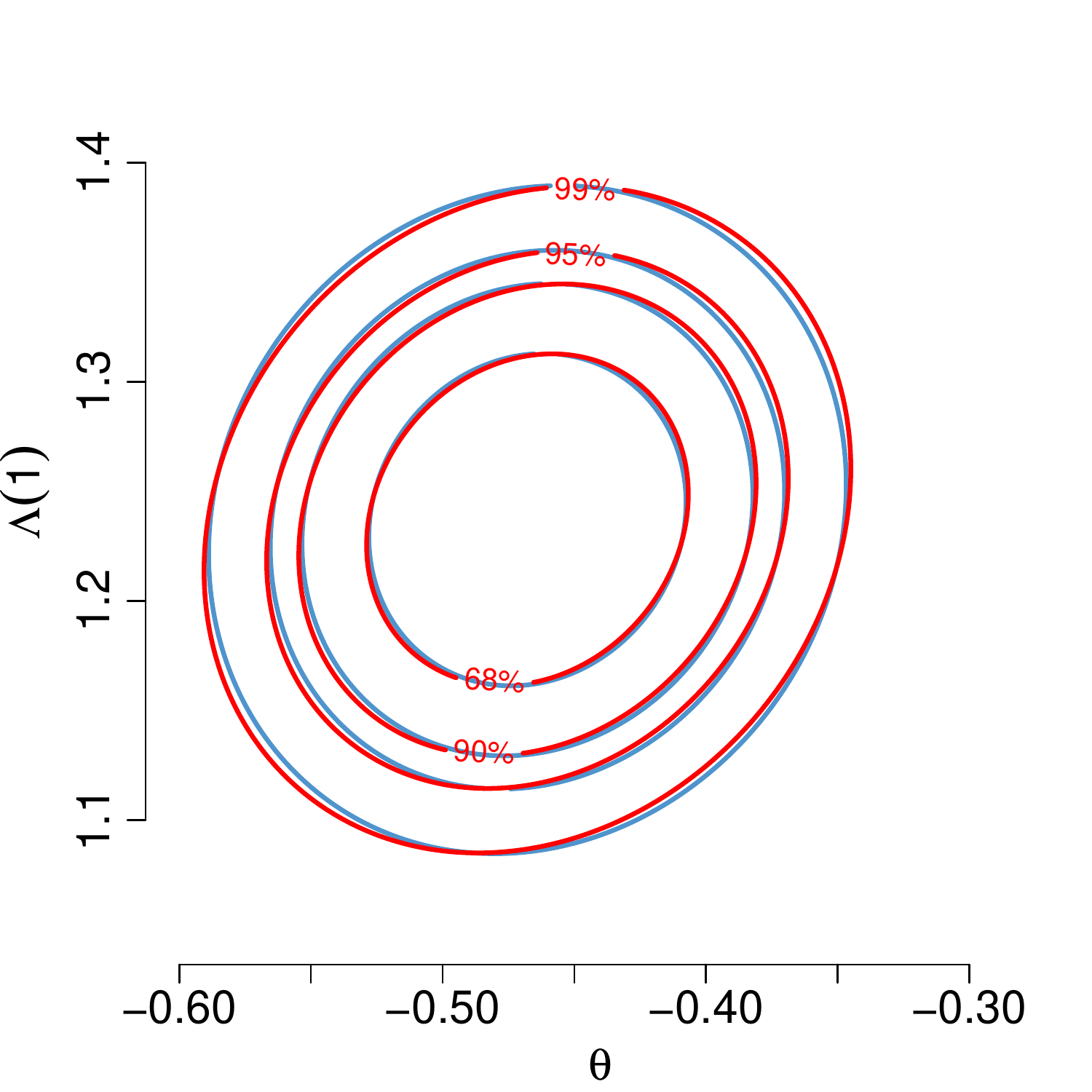}}
\caption{Contour plots of the elliptical credible sets at the 68\%, 90\%, 95\%, and 99\% quantiles. (a) is obtained using the joint limiting distribution given in Corollary \ref{cor-1}, (b) is obtained using the 1,000 independent draws of the pair $(\theta, \Ld(1))$ from the MCMC output. In (c), we overlay the credible sets in (a) and (b). The 1,000 draws are plotted with gray color.}
\label{fig:joint-density}
\end{figure}

Next, we study the joint posterior distribution of $\theta$ and $\Ld(1)$, which involves the correlation between the two quantities.
In Figure \ref{fig:joint-density}, we give three plots. 
In (a), we plot the 86\%, 90\%, 95\%, and 99\% contour plots of the limiting joint distribution in Corollary \ref{cor-1}. In (b), we plot the contours with the same four quantiles for a bivariate normal distribution, which its mean, variances, and correlations are estimated from the 1,000 draws. 
In (c), we found that the two sets of contour plots in (a) and (b) indeed align quite well, 
which suggests that the empirical distribution matches with the theoretical limiting distribution in the corollary. 
Our calculation reveals that the correlation between $\theta$ and $\Ld(1)$ in (a) is 0.15 and that in (b) is 0.10. 
A benefit of studying the joint posterior distribution with the correlation is that one can obtain the elliptical credible sets instead of rectangular credible sets. The length and the width of the rectangular credible sets are the 97.5\% credible intervals of $\theta$ and $\Ld(1)$ respectively. 
Therefore, the area of a rectangular credible set is typically larger than that of an elliptical credible set. For example, in (b), the area of the 95\% elliptical credible set is 1.07 and that of the 95\% rectangular credible set is 1.76, which is 64\% bigger.

\subsection{Study II: Comparing the coverage and the area between the credible bands and the confidence bands}
\label{study-II}

The study in the last section is based on a single simulated dataset. We now provide a more thorough study to compare the coverage and the area of the credible (or confidence) bands under various settings. Specially, we want to compare: 1) the two Bayesian methods using the independent gamma  and the dependent gamma priors {(while the prior for $\theta$ is chosen to be the standard normal distribution)}; 2) datasets with two different sample sizes $n = 200$ and $n = 1,000$; 3) data with two different types of censoring: administrative censoring only and administrative censoring with additional uniform censoring; 4) coverages of the baseline survival function and {survival function conditional on $z$} and 5) data are generated from the continuous function in \ref{chzf} and that from the piecewise constant function in \ref{pwhzf}.

{
Two different baseline hazard functions are used to generate the data:}
\begin{enumerate}[label=\textbf{(1)}]
\item \label{chzf} $\lambda_0(t) = 0.8\times \sin(2\pi (t+0.05))+ (t + 0.05)^4 - 1.8 \times (t + 0.05)^2 + 2, \ t \in [0,1]$,
\end{enumerate}
\begin{enumerate}[label=\textbf{(2)}]
\item \label{pwhzf} $\lambda_0(t) = 3 \times \mathbbm{1}_{[0,0.4)}(t) + 1.5 \times \mathbbm{1}_{[0.4, 0.6)}(t) + 2 \times \mathbbm{1}_{[0.6, 1]}(t)$,
\end{enumerate}

\begin{figure}[!h]
\centering
\begin{minipage}{0.45\textwidth}
\centering
\includegraphics[width=1\textwidth]{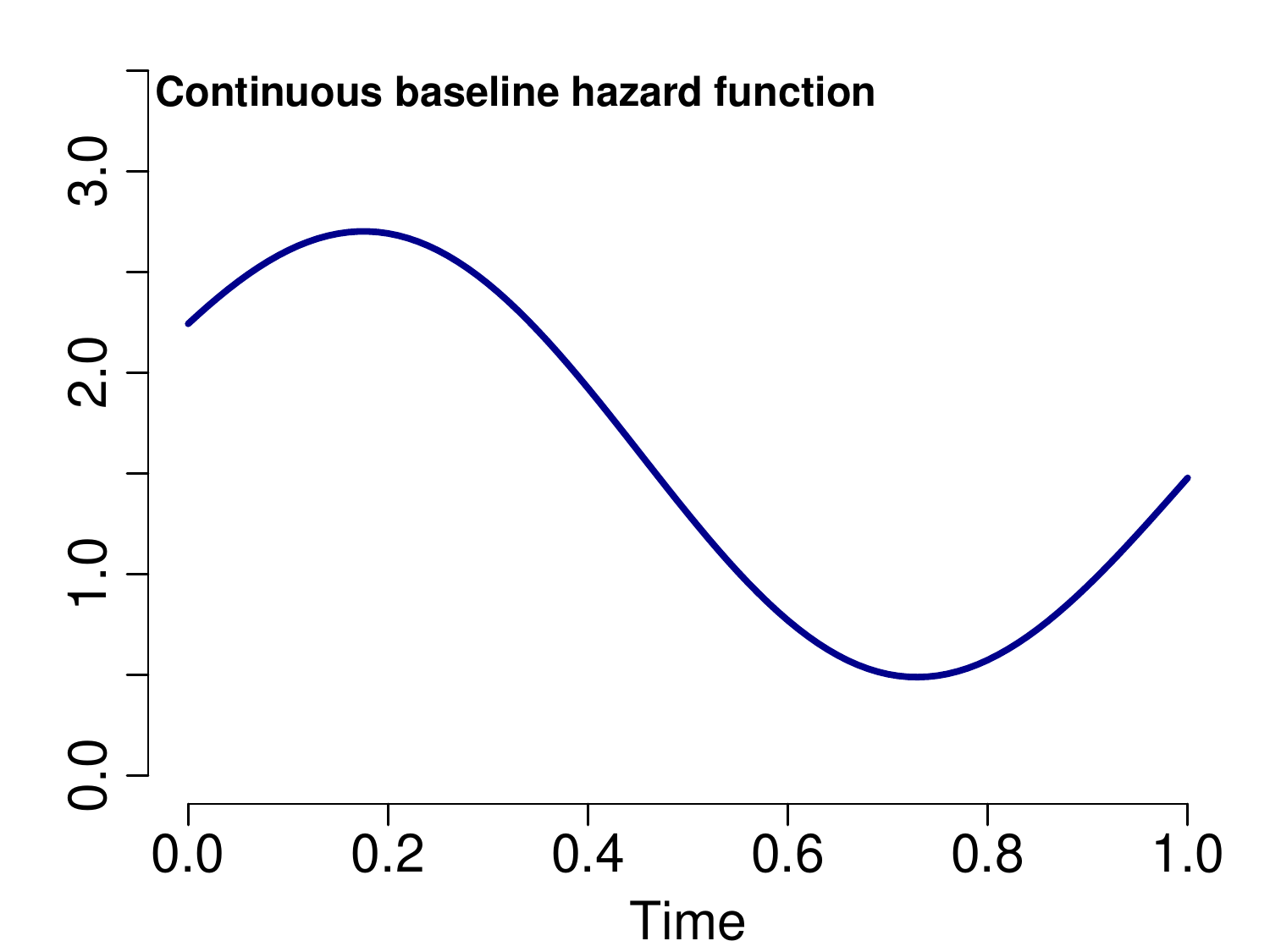}
\caption{Plot of the continuous baseline hazard function in {\rm \ref{chzf}}.}
\label{fig:chzf}
\end{minipage}%
\hspace{0.04\textwidth}
\begin{minipage}{0.45\textwidth}
\includegraphics[width = 1\textwidth]{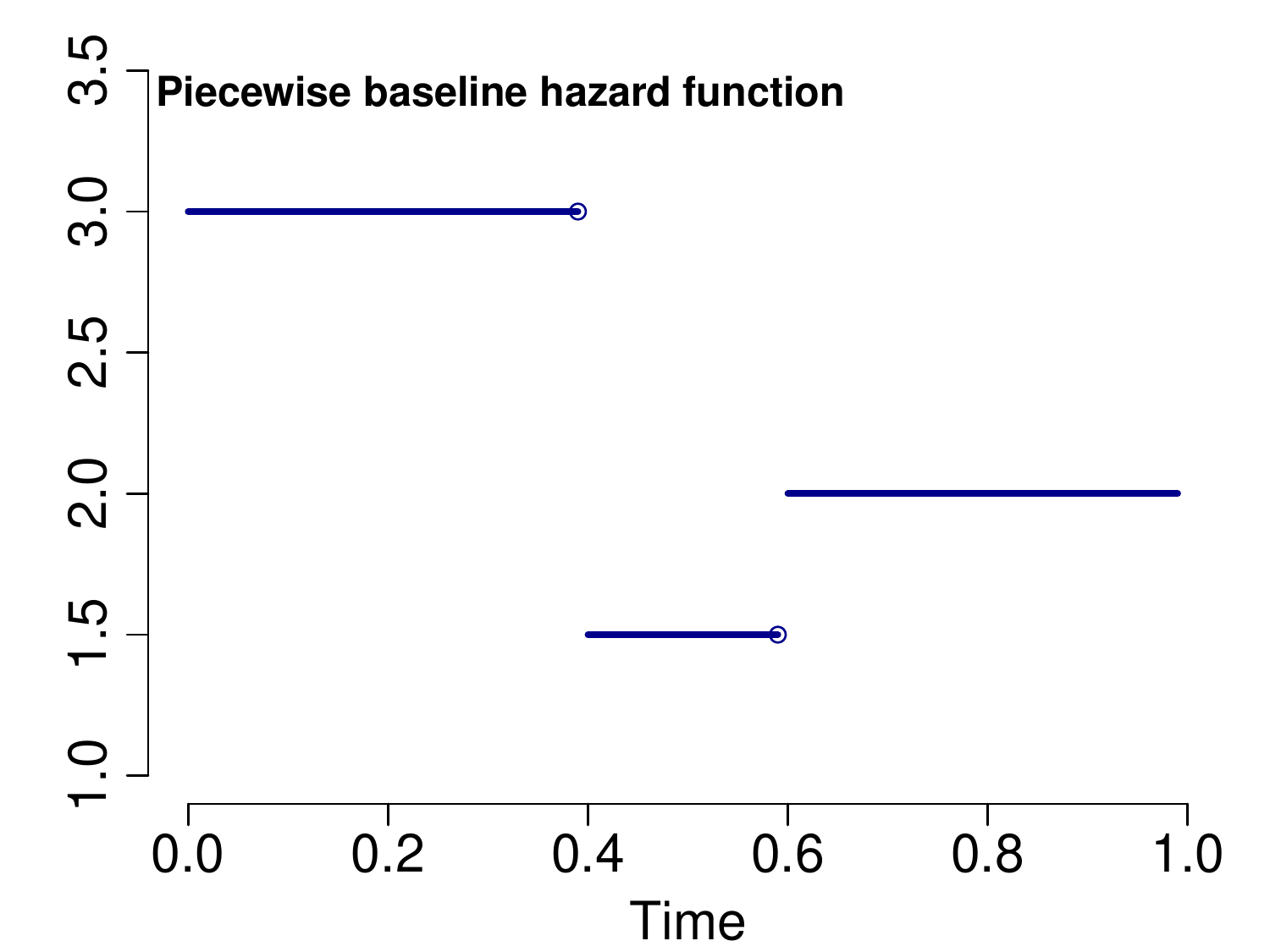}
\caption{Plot of the piecewise constant baseline hazard function in {\rm \ref{pwhzf}}.}
\label{fig:pwhzf}
\end{minipage}
\end{figure}

{The first one is a smooth function and the second one is piecewise constant. 
Plots of the two functions are given in Figures \ref{fig:chzf} and \ref{fig:pwhzf} respectively. These numerical choices are for illustration purposes and otherwise fairly arbitrary. Similar simulation results would hold if choosing other either smoothly varying or piecewise constant functions (we note once again that, although our theoretical results assume H\"older smoothness of the true log-hazard, the techniques go through for histogram true hazards as well).}

For each setting, we generate 1,000 datasets. 
For each dataset, we run the MCMC sampler to obtain the 95\% credible (or confidence) band. 
The coverage is the percentage of the credible (or confidence) bands encompassing the true function. The area estimated by taking the average of 1,000 areas of the credible (or confidence) bands. Results using the continuous baseline hazard function are given in Table \ref{tab:coverage-1} and those using the piecewise constant baseline hazard function are given in Table \ref{tab:coverage-2}.

\begin{table}[ht]
\caption{{Coverages and areas of the 95\% Bayesian credible bands and of the 95\% confidence bands with $\theta_0 \in \mathbbm{R}^5$ and $\ld_0$ is chosen as the continuous function in {\rm \ref{chzf}}.}}
\centering 
\begin{tabular}{ccccccccc} 
\hline
 & \multicolumn{4}{c}{\bf Adm. censoring only} & \multicolumn{4}{c}{\bf Adm. $+$ Unif. censoring}\\
\hline
& \multicolumn{2}{c}{baseline survival} & \multicolumn{2}{c}{cond. survival} 
& \multicolumn{2}{c}{baseline survival} & \multicolumn{2}{c}{cond. survival}\\
 & coverage & area & coverage & area & coverage & area & coverage & area\\[2pt]
\hline
ind. 200 & 0.96 & 0.16 & 0.94 & 0.16 & 0.97 & 0.20 & 0.96 & 0.21\\
dep. 200 & 0.93 & 0.16 & 0.93 & 0.16 & 0.92 & 0.19 & 0.91 & 0.19 \\
freq. 200 & 0.93 & 0.18 & 0.93 & 0.18 & 0.92 & 0.22 & 0.91 & 0.22 \\[2pt]
\hline 
ind. 1000 & 0.95 & 0.08 & 0.93 & 0.08 & 0.98 & 0.10 & 0.96 & 0.10 \\
dep. 1000 & 0.93 & 0.08  & 0.92 & 0.08 &  0.92 & 0.09 & 0.92 & 0.09 \\
freq. 1000 & 0.94 & 0.08 & 0.93 & 0.08 & 0.94 & 0.10 & 0.94 & 0.10 \\[2pt]
\hline
\end{tabular}
\label{tab:coverage-1}
\end{table}

From Table \ref{tab:coverage-1}, in which $\ld_0$ is chosen to be a continuous function, first, we found that the two Bayesian methods, either using the independent or the dependent gamma prior, produce similar converge results and areas for the credible bands. 
Second, when $n = 200$, 
%\sout{not only the coverage of the two Bayesian methods is better than that of the frequentist method, but the area is also smaller.} 
{the two Bayesian methods yield a smaller area than the frequentist method. } 
The coverages and the areas of the three methods become similar when $n = 1,000$.
As approximation becomes more accurate when the sample size increases, 
%\sout{this leads to higher coverage and a smaller area of the confidence band.}
{it is well expected that both the Bayesian approach and the frequentist approach give comparable coverage and area of the confidence band.} 
Third, we found that when data are administratively and uniformly censored, the area of the credible bands is larger than those only administratively censored. 
Such a result is expected, as we found that in a typical simulated dataset, a former has %{\sout{$\sim$45\%} 
{$\sim$ 40\%} data are censored, and the latter 
%\sout{only} 
has 
%\sout{$\sim$10\%} 
{$\sim$20\%} data are censored.  
Last, there is no significant difference between the coverage and the area of the baseline survival function and the {survival function conditional on $z$}, even though the latter accounts for the uncertainty for estimating the regression coefficients. 

\begin{table}[ht]
\caption{Coverages and areas of the 95\% Bayesian credible bands and of the 95\% confidence bands with $\theta_0 \in \mathbbm{R}^5$ and $\ld_0$ is chosen as the piecewise constant function in {\rm \ref{pwhzf}}.}

\centering 
\begin{tabular}{ccccccccc} 
\hline
 & \multicolumn{4}{c}{\bf Adm. censoring only} & \multicolumn{4}{c}{\bf Adm. $+$ Unif. censoring}\\
\hline
& \multicolumn{2}{c}{baseline survival} & \multicolumn{2}{c}{cond. survival} 
& \multicolumn{2}{c}{baseline survival} & \multicolumn{2}{c}{cond. survival}\\
 & coverage & area & coverage & area & coverage & area & coverage & area\\[2pt]
\hline
ind. 200 & 0.94 & 0.15 & 0.93 & 0.15 & 0.96 & 0.08 & 0.95 & 0.18    \\
dep. 200 & 0.95 & 0.15 & 0.94 & 0.15 & 0.93 & 0.18 & 0.93 & 0.18  \\
freq. 200 & 0.94 & 0.17 & 0.94 & 0.17 & 0.90 & 0.21 & 0.90 & 0.21 \\[2pt]
\hline 
ind. 1000 & 0.95 & 0.09 & 0.93 & 0.07 & 0.95 & 0.09 & 0.94 & 0.09 \\
dep. 1000 & 0.93 & 0.07 & 0.93 & 0.07 & 0.93 & 0.09 & 0.92 & 0.09 \\
freq. 1000 & 0.93 & 0.07 & 0.93 & 0.08 & 0.92 & 0.10 & 0.92 & 0.10 \\[2pt]
\hline
\end{tabular}
\label{tab:coverage-2} 
\end{table}

Table \ref{tab:coverage-2} gives the results for data are generated from the baseline hazard that is the piecewise constant function in \ref{pwhzf}. 
%\sout{We notice that the coverage of the frequentist confidence bands improved significantly when $n = 200$ compared to that in Table \ref{tab:coverage-1}. This suggests that the frequentist approach is more accurate for estimating a piecewise constant function. Indeed, the Breslow estimator treats the $\ld_0$ as piecewise constant between uncensored failure times} \citep[see Page 473 of][]{lin07}.
{We observe similar results for the coverage of the frequentist confidence bands as in Table \ref{tab:coverage-1}. We also}
%\sout{However, we} 
observe that the area of the credible bands provided by the two Bayesian methods is smaller than that of the confidence bands. 
We also found that both of the two Bayesian methods provide similar coverage results whether $\ld_0$ is chosen to be the continuous function or the piecewise constant function.

In summary, using the Bayesian method can be attractive for estimating data with a relatively small sample size, 
%\sout{especially when the true baseline hazard function is not piecewise constant} 
{as it gives a smaller area}. 
The frequentist method needs to apply an asymptotical approximation to obtain the confidence band, and the approximation can perform {slightly} poorly when the sample size is relatively small. On the other hand, the proposed Bayesian method provides a credible band without using any asymptotical approximation. Notice that the width for the credible band is constant over time. It should be possible to build a varying-width credible band whose overall area is smaller, both in small samples and asymptotically, but 
the construction and analysis of such a band is outside the scope of the paper. 
Yet, the considered fixed--width credible band performs already remarkably well, in particular in finite samples, and achieves the asymptotic limits expected from the Donsker theorem for large sample sizes.
%As we found from Table \ref{tab:coverage-1} that the 95\% Bayesian credible band is comparable to the 95\% confidence band, particularly when the sample size is large, we conclude that our empirical findings match with the Donsker result for the conditional survival function in Section \ref{sec:npBvM}. 

\section{Discussion}
\label{sec:diss}

We provide three new exciting results for the study of the Bayesian Cox model: 1) a joint Bernstein--von Mises theorem for the linear functionals of $\theta$ and $\ld$; in particular, the correlation between the two functionals is captured by the results; 2) a Bayesian Donsker theorem for {the hazard function conditional on $z$} and the {survival function conditional on $z$}; 3) a supremum-norm posterior contraction rate for {the hazard function conditional on $z$}. 

The paper makes major advances on two fronts: on the one hand, it provides new results on optimal posterior convergence rates both in $L^1$-- and $L^\infty$--sense for the hazard; uncertainty quantification is considered for finite dimensional functionals as well as for the posterior cumulative hazard process: those are the first results of this kind for non--conjugate priors (in particular priors for which explicit posterior expressions are not accessible) in this model. On the other hand, the paper provides validation for several classes of practically used histogram priors (see e.g. \cite{ibrahim01}), both for dependent and independent histogram heights.

As a comparison, the results from \cite{cast12} (Theorem 5) and \cite{ghosal17} (Theorem 12.12) require a fast enough preliminary posterior contraction rate of $n^{-3/8}$ in terms of the Hellinger distance. This effectively rules out the use of regular histogram priors, which are limited in terms of rate by $n^{-1/3}$ (corresponding to the optimal minimax rate for Lipschitz functions). Two key novelties here are that {\em a)} we only require a preliminary rate of an order faster than $n^{-1/4}$ (corresponding to $\beta=1/2$ in \eqref{rateveps}) {\em b)} the use of the multiscale approach introduced in \cite{cast14a} enables one to provide both optimal supremum norm contraction rates for the conditional hazard, justifying practically the visual closeness of estimated hazard curves to the true curve, and uncertainty quantification for the conditional cumulative hazard, which follows a BvM for $\Ld(\cdot) e^{\theta'z}$.

{Comparing to \citet{cast20survival} which studied the nonparametric right-censoring model, we would like to highlight the challenges that are unique to our study of the Cox model:
First, deriving the joint Bernstein-von Mises (BvM) results for the Cox model is more challenging than for the right-censoring model as one needs to construct local paths for both $\theta$ and $\ld$ in \eqref{path-theta}-\eqref{path-r}. Our construction of these local paths in multi-dimensions, and jointly with linear functionals of the hazard, is new. }

{Second, one needs to invoke Proposition \ref{prop-tightness} in \citet{ning22supp} to obtain the joint nonparametric BvM theorem. This is in contrast with the study of the right-censoring model, for which Proposition 6 of \citet{cast14a} can be used directly for handling a single non-parametric quantity. Hence, the study of semiparametric models, including the Cox model, requires an extended argument. It is worth mentioning that these results can be useful for studying other semiparametric models in the future. In particular the joint and multidimensional BvMs obtained here are new and could be obtained elsewhere by following our arguments.}

{Third, controlling the LAN reminder terms and the semiparametric bias are significantly more complex tasks for the Cox model. This essentially leads to the requirement of the regularity $\beta > 1/2$. This requirement seems to be unavoidable with the current proof techniques.}

{Last, in practice, confidence bands for the survival function of the Cox model are not commonplace; however, for the right-censoring method, confidence bands for survival functions are routinely used. We are not aware of an easy-to-implement computational algorithm for obtaining the confidence band. An additional contribution of our paper is that we showed that this band could be easily obtained by using the Bayesian approach introduced here. }

We also underline that although not investigated here in details for reasons of space, the results extend to smoother dictionaries than histograms: for instance, if the true hazard is very smooth, one can derive correspondingly very fast posterior rates (obtaining optimal rates $n^{-\beta/(2\beta+1)}$ for any $\beta>1/2$, up to log factors) if one chooses the basis $(\psi_{lk})$ to be a suitably smooth wavelet basis. We refer the interested reader to \cite{cast20survival} for more on how to effectively obtain this. 
{For the frequentist approach considered in this article, the Breslow estimator, which treats the $\ld_0$ as piecewise constant between uncensored failure times \citep[see Page 473 of][]{lin07}, can also be replaced by a smoother estimate, e.g. by taking the kernel-smoothing approach  as in \citet{ramlau83} and \citet{guilloux16} to model smooth baseline hazard functions.}

The present work studies the classical Cox model.
Many extensions of the model have been proposed, such as the Cox model with time-varying covariates \citep{fisher99}, the nonproportional hazards model \citep{schemper02}, the Cox-Aalan model \citep{scheike02} to name a few. The Bayesian nonparametric perspective is particularly appealing in these more complex settings; let us cite two recent practical success stories of the approach in settings going beyond the Cox model (in particular enabling more complex dependencies in terms of covariates and hazard, and/or time dependence):  one is the use of  BART ({Bayesian} additive regression trees) priors in \citet{sparapani16}, another is the use of dependent Dirichlet process priors in \citet{xu19}. 
It would be very desirable to obtain theory and validation for these more complex settings: the present work can be seen as a first step towards this aim.

\section{Proof of Theorem \ref{sec:jointBvM}}
\label{proof:thm1}

{In this section, we prove Theorem \ref{sec:jointBvM}. 
Section \ref{sec:bg} provides the necessary background for studying the Cox model, including the expressions of the LAN-norm, the LAN expansion, the log-likelihood ratio, the squared Hellinger distance between $f_\eta$ and $f_{\eta_0}$, and the relation between the random histogram and Haar wavelets. The relation is useful for studying the histogram priors in Section \ref{sec:priors} as one could invoke results based on the use of Haar wavelets priors directly.
The main proof of Theorem \ref{sec:jointBvM} is given in Section \ref{proof:main-thm1}.
The proofs of the remaining theorems and lemmas in Section \ref{sec:main-results} are given in the Supplemental Materials \citep{ning22supp}. }

\subsection{Background}
\label{sec:bg}
%\subsection{Log-likelihood ratio, LAN expansion, and the least-favorable direction}
We first review several properties of the Cox model that will be frequently used in the proofs of the theorems:
\begin{enumerate}
\item Let us introduce the Hilbert inner product between $(\vartheta_1, g_1)$ and $(\vartheta_2, g_2)$, for any $\vartheta_1, \vartheta_2 \in \mathbb{R}^p$ and any $g_1, g_2 \in L^2(\Lambda_0)$, as
\begin{align}
\langle
(\vartheta_1, g_1), (\vartheta_2, g_2)
\rangle_L
= \Lambda_0 \left\{
\vartheta_1' M_2(\cdot) \vartheta_2 + (g_1(\cdot) \vartheta_2' + g_2(\cdot) \vartheta_1') M_1(\cdot)
+ g_1(\cdot) g_2(\cdot) M_0(\cdot)
\right\}.
\label{lan-norm-inner-product}
\end{align}
This inner-product features in the likelihood expansion (also called Locally Asymptotically Normal expansion or LAN) in the next point and is simply called LAN norm. 

\item For the log-likelihood ratio given in (\ref{log-likelihood}), 
the LAN expansion for this log-likelihood ratio can be written as
\begin{align}
\label{lan}
    \ell_n(\eta) - \ell_n(\eta_0) 
    = - \frac{n}{2}\|\theta-\theta_0, r-r_0\|^2_L + \sqrt{n} W_n(\theta - \theta_0, r-r_0) + R_n(\eta, \eta_0),
\end{align}
where
\begin{equation}
    \label{LAN-norm}
\begin{split}
    \|\theta-\theta_0, r - r_0\|_L^2
    = \Lambda_0
    \big\{
    & (\theta - \theta_0)'M_2(u) (\theta - \theta_0)
    + 2 (r-r_0)(u) (\theta - \theta_0)'M_1(u) \\
    & + (r-r_0)^2(u) M_0(u)
    \big\},
  \end{split}  
 \end{equation}   
 is the LAN-norm part,
 \begin{equation}
 \begin{split}
    W_n(\theta - \theta_0, r-r_0)
    =
    \frac{1}{\sqrt{n}}
    \sum_{i=1}^n \big\{
    & \delta_i \big(
    (\theta - \theta_0)'Z_i + (r - r_0)(Y_i)
    \big)
    - e^{\theta_0'Z_i}
    \big(
    (\theta - \theta_0)'Z_i \Lambda_0(Y_i)\\
    & + (\Lambda_0 (r- r_0))(Y_i)
    \big)
    \big\},
\end{split}    
\label{W}
\end{equation}
is the stochastic part,
and $R_n(\eta, \eta_0)$ is the remainder part, which can be further written as
\begin{align}
\label{remainder}
R_n(\eta, \eta_0) = R_{n,1}(\eta, \eta_0)  + R_{n,2}(\eta, \eta_0),
\end{align}
where 
\begin{align}
\label{R1}
R_{n,1}(\eta, \eta_0) = - \mathbb{G}_n\Psi_n(\eta).
\end{align} 
For a measurable function $f$, 
$
\mathbb{G}_n(f) = \frac{1}{\sqrt{n}}\sum_{i=1}^n (f(X_i) - P_{\eta_0} f),
$
which is the centered and scaled version of the empirical measure,
and
\begin{equation*}
\begin{split}
    \Psi_n(\eta)(X_i) 
     =
    \sqrt{n}
    \Big\{
    & e^{\theta'Z_i}\Ld_0\{e^{r - r_0}\}(Y_i) - e^{\theta_0'Z_i} \Ld_0(Y_i)
	- e^{\theta_0'Z_i} (\theta- \theta_0)'Z_i \Ld_0(Y_i) \\
	& - e^{\theta_0'Z_i}\Ld_0\{r - r_0\}(Y_i)   
	\Big\}.
\end{split}	
\end{equation*}
Let
\begin{align}
\label{M0-theta}
M_0(\theta)(\cdot) = \mathbb{E}_{\eta_0} (\mathbbm{1}_{u \leq Y} e^{\theta'Z})
 = \int \bar G_z(u) e^{\theta'z} e^{-\Ld_0(u) e^{\theta_0'z}} f_Z(z) dz,
\end{align}
and $M_0(\cdot)$ and $M_1(\cdot)$ in (\ref{M0}) and (\ref{M1}) respectively,
\begin{equation}
\begin{split}	
    R_{n,2}(\eta, \eta_0) 
    = 
    & - n\Lambda_0 \Big\{
      M_0(\theta)(\cdot) e^{(r - r_0)(\cdot)} - M_0(\cdot) - (\theta - \theta_0)'M_1(\cdot) \\
    &   - (r-r_0)(\cdot)M_0(\cdot)
    \Big\} + \frac{n}{2}\| \theta - \theta_0, r - r_0\|_L^2.
\end{split}    
    \label{R2}    
\end{equation}

\item Recall that $\gamm = M_1/M_0$ is the least favorable direction and $\tilde I_{\eta_0}$ is the efficient information matrix in (\ref{eff-info-mtx}), the LAN-norm in (\ref{LAN-norm}) can be also written as 
\begin{align}
\label{LAN-norm-v2}
\|\theta - \theta_0, r-r_0\|_L^2 = (\theta - \theta_0)'\tilde I_{\eta_0} (\theta - \theta_0) + \|0, r - r_0 + \gamm' (\theta - \theta_0)\|_L^2.
\end{align}

\item For the density function in (\ref{eqn:density-function}), the squared Hellinger distance between $f_{\eta}$ and $f_{\eta_0}$ is
\begin{align}
 h^2(f_{\eta}, f_{\eta_0}) 
    & = 
    \int \int_0^1 
    \left[
        \sqrt{S_\eta} - \sqrt{S_{\eta_0}}
    \right]^2(u, z) g_{z}(u) f_Z(z) du dz
    \label{hellinger-1}
    \\
    & \quad +
    \int \int_0^1 \bar{G}_{z} (u) 
    \left[
        \sqrt{\lambda S_\eta e^{\theta' z}}
        - 
        \sqrt{\lambda_0 S_{\eta_0} e^{\theta_0' z}}
    \right]^2 (u, z) f_Z(z) du dz
    \label{hellinger-2}\\
    & \quad + 
    \int \bar{G}_{z} (1) 
    \left[
        \sqrt{S_\eta} - \sqrt{S_{\eta_0}}
    \right]^2(1, z) f_Z(z) dz,
    \label{hellinger-3}
\end{align}
where we slightly abuse the notation by denoting $S_{\eta}(u, z) = e^{-\Ld_0(u) e^{\theta'z}}$. 
A similar expression appears on page 34 of \citet{cast12} (up to a typo in the third term in his expression, fixed in the last display) 

\end{enumerate}

We also review the relation between the random histogram and Haar wavelets.
Let's denote $r_H = (r_1, \dots, r_{2^{L} + 1})'$, where $r_H = \log \ld_H$,  as the step heights of the random histogram and $r_S = (r_{-1}, r_{00}, r_{01}, \dots, r_{L(2^L - 1)})'$ as the coefficients in the Haar wavelet prior, then through Haar transformation,
\begin{align}
\label{haar-transf}
r_S = \Psi r_H
\end{align}
where for $I_k^{l} = (k 2^{-l}, (k+1)2^{-l}]$,
$$
\Psi_{-1,j} = 2^{-(L+1)}, \quad 
\Psi_{lk,j} = 2^{-(L+1)+l/2}
[
\mathbbm{1}_{I_{j-1}^{L+1} \subset I_{2k}^{l+1}} - \mathbbm{1}_{I_{j-1}^{L+1} \subset I_{2k+1}^{l+1}}
],
$$
and $2^{(L+1)/2} \Psi$ is an orthogonal matrix.

Last, as the posterior concentrates on the set $\|\ld - \ld_0\|_\infty \leq \zeta_n$ in Lemma \ref{lemma1}, we use the fact that
$$
\ld - \ld_0 = e^{r} - e^{r_0} = e^{r_0} \left(e^{r - r_0} - 1\right);
$$
therefore, as long as $\zeta_n = o(1)$ which is the case since $\beta > 1/2$, by Taylor's theorem and assumption \ref{asp:iii} such that $\|\ld_0\|_\infty$ by some constant, $\|\ld - \ld_0\|_\infty = O(\zeta_n)$ automatically translates into the same rate for $\|r - r_0\|_\infty$. This fact will be automatically applied in our proofs.  

\subsection{Proof of the main theorem}
\label{proof:main-thm1}

{We now prove Theorem \ref{sec:jointBvM}. We follow \citet{cast15BvM} and show that the Laplace transform of the induced posterior distribution on the functional of interest converges to the corresponding Laplace transform of the optimal (efficient) Gaussian limit.}
From Lemma 1 and Lemma 2 of \citet{cast15BvMSupp}, it is sufficient to show that the Laplace transform of the induced posterior distribution on the functional of interest in (\ref{thm1-1}) converges to the corresponding Laplace transform of the optimal (efficient) Gaussian limit.
That is,
define $\varphi_a(\theta) = \theta'a$,
$\varphi_b(\lambda) = \int \lambda b$
and 
let $\hat \varphi_a = 
\varphi_a(\theta_0) + W_n^{(1)}(a)/\rn$
and
$\hat \varphi_b = 
\varphi_b(\ld_0) + W_n^{(2)}(b)/\rn$,
where 
\begin{align*}
& W_n^{(1)}(a) = W_n(\tIinv a,\ - \gamm' \tIinv a), \\
& W_n^{(2)}(b) = W_n\left(- \tIinv \Ld_0\{b\gamm\},\ \gamb + \gamm' \tIinv \Ld_0\{b\gamm\}\right),
\end{align*}
and $W_n(\cdot, \cdot)$ is given in (\ref{W}), 
our goal is to show that for any $h = (t, s) \in \mathbb{R}^2$, 
\begin{align}
\label{pf-thm1-1}
\mathbb{E}
\left[
e^{\rn h \left( \varphi_a(\theta) - \hat \varphi_a, \ \varphi_b (\lambda) - \hat \varphi_b \right)'} \given X, A_n
\right] 
\stackrel{P_{\eta_0}}{\to}
e^{h \Sigma_{a,b} h'/2},
\end{align}
where 
\begin{equation}
\label{Sigma-ab}
\Sigma_{a,b} = 
\begin{pmatrix}
a' \tIinv a & - a' \tIinv \Lambda_0\{ b \gamm\}\\
- a' \tIinv \Lambda_0\{ b \gamm\} & 
\Lambda_0\{b \gamb\} + \Lambda_0\{b \gamm'\} \tIinv \Lambda_0\{b\gamm\}
\end{pmatrix}.
\end{equation} 
%Recall that 
%$$
%W_n(g,h) = \frac{1}{\sqrt{n}} \sum_{i=1}^n \left\{
%\delta_i (g'Z_i + h(Y_i)) - e^{\theta_0'Z_i} (g'Z_i \Ld_0(Y_i) + (\Ld_0h)(Y_i)) 
%\right\}.
%$$
By applying Bayes' formula and dividing the expression at the right hand side on both side of (\ref{pf-thm1-1}), the display in (\ref{pf-thm1-1}) can be written as 
\begin{align}
\label{pf-thm1-2}
\frac{
\int_{A_n}
e^{\sqrt{n}h(\vpa(\theta) - \hat\varphi_a, \vpb(\ld) - \hat \varphi_b)' + \ell_n(\eta) - \ell_n(\eta_0) - h'\Sigma_{a,b} h/2}
d\Pi(\eta)
}{
\Pi(A_n \given X) \int e^{\ell_n(\eta) - \ell_n(\eta_0)} d\Pi(\eta)
} \stackrel{P_{\eta_0}}{\to}1.
\end{align}
Below, we will provide the key steps for proving (\ref{pf-thm1-2}). Intermediate lemmata along with their proofs are left to Section \ref{sec:support-lemma-thm1} in the supplemental material.

To bound the numerator at the left hand side of (\ref{pf-thm1-2}),
an important step is to show that
\begin{equation}
\label{pf-thm1-2-numerator}
\begin{split}
& \sup_{\eta \in A_n} |\rn h(\vpa(\theta) - \hat\varphi_a, \vpb(\ld) - \hat \varphi_b)' + \ell_n(\eta) - \ell_n(\eta_0) - h'\Sigma_{a,b} h/2| \\
& \quad \leq
\sup_{\eta \in A_n} |\ell_n(\eta_h) - \ell_n(\eta_0)| + o(1) + o_{P_{\eta_0}}(1),
\end{split}
\end{equation}
where $\eta_h = (\theta_h, r_h)$, and $\theta_h$ and $r_h$ are given in (\ref{path-theta}) and (\ref{path-r}) respectively.

To prove (\ref{pf-thm1-2-numerator}), using the expression of the LAN-norm given in (\ref{lan}) and note that one can write
$
\ell_n(\eta) - \ell_n(\eta_0) = \ell_n(\eta) - \ell_n(\eta_0) - [\ell_n(\eta_h) - \ell_n(\eta_0)] + \ell_n(\eta_h) - \ell_n(\eta_0)
$ 
and then obtain
\begin{align}
& \ell_n(\eta) - \ell_n(\eta_0) - [\ell_n(\eta_h) - \ell_n(\eta_0)] \nonumber \\
& \quad = 
-\frac{n}{2}\left(\|\theta - \theta_0, r - r_0\|_L^2 - \|\theta_h - \theta_0, r_h - r_0\|_L^2 \right)
\label{pf-thm1-5}\\
& \qquad + \sqrt{n}\left(W_n(\theta - \theta_0, r- r_0) - W_n(\theta_h - \theta_0, r_h - r_0)\right)
\label{pf-thm1-6}\\
& \qquad + R_n(\eta, \eta_0) - R_n(\eta_h, \eta_0),
\label{pf-thm1-7}
\end{align}
where $\|\cdot, \cdot\|_L$, $W_n(\cdot, \cdot)$, and $R_n(\cdot, \cdot)$ are defined in (\ref{LAN-norm}), (\ref{W}), and (\ref{remainder}) respectively. 

On the other hand, we write $(\theta - \theta_0)'a$ and $\Ld\{b\} - \Ld_0\{b\}$ as their LAN-norm Hilbert inner product forms, i.e., 
\begin{align}
(\theta & - \theta_0)'a
=  \left\langle 
(\theta - \theta_0, r - r_0), 
\left(\tIinv a, - \gamm' \tIinv a\right)
\right\rangle_L,
\label{theta-lan-norm}\\
\Ld\{b\} -  \Ld_0 \{b\}  
& = \left\langle
\left(\theta - \theta_0, \frac{\ld - \ld_0}{\ld_0}\right),
\left( - \tIinv \Ld_0\{b\gamm\}, \gamb + \gamm'\tIinv \Ld_0\{b \gamm\} \right)
\right\rangle_L.
\label{lambda-lan-norm}
\end{align}
by using fact that (also, see (\ref{lan-norm-inner-product})) the LAN-norm Hilbert inner product between $(\vartheta_1, g_1)$ and $(\vartheta_2, g_2)$, for any $\vartheta_1, \vartheta_2 \in \mathbb{R}^p$ and any $g_1, g_2 \in L^2(\Lambda_0)$, is defined as
\begin{align*}
\langle
(\vartheta_1, g_1), (\vartheta_2, g_2)
\rangle_L
= \Lambda_0 \left\{
\vartheta_1' M_2(\cdot) \vartheta_2 + (g_1(\cdot) \vartheta_2' + g_2(\cdot) \vartheta_1') M_1(\cdot)
+ g_1(\cdot) g_2(\cdot) M_0(\cdot)
\right\}.
\end{align*}
The right hand side of (\ref{lambda-lan-norm}) can be further decomposed into three parts such that
\begin{align}
\label{lambda-decomp}
\Ld\{b\} - \Ld_0\{b\} = B_1(\eta, \eta_0) + B_2(\eta, \eta_0) - B_3(\eta, \eta_0),
\end{align}
where
\begin{align}
& B_1(\eta, \eta_0) = \left\langle
\left(\theta - \theta_0, r - r_0 \right), 
\left(- \tIinv \Ld_0\{b\gamm\}, \gambln + \gammln' \tIinv \Ld_0\{b \gamm\}\right)
\right\rangle_L,
\label{B1}\\
& B_2(\eta, \eta_0) = 
\left\langle
\left(0, \frac{\lambda - \lambda_0}{\lambda_0}\right), 
\left(0, \gamb - \gambln + (\gamm - \gammln)' \tIinv \Ld_0\{b \gamm\}\right)
\right\rangle_L,
\label{B2}\\
& B_3(\eta, \eta_0) = 
\left\langle
\left(0, r - r_0 - \frac{\lambda - \lambda_0}{\lambda_0}\right), 
\left(- \tIinv \Ld_0\{b\gamm\}, \gambln + \gammln' \tIinv \Ld_0\{b \gamm\}\right)
\right\rangle_L.
\label{B3}
\end{align}
The third term (\ref{B3}) is a semiparametric bias. 

Now we plug (\ref{pf-thm1-5})-(\ref{pf-thm1-7}) and (\ref{lambda-decomp}) into (\ref{pf-thm1-2-numerator}) and note that $\hat \varphi_a = \theta_0'a + W_n^{(1)}(a)/\rn$ and $\hat \varphi_b = \Ld_0\{b\} + W_n^{(2)}(b)/\rn$, then the left hand side of (\ref{pf-thm1-2-numerator}) can be written as 
\begin{align}
& \rn h (\varphi_a(\theta) - \hat \varphi_a,  \varphi_b(\ld) - \hat\varphi_b)'
+ \ell_n(\eta) - \ell_n(\eta_0) - [\ell_n(\eta_h) - \ell_n(\eta_0)] - h'\Sigma_{a,b}h/2
\nonumber \\
& \quad
= 
t \rn (\theta - \theta_0)'a 
+ s\rn B_1(\eta, \eta_0) - \frac{n}{2}(\|\theta - \theta_0, r - r_0\|_L^2 
- \|\theta_h - \theta_0, r_h - r_0\|_L^2)
\label{pf-thm1-part1}
\\
& \qquad 
+ s\rn B_2(\eta, \eta_0) - {h'\Sigma_{a,b}h}/{2}
+ R_n(\eta, \eta_0) - R_n(\eta_h, \eta_0) - s \rn B_3(\eta, \eta_0)
\label{pf-thm1-part2} \\
& \qquad
- tW_n^{(1)}(a) - sW_n^{(2)}(b) + \rn W_n(\theta - \theta_0, r - r_0) - \rn W_n(\theta_h - \theta_0, r_h - r_0).
\label{pf-thm1-part3}
\end{align}
Thus to proof (\ref{pf-thm1-2-numerator}),
one needs to show the last display is bounded by $o(1) + o_{P_{\eta_0}}(1)$ uniformly on the set $\eta \in A_n$. 
We will bound each line in the last display:

\begin{enumerate}
\item To bound (\ref{pf-thm1-part1}), by plugging-in the expressions of $\theta_h$, $r_h$, and (\ref{theta-lan-norm}), one can check that  
$$
t \rn (\theta - \theta_0)'a + s\rn B_1(\eta, \eta_0) = n \langle (\theta - \theta_0, r - r_0), (\theta - \theta_h, r - r_h) \rangle_L.
$$
Also, note that
$
\langle (\theta_h - \theta_0, r_h - r_0) \rangle_L
= 
\langle (\theta_h - \theta, r_h - r) \rangle_L
+ 
\langle (\theta - \theta_0, r - r_0) \rangle_L,
$
by expanding the two squared LAN-norms in (\ref{pf-thm1-part1}), we have
\begin{align*}
& n \langle (\theta-\theta_0, r- r_0), (\theta - \theta_h, r - r_h)\rangle_L 
- \frac{n}{2} \|\theta - \theta_0, r - r_0\|_L^2 + \frac{n}{2} \|\theta_h - \theta_0, r_h - r_0\|_L^2\\
& \quad = \frac{n}{2} \|\theta - \theta_h, r - r_h\|_L^2.
\end{align*}
Define
\begin{align}
\label{Dn}
D_n & = \frac{n}{2} \|\theta - \theta_h, r - r_h\|_L^2  
- \frac{h'\Sigma_{a,b}h}{2}.
\end{align}
Note that the expression in (\ref{pf-thm1-part1}) equals to $D_n$.
By invoking Lemma \ref{lemma:bounding-Dn} in the supplemental material, we obtain
\begin{align*}
\sup_{\eta \in A_n}|D_n| 
&  \leq s^2\|\gamb^2 - \gambln^2\|_\Lone + (t^2 + s^2) (p^2 2^{-L_n} \|\gambln\|_\Lone + p^4 2^{-2L_n}).
\end{align*}
To bound the last display, we first invoke Lemma \ref{lemma:psib-1} in the supplemental material to obtain $\|\gambln\|_\Lone \leq \|\gambln\|_\Ltwo \lesssim \|b\|_\Ltwo$.
Then applying the inequality $\|fg\|_\Lone \leq \|f\|_\Lone \|g\|_\Linfty$ and using {\ref{cond:B}} to obtain the bound $\|\gamb^2 - \gambln^2\|_\Lone \leq 
(\|\gamb\|_\Lone + \|\gambln\|_\Lone)(\|\gamb - \gambln\|_\Linfty) \leq (\|\gamb\|_\Lone + \|\gambln\|_\Lone)/(\rn \epsilon_n)$. 
Since $b \in L^\infty([0,1])$ and $p, t, s$ are constants, by plugging-in the two upper bounds, the last display is bounded by 
$
s^2 \|b\|_\Ltwo / (\rn \epsilon_n)+ p^2 2^{-L_n} \|b\|_\Ltwo + o(1),
$
which is $o(1)$, as $\rn \epsilon_n \to \infty$ and $L_n \to \infty$ as $n \to \infty$.

\item To bound (\ref{pf-thm1-part2}), we first deal with the first term.
By (\ref{bound-B2}) in Lemma \ref{lemma:bounding-B2},
$$
\rn \sup_{\eta \in A_n} |sB_2(\eta, \eta_0)| 
\lesssim |s| \rn \epsilon_n \|\gamb - \gambln\|_\Linfty + |s| p^2 \rn \epsilon_n 2^{-L_n} \|b\|_\Lone.
$$
Using \ref{cond:B} and the assumption $\rn \epsilon_n 2^{-L_n} = o(1)$, note that $b \in L^{\infty}([0,1])$, the last display is $o(1)$.

To bound the last three terms in (\ref{pf-thm1-part2}),
from (\ref{remainder}), $R_n(\eta, \eta_0) = R_{n,1}(\eta, \eta_0) + R_{n,2}(\eta, \eta_0)$,
where the expressions of $R_{n,1}(\eta, \eta_0)$ and $R_{n,2}(\eta, \eta_0)$ are given in (\ref{R1}) and (\ref{R2}) respectively. 
Then, 
\begin{equation}
\begin{split}
& \sup_{\eta \in A_n}|R_n(\eta, \eta_0) - R_n(\eta_h, \eta_0) - s\rn B_3(\eta, \eta_0)|  
\\
& \quad
\leq 
\sup_{\eta \in A_n}|R_{n,1}(\eta, \eta_0) - R_{n,1}(\eta_h, \eta_0)| 
\\
& \quad \quad 
+ 
\sup_{\eta \in A_n}|R_{n,2}(\eta, \eta_0) - R_{n,2}(\eta_h, \eta_0) - s\rn B_3(\eta, \eta_0)|.
\label{bounding-remainder}
\end{split}
\end{equation}
%\end{align}

We apply Lemma \ref{lemma-R1} in the supplemental material to bound the first term in the last display. 
To verify the conditions in Lemma \ref{lemma-R1}, since $\|a\|_\infty$ is bounded and $b \in L^{\infty}([0,1])$, for $\Delta_1$ and $\Delta_{2,L_n}$ defined in (\ref{K1}) and (\ref{K2})  in the supplemental material, it is easy to check that both $\Delta_1/\rn = o(1)$ and $\Delta_{2,L_n}/\rn = o(1)$ by applying triangular inequalities.  
Thus we can invoke Lemma \ref{lemma-R1} to obtain 
$$
\sup_{\eta \in A_n} |R_{n,1}(\eta, \eta_0) - R_{n,1}(\eta_h, \eta_0)|
\leq O_{P_{\eta_0}}\left(L_n^2/\rn + \epsilon_n L_n\right) = o_{P_{\eta_0}}(1),
$$
as $L_n^2/\rn = o(1)$ and $\epsilon_n L_n= o(1)$ by assumptions.

To bound the last line in (\ref{bounding-remainder}),
define $K_{a,b,t,s} = p^2(|t| \|a\|_\infty + |s|\|b\|_\Lone) + |s|\|b\|_\Linfty$
and $\tilde K_{a,b,t,s} = p^2(|t|\|a\|_\infty + |s|\|b\|_\Lone) + |s|\|b\|_\Ltwo$. Since $t, s, p$ are all fixed constants, $\|a\|_\infty$ is bounded, and $b \in L^\infty([0,1])$, 
then $K_{a,b,t,s}=O(1)$ and $\tilde K_{a,b,t,s} = O(1)$. By invoking Lemma \ref{thm1-lemma2} in the supplemental material, we have
\begin{align*}
& \sup_{\eta \in A_n}|R_{n,2}(\eta, \eta_0) - R_{n,2}(\eta_h, \eta_0) - s\rn B_3(\eta, \eta_0)|\\
& \quad \lesssim \tilde K_{a,b,t,s}^3/\rn + K^2_{a,b,h,p}L_n^2 \epsilon_n + |s| p^2 \rn \epsilon_n2^{-L_n} + \rn \epsilon_n^2 L_n K_{a,b,t,s}\\
& \quad \lesssim L_n^2 \epsilon_n + \rn \epsilon_n 2^{-L_n} + \rn \epsilon_n^2 L_n,
\end{align*}
which is $o(1)$ as $\epsilon_n L_n = o(1)$ and $\rn \epsilon_n^2 L_n = o(1)$ by assumption.
We thus showed that (\ref{pf-thm1-part2}) is $o_{P_{\eta_0}}(1)$ for $\eta \in A_n$.

\item To bound (\ref{pf-thm1-part3}),
we directly plug-in the expressions of $\theta_h$ and $r_h$ to obtain
\begin{align*}
& \rn W_n(\theta - \theta_0, r- r_0) - \rn W_n(\theta_h - \theta_0, r_h - r_0)\\
& \quad = W_n(t \tIinv a - s\tIinv \Ld_0\{b\gamm\}, \ - t \gammln' \tIinv a + s\gambln + s \gammln' \tIinv \Ld_0\{b\gamm\}).
\end{align*}
Due to the linearity of $W_n(\cdot, \cdot)$,
(\ref{pf-thm1-part3}) can be also written as
$$
{W_n(0, t(\gamm - \gammln)' \tIinv a)} - {W_n(0, s(\gamb - \gambln))} - {W_n(0, s(\gamm - \gammln)' \tIinv \Ld_0(b \gamm))}.
$$ 
Applying the fourth point of Lemma \ref{lemma:psib-2} in the supplemental material and using the fact that $2^{-L_n/2} = o(1)$ as $n\to \infty$, the second term in the last display is $o_{P_{\eta_0}}(1)$. 
To bound the first term in the last display,
since $p, t, s$ are fixed constants and $\|a\|_\infty$ is bounded, by \ref{asp:v}, we have
\begin{align*}
|t(\gamm - \gammln)' \tIinv a| & \leq |t| p^2 \|a\|_\infty \|\tIinv\|_{(\infty, \infty)} \max_j \|\gamm^j - \gammln^j\|_\Linfty  \\
& \lesssim 
\max_j \|\gamm^j - \gammln^j\|_\Linfty.
\end{align*}
By the fourth point of Lemma \ref{lemma:psib-2}, then $W_n(0, t(\gamm - \gammln)' \tIinv a) = o_{P_{\eta_0}}(1)$. 
The bound for the third term in the last display can be obtained similarly, which is also $o_{P_{\eta_0}}(1)$. 
Therefore, we showed that (\ref{pf-thm1-part3}) is $o_{P_{\eta_0}}(1)$. 
\end{enumerate}
Now by collecting the bounds derived above for (\ref{pf-thm1-part1})-(\ref{pf-thm1-part3}), we thus verified (\ref{pf-thm1-2-numerator}). 

With (\ref{pf-thm1-2-numerator}), the expression at the left hand side in (\ref{pf-thm1-2}) is bounded by
$$
\frac{
\int_{A_n} e^{\ell_n(\eta_h) - \ell_n(\eta_0) + o_{P_{\eta_0}(1)}}
d \Pi(\eta)
}{
\Pi(A_n \given X) \int e^{\ell_n(\eta) - \ell_n(\eta_0)} d\Pi(\eta)
}.
$$
By \ref{cond:P}, i.e., $\Pi(A_n \given X) = 1 + o_{P_{\eta_0}}(1)$, and then by \ref{cond:C1}, the change of variables condition, using the inequality $e^{x} = 1 + o(x)$ if $x = o(1)$, 
the last display is bounded by $1 + o_{P_{\eta_0}}(1)$. We thus complete the proof. 

\section*{Acknowledgement}

The authors would like to thank St\'ephanie van der Pas for helpful discussions with the $\mathsf{R}$ code for simulation studies. The authors would also like to thank the Associate
Editor and two referees for insightful comments.

\begin{supplement}
The supplement \citet{ning22supp} includes the proofs of the results stated in this paper. 
The $\mathsf{R}$ code of the Bayesian method in the paper is available on the website \url{https://github.com/Bo-Ning/Bayesian-Cox-Piecewise-Constant-Hazard-Model}. 
\end{supplement}

%%%%%%%%%%%%%%%%%%%%%%%%%%%%%%%%%%%%%%%%%%%%%%%%%%%%%%%%%%%%%%

%%%%%%%%%%%%%%%%%%%%%%%%%%%%%%%%%%%%%%%%
%%%%%%%%%%%%%%%% Bibliography %%%%%%%%%%%%%%%%%
%%%%%%%%%%%%%%%%%%%%%%%%%%%%%%%%%%%%%%%
\bibliographystyle{chicago}
\bibliography{citation}

%%%%%%%%%% Merge with supplemental materials %%%%%%%%%%
%% \section{Supplementary material}
%\section{Supplementary Material}
%\beginsupplement
%
%This supplementary materials contains the proofs of the main results in Section \ref{main-results}. 

%%%%%%%%%% Merge with supplemental materials %%%%%%%%%%
%\widetext
\clearpage
\setcounter{equation}{0}
\setcounter{figure}{0}
\setcounter{table}{0}
\setcounter{page}{1}
\makeatletter
\renewcommand{\theequation}{S\arabic{equation}}
\renewcommand{\thefigure}{S\arabic{figure}}
\setcounter{lemma}{0}
\renewcommand{\thelemma}{S\arabic{lemma}}
\setcounter{theorem}{0}
\renewcommand{\thetheorem}{S\arabic{theorem}}
\setcounter{remark}{0}
\renewcommand{\theremark}{S\arabic{remark}}
\setcounter{corollary}{0}
\renewcommand{\thecorollary}{S\arabic{corollary}}
\setcounter{prop}{0}
\renewcommand{\theprop}{S\arabic{prop}}
\renewcommand{\bibnumfmt}[1]{[S#1]}
\renewcommand{\citenumfont}[1]{S#1}
\setcounter{section}{0}
\renewcommand{\thesection}{S\arabic{section}}

\begin{frontmatter}

\title{Supplement to ``Bayesian multiscale analysis of the Cox model''}

%%%%%%%%%% Merge with supplemental materials %%%%%%%%%%
%%%%%%%%%% Prefix a "S" to all equations, figures, tables and reset the counter %%%%%%%%%%

%\tableofcontents
%\contentsmargin{1em}
%\dottedcontents{section}[]{}{4em}{}

\author{\fnms{Bo Y.-C.} \snm{Ning}\thanksref{t1}\ead[label=e2]{bycning@ucdavis.edu}}
\and
\author{\fnms{Isma\"el} \snm{Castillo}\thanksref{t2}\ead[label=e1]{ismael.castillo@upmc.fr}}

 \thankstext{t1}{The author gratefully acknowledges the support from the Fondation Sciences Math\'ematiques de Paris postdoctoral fellowship.}
\thankstext{t2}{The author gratefully acknowledges support from the Institut Universitaire de France (IUF) and from the ANR grant ANR-17-CE40-0001 (BASICS).}
\affiliation{Sorbonne University \& UC Davis\thanksmark{m1} and Sorbonne University \& IUF\thanksmark{m2}}

\address{University of California, Davis \\
Department of Statistics,\\
1227 Mathematical Science Building,\\
One Shields Avenue, 
Davis, CA 95616 United States\\
\printead{e2}}

\address{Sorbonne Universit\'e \&
Institut Universitaire de France\\
  Laboratoire de Probabilit\'es, 
Statistique et Mod\'elisation\\ 4, Place Jussieu, 75252, Paris cedex 05, France\\
  \printead{e1}}
\runauthor{Ning \& Castillo}

\begin{abstract}
In this supplemental, we include the proofs for the results stated in the main paper. We also provide a summary of contents and the background on the Cox model.
\end{abstract}

\end{frontmatter}

%\ma{I would regroup some sections to slightly reduce their numbers, e.g. S4 and S5, and S6 and S7. I would also put the Hellinger rate, verification of (P) and also possibly (B) before, just before the proof of Theorem 1, as a 'warm up' to that proof, since these preliminary results are in a sense more `elementary'}

%\startcontents
%\tableofcontents
\section*{Table of Contents}
\vspace{-0.2cm}
{\linespread{-1}\selectfont
\startcontents
\printcontents{ }{1}{}
}

\section{Summary of contents}
%Section \ref{sec:bg} 

Section \ref{proof-hellinger} focuses on deriving a preliminary Hellinger contraction rate, ${{\nu}}_n$ given in \eqref{rateveps}, for the posterior and verifying condition \ref{cond:P} for the specific priors considered in Section \ref{sec:priors}. 
Comparing with previous studies by \citet{cast12} and \citet{ghosal17}, we adopt a novel argument that enables us to obtain a slower $\ell_1$-rate for the hazard $\lambda$. 
This rate enables us to work with the piecewise constant prior as, for example, in Condition \ref{cond:P}, if only the $\ell_\infty$-rate is available, the assumption $\rn \epsilon_n^2 L_n = o(1)$ shall be replaced by $\rn \epsilon_n\zeta_n L_n = o(1)$ instead (one can check this in Lemma \ref{sec:bounding-R2}). 
Then by plugging-in $L_n$ in \eqref{choicel} and the rate in \eqref{rateveps},  
one easily checks that the former assumption implies $\beta > 1/2$ but the latter implies $\beta > 1$.

{{Section \ref{proof:thm1} provides the proof of those intermediate lemmata used to prove Theorem \ref{thm:bvm-linear-fns}.
More specifically, we} obtain upper bounds for the ``semiparametric bias'' and the two remainder terms from the LAN expansion. Since those bounds are not only used in the proof in this section but also in the proofs of the nonparametric BvM results and the supremum-norm rate later on, we include these as separate lemmas.

Two nonparametric BvM theorems are established in Section \ref{joint-npBvM}. 
The first theorem concerns the joint posterior distribution of $\theta$ and $\ld$, and the second one concerns {the hazard function conditional on $z$}. 
To prove the first theorem, a key step (Proposition \ref{prop-tightness}) consists in establishing a parametric $1/\rn$--rate for the baseline hazard function through weakening the $L^2$ norm to the multiscale norm defined in \eqref{multiscale-space}. 
This extends the tightness condition for nonparametric models used in \citet{cast14a} (Proposition 6 therein) to semiparametric models. 
It can be useful for studying other semiparametric models as well. 
Our nonparametric BvM result not only can be of independent interest but also serves as an intermediate step for obtaining the joint Bayesian Donsker result in Theorem \ref{thm:donsker-joint}.

Using the nonparametric BvM result for the joint posterior distribution, one can apply the functional delta method to derive the Bayesian Donsker theorem for the cumulative hazard function. Its proof is given in Section \ref{pf-donsker}.
One can further obtain the supremum-rate by following the approach developed by \citet{cast14b}. The proof is given in Section \ref{pf-supnorm}.
In Section \ref{sec:lower-bound}, we derive a lower bound for the supremum-norm rate and show this lower bound matches with the supremum-norm rate, implying that the obtained rate is optimal.

In Section \ref{verify-conditions}, we verify conditions \ref{cond:B}, \ref{cond:C1}, and \ref{cond:C2} for the specific priors considered in the paper.
In Section \ref{auxiliary}, we gather the remaining lemmas used in our proofs.
In particular, it includes the key proposition we mentioned earlier for establishing the tightness condition and two lemmas (Lemmas \ref{centering1} and \ref{lemma-centering-Lambda}) on centering and efficiency that imply that the joint posterior can be centered at efficient frequentist estimators.

%%%%%%%%%%%%%%%%%%%%%%%%%%%%%%%%%%%%%%%%%%%%%%%%%
%%%%%%%%%%%%%%%%%%%%%%%%%%%%%%%%%%%%%%%%%%%%%%%%%
\section{Proof of Lemma \ref{lemma1}}
\label{proof-hellinger}

In this section, we obtain the Hellinger rate ${{\nu}}_n = \left({\log n}/{n}\right)^{\frac{\beta}{2\beta + 1}}$ by invoking the general theory of posterior contraction proposed by \citet{ghosal00} (see also Theorem 8.9 of \citet{ghosal17}).
Let's define the Kullback-Leibler divergence and the Kullback-Leibler variation between densities $f$ and $g$ as $K(f, g) = \int f \log (f/g)$ and $V(f, g) = \int f (\log (f/g) - K(f/g))^2$, and let $N(\epsilon, \mathcal{F}, \rho)$ stands for the $\epsilon$-covering number of a set $\mathcal{F}$ with respect to a metric $\rho$, which is the minimal number of $\epsilon$-balls in $\rho$-metric needed to cover the set $\mathcal{F}$. 

\subsection{Auxiliary lemmata for proving Lemma \ref{lemma1}}
\label{aux-lemma1}

The following two lemmas are useful for verifying the prior mass condition in the general theory.
\begin{lemma}
\label{bounding-hellinger}
For fixed $r_1$, $r_2$, $\theta_1$, and $\theta_2$, 
let $p_1, p_2$ be the distributions associates to $(\theta_1, r_1)$ and $(\theta_2, r_2)$.
Assuming that there exists a constant $0 \leq Q \leq 1/4$ such that $\|r_1 - r_2\|_\Linfty + \|\theta_1 - \theta_2\|_1 C \leq Q$ for some constant $C$, then there exist a constant $c$ depending on $C$ only such that 
$$
h^2(p_1, p_2) \leq c Q^2 e^{2Q}.
$$
\end{lemma}

\begin{proof}
The proof is similar to that of Lemma 7 of \citet{cast12}. Although in his proof, $\theta$ is assumed to be a scalar, the same proof carries out for multivariate $\theta$ without difficulty.  
\end{proof}

\begin{lemma}
\label{bounding-KL}
Under the same setting as in Lemma \ref{bounding-hellinger}, assuming that $\|r_1 - r_2\|_\Linfty + \|\theta_1 - \theta_2\|_1 C \leq Q$ for some constant $Q$, then 
$$
K(p_1, p_2) = p_1 \log (p_1/p_2) \lesssim h^2(p_1, p_2), 
\quad 
V(p_1, p_2) = p_1 \log^2 (p_1/p_2) \lesssim h^2(p_1, p_2).
$$
\end{lemma}

\begin{proof}
The proof is similar to that of Lemma 8 of \citet{cast12}. 
\end{proof}

\subsection{Proof of Lemma \ref{lemma1}}
\label{proof-lemma1}

By invoking the general theory of posterior contraction in \citet{ghosal00}, 
we need to verify the following three conditions,
\begin{align}
    &\Pi(\mathcal{F}_n^c) \lesssim \exp(-(C_1+4) n\epsilon_n^2),
    \label{cond-1} \\
    &\Pi(B_{\KL}(\eta_0, C_2\epsilon_n)) \gtrsim \exp(-C_1 n \epsilon_n^2),
    \label{cond-2} \\
    &\log N(\epsilon_n, \mathcal{F}_n, h) \leq C_3 n\epsilon_n^2,
    \label{cond-3}
\end{align}
where $B_{\KL}(\eta_0, \epsilon) = \left\{\eta: K(f_{\eta_0}, f_\eta) \leq \epsilon^2,\ V(f_{\eta_0}, f_\eta) \leq \epsilon^2\right\}$ and $C_1, \dots, C_3$ are positive constants.

We first verify the above three conditions for the Haar wavelet prior \ref{prior:W}. 
Consider either independent standard normal prior or independent standard Laplace prior on $Z_{lk}$, then, $r_{lk} \sim N(0, \sigma_l^2)$ or $r_{lk} \sim \text{Laplace}(0, \sigma_l)$. 
To verify (\ref{cond-1}), 
let $\mathcal{F}_n = \{(\theta, r): \|\theta\|_\infty \leq C, \ r_{lk} = 0 \ (\forall \ l > L_n, k), \ |r_{lk}| \leq n \ (\forall \ l \leq L_n, k)\}$, where $C$ is the same as it in {\rm \ref{prior:T}}.
Then, note that the prior of $\theta$ is truncated at $-C$ and $C$, by following from a union bound, we obtain $\Pi(\mathcal{F}_n^c) \leq \Pi(\theta: \|\theta\|_\infty > C) + \sum_{l,k} \Pi(r: |r_{lk}| > n)
\leq \sum_{l,k} \Pi(|r_{lk}| > n)$.
Using $\sigma_l \leq 1$ (either choosing $\sigma_l = 1$ or $\sigma_l = 2^{-l}$) and the Gaussian or Laplace tail bound, for each $0 \leq l \leq L_n$ and $0 \leq k < 2^{l}$, $\Pi(|r_{lk}| > n) \lesssim e^{-n}$. Hence, $\Pi(\mathcal{F}_n^c) \lesssim 2^{L_n} e^{-n} \leq \exp(-(C_1+4) n\epsilon_n^2)$ for a sufficiently large $C_1$. 

Next, to verify (\ref{cond-2}). From Lemma \ref{bounding-KL}, $K(f_{\eta_0}, f_\eta) \lesssim h^2(f_{\eta_0}, f_\eta)$ and $V(f_{\eta_0}, f_\eta) \lesssim h^2(f_{\eta_0}, f_\eta)$, thus, 
$\Pi(B_{\KL}(\eta_0, \epsilon_n)) \geq \Pi(\eta: h^2(f_{\eta_0}, f_\eta) \lesssim \epsilon_n^{{2}})$.
Moreover, from Lemma \ref{bounding-hellinger}, if $\|r - r_0\|_\Linfty + |(\theta - \theta_0)'z| \leq Q$, for $Q \leq 1/4$, then $h^2(f_{\eta_0}, f_\eta) \lesssim Q^{{2}}e^{2Q}$. Since $Q^{{2}}e^{2Q} {{\asymp}} \epsilon_n^{{2}}$ implies $Q {{\lesssim}} \epsilon_n$, we obtain the following lower bound:
\begin{align}
\label{cond-2-0.1}
\Pi(B_{\KL}(\eta_0, \epsilon_n)) \geq \Pi(\{\theta: \|\theta - \theta_0\|_1 \lesssim \epsilon_n\}) \times 
\Pi(\{r: \|r - r_0\|_\Linfty \lesssim \epsilon_n\}).
\end{align}
From the second paragraph in Section S-6.1 of \citet{cast21survivalSupp}, we immediately obtain
$\Pi(r: \|r - r_0\|_\Linfty \lesssim \epsilon_n) \gtrsim \exp(-C_1'n\epsilon_n^2)$ for some constant $C_1'$. 
To bound the first term in the last display,
denote $\vartheta_j = \theta_j - \theta_{0,j}$ and change variables from $\theta_j$ to $\vartheta_j$ for $\forall j \in \{1, \dots, p\}$, we have
\begin{align}
    \Pi(\theta: \|\theta - \theta_0\|_1 \leq c\epsilon_n)
    & \geq \prod_{j=1}^p \Pi(|\theta_j - \theta_{0,j}| \leq c\epsilon_n/p)
    \gtrsim \prod_{j=1}^p \Pi(|\vartheta_j| \leq c\epsilon_n/p).
    \label{cond-2-1}
\end{align}
The second inequality is obtained by using the fact that $\|\theta_{0}\|_\infty$ is bounded in \ref{asp:ii} and $p$ is a fixed constant.
We can further obtain $\Pi(|\vartheta_j| \leq c\epsilon_n/p) \geq b_1\epsilon_n/p$ for some constant $b_1$.
Therefore,  
$$
\Pi(B_{\KL}(\eta_0, \epsilon_n))
\gtrsim
\left(
\frac{b_1 \epsilon_n}{p}
\right)^p
\times e^{- C_1' n\epsilon_n^2} \geq e^{-C_1 n\epsilon_n^2},
$$
by choosing a sufficiently large $C_1 > C_1'$.

Last, we verify (\ref{cond-3}). Given that $\mathcal{F}_n$ as above and define $\mathcal{A}_n = \{r: \|r - r_0\|_\infty \leq \epsilon_n\}$, we have
$$
 \log N(\epsilon_n, \mathcal{F}_n, h) \leq \log N(\epsilon_n, \{\theta: \|\theta \|_\infty \geq C\}, \|\cdot\|_1) + \log N(\epsilon_n, \mathcal{A}_n, \|\cdot\|_2).
$$
where $\|\cdot\|_1$ and $\|\cdot\|_2$ stand for the $\ell_1$- and $\ell_2$-norm of a vector respectively.
Again by following the argument in the third paragraph of Section S-6.1 of \citet{cast21survivalSupp} for bounding the entropy, the second term in the last display is bounded by $C_3n\epsilon_n^2$.
By Proposition C.2 of \citet{ghosal17}, the first term in the last display can be bounded by 
\begin{align*}
\log N(\epsilon_n, \{\theta \in \mathbb{R}^p: \|\theta \|_\infty \geq C\}, \|\cdot\|_1) 
& \leq 
\log N(\epsilon_n, \{\theta \in \mathbb{R}^p: \|\theta\|_1 \geq pC\}, \|\cdot\|_1) \\
& \leq p \log \left(
\frac{3pC}{\epsilon_n}
\right) \leq  C_3 n \epsilon_n^2,
\end{align*}
for some constant $C_3$.
Therefore, by combining the upper bounds for 
$\log N(\epsilon_n, \{\theta: \|\theta \|_\infty \geq C\}, \|\cdot\|_1)$
and $\log N(\epsilon_n, \mathcal{A}_n, \|\cdot\|_2)$ derived above,
we obtain (\ref{cond-3}). 
We thus verified all three conditions.

We now consider the use of the random histogram prior {\ref{prior:H}}.
We choose the dependent Gamma prior on $(\lambda_k)$. 
Proving the independent Gamma prior case is simpler, thus we omit further details for brevity.  
Introducing the set $\mathcal{A}_n' = \{r: |r_k| \leq n^2, \ 0 \leq k < 2^{L_n+1}\}$, and then define 
$\mathcal{F}'_n = \{(\theta, r): \|\theta\|_\infty \leq C,  \ |r_{k}| \leq n^2, \ \forall k\}$, which is similar to $\mathcal{F}_n$ with $(r_{lk})$ is replaced by the histogram heights $(r_k)$. 
Then 
$\Pi(\mathcal{F}_n^c) \leq \Pi(\|\theta\|_\infty >C ) + \Pi(r: |r_k| > n^2)$.
We only need to bound the second term as the first term is 0 since the prior of $\theta$ is truncated as $C$. 
By following the proof in Section S-6.2 of \citet{cast21survivalSupp}, the second term is bounded by $e^{-n^2} + 2^{L_n+1} \exp(-2^{L_n})$ for $L_n$ is chosen as in (\ref{choicel}), this term is smaller than $\exp(-(C_1+4)n\epsilon_n^2)$ for a sufficiently large $C_1$. 

To verify (\ref{cond-2}), from (\ref{cond-2-0.1}) and (\ref{cond-2-1}), what left is to lower bound the prior probability $\Pi(r: \|r - r_0\|_\infty \lesssim \epsilon_n)$, which is bounded below by 
$\exp(2^{L_n+1} \log \tilde \epsilon_n)$, where 
$$
\tilde \epsilon_n = \frac{\epsilon_n \alpha^\alpha}{\Gamma(\alpha)}
\exp \left(\alpha
\left[
D2^{-(L_n+1) \beta }+ \epsilon_n - \exp(D2^{-(L_n + 1)\beta} + \epsilon_n)
\right]
\right),
$$
which is approximately $\epsilon_n \alpha^{\alpha} e^{-\alpha}/\Gamma(\alpha)$ as 
 $\alpha, D$ are constants, $L_n$ in (\ref{choicel}), and $\epsilon_n = o(1)$.
Since, $\exp(2^{L_n+1} \log \tilde \epsilon_n) \geq \exp(-C_1 n\epsilon_n^2)$ for a sufficiently large $C_1$. We thus verified (\ref{cond-2}).

To verify (\ref{cond-3}), we only need to bound $\log N(\epsilon_n, \mathcal{A}_n', \|\cdot\|_2)$. This quantity is bounded by $cL_n \log L_n \leq C_3 n \epsilon_n^2$, for some constant $C_3, c$, by following the same argument as on Page 26 of \citet{cast21survivalSupp}. 

We have verified all three conditions for the Haar wavelet prior as well as the random histogram prior. For both priors, the same Hellinger rate $\epsilon_n$ is obtained. By choosing $\epsilon_n = {{\nu}}_n$, the same result holds. 

\subsection{Verifying {\rm \ref{cond:P}}}
\label{sec:boundparameters}
\rm 
Recall that the squared Hellinger distance between $f_{\eta}$ and $f_{\eta_0}$ is given by
\begin{align}
    h^2(f_{\eta}, f_{\eta_0}) 
    & = 
    \int \int_0^1 
    \left[
        \sqrt{S_\eta} - \sqrt{S_{\eta_0}}
    \right]^2(u,z) g_{z}(u) f_Z(z) du dz
    \label{hellinger-1}
    \\
    & \quad +
    \int \int_0^1 \bar{G}_{z} (u) 
    \left[
        \sqrt{\lambda S_\eta e^{\theta' z}}
        - 
        \sqrt{\lambda_0 S_{\eta_0} e^{\theta_0' z}}
    \right]^2 (u, z) f_Z(z) du dz
    \label{hellinger-2}\\
    & \quad + 
    \int \bar{G}_{z} (1) 
    \left[
        \sqrt{S_\eta} - \sqrt{S_{\eta_0}}
    \right]^2(1, z) f_Z(z) dz,
    \label{hellinger-3}
\end{align}
where $S_\eta(u,z) = \exp(-\Ld(u) e^{\theta'z})$.
Denote the $L_1$ distance between the two functions 
$s_{\eta_1} := {s_{\eta_1}(u, z)}$ and $s_{\eta_2} := {s_{\eta_2}(u, z)}$ as
$$
H_1\left({s_{\eta_1}}, {s_{\eta_2}}\right)
= 
\int \int_0^1 
\left|
\sqrt{s_{\eta_1}} - \sqrt{s_{\eta_2}}
\right|(u, z) du dF_Z(z),
$$
where $dF_Z(z) = f_Z(z) dz$,
and the squared ``pseudo-Hellinger'' distance between the same two functions as
$$
H_2^2\left(s_{\eta_1}, s_{\eta_2}\right)
= 
\int \int_0^1 
\left(
\sqrt{s_{\eta_1}} - \sqrt{s_{\eta_2}}
\right)^2 (u, z) du dF_Z(z).
$$

\begin{lemma}
\label{lemma-1}
Suppose assumptions {\rm \ref{asp:i}-\ref{asp:v}} hold, if $h^2(f_\eta, f_{\eta_0}) \leq \epsilon_n^2$ and $\|\theta\|_\infty \leq C$, then there exist constant $C_1 > 0$ such that $\Ld(1)\leq C_1$. 
\end{lemma}
\begin{proof}
From the definition of $h^2(f_\eta, f_{\eta_0})$ and by \ref{asp:i} and \ref{asp:iv}, one can deduce that
\begin{align*}
\epsilon_n^2 \geq h^2(f_\eta, f_{\eta_0})
& \geq \int \bar G_z(1) \left(
e^{-\Lambda(1) e^{\theta'z}/2} - e^{-\Lambda_0(1) e^{\theta_0'z}/2}
\right)^2 f_Z(z) dz \\
& \gtrsim
\int_{z \in [-c_1, c_1]^p} 
\left(
e^{-\Lambda(1) e^{\theta'z}/2} - e^{-\Lambda_0(1) e^{\theta_0'z}/2}
\right)^2 f_Z(z) dz.
\end{align*}
By the mean-value theorem and the fact that $\int_{z \in [-c_1, c_1]^p} f_Z(z) = 1 > 0$, the last display implies that there exists a $z^\star \in [-c_1, c_1]^p$ such that 
$$
\left|\sqrt{e}^{-\Lambda(1) e^{\theta'z^\star}} - \sqrt{e}^{-\Lambda_0(1) e^{\theta_0'z^\star}}
\right| \lesssim \epsilon_n,
$$
which further implies that there exists a constant $d_1$ such that
$$
\sqrt{e}^{-\Lambda_0(1) e^{\theta_0'z^\star}} - d_1 \epsilon_n \leq 
\sqrt{e}^{-\Lambda(1) e^{\theta'z^\star}} \leq d_1 \epsilon_n + \sqrt{e}^{-\Lambda_0(1) e^{\theta_0'z^\star}}.
$$
Using assumptions \ref{asp:i}-\ref{asp:iii}, $\sqrt{e}^{-\Lambda_0(1) e^{\theta_0'z^\star}}$ is bounded both from above and in below by some constants, therefore, 
the last display implies that $\sqrt{e}^{-\Lambda(1) e^{\theta'z^\star}}$ is also bounded both from above and in below. 
Since $\|\theta\|_\infty \leq C$, we conclude that $\Ld(1)$ is bounded.
\end{proof}

\begin{lemma}
\label{lemma-2}
Suppose assumptions {\rm \ref{asp:i}-\ref{asp:iv}} hold, if $h^2(f_\eta, f_{\eta_0}) \leq \epsilon_n^2$ and 
$\|\theta\|_\infty \leq C$, then
$H_2(\lambda e^{\theta'z}, \lambda_0 e^{\theta_0'z}) \lesssim \epsilon_n$.
\end{lemma}

\begin{proof}
By the definition of the ``pseudo-Hellinger'' distance, 
\begin{align*}
    H_2^2(\lambda e^{\theta'z}, \lambda_0 e^{\theta_0'z}) 
    & = \int \int_0^1 \left(
    \sqrt{\lambda e^{\theta'z}} - \sqrt{\lambda_0 e^{\theta_0'z}}
    \right)^2(u, z) dudF_Z(z).
\end{align*}
Using the fact that $S_\eta(u,z) \geq S_\eta(1, z)$ and by Lemma \ref{lemma-1}, $\Ld(1) \leq C_1$, the last display can be bounded by 
\begin{align}    
    & \int \frac{1}{S_{\eta}(1, z)} 
    \int_0^1 
    \left(
      \sqrt{\lambda e^{\theta'z}S_{\eta}} - \sqrt{\lambda_0 e^{\theta'z}S_{\eta}} 
    \right)^2 (u, z) du dF_Z(z) \nonumber \\
    &\quad \leq e^{\Ld(1) e^{\max_z |\theta'z|}} \int \int_0^1 \left(
      \sqrt{\lambda e^{\theta'z}S_{\eta}} - \sqrt{\lambda_0 e^{\theta'z}S_{\eta}} 
    \right)^2 (u, z) du dF_Z(z),
    \label{pf-lemma1-1}
\end{align}   
By invoking \ref{asp:i}, we have $e^{\Ld(1) e^{\max_z|\theta'z|}} \leq e^{C_1 e^{pCc_1}}: = C_2$,
thus, (\ref{pf-lemma1-1}) is bounded by 
\begin{align}
\label{pf-lemma1-1.1}
C_2 \int \int_0^1 
    \left(
    \sqrt{\lambda e^{\theta'z} S_{\eta}} - \sqrt{\lambda_0 e^{\theta_0'z} S_{\eta}}
    \right)^2 (u, z) du dF_Z(z).
\end{align}
Applying the inequality $(a+b)^2 \leq 2 a^2 + 2b^2$,
(\ref{pf-lemma1-1.1}) can be further bounded by
\begin{align}
  C_2  \int \int_0^1 
    \left[   
    \left(
    \sqrt{\lambda e^{\theta'z} S_{\eta}} - \sqrt{\lambda_0 e^{\theta_0'z} S_{\eta_0}}
    \right)^2 
    + 
    \left(
    \sqrt{\lambda_0 e^{\theta_0'z} S_{\eta_0}} - \sqrt{\lambda_0 e^{\theta_0'z} S_{\eta}}
    \right)^2
    \right] (u, z) du dF_Z(z)
    \label{pf-lemma1-2}
\end{align}
To bound the last display, 
first, using that $h^2(f_{\eta}, f_{\eta_0}) \leq \epsilon_n^2$ and \ref{asp:iv}, 
we have (\ref{hellinger-2}) $\leq \epsilon_n^2$.
Thus, the first term in the sum of (\ref{pf-lemma1-2}) is bounded by $\epsilon_n^2$ up to some constant.
Next, by \ref{asp:i}-\ref{asp:iv}, $\lambda_0$, $\theta_0'z$, $g_z(u)$ are all bounded.
The second term in the sum of (\ref{pf-lemma1-2}) is bounded by $\int \int_0^1 [\sqrt{S_{\eta_0}} - \sqrt{S_\eta}]^2(u,z) du f(z) dz$.
Since $h^2(f_\eta, \eta_0) \lesssim \epsilon_n^2$, (\ref{hellinger-1}) is $\lesssim \epsilon_n^2$,
which implies that $\int\int_0^1 \left[\sqrt{S_\eta} - \sqrt{S_{\eta_0}} \right]^2 (u,z) du f(z) dz \lesssim \epsilon_n^2$
Therefore, we obtain $H_2^2(\ld e^{\theta'z}, \ld_0 e^{\theta_0'z}) \lesssim \epsilon_n^2$. 
\end{proof}

\begin{lemma}
\label{lemma-3}
Suppose assumptions {\rm \ref{asp:i}-\ref{asp:iv}} hold, if $h^2(f_{\eta}, f_{\eta_0}) \leq \epsilon_n^2$, and $\|\theta\| \leq C$, define $\bar{\lambda} = \lambda/\Lambda(1)$ and $\bar{\lambda}_0 = \lambda_0/\Lambda_0(1)$, then
$H_1({\bar{\lambda}}, {\bar{\lambda}_0}) \lesssim \epsilon_n$ and 
$\int \int_0^1 |\lambda e^{\theta'z} - \lambda_0 e^{\theta_0'z}| du dF_Z(z) \lesssim \epsilon_n$.
\end{lemma}

\begin{proof}
Recall the definition of 
the $L_1$ distance, we have
\begin{align*}
H_1({\bar{\lambda}}, {\bar{\lambda}_0}) 
& = \int \int_0^1 \left|
\sqrt{\frac{\lambda}{\Lambda(1)}} - 
\sqrt{\frac{\lambda_0}{\Lambda_0(1)}}\right| du dF_Z(z)  
= 
\int \int_0^1 
\left|
\sqrt{\frac{\lambda e^{\theta'z}}{\Lambda(1) e^{\theta'z}}} - 
\sqrt{\frac{\lambda_0 e^{\theta_0'z}}{\Lambda_0(1) e^{\theta_0'z}}}
\right| du dF_Z(z).
\end{align*}
Applying the triangle inequality, the last display is bounded by 
\begin{align}
    & \int \int_0^1 
    \left[
    \left|
    \sqrt{\frac{\lambda e^{\theta'z}}{\Lambda(1) e^{\theta'z}}}
     - 
    \sqrt{\frac{\lambda e^{\theta'z}}{\Lambda_0(1) e^{\theta_0'z}}}
    \right|
    + 
     \left|
    \sqrt{\frac{\lambda e^{\theta'z}}{\Lambda_0(1) e^{\theta_0'z}}} - 
    \sqrt{\frac{\lambda_0 e^{\theta_0'z}}{\Lambda_0(1) e^{\theta_0'z}}}
    \right|
        \right]
    %\right]
    du dF_Z(z) 
    \nonumber \\
    & \quad = 
    \int \sqrt{\Lambda(1) e^{\theta'z}}
    \left|
      (\Lambda(1) e^{\theta'z})^{-1/2} - (\Lambda_0(1) e^{\theta_0'z})^{-1/2}
    \right| dF_Z(z) 
    \label{pf-lemma2-1}\\
   &\qquad + \int \frac{1}{\Ld_0(1) e^{\theta_0'z}} 
   \int_0^1 
   \left|
   \sqrt{\lambda e^{\theta'z}} - \sqrt{\ld_0 e^{\theta_0'z}}
   \right| du dF_z(z)
      \label{pf-lemma2-2}
\end{align}
First, we bound (\ref{pf-lemma2-2}). 
Using \ref{asp:i}-\ref{asp:iii}, (\ref{pf-lemma2-2}) can be bounded by 
$$
\frac{1}{\Ld_0(1) e^{\min_z (\theta_0'z)}} \int 
\int_0^1 
\left|
\sqrt{\lambda e^{\theta'z}} - \sqrt{\ld_0 e^{\theta_0'z}}
\right| du dF_z(z)
\lesssim H_1(\ld e^{\theta'z}, \ld_0 e^{\theta_0'z})
\leq H_2(\ld e^{\theta'z}, \ld_0 e^{\theta_0'z}),
$$
where $H_2(\cdot, \cdot)$ is the ``pseudo-Hellinger'' distance.
Then, by Lemma \ref{lemma-2},  the last display is bounded by $\epsilon_n$ up to some constant. 

Next, we bound (\ref{pf-lemma2-1}), which can be written as
\begin{align*}
   & \int \sqrt{\Lambda(1) e^{\theta'z}}
    \left|
      (\Lambda(1) e^{\theta'z})^{-1/2} - (\Lambda_0(1) e^{\theta_0'z})^{-1/2}
    \right| dF(z)  \\
    & \quad = \int \frac{1}{\sqrt{\Lambda_0(1) e^{\theta_0'z}}} 
    \left|
    \frac{
        \Lambda_0(1) e^{\theta_0'z} - \Lambda(1) e^{\theta'z}
    }{
        \sqrt{\Lambda_0(1) e^{\theta_0'z}} + \sqrt{\Lambda(1) e^{\theta'z}}
    }
    \right| d F_Z(z).
\end{align*}
Using assumptions \ref{asp:i}-\ref{asp:iii}, the last display can be bounded by 
\begin{align*} 
    & \frac{1}{\Ld_0(1) e^{\min_z(\theta_0'z)}}
    \int \left|\Lambda_0(1) e^{\theta_0'z} - \Lambda(1) e^{\theta'z}\right| dF_Z(z) \\
    & \quad \lesssim
    \int \int_0^1 |\lambda e^{\theta'z} - \lambda_0 e^{\theta_0'z}| du  dF_Z(z)\\
    &  \quad = \int \int_0^1 \left|\sqrt{\lambda e^{\theta'z}} - \sqrt{\lambda_0 e^{\theta_0'z}}\right|
       \left|\sqrt{\lambda e^{\theta'z}} + \sqrt{\lambda_0 e^{\theta_0'z}}\right|
       du  dF_Z(z).
  \end{align*}
By applying the Cauchy-Schwartz inequality, the last display be bounded by 
\begin{align}       
    \left(\int 
     \int_0^1 
    \left(\sqrt{\lambda e^{\theta'z}} - \sqrt{\lambda_0 e^{\theta_0'z}}
    \right)^2 du dF_Z(z) \right)^{1/2} 
    \left( \int \int_0^1 \left(\sqrt{\lambda e^{\theta'z}} + \sqrt{\lambda_0 e^{\theta_0'z}}\right)^2 du
     dF_Z(z) \right)^{1/2}
    \label{pf-lemma2-3}
\end{align}
By Lemma \ref{lemma-2}, the first term in the product of (\ref{pf-lemma2-3}) is $\lesssim \epsilon_n$.
Using the inequality $(a+b)^2 \leq 2a^2 + 2b^2$, the second term of (\ref{pf-lemma2-3}) is bounded by 
\begin{align*}
2 \int \int_0^1 \lambda e^{\theta'z} du dF_Z(z)
+  2 \int \int_0^1 \lambda_0 e^{\theta_0'z} du dF_Z(z).
\end{align*}
Using \ref{asp:i}-\ref{asp:iii} and 
by Lemma \ref{lemma-1}, the last display is bounded by a constant. 
Therefore, (\ref{pf-lemma2-3}) is $\epsilon_n$ times some constant.
By combining the upper bounds for (\ref{pf-lemma2-1}) and (\ref{pf-lemma2-2}), we obtain  $H_1({\bar{\lambda}}, {\bar{\lambda}_0}) \lesssim \epsilon_n$.

To bound 
$
\int \int_0^1 \left|\ld e^{\theta'z} - \ld_0 e^{\theta_0'z}\right| du dF_Z(z),
$ 
one can write it as 
$$
\int \int_0^1 
\left| \left(\sqrt{\lambda e^{\theta'z}} - \sqrt{\lambda_0 e^{\theta_0'z}} \right)
\left(\sqrt{\lambda e^{\theta'z}} + \sqrt{\lambda_0 e^{\theta_0'z}}\right)\right| du dF_Z(z).
$$ 
Then, by the Cauchy-Schwarz inequality, 
the last display is bounded by (\ref{pf-lemma2-3}). The remaining proof is the same as above. 
\end{proof}

%%%%%%%%%%%%%%%%%%%%%%%%%%%%%%%%%%%%%%%%%%%%

\begin{lemma}
\label{lemma-4}
Suppose $r_0 \in \mathcal{H}(\beta, D)$ with $\beta > 1/2$ and $D > 0$ and assumptions {\rm \ref{asp:i}, \ref{asp:ii}, \ref{asp:iii}}, and {\rm \ref{asp:v}} hold. If $h^2(f_{\eta_0}, f_\eta) \lesssim \epsilon_n^2$ and $\|\theta\|_\infty \leq C$, then 
\begin{align*}
    & \|\theta - \theta_0\|^2 + \|r - r_0\|_\Ltwo^2 \lesssim 
    \|\theta - \theta_0, r - r_0\|_L^2 \lesssim \epsilon_n^2, \\
    & \|r - r_0\|_\Linfty \lesssim 2^{L_n/2} \epsilon_n + 2^{-\beta L_n} = o(1).
\end{align*}
In particular, if choosing $L_n$ such that $2^{L_n} = (n/\log n)^{1/(2\beta + 1)}$, then 
$\|r - r_0\|_\Linfty \lesssim 2^{L_n/2} \epsilon_n$.
\end{lemma}

\begin{proof}
By following the proof of Lemma 10 of \citet{cast12}, we have
\begin{align}
\label{lemma4-1}
\|\theta - \theta_0, r - r_0\|_L^2 \lesssim \int f_{\eta_0} \log^2\left(\frac{f_{\eta_0}}{f_\eta}\right) .
\end{align}
By invoking Lemma 8 of \citet{ghosal07DM}, the last display can be further bounded by  
\begin{align}
\label{lemma4-2}
h^2(f_{\eta_0}, f_\eta) (1 + \log \|f_{\eta_0}/f_\eta\|_\Linfty)^2.
\end{align}
Simply calculations reveals that 
$
\log (f_0/f) = \delta \left( (r - r_0) + (\theta - \theta_0)'z \right) - \Lambda e^{\theta'z} + \Lambda_0 e^{\theta_0'z}
$ and thus,
\begin{align}
\label{lemma4-3}
\log \|f_0/f\|_\Linfty \leq \| \log (f_0/f)\|_\Linfty 
\lesssim 
\|r - r_0\|_\Linfty + \|\theta - \theta_0\|_1 \|z\|_\infty + \max_z\|\Lambda e^{\theta'z} - \Lambda_0 e^{\theta_0'z}\|_\Linfty. 
\end{align}
By applying the triangle inequality, 
since $\Ld_0$, $z$, and $\theta_0$ are all bounded quantities from assumptions \ref{asp:i}, \ref{asp:ii}, \ref{asp:iii},
the third term on the right hand side of (\ref{lemma4-3}) can be bounded by 
\begin{align*}
\max_z \|(\Ld - \Ld_0) e^{\theta'z}\|_\Linfty + \max_z\|\Ld_0 (e^{\theta'z} - e^{\theta_0'z})\|_\Linfty
\lesssim \|\Ld - \Ld_0\|_\Linfty + \|\theta - \theta_0\|_1 \|z\|_\infty.
\end{align*}
One can further bound $\|\Ld - \Ld_0\|_\Linfty \leq \|\ld - \ld_0\|_\Lone = \int_0^1 e^{r_0}|e^{r - r_0} - 1| \lesssim \int_0^1 |e^{r-r_0}-1| \leq \|r - r_0\|_\Lone \leq \|r - r_0\|_\Ltwo$.
Thus, 
by plugging the above upper bound back into (\ref{lemma4-3}) and then into (\ref{lemma4-2}), 
also, using the inequality $(a+b)^2 \leq a^2 + b^2$, we obtain
\begin{align}
\label{lemma4-4}
\|\theta - \theta_0, r - r_0\|_L^2  \lesssim \epsilon_n^2 (1 + \|r - r_0\|_\Ltwo^2 + \|\theta - \theta_0\|_1^2).
\end{align}
On the other hand, using \ref{asp:i}-\ref{asp:iii} and \ref{asp:v}, one obtains 
\begin{align}
\label{lemma4-5}
\|\theta - \theta_0, r-r_0\|_L^2 \gtrsim \|\theta - \theta_0\|^2 + \|r - r_0\|_\Ltwo^2. 
\end{align}
Due to $\|\theta - \theta_0\|_1 \leq \sqrt{p}\|\theta - \theta_0\|$, $p$ a fixed constant, 
by combing (\ref{lemma4-4}) with (\ref{lemma4-5}), we have 
$$
\|\theta - \theta_0\|^2 + \|r - r_0\|_\Ltwo^2 \lesssim \epsilon_n^2 (1 + \|\theta - \theta_0\|^2 + \|r - r_0\|_\Ltwo^2),
$$
which implies $\|\theta - \theta_0\|^2 + \|r - r_0\|_\Ltwo^2 \lesssim \epsilon_n^2$. 

To obtain the upper bound for $\|r - r_0\|_\Linfty$, using that $\|r - r_0\|_\Ltwo^2 \lesssim \epsilon_n^2$, one can directly apply the proof of Lemma S-7 in \citet{cast21survivalSupp} to obtain the upper bound.  Note that when using their argument, we let $\gamma$ in their proof equals to $\beta$. 

If $2^{L_n} = (n/\log n)^{1/(2\beta + 1)}$, note that $\epsilon_n = {{\nu}}_n$ in (\ref{rateveps}), then  $2^{L_n/2}\epsilon_n = (n/\log n)^{\frac{1-2\beta}{2(2\beta + 1)}}$, which is larger than $2^{-\beta L_n}$. We thus obtain the result. Note that since we assume $\beta > 1/2$, $2^{L_n/2}\epsilon_n = o(1)$ still holds.  
\end{proof}

%%%%%%%%%%%%%%%%%%%%%%%%%%%%%%%%%%%%%%%%%%%%
\begin{lemma}
\label{lemma-5}
If assumption {\rm \ref{asp:iv}} and the same conditions as in Lemma \ref{lemma-4} hold, then
$\|\ld - \ld_0\|_\Lone \lesssim \epsilon_n$. 
\end{lemma}

\begin{proof}
By Lemma \ref{lemma-3}, we have
$
\int \int_0^1 |\ld e^{\theta'z} - \ld_0 e^{\theta_0'z}|du dF_Z(z) \lesssim \epsilon_n. 
$
With \ref{asp:i}-\ref{asp:iii}, this inequality implies that 
$$
\int \int_0^1 |e^{r - r_0 + (\theta - \theta_0)'z} - 1| du dF_Z(z)
\lesssim 
\int \int_0^1 e^{r_0 + \theta_0'z} |e^{r - r_0 + (\theta - \theta_0)'z} - 1| du dF_Z(z)
\lesssim \epsilon_n^2. 
$$
On the other hand,
by \ref{lemma-4}, $\|\theta - \theta_0\| \lesssim \epsilon_n$ and $\|r - r_0\|_\Linfty = o(1)$, thus, the last display implies $\int \int_0^1 |(\theta - \theta_0)'z + (r - r_0)| du dF_Z(z) \lesssim \epsilon_n$. 
Since $Z$ is bounded by assumption \ref{asp:i}, then using that $\|\theta - \theta_0\| \lesssim \epsilon_n$, we have $\|r - r_0\|_\Lone \lesssim \epsilon_n$. 
To show $\|\ld - \ld_0\|_\Lone \lesssim \epsilon_n$, we first write 
$\ld - \ld_0 = e^{r_0}(e^{r - r_0} - 1)$. Since $\ld_0$ is bounded by \ref{asp:iii} and $\|r - r_0\|_\Linfty = o(1)$, we then apply Taylor's theorem and obtain that $\|\ld - \ld_0\|_\Lone \lesssim \|r - r_0\|_\Lone \lesssim \epsilon_n$. 
\end{proof}

%%%%%%%%%%%%%%%%%%%%%%%%%%%%%%%%%%%%%%%%%%%%%%%%%
%\section{Proof of Theorem \ref{thm:bvm-linear-fns}} 
%%\re{why the plural in proofs? (also above for Lemma 1)}
%\label{proof:thm1}

%%%%%%%%%%%%%%%%%%%%%%%%%%%%%%%%%%%%%%%%%%%%%%%%%
\section{Supporting lemmata for Theorem \ref{thm:bvm-linear-fns}}
\label{sec:support-lemma-thm1}

We prove the intermediate steps in Section \ref{proof:thm1} for Theorem \ref{thm:bvm-linear-fns} in this section.
There are three subsections:
Section \ref{sec:bounding-LAN-norm} obtains bounds for 
$D_n$ in (\ref{Dn}) and $\sup_{\eta \in A_n} | B_2(\eta, \eta_0)|$ for $B_2(\eta, \eta_0)$ in (\ref{B2}),
Section \ref{sec:bounding-R1} provides an upper bound for $\sup_{\eta \in A_n}|R_{n,1}(\eta, \eta_0) - R_{n,1}(\eta_h, \eta_0)|$, $R_{n,1}(\eta, \eta_0)$ in (\ref{R1}), and Section \ref{sec:bounding-R2} derives an upper bound for $\sup_{\eta \in A_n} |R_{n,2}(\eta, \eta_0) - R_{n,2}(\eta_h, \eta_0) - s \rn B_3(\eta, \eta_0)|$, where $R_{n,2}(\eta,\eta_0)$ and $B_3(\eta, \eta_0)$ are given in (\ref{R2}) and (\ref{B3}) respectively. 
%
%The first one is about bounding $D_n$ and the rest two lemmata are on controlling the two remainder terms given in (\ref{R1}) and (\ref{R2}).   
%Before giving the lemmas and their proofs, recall that 
%\begin{align*}
%& \qquad \Delta_1 := \Delta_1(a, b, h) = t\tIinv a + s \tIinv \Ld_0 \{b\gamm\},\\
%%\label{delta1}\\
%\Delta_{2,L_n} &: = \Delta_{2,L_n}(a,b,h) = t \gammln \tIinv a - s\gambln - s\gammln' \tIinv \Ld_0\{b \gamm\}.
%%\label{delta2}
%\end{align*}
%It is easy to check that $\theta_h = \theta + \Delta_1/\rn$ and $r_h = r+ \Delta_{2,L_n}/\rn$, where $\theta_h$ and $r_h$ are defined in (\ref{path-theta}) and (\ref{path-r}).

\subsection{Bounding $\sup_{\eta \in A_n} |B_2(\eta, \eta_0)|$ and $D_n$}
\label{sec:bounding-LAN-norm}

\begin{lemma}
\label{lemma:bounding-B2}
Suppose assumptions {\rm \ref{cond:P}}, {\rm \ref{asp:i}, \ref{asp:iii}, {\it and} \ref{asp:v}} hold, then 
\begin{align}
%& \sup_{\eta \in A_n} |B_2(\eta, \eta_0)| 
%\lesssim 
%p^2 \epsilon_n^2 2^{-L_n} + \epsilon_n^2 \|\gamb - \gambln\|_\Linfty  + \epsilon_n^2 \|\gamb\|_\Linfty
%\label{bound-B2}\\
&\sup_{\eta \in A_n} |B_2(\eta, \eta_0)| 
\lesssim 
\epsilon_n \|\gamb - \gambln\|_\Linfty + p^2 \epsilon_n 2^{-L_n} \|b\|_\Lone
\label{bound-B2}
\end{align}
where 
%$$
%B_2(\eta, \eta_0) = 
%\left\langle
%\left(0, r - r_0 - \frac{\ld - \ld_0}{\ld_0}\right),
%\left(- \tIinv \Ld_0\{b\gamm\}, \gambln + \gammln' \tIinv \Ld_0\{b \gamm\} \right)
%\right\rangle_L
%$$
%and 
$$
B_2(\eta, \eta_0) = 
\left\langle
\left(0, \frac{\ld - \ld_0}{\ld_0}\right),
\left(0, \gamb - \gambln + (\gamm - \gammln)' \tIinv \Ld_0\{b \gamm\} \right)
\right\rangle_L.
$$
\end{lemma}
\begin{proof}

By the definition of the LAN-norm Hilbert inner product in ({\ref{lan-norm-inner-product}}), $B_2(\eta, \eta_0)$ can be re-written as
\begin{align*}
& B_2(\eta, \eta_0) 
= \Ld_0\left\{
\frac{\ld - \ld_0}{\ld_0} \left(
\gamb - \gambln + (\gamm - \gammln)' \tIinv \Ld_0\{b\gamm\}
\right) M_0
\right\}. 
\end{align*}
Applying the inequalities $\Ld_0\{f(\cdot) g(\cdot)\} \leq \|f(\cdot)\|_\Lone \|g(\cdot)\|_\Linfty \|\Ld_0\|_\Linfty$ and $\|f + g\|_\Linfty \leq \|f\|_\Linfty + \|g\|_\Linfty$ for any $f, g \in L^2\{\Ld_0\}$, the last display is bounded by 
\begin{align}
\|\ld - \ld_0\|_\Lone \|M_0\|_\Linfty
\left(
\|\gamb - \gambln\|_\Linfty + 
\|(\gamm - \gammln)' \tIinv \Ld_0\{b\gamm\}\|_\Linfty
\right).
\label{pf-thm1-lem1-3}
\end{align}
Moreover,
\begin{align*}
\|(\gamm - \gammln)' \tIinv \Ld_0\{b\gamm\}\|_\Linfty
\leq p^2 \max_j \|\gamm^j - \gammln^j\|_\Linfty \|\tIinv\|_{(\infty, \infty)} 
 \|\Ld_0 \{b \gamm\}\|_\infty.
\end{align*}
To bound the last display, 
first, using the third point of Lemma \ref{lemma:psib-2}, then
$\max_j \|\gamm^j - \gammln^j\|_\Linfty \lesssim 2^{-L_n}$;
next, by \ref{asp:v}, $\|\tIinv\|_{(\infty, \infty)}$ is bounded by a constant;
last, by \ref{asp:i} and \ref{asp:iii},
$\|\Ld_0\{b \gamm\}\|_\infty  \leq \|\ld_0\|_\Linfty \|b\|_\Lone \|\gamm\|_\infty
\leq \|\ld_0\|_\Linfty \|b\|_\Lone \|z\|_\infty
\lesssim \|b\|_\Lone$.
Thus,
the last display is bounded by $C_3p^2 2^{-L_n} \|b\|_\Lone$ for some constant $C_3$,
and hence (\ref{pf-thm1-lem1-3}) is bounded by 
$$
\|\ld - \ld_0\|_1\|M_0\|_\infty \left(\|\gamb - \gambln\|_\infty + C_3p^2 2^{-L_n}\|b\|_\Lone \right).
$$
Using \ref{cond:P}, $\eta \in A_n$, and the fact that $M_0(\cdot)$ is bounded, we obtain (\ref{bound-B2}).
\end{proof}

\begin{lemma}
\label{lemma:bounding-Dn}
Suppose assumptions {\rm \ref{asp:i}, \ref{asp:iii}, \ref{asp:v}}, and  {\rm \ref{cond:P}} hold, define 
\begin{align*}
D_n = \frac{n}{2} \|\theta - \theta_h, r-r_h\|_L^2 - \frac{h\Sigma_{a,b} h'}{2},
\end{align*}
where the LAN-norm
$\|\cdot, \cdot\|_L$ is given in (\ref{LAN-norm}), 
$\theta_h$ and $r_h$ are given in (\ref{path-theta}) and (\ref{path-r}),
$h=(t,s)$, and 
\begin{equation*}
\Sigma_{a,b} 
= 
\begin{pmatrix}
a' \tIinv a & - a' \tIinv \Ld_0\{b\gamm\} \\
- a' \tIinv \Ld_0\{b\gamm\} & \Ld_0\{b\gamb\} + \Ld_0\{b\gamm'\} \tIinv \Ld_0\{b\gamm\}
\end{pmatrix},
\end{equation*}
then,
$
|D_n| 
\lesssim 
s^2 \|\gamb^2-\gambln^2\|_\Lone + (t^2 + s^2) (p^2 2^{-L_n} \|\gambln\|_\Lone + p^4 2^{-2L_n}).
$
\end{lemma}

\begin{proof}
By (\ref{LAN-norm-v2}), one can write the squared LAN-norm as $\|\theta - \theta_h, r - r_h\|_L^2 = (\theta - \theta_h)' \tilde I_{\eta_0} (\theta - \theta_h) + \|0, r - r_h + \gamm'(\theta - \theta_h)\|_L^2$,
where
\begin{align*}
\theta_h & = \theta - \frac{t \tIinv a}{\rn} + \frac{s \tIinv \Ld_0\{b \gamm\}}{\rn},\\
r_h = r + & \frac{t \gammln' \tIinv a}{\rn} - \frac{s \gambln}{\rn} - \frac{s \gammln' \tIinv \Ld_0\{b\gamm\}}{\rn}.
\end{align*}
By plugging-in the expressions of $\theta_h$ and $r_h$, we obtain
$$
n(\theta - \theta_h)' \tilde I_{\eta_0} (\theta - \theta_h) 
= t^2 a' \tIinv a - 2ts a' \tIinv \Ld_0\{b \gamm\} + s^2 \Ld_0\{b \gamm'\} \tIinv \Ld_0\{b \gamm\}
$$
and
\begin{align*}
&n \|0, r - r_h + \gamm'(\theta - \theta_h)\|_L^2 \\
& \quad = 
\Ld_0\left\{
\left[t (\gamm - \gammln)' \tIinv a + s\gambln - s(\gamm - \gammln)' \tIinv \Ld_0\{b\gamm\} \right]^2 M_0
\right\}.
\end{align*}

One the other hand, by plugging-in the expression of $\Sigma_{a,b}$, we can write the second term in $D_n$ as
$$
h'\Sigma_{a,b}h/2
= 
t^2 a' \tIinv a/2 - ts a' \tIinv \Ld_0\{b\gamm\} + s^2 \Ld_0\{b\gamb\}/2 + s^2 \Ld_0\{b\gamm'\} \tIinv \Ld_0\{b\gamm\}/2.
$$
By collecting all the relevant terms and letting $\Delta_n = (\gamm - \gammln)' \tIinv (ta - s\Ld_0\{b\gamm\})$, we obtain
\begin{align}
\label{pf-thm1-lem2-1}
D_n = 
\frac{n}{2} \| \theta - \theta_h, r - r_h\|_L^2 - \frac{h \Sigma_{a,b} h'}{2}
= \Ld_0\{ [(\Delta_n + s\gambln)^2 - s^2 \gamb^2] M_0 \}.
\end{align}
The last display is bounded by 
\begin{align*}
\|M_0\|_\Linfty \|(\Delta_n + s\gambln)^2 - s^2 \gamb^2\|_\Lone \|\ld_0\|_\Linfty
\lesssim \|(\Delta_n + s\gambln)^2 - s^2 \gamb^2\|_\Lone,
\end{align*}
where we used the inequality $\Ld_0\{f(\cdot) g(\cdot)\} \leq \|f(\cdot)\|_\Linfty \|g(\cdot)\|_\Linfty \|\ld_0\|_\Lone$ for any $f, g \in L^2\{\Ld_0\}$, Lemma \ref{lemma:M0M1}, and assumption \ref{asp:iii}.
Applying the triangle inequality for the supremum metric, then the upper bound in the last display can be further bounded by 
\begin{align}
\label{pf-thm1-lem2-2}
s^2\|\gamb^2 - \gambln^2\|_\Lone + 2\|s \gambln \Delta_n\|_\Lone + \|\Delta_n\|_\Lone^2. 
\end{align}
%The first term in (\ref{lem1-thm1-6}) is bounded by 
%$s^2 \|\gamb - \gambln\|_\Linfty \|\gamb + \gambln\|_\Linfty \leq 2 s^2 \|\gamb - \gambln\|_\Linfty \|b\|_\Linfty$ as $\|\gambln\|_\Linfty \leq \|\gamb\|_\Linfty \lesssim \|b\|_\Linfty$. 
Using the third point of Lemma \ref{lemma:psib-2}, \ref{asp:i}, and \ref{asp:v}, we have
\begin{align*}
\|\Delta_n\|_\Linfty 
& \leq p^2 \max_j\|\gamm^j - \gammln^j\|_\Linfty \|\tIinv\|_{(\infty, \infty)} \|ta - s\Ld_0\{b\gamm\}\|_\infty\\
& \lesssim p^2 2^{-L_n} (|t| + |s|).
\end{align*}
Then, the second term in (\ref{pf-thm1-lem2-2}) is bounded by 
$2|s| \|\gambln\|_\Lone \|\Delta_n\|_\Linfty \lesssim p^2 (|ts| + s^2) 2^{-L_n} \|\gambln\|_\Lone$
and the third term is bounded by 
$\|\Delta_n\|_\infty^2 \lesssim p^4 (t^2 + s^2) 2^{-2L_n}$, as $\|\Delta_n\|_1 \leq \|\Delta \|_\infty$.
Thus, (\ref{pf-thm1-lem2-2}) is bounded by $C(s^2\|\gamb^2 - \gambln^2\|_\Lone + 
 p^2 (|ts| + s^2) 2^{-L_n} \|\gambln\|_\Lone + p^4 (t^2 + s^2) 2^{-2L_n})$ for some constant $C$.
Last, using the inequality $|ts| \leq t^2 + s^2$ to complete the proof.
\end{proof}

%%%%%%%%%%%%%%%%%%%%%%%%%%%%%%%%%%%%%%%%%%%%%%%%%%%
%%%%%%%%%%%%%%%%%%%%%%%%%%%%%%%%%%%%%%%%%%%%%%%%%%%
\subsection{Bounding $\sup_{\eta \in A_n}|R_{n,1}(\eta, \eta_0) - R_{n,1}(\eta_h, \eta_0)|$}
\label{sec:bounding-R1}

%Let $\mathcal{Q}_n$ be a class of functions such that for some constant $D > 0$ and sequences of $(\mu_n)$, $(v_n)$, and $(\epsilon_n)$ of positive real numbers, we define

From the definition of $R_{n,1}$ given in (\ref{R1}), let $g_n(\eta): = g_n(\eta)(y,z)$ for
$$
g_n(\eta)(y,z) = - \rn \left(e^{\theta'z}\Ld_0\{e^{r - r_0}\}(y) - e^{\theta_0'z} \Ld_0(y)
- e^{\theta_0'z} (\theta- \theta_0)'z \Ld_0(y) - e^{\theta_0'z}\Ld_0\{r - r_0\}(y) \right),
$$
we can write
$R_{n,1}(\eta, \eta_0) - R_{n,1}(\eta_h, \eta_0) = \mathbb{G}_n(g_n(\eta) - g_n(\eta_h))$. 
Furthermore, 
let $\Ld_h(\cdot) = \int_0^\cdot e^{r_h}$, then 
\begin{align*}
& g_n(\eta) - g_n(\eta_h)\\
&= \rn e^{\theta_0'z} 
e^{(\theta - \theta_0)'z}
\Ld_0
\left\{
e^{(r - r_0)(y)} 
\left[
e^{(\theta_h - \theta)'z + (r_h - r)(y)} - (\theta_h - \theta)'z - (r_h - r)(y) - 1
\right]
\right\} \\
& \quad + \rn e^{\theta_0'z} 
\Ld_0\left\{
\left[e^{(\theta - \theta_0)'z + (r - r_0)(y)} - 1\right]
\left((\theta_h - \theta)'z + (r_h - r)(y) \right)
\right\}
\end{align*}

Denote the set  
\begin{align}
\label{L_n}
\mathcal{L}_n = \{(\theta, \lambda): \theta \in \mathbb{R}^p, \ld \in L^{\infty}[0,1], \ \|\theta - \theta_0\| \leq \epsilon_n, 
\ \|\ld - \ld_0\|_\Lone \leq \epsilon_n\},
\end{align}
where $(\epsilon_n)$ is a sequence of positive real numbers which is typically reduce to the posterior rate in Hellinger distance. Note that $\mathcal{L}_n \subseteq A_n$.

For $\eta \in \mathcal{L}_n$, let's define the following two classes of functions:
\begin{align}
\mathcal{F}_{n,1} & = 
\left\{
\rn
\int_0^\cdot
e^{\theta'z + r(\cdot)} 
\left[
e^{(\theta_h - \theta)'z + (r_h - r)(\cdot)} - (\theta_h - \theta)'z - (r_h - r)(\cdot) - 1
\right], \ \eta \in \mathcal{L}_n \right\}, 
\label{F-n1}\\
\mathcal{F}_{n,2} 
& =
\left\{ \rn 
\int_0^\cdot 
e^{\theta_0'z + r_0(\cdot)}
\left(
\left[e^{(\theta - \theta_0)'z + (r-r_0)(\cdot)} - 1\right]
\left[(\theta_h - \theta)'z + (r_h - r)(\cdot) \right]
\right),\ \eta \in \mathcal{L}_n
\right\}.
\label{F-n2}
\end{align}
One can easily verify that for $f_{n,i} \in \mathcal{F}_{n,i}$, $i= 1,2$,
\begin{align}
\label{eqn:gn}
\mathbb{G}_n(g_n(\eta) - g_n(\eta_h)) =  \mathbb{G}_n(f_{n,1}) + \mathbb{G}(f_{n,2}).
\end{align}
In this way, the empirical process is decomposed into two parts, $\mathbb{G}_n(f_{n,1})$ and $\mathbb{G}_n(f_{n,2})$. 
In the next lemma, we derive upper bounds for both parts. 

\bigskip
\begin{lemma}
\label{R1-version1}
For $\mathcal{F}_{n,1}$ and $\mathcal{F}_{n,2}$ defined in (\ref{F-n1}) and (\ref{F-n2}) respectively
and $\mathcal{L}_n$ in (\ref{L_n}), define
\begin{align}
\Delta_1 := \Delta_1(a,b,h) & = - t \tIinv a + s \tIinv \Ld_0\{b\gamm\},
\label{K1}\\
\Delta_{2,L_n}:=\Delta_{2,L_n}(a,b,h) &= t \gammln' \tIinv a - s\gambln + s\gammln \tIinv \Ld_0\{b\gamm\},
\label{K2}
\end{align}
if $\|\Delta_1\|_\infty/\rn \leq d_1$ and $\|\Delta_{2,L_n}\|_\Linfty/\rn \leq d_2$, $d_1 + d_2 < 1$,
denote $\|\mathbb{G}_n\|_{\mathcal{F}_{n,i}} = \sup_{f_{n,i} \in \mathcal{F}_{n,i}} |\mathbb{G}_n(f_{n,i})|$
for $i = 1, 2$,
then,
\begin{align*}
& \mathbb{E}_{\eta_0}^\star \left[ \|\mathbb{G}_n\|_{\mathcal{F}_{n,1}} \right] 
\lesssim (\|\Delta_1\|_\infty + \|\Delta_{2,L_n}\|_\Linfty)^2/\rn,\\
& \mathbb{E}_{\eta_0}^\star \left[
\|\mathbb{G}_n\|_{\mathcal{F}_{n,2}} 
\right]
\lesssim \epsilon_n (\|\Delta_1\|_\infty + \|\Delta_{2,L_n}\|_\Linfty).
\end{align*}
\end{lemma}

\begin{proof}
For $\theta_h$ and $r_h$ given in (\ref{path-theta}) and (\ref{path-r}), we have $\theta_h - \theta = \Delta_1/\rn$ and $r_h - r = \Delta_{2,L_n}/\rn$. Let us further denote $\tilde \Delta_n = \Delta_1'z + \Delta_{2,L_n}$,
then $(\theta_h - \theta)'z + (r_h - r)  = \tilde \Delta_n$. 

To bound $\mathbb{E}_{\eta_0}^\star \left[ \|\mathbb{G}_n\|_{\mathcal{F}_{n,1}} \right]$,
let $f_{n,1} = f_{n,11}  + f_{n,22}$, where
\begin{align*}
f_{n,11} & = \rn e^{\theta_0'z} \int_0^\cdot  e^{r_0} e^{(\theta - \theta_0)'z} (e^{r-r_0}-1) (e^{\tilde \Delta_n/\rn} - \tilde\Delta_n/\rn - 1),\\
f_{n,12} & = \rn e^{\theta_0'z} \int_0^\cdot e^{r_0} e^{(\theta - \theta_0)'z} (e^{\tilde\Delta_n/\rn} - \tilde\Delta_n/\rn - 1).
\end{align*}
Let $\mathcal{F}_{n,11}$ and $\mathcal{F}_{n,12}$ be classes of functions such that $f_{n,11} \in \mathcal{F}_{n,11}$ and $f_{n,12} \in \mathcal{F}_{n,12}$, where
\begin{align*}
\mathcal{F}_{n,11} & = 
\left\{\rn e^{\theta_0'z} \int_0^\cdot e^{r_0} e^{(\theta - \theta_0)'z}  (e^{r-r_0}-1) (e^{\tilde\Delta_n/\rn} - \tilde\Delta_n/\rn - 1), \ \eta \in \mathcal{L}_n\right\}, \\
\mathcal{F}_{n,12} &= 
\left\{\rn e^{\theta_0'z}  \int_0^\cdot e^{r_0} e^{(\theta - \theta_0)'z} (e^{\tilde\Delta_n/\rn} - \tilde\Delta_n/\rn - 1), \ \eta \in \mathcal{L}_n\right\},
\end{align*}
then, we can further bound
\begin{align}
\label{gn-F1}
\mathbb{E}_{\eta_0}^\star \left[ \|\mathbb{G}_n\|_{\mathcal{F}_{n,1}} \right]
\leq 
\mathbb{E}_{\eta_0}^\star \left[ \|\mathbb{G}_n\|_{\mathcal{F}_{n,11}} \right]
+ \mathbb{E}_{\eta_0}^\star \left[ \|\mathbb{G}_n\|_{\mathcal{F}_{n,12}} \right].
\end{align}

To bound the first term in the last display, we use Lemma \ref{lemma-ep-3}. 
We shall check the conditions in Lemma \ref{lemma-ep-3} first.
For $\|\theta_1 - \theta_0\| \leq \epsilon_n$ and $\|\theta_2 - \theta_0\| \leq \epsilon_n$, $\theta_1, \theta_2 \in \mathcal{L}_n$, 
using the fact that $e^x$ is a Lipschitz continuous function for a bounded $x$, by \ref{asp:i} and \ref{asp:ii}, we have
$|e^{\theta_1 'z} - e^{\theta_2'z} |
= e^{\theta_0'z} |e^{(\theta_1 -\theta_0)'z} - e^{(\theta_2 - \theta_0)'z} |
 \lesssim \|z\| \|\theta_1 - \theta_2\| \lesssim \|\theta_1 - \theta_2\|$. 
Thus $|e^{\theta_1'z} - e^{\theta_2'z} | \lesssim \|\theta_1 - \theta_2\|$. 
Next, by applying Taylor's theorem and using \ref{asp:ii},
$e^{\theta_0'z} |e^{(\theta - \theta_0)'z}| \lesssim 1 + o(1)$. 
Thus the condition for $g_\theta$ part in Lemma \ref{lemma-ep-3} is satisfied. 

Next, we check the condition for $h$ part, here $h := h_{n,11}$, where 
$$
h_{n,11} = \rn \int_0^\cdot e^{r_0} (e^{r - r_0} - 1) (e^{\tilde \Delta_n/\rn} - \tilde \Delta_n/\rn - 1).
$$
From the last display, one immediately has $h_{n,11}(0) = 0$. 
One also needs derive an upper bound for $\|h_{n,11}\|_{BV}$.
As $\tilde \Delta_n/\rn \leq \|\Delta_1\|_\infty/\rn + \|\Delta_{2, L_n}\|_\Linfty/\rn \leq d_1 + d_2 < 1$ by assumption, using $\lambda \in \mathcal{L}_n$, \ref{asp:iii}, and Taylor's theorem, we have
$$
\|h_{n,11}\|_{BV} \lesssim \|\ld - \ld_0\|_\Lone (\|\Delta_1\|_\infty + \|\Delta_{2,L_n}\|_\Linfty)^2/\rn
\leq \epsilon_n (\|\Delta_1\|_\infty + \|\Delta_{2,L_n}\|_\Linfty)^2/\rn.
$$ 

With all the conditions in Lemma \ref{lemma-ep-3} are verified, applying this lemma, we obtain
$$
\mathbb{E}_{\eta_0}^\star \|\mathbb{G}_n\|_{\mathcal{F}_{n,11}} 
\lesssim \epsilon_n (\|\Delta_1\|_\infty + \|\Delta_{2,L_n}\|_\Linfty)^2/\rn.
$$

Bounding the second term in (\ref{gn-F1}) is similar. 
Again, we use Lemma \ref{lemma-ep-3}.
Since the $g_\theta$ part is the same as in $f_{n,12}$, by following the same argument, the first condition in Lemma \ref{lemma-ep-3} is verified. 
To verify the second condition, let $h := h_{n, 12} = \rn \int_0^\cdot e^{r_0} (e^{\tilde \Delta_n/\rn} - \tilde \Delta_n/\rn - 1)$ instead, 
it is also clear that $h_{n, 12}(0) = 0$. Since $\tilde \Delta_n/\rn < 1$, by \ref{asp:iii} and Taylor's theorem, we have
$
\|h_{n,12} \|_{BV} \lesssim (\|\Delta_1\|_\infty + \|\Delta_{2,L_n}\|_{\Linfty})^2 /\rn.
$
Thus, we obtain 
$$
\mathbb{E}_{\eta_0}^\star \|\mathbb{G}_n\|_{\mathcal{F}_{n,12}} \lesssim  (\|\Delta_1\|_\infty + \|\Delta_{2,L_n}\|_\Linfty)^2/\rn.
$$

Now by combining the two upper bounds derived above, we have
$$
\mathbb{E}_{\eta_0}^\star \|\mathbb{G}_n\|_{\mathcal{F}_{n,1}} \lesssim  (1+\epsilon_n)(\|\Delta_1\|_\infty + \|\Delta_{2,L_n}\|_\Linfty)^2/\rn \lesssim (\|\Delta_1\|_\infty + \|\Delta_{2,L_n}\|_\Linfty)^2/\rn.
$$

To bound an upper bound for $\mathbb{E}_{\eta_0}\|\mathbb{G}_n\|_{f_{n,2}}$,
let 
$
f_{n,2} = f_{n,21} + f_{n,22},
$
where
\begin{align*}
f_{n,21} & = \int_0^\cdot e^{\theta_0'z + r_0(\cdot)}  e^{(\theta-\theta_0)'z} \left(e^{(r- r_0)(\cdot)} - 1\right) \tilde \Delta_n,\\
f_{n,22} & = \int_0^\cdot e^{\theta_0'z + r_0(\cdot)} ( e^{(\theta-\theta_0)'z} - 1) \tilde \Delta_n,
\end{align*}
and denote
\begin{align*}
\mathcal{F}_{n,21} 
& = \left\{
 \int_0^\cdot e^{\theta_0'z + r_0(\cdot)}  e^{(\theta-\theta_0)'z} \left(e^{(r- r_0)(\cdot)} - 1\right) \tilde \Delta_n,\  \eta \in \mathcal{L}_n
\right\}, \\
\mathcal{F}_{n,22} 
& = \left\{
\int_0^\cdot e^{\theta_0'z + r_0(\cdot)} ( e^{(\theta-\theta_0)'z} - 1) \tilde \Delta_n, \ \eta \in \mathcal{L}_n
\right\}
\end{align*}
such that $f_{n,21} \in \mathcal{F}_{n,21}$ and $f_{n,22} \in \mathcal{F}_{n,22}$.
Then,
\begin{align}
\label{gn-F2}
\mathbb{E}_{\eta_0}^\star \|\mathbb{G}_n\|_{\mathcal{F}_{n,2}} \leq
\mathbb{E}_{\eta_0}^\star \|\mathbb{G}_n\|_{\mathcal{F}_{n,21}} + 
\mathbb{E}_{\eta_0}^\star \|\mathbb{G}_n\|_{\mathcal{F}_{n,22}}.
\end{align}

Bounding (\ref{gn-F2}) is similar to bounding (\ref{gn-F1}). 
We need verify the conditions in Lemma \ref{lemma-ep-3} for each term in the upper bound.

For the first term, consider $g_{n,21} = e^{\theta'z}$ and 
$h_{n,21} =  \int_0^\cdot \left(e^{r-r_0} - 1\right) \tilde \Delta_n$. 
We already showed that $g_{n,21}$ is bounded by a constant as it is the same as the $g_\theta$ part in the function $f_{n,11}$. 
From the expression of $h_{n,21}$, it is easy to check that $h_{n,21}(0) = 0$. Also, 
\begin{align}
\label{h-n21}
\|h_{n,21}\|_{BV} \leq \epsilon_n \|\tilde \Delta_n\|_\infty \leq \epsilon_n (\|\Delta_1\|_\infty \|z\|_1 + \|\Delta_{2, L_n}\|_\Linfty).
\end{align}
Thus by \ref{asp:i} and Lemma \ref{lemma-ep-3},
$$
\mathbb{E}_{\eta_0}^\star \|\mathbb{G}_n\|_{\mathcal{F}_{n,21}}  
\lesssim \epsilon_n (\|\Delta_1\|_\infty + \|\Delta_{2,L_n}\|_\Linfty). 
$$

To bound the second term in (\ref{gn-F2}), 
define $g_{n, 22} = \int_0^1 e^{\theta_0'z + r_0(\cdot)} ( e^{(\theta-\theta_0)'z} - 1) \tilde \Delta_n$
and $h_{n,22}(u) = \mathbbm{1}_{[0, u]}$, $u \in [0, 1]$.
Then, $h_{n,22}(0) = 0$ and $\|h_{n,22}(u) \|_{BV} = 1$. 
Before applying Lemma \ref{lemma-ep-3}, we also need to bound $|g_{n,22}(\theta_1) - g_{n,22}(\theta_2)|$ for $\theta_1, \theta_2 \in \mathcal{L}_n$.
We have
$$
|g_{n,22}(\theta_1) - g_{n,22}(\theta_2)|
\leq e^{\theta_0'z} \|\Ld_0\|_1 \|\tilde \Delta_n\|_\Linfty 
|e^{(\theta_1 - \theta_0)'z} - e^{(\theta_2 - \theta_0)'z}|.
$$
By Taylor's theorem and \ref{asp:i}-\ref{asp:iii}, the last display is bounded by 
some constant times 
$$
\|\theta_1 - \theta_2\| \|\tilde \Delta_n\|_\Linfty \leq \epsilon_n \|\tilde \Delta_n\|
\leq \epsilon_n \left( \|\Delta_1\|_\infty + \|\Delta_{2,L_n}\|_\infty \right).
$$
Therefore, 
$\mathbb{E}_{\eta_0}^\star \|\mathbb{G}_n\|_{\mathcal{F}_{n,22}} \lesssim \epsilon_n (\|\Delta_1\|_\infty + \|\Delta_{2,L_n}\|_\Linfty)$. 
Thus, we obtain
$$
\mathbb{E}_{\eta_0}^\star \|\mathbb{G}_n\|_{\mathcal{F}_{n,2}}  
\leq 
\mathbb{E}_{\eta_0}^\star \|\mathbb{G}_n\|_{\mathcal{F}_{n,21}}  
+
\mathbb{E}_{\eta_0}^\star \|\mathbb{G}_n\|_{\mathcal{F}_{n,22}}  
\lesssim \epsilon_n (\|\Delta_1\|_\infty + \|\Delta_{2,L_n}\|_\Linfty).
$$
\end{proof}

%%%%%%%%%%%%%%%%%%%%%%%%%%%%%%%%%%%%%%%%%%%%%%
%%%%%%%%%%%%%%%%%%%%%%%%%%%%%%%%%%%%%%%%%%%%%%
We have obtained an upper bound for $\mathbb{E}_{\eta_0}^\star \|\mathbb{G}_n\|_{\mathcal{F}_{n,1}}$ and $\mathbb{E}_{\eta_0}^\star \|\mathbb{G}_n\|_{\mathcal{F}_{n,2}}$ respectively.
The derivation at the beginning of this subsection shows that $\sup_{\eta \in A_n} |R_{n,1}(\eta, \eta_0) - R_{n,1}(\eta_h, \eta_0)|$ can be bounded by the summation of these two upper bounds. The bound is given in the next lemma. 

\begin{lemma}
\label{lemma-R1}
Under the same condition as in Lemma \ref{R1-version1}, for 
$\Delta_1$ and $\Delta_{2,L_n}$ defined in (\ref{K1}) and (\ref{K2}) respectively, 
if $\|a\|_\infty$ is bounded and $b \in L^{\infty}([0,1])$, then under assumptions {\rm \ref{asp:i}-\ref{asp:v}},
$$
\sup_{\eta \in \mathcal{L}_n} |R_{n,1}(\eta, \eta_0) - R_{n,1} (\eta_h, \eta_0)| 
= O_{P_{\eta_0}}(L_n^2/\rn + \epsilon_n L_n).
$$
\end{lemma}
\begin{proof}
From (\ref{eqn:gn}), we have
$
R_{n,1}(\eta, \eta_0) - R_{n,1}(\eta_h, \eta_0)
= \mathbb{G}_n (f_{n,1}) + \mathbb{G}_n (f_{n,2}),
$
where $f_{n,1} \in \mathcal{F}_{n,1}$ in (\ref{F-n1}) and $f_{n,2} \in \mathcal{F}_{n,2}$ in (\ref{F-n2}). 
Thus
\begin{align}
\label{rn1-ub}
\sup_{\eta \in \mathcal{L}_n} |R_{n,1}(\eta, \eta_0) - R_{n,1}(\eta_h, \eta_0)| 
\leq \| \mathbb{G}_n\|_{\mathcal{F}_{n,1}} + \| \mathbb{G}_n\|_{\mathcal{F}_{n,1}}.
\end{align}
We apply Lemma \ref{R1} to bound each term in the last display.
The upper bound in Lemma \ref{R1} involves $\|\Dl_1\|_\infty$ and $\|\Dl_{2,L_n}\|_\Linfty$,
where 
\begin{align}
& \qquad \Delta_1 := \Delta_1(a, b, h) = t\tIinv a + s \tIinv \Ld_0 \{b\gamm\},
\label{delta1}\\
\Delta_{2,L_n} &: = \Delta_{2,L_n}(a,b,h) = t \gammln \tIinv a - s\gambln - s\gammln' \tIinv \Ld_0\{b \gamm\},
\label{delta2}
\end{align}
we need to bound $\|\Delta_1\|_\infty/\rn$ and $\|\Delta_{2,L_n}\|_\infty/\rn$ first. 

By triangular inequality and $\|fg\|_\infty \leq \|f\|_\infty \|g\|_\infty$,
\begin{align*}
\|\Dl_1\|_\infty & = \|t a'\tIinv  + s\tIinv \Ld_0 \{b\gamm\}\|_\infty \\
&  \leq |t| p^2 \|a\|_\infty \|\tIinv\|_{(\infty, \infty)} + |s| p^2 \|\Ld_0\{b\gamm\}\|_\infty \|\tIinv\|_{(\infty, \infty)}.
\end{align*}
By \ref{asp:i}-\ref{asp:v}, since $t,s,p$ are fixed and $\|a\|_\infty$ and $\|b\|_\infty$ are both bounded by assumption, and $\|\Ld_0\{b\gamm\}\|_\infty \leq \|\Ld_0\|_\infty \|b\|_1 \|\gamma\|_\infty$, the last display is $O(1/\rn) = o(1)$. 

Again, by triangular inequality, 
\begin{align*}
\|\Dl_{2,L_n}\|_\Linfty 
& \leq \|t \gammln ' \tIinv a\|_\Linfty + \|s \gambln\|_\Linfty + \|s\gammln' \tIinv \Ld_0\{b \gamm\}\|_\Linfty\\
& \lesssim p^2 (|t| + |s|) L_n + |s| L_n,
\end{align*}
Using \ref{asp:i}-\ref{asp:v}, the first point in Lemma \ref{lemma:psib-1}, and the third point in Lemma \ref{lemma:psib-2}, the last display is bounded by $O(L_n /\rn) = o(1)$.

We now apply Lemma \ref{R1-version1} and obtain
\begin{align*}
\sup_{f_{n,1} \in \mathcal{F}_{n,1}}|\mathbb{G}_n (f_{n,1})| 
& = O_{p_{\eta_0}}(L_n^2/\rn),\\
\sup_{f_{n,2} \in \mathcal{F}_{n,2}}|\mathbb{G}_n (f_{n,2})| 
& = O_{p_{\eta_0}}( \epsilon_n L_n ).
\end{align*}
Thus, 
$
\sup_{\eta \in \mathcal{L}_n} |R_{n,1}(\eta, \eta_0) - R_{n,1} (\eta_h, \eta_0)| 
= O_{P_{\eta_0}}(L_n^2/\rn + \epsilon_n L_n).$
\end{proof}

%%%%%%%%%%%%%%%%%%%%%%%%%%%%%%%%%%%%%%%%%%%%%%%%
\subsection{Bounding $\sup_{\eta \in A_n} |R_{n,2}(\eta, \eta_0) - R_{n,2}(\eta_h, \eta_0) - s\rn B_3(\eta, \eta_0)|$}
\label{sec:bounding-R2}

\begin{lemma}
\label{thm1-lemma2}
Suppose {\rm \ref{asp:i}-\ref{asp:v}} and {\rm\ref{cond:P}} hold, 
let $K_{a,b,t,s} = p^2 (|t|\|a\|_\infty + |s| \|b\|_\Lone) + |s| \|b\|_\Linfty$
and $\tilde K_{a,b,t,s} = p^2 (|t|\|a\|_\infty + |s| \|b\|_\Lone) + |s|\|b\|_\Ltwo$,
if $K_{a,b,t,s} L_n/\rn = o(1)$, then 
\begin{align*}
& \sup_{\eta \in A_n} |R_{n,2}(\eta_h, \eta_0) - R_{n,2}(\eta, \eta_0) - s\rn B_3(\eta, \eta_0)| \\
& \quad \lesssim 
\tilde K_{a,b,t,s}^3/\rn + K^2_{a,b,h,p}L_n^2 \epsilon_n+ |s|p^2 \rn \epsilon_n 2^{-L_n} + \rn \epsilon_n^2 L_n K_{a,b,t,s}.
\end{align*}
\end{lemma}

\begin{proof}
By plugging-in the expression for $R_{n,2}(\eta,\eta_0)$ in (\ref{R2}), 
\begin{align}
& R_{n,2}(\eta, \eta_0) - R_{n,2}(\eta_h, \eta_0)
\nonumber \\
& \quad = 
n \Ld_0 \Big\{
M_0(\theta_h) e^{r_h-r_0} - M_0(\theta) e^{r-r_0}
- (\theta_h - \theta)'M_1 - (r_h - r) M_0 \Big\} 
\label{lem6-thm1-1} \\
& \qquad - \frac{n}{2} \left(\|\theta_h - \theta_0, r - r_0\|_L^2 - \|\theta - \theta_0, r_h - r_0\|_L^2\right)
\label{lem6-thm1-2}
\end{align}
Denote $m(u,z) = e^{\theta_0'z} e^{-\Ld_0(u) e^{\theta_0'z}}$,
for $M_0(\theta)(\cdot)$ in (\ref{M0-theta}), we can write
$$
M_0(\theta) = \int \bar G_z(u) e^{(\theta - \theta_0)'z} m(u,z)f_Z(z) dz,
$$
for the expression in (\ref{lem6-thm1-1}), we can write 
\begin{align}
& n\Lambda_0 \Bigg\{ 
\int \bar G_z(u) e^{(\theta - \theta_0)'z + (r-r_0)} \left[e^{(\theta_h - \theta)'z + r_h - r} - 1\right] m(u, z) f_Z(z) dz 
\label{lem6-thm1-3}\\
& \qquad - (\theta_h - \theta)'M_1 - (r_h - r)M_0 
\Bigg\},
\label{lem6-thm1-4}
\end{align}
and for the expression in (\ref{lem6-thm1-2}) and the LAN-norm in (\ref{LAN-norm-v2}), we can write
\begin{align*}
& \|\theta_h - \theta_0, r_h - r_0\|_L^2 - \|\theta - \theta_0, r - r_0\|_L^2
= \|\theta_h - \theta, r_h - r\|_L^2 \\
& \qquad + 2 \Ld_0 
\left\{
\int \bar G_z(u) 
\left[(\theta_h - \theta)'z + r_h - r\right]
\left[(\theta - \theta_0)'z + r- r_0\right]
m(u,z) f_Z(z) dz
\right\}.
\end{align*}
By plugging-in the above expressions and $B_3(\eta, \eta_0)$ in (\ref{B3}),
\begin{align}
& R_{n,2}(\eta, \eta_0) - R_{n,2}(\eta_h, \eta_0) - s \rn B_3(\eta, \eta_0) 
\label{Rn2-diff}\\
& \quad = 
n\Ld_0
\Big\{
\int \bar G_z(u) 
e^{(\theta - \theta_0)'z + r-r_0}
\left[ e^{(\theta_h - \theta)'z + r_h - r} - (\theta_h - \theta)'z - (r_h - r) - 1\right] 
\label{lem6-thm1-5} \\
& \qquad \qquad \ \times m(u,z) f_Z(z) dz
\Big\} 
 - \frac{n}{2} \|\theta_h - \theta, r_h - r\|_L^2
\label{lem6-thm1-6} \\
& \qquad 
+ B^h_n (\eta, \eta_0) - s \rn B_3(\eta, \eta_0),
\label{lem6-thm1-7}
\end{align}
where
\begin{align*}
B^h_n(\eta, \eta_0)
& = n\Ld_0
\Bigg\{ \int \bar G_z(u)
\left[
e^{(\theta - \theta_0)'z+ r-r_0} - (\theta - \theta_0)'z - (r - r_0) - 1
\right] \\
& \quad \times ((\theta_h - \theta)'z + r_h - r)
m(u,z) f_Z(z) dz
\Bigg\}.
\end{align*}
To bound (\ref{Rn2-diff}), we first bound (\ref{lem6-thm1-5}) and (\ref{lem6-thm1-6}).
We first rewrite (\ref{lem6-thm1-5}) as follows:
\begin{align*}
& n\Ld_0
\left\{
\int \bar G_z(u) 
(e^{(\theta - \theta_0)'z + r-r_0} - 1)
\left[ e^{(\theta_h - \theta)'z + r_h - r} - (\theta_h - \theta)'z - (r_h - r) - 1\right]
m(u,z) f_Z(z) dz
\right\} \\
& \quad + 
n\Ld_0
\left\{\int \bar G_z(u) 
\left[ e^{(\theta_h - \theta)'z + r_h - r} - (\theta_h - \theta)'z - (r_h - r) - 1\right]
m(u,z) f_Z(z) dz
\right\} \\
& \quad = (I) + (II). 
\end{align*}
We then bound $(I)$ and $(II)-\frac{n}{2} \|\theta_h - \theta, r_h - r\|_L^2$ separately. 

We first bound $(II)-\frac{n}{2} \|\theta_h - \theta, r_h - r\|_L^2$. In order to apply Taylor's theorem to the exponential part in the expression,
one shall check $\|(\theta_h - \theta)'z + r_h - r\|_\infty = o(1)$. 
From (\ref{thm1-lemma7-1}) in Lemma \ref{thm1-lemma7}, 
\begin{align*}
\max_z |(\theta_h - \theta)'z| + \|r_h - r\|_\Linfty 
\lesssim p^2 L_n (|t|\|a\|_\infty + |s| \|b\|_\Lone)/\rn + |s|\|\gambln\|_\Linfty/\rn.
\end{align*}
By Lemma \ref{lemma:psib-1}, we have $\|\gambln\|_\Linfty \lesssim L_n \|b\|_\Linfty$, 
then the last display is bounded by a constant times $L_n K_{a,b,t,s}/\rn$
for $K_{a,b,t,s} = p^2 (|t|\|a\|_\infty + |s| \|b\|_\Lone) + |s| \|b\|_\Linfty$.
Using the assumption that $L_n K_{a,b,t,s}/\rn = o(1)$, we obtain $\|(\theta_h - \theta)'z + r_h - r\|_\infty = o(1)$.

We now apply Taylor's theorem for $e^{(\theta_h - \theta)'z + r_h - r}$,
then 
\begin{align}
\label{lem6-thm1-7-1}
 &(II) - \frac{n}{2} \|\theta_h - \theta, r_h - r\|_L^2
 \leq n\|(\theta_h - \theta)'z +r_h - r\|_\Lone^3 \|m(u,z)\|_\Linfty \|f_Z(z)\|_\Linfty \|\ld_0\|_\Linfty.
\end{align}
From (\ref{thm1-lemma7-2}), 
$$\|(\theta - \theta_h)'z + r - r_h\|_\Lone \lesssim p^2 (|t|\|a\|_\infty + |s|\|b\|_\Lone)/\rn + |s| \|\gambln\|_\Ltwo/\rn.$$
Let $\tilde K_{a,b,t,s} = p^2 (|t| \|a\|_\infty + |s|\|b\|_\Lone) + |s|\|b\|_\Ltwo$, 
$\|\gambln\|_\Ltwo \lesssim \|b\|_\Ltwo$ by Lemma \ref{lemma:psib-1}, 
and \ref{asp:i}--\ref{asp:iv},
(\ref{lem6-thm1-7-1}) is bounded by 
$C_1\tilde K_{a,b,t,s}^3 /\rn$ for some constant $C_1 > 0$.

Next, we bound $(I)$.
Since $\|(\theta - \theta_0)'z + r-r_0\|_\Linfty \lesssim 2^{L_n/2}\epsilon_n + 2^{-\beta L_n} = o(1)$ for $\beta > 1/2$, apply Taylor's theorem and $\|\ld - \ld_0\|_\Lone \leq \epsilon_n$ as $\eta \in A_n$, we obtain  
\begin{align*}
\|e^{(\theta - \theta_0)'z + r - r_0} - 1\|_\Lone K_{a,b,t,s}^2 L_n^2 
& \leq \tilde K_{a,b,t,s}^2 L_n^2 
	\| e^{(\theta - \theta_0)'z}(e^{r-r_0} - 1) + (e^{(\theta - \theta_0)'z} - 1)\|_\Lone\\
& \lesssim \tilde K_{a,b,t,s}^2 L_n^2  (\|\theta - \theta_0\|_1 + \|\lambda - \lambda_0\|_\Lone)\\
& \lesssim \tilde K_{a,b,t,s}^2 L_n^2 \epsilon_n.
\end{align*}
By combining the above upper bounds, we obtain 
$$
\sup_{\eta \in A_n} ((I) + (II) - \frac{n}{2}\|\theta_h - \theta, r_h - r\|_L^2) \lesssim 
\tilde K^2_{a,b,h,p}L_n^2 \epsilon_n + \tilde K_{a,b,t,s}^3/\rn.
$$ 
To bound (\ref{lem6-thm1-7}),
we rewrite $B_3(\eta, \eta_0)$ in (\ref{B3}) as
\begin{align*}
\rn \left\langle
\left(0, r - r_0 - \frac{\ld - \ld_0}{\ld_0} \right),
\left( \theta - \theta_h - \frac{r a'\tIinv}{\rn}, 
r - r_h + \frac{t\gammln' \tIinv a}{\rn} \right)
\right\rangle_L.
\end{align*}
Then,
\begin{align}
& B^h_n (\eta, \eta_0) - s\rn B_3(\eta, \eta_0) \nonumber \\
& \quad =
n \Ld_0 \left\{
\int \bar G_z(u) \left[ (\theta - \theta_h)'z + r-r_h\right] 
\left[
(\theta - \theta_0)'z - e^{(\theta - \theta_0)'z + r- r_0} + e^{r - r_0} 
\right] m(u,z) f_Z(z) dz
\right\}
\label{lem6-thm1-7}
\\
& \quad \quad 
+ \rn s \Ld_0\left\{
\left(r - r_0 -e^{r- r_0} +1\right)
(\gammln - \gamm)' \tIinv a M_0
\right\}.
\label{lem6-thm1-8}
\end{align}
Using the fact that $e^{r-r_0} -1 = (\ld - \ld_0)/\ld_0$, (\ref{lem6-thm1-8}) can be bounded by 
\begin{align*}
& \rn |s| p^2 (\|\log\ld - \log \ld_0\|_\Lone + \|\ld - \ld_0\|_\Lone \|\ld_0^{-1}\|_\Linfty)  \\
& \quad \times \max_j\|\gammln^j - \gamm^j\|_\Linfty \|\tIinv\|_{(\infty, \infty)}\|a\|\|M_0\|_\Linfty \|\ld_0\|_\Linfty.
\end{align*}
Using \ref{asp:i}-\ref{asp:v} and the fact that $\log(\cdot)$ is Lipschitz, 
$\|\log\ld - \log \ld_0\|_\Lone \lesssim \|\ld - \ld_0\|_\Lone$ as $\|\ld - \ld_0\|_\Linfty = o(1)$ due to $\beta > 1/2$, and the third point of Lemma \ref{lemma:psib-2}, the last display is bounded by a constant times $|s| p^2 \rn \epsilon_n 2^{-L_n}$.

To bound (\ref{lem6-thm1-7}), 
write 
$$
(\theta - \theta_0)'z - e^{(\theta - \theta_0)'z + r- r_0} + e^{r - r_0}
= 
(e^{r-r_0} - 1)(1 - e^{(\theta - \theta_0)}) - e^{(\theta - \theta_0)'z} + (\theta - \theta_0)'z + 1
$$
Since $\|\theta - \theta_0\| \leq \epsilon_n$ and $\|e^{r- r_0} - 1\|_\Lone \lesssim \|\ld - \ld_0\|_\Lone \leq \epsilon_n$ as $\|r - r_0\|_\infty = o(1)$ due to $\beta > 1/2$, we obtain 
\begin{align*}
& \|(e^{r-r_0} - 1)(1 - e^{(\theta - \theta_0)'z})\|_\Lone \lesssim \epsilon_n^2, 
\quad  \max_z \|e^{(\theta - \theta_0)'z} - (\theta - \theta_0)'z - 1\|\lesssim \epsilon_n^2.
\end{align*}
Also, 
$\| (\theta - \theta_h)'z + r-r_h\|_\Linfty \leq L_n \tilde K_{a,b,t,s}/\rn$ as we argued above,
we can bound (\ref{lem6-thm1-7}) by a constant times
$
\rn \epsilon_n^2 L_n K_{a,b,t,s}.
$
By adding the upper bounds of (\ref{lem6-thm1-7}) and (\ref{lem6-thm1-8}), we obtain
\begin{align}
\label{lem6-thm1-9}
\sup_{\eta \in A_n} |B_n^h(\eta, \eta_0) - s\rn B_3(\eta, \eta_0)| \lesssim 
|s| p^2 \rn \epsilon_n 2^{-L_n} + \rn \epsilon_n^2 L_n \tilde K_{a,b,t,s}.
\end{align}
By combining the upper bound for (\ref{lem6-thm1-7-1}), which is $C_1\tilde K^3_{a,b,t,s}/\rn$, and (\ref{lem6-thm1-9}), we thus complete the proof.
\end{proof}

%%%%%%%%%%%%%%%%%%%%%%%%%%%%%%%%%%%%%%%
\begin{lemma}
\label{thm1-lemma7}
For $\theta_h$ defined in (\ref{path-theta}) and $r_h$ in (\ref{path-r}), suppose assumptions
{\rm \ref{asp:i}-\ref{asp:v}}
hold,
then,
\begin{align}
\label{thm1-lemma7-1}
\max_z |(\theta_h-\theta)'z| + \|r_h-r\|_\Linfty 
& \lesssim  p^2 L_n (|t|\|a\|_\infty + |s|\|b\|_\Lone)/\rn + |s| \|\gambln\|_\Linfty/\rn,\\
\label{thm1-lemma7-2}
\max_z |(\theta_h-\theta)'z| + \|r_h-r\|_\Lone
& \lesssim p^2 (|t|\|a\|_\infty + |s|\|b\|_\Lone)/\rn + |s| \|\gambln\|_\Ltwo/\rn,
\end{align}
\end{lemma}

\begin{proof}
Recall the definitions of $\theta_h$ and $r_h$, we immediately obtain
\begin{align*}
& (\theta_h- \theta) 'z = - \frac{t a'\tIinv z}{\rn} + \frac{s  \Ld_0\{b\gamm'\} \tIinv z}{\rn},\\
r_h - & r = \frac{t \gammln' \tIinv a}{\rn} - \frac{s \gambln}{\rn} - \frac{s \gammln' \tIinv \Ld_0\{b \gamm\}}{\rn}.
\end{align*}
First, by applying 
the inequality $a'Hb \leq p^2 \|a\|_\infty \|H\|_{(\infty, \infty)} \|b\|_\infty$ for any $a, b \in \mathbb{R}^p$ and $H \in \mathbb{R}^{p \times p}$, one obtains
$$
\rn \max_z |(\theta_h - \theta)'z| 
\leq |t| p^2 \|a\|_\infty \|\tIinv\|_{(\infty, \infty)} \|z\|_\infty + |s| p^2 \|\Ld_0\{b\gamm\}\|_\infty \|\tIinv\|_{(\infty, \infty)} \|z\|_\infty.
$$
By \ref{asp:i}, \ref{asp:iii}, \ref{asp:v}, and $\|\Ld_0\{b \gamm\}\|_\infty \leq \max_j \|\gamma_{M_{j}}\|_\Linfty \|b\|_\Lone \|\lambda_0\|_\Linfty
\leq \|z\|_\infty \|b\|_\Lone \|\lambda_0\|_\Linfty \leq C_1 \|b\|_\Lone$ for some constant $C_1 > 0$, 
the last display is thus bounded by 
\begin{align}
\label{thm1-lemma7-3}
\max_z |(\theta_h - \theta)'z|  \lesssim p^2 (|t|\|a\|_\infty + |s| \|b\|_\Lone) /\rn. 
\end{align}
Next, applying the same triangle inequality again, one obtains 
\begin{align*}
\rn \|r_h  - r\|_\Linfty
& \leq |t| p^2 \|a\|_\infty \|\tIinv\|_{(\infty, \infty)}\max_j \|\gammln^j\|_\Linfty + |s| \|\gambln\|_\Linfty \\
& \quad + |s| p^2 \|\Ld_0\{b \gamm\}\|_\infty \|\tIinv\|_{(\infty, \infty)} \max_j\|\gammln^j\|_\Linfty,
\end{align*}
where $\gammln^j$ is the $j$-th coordinate of $\gammln$. 
Then by \ref{asp:i}, \ref{asp:iii}, \ref{asp:v}, and $\|\Ld_0\{b \gamm\}\|_\infty \lesssim \|b\|_\Lone$, using  the third point of Lemma \ref{lemma:psib-2}, $\max_j\|\gammln\|_\Linfty \lesssim L_n$, we obtain
\begin{align}
\label{thm1-lemma7-4}
\|r_h  - r\|_\Linfty \lesssim p^2 L_n (|t|\|a\|_\infty + |s|\|b\|_\Lone)/\rn + |s| \|\gambln\|_\Linfty/\rn.
\end{align}
Now, combining the bounds in (\ref{thm1-lemma7-3}) and (\ref{thm1-lemma7-4}), we obtain
(\ref{thm1-lemma7-1}). 

Proving (\ref{thm1-lemma7-2}) is similar. Since $\|\cdot\|_1 \leq \|\cdot\|_2$, we have
\begin{align}
\rn \|r_h  - r\|_\Ltwo
& \leq |t| p^2 \|a\|_\infty \|\tIinv\|_{(\infty, \infty)} \max_j \|\gammln^j\|_\Ltwo + |s| \|\gambln\|_\Ltwo
\nonumber \\
& \quad + |s| p^2 \|\Ld_0\{b\gamm\}\|_\infty \|\tIinv\|_{(\infty, \infty)} \max_j \|\gammln^j\|_\Ltwo
\nonumber \\
& \lesssim p^2|t|\|a\|_\infty + |s| \|\gambln\|_\Ltwo + p^2 |s|\|b\|_\Lone,
\label{thm1-lemma7-5}
\end{align}
where we used triangular inequality and the first inequality in the third point of Lemma \ref{lemma:psib-2}.
By combing the bounds in (\ref{thm1-lemma7-3}) and (\ref{thm1-lemma7-5}), we proved (\ref{thm1-lemma7-2}).
\end{proof}

%%%%%%%%%%%%%%%%%%%%%%%%%%%%%%%%%%%%%%%%%%%%%%%%%%%%%%%%%%%%%%

\section{Joint nonparametric BvM for $\eta = (\theta, \ld)$}
\label{joint-npBvM}
In this section, we establish two nonparametric BvM theorems: first, the joint BvM theorem for the regression coefficients $\theta$ and the nonparametric part $\lambda$ and next, the BvM theorem for {the hazard function conditional on $z$}. 
The second BvM theorem serves as an important step for obtaining the Donsker theorem for the conditional cumulative hazard function and the {survival function conditional on $z$}. 

\subsection{Nonparametric BvM theorems}
\label{nonparametric-BvM}

Since the baseline hazard function $\lambda$ is a nonparametric quantity, it is well known that $\lambda$ is only estimable with a slower rate than $1/\rn$ in $L^2$ (hence, $L^\infty$-losses).
In order to obtain a rate of the order of $1/\rn$, one has to choose some larger spaces than $L^2$; e.g., the {\it Sobolev spaces} or the {\it `logarithmic' Sobolev spaces} with an order $s \leq -1/2$ introduced by  \citet{cast13} and {\it the multiscale space} proposed by \citet{cast14a}. Here, we work with the {\it multiscale space} as it contains the order $s = -1/2$ Sobolev space and is more adapted to obtain supremum-norm contraction rates. 

Let us define $Q_0$ the probability measure on $[0,1]$ with density $q_0 = \lambda_0/M_0$ with respect to Lebesgue's measure, i.e., $dQ_0(z) = q_0(x) dx$. 
Denote by $\mathbb{Z}_{Q_0}$ the $Q_0$-white noise process indexed by the Hilbert space 
$L^2(Q_0) = \{f: \int_0^1 f^2 dQ_0 < \infty\}$: that is, the zero-mean Gaussian process with its covariance function given by
\begin{align}
\label{WN-cov}
\mathbb{E}(\mathbb{Z}_{Q_0}(g) \mathbb{Z}_{Q_0}(h)) = \int_0^1 gh dQ_0.
\end{align} 
Let $(w_l)$ be a sequence $w_l/\sqrt{l} \to \infty$ as $l \to\infty$. 
We require $w_l \geq 1$ so that $\|x\|_{\mathcal{M}} \leq \|x\|_{L^2}$, $x \in \mathcal{M}$. We follow Definition 1 in \citet{cast14a} and call $(w_l)$ an admissible sequence.
Let $(\psi_{lk})$ be the Haar wavelet basis, the {\it multiscale space} is defined as
\begin{align}
\label{multiscale-space}
\mathcal{M}: = \mathcal{M}(w) = \left\{
\lambda = \{\langle \lambda, \psi_{lk} \rangle\},
\ \sup_{l \leq L} \max_{0 \leq k \leq 2^l} 
\frac{|\langle \lambda, \psi_{lk} \rangle|}
{w_l} < \infty
\right\}.
\end{align}
Furthermore, a {\it separable multiscale subspace} thereof is defined as
\begin{align}
\label{multiscale-subspace}
\mathcal{M}_0: = \mathcal{M}_0(w) = \left\{
\lambda = \{\langle \lambda, \psi_{lk} \rangle \},
\ \sup_{l \to \infty} \max_{0 \leq k \leq 2^l} 
\frac{|\langle \lambda, \psi_{lk} \rangle|}
{w_l} = 0
\right\}.
\end{align}

Let us denote 
$W_n^{(1)} = W_n^\star(\tIinv, \ - \gamm' \tIinv)$ as in (\ref{Wn-1})
and $W_n^{(2)}(b) = W_n(-\tIinv \Ld_0\{b \gamm\}, \ \gamb + \gamm' \tIinv \Ld_0 \{b \gamm\})$ as in (\ref{Wn-2}), and define the centering sequences for $\theta$ and $\ld$ as
$$
{T_n^\theta} = \theta_0 + W_n^{(1)}/\rn
$$
and 
\begin{align}
\label{centering-Tn}
\langle T_{n}^\ld, \psi_{lk} \rangle 
= 
\begin{cases}
\langle \lambda_0, \psi_{lk} \rangle 
+
W_n^{(2)}(\psi_{lk})/\rn & \quad \text{if} \ l\leq L_n,\\
0 & \quad \text{if} \ l > L_n.
\end{cases}
\end{align}

Let $\tau_{T_n}$ be the map 
$$
\tau_{T_n}: \eta \to \rn (\theta - T_n^{\theta}, \ \langle \ld - T_n^\ld, b \rangle),
$$
and $\Pi(\cdot \given X) \circ \tau_{T_n}^{-1}$ be the distribution induced on $\rn (\theta - {T_n^{\theta}}, \ \langle \ld - T_n^\ld, b \rangle)$. 

In order to obtain the nonparametric BvM, one needs to assume a stronger version of the change of variables condition than \ref{cond:C1}, as $t$ and $s$ can increase with $n$:
\begin{enumerate}[label=\textbf{(C2)}]
\item \label{cond:C2} (Change of variables condition, version 2) with the same setting as in \ref{cond:C1}, 
for any $|t|, |s| \leq \log n$, one assumes, for $A_n$ as in \ref{cond:P} and some constant $C_1 > 0$,
$$
\frac{
\int_{A_n} e^{\ell_n(\eta_h) - \ell_n(\eta_0)} d\Pi(\eta)
}
{
\int e^{\ell_n(\eta) - \ell_n(\eta_0)} d\Pi(\eta)
}
\leq e^{C_1(1 +t^2 + s^2)},
$$
for some $\eta_h = (\theta_h, r_h)$, $\theta_h$ and $r_h$ as in \eqref{path-theta} and \eqref{path-r}, for a fixed $a = z\in \mathbb{R}^p$ and a collection of functions $b$ to be specified below. 
\end{enumerate}

Also, for $\epsilon_n$ and $\zeta_n$ in $A_n$ and for $L_n$ in (\ref{choicel}), we assume 
\begin{align}
\label{rate-condition}
\rn \epsilon_n^2 L_n = o(1), \quad 
\rn \epsilon_n 2^{-L_n} L_n = o(1), \quad
\zeta_n L_n^2 = o(1).
\end{align}

\begin{theorem}[Joint nonparametric BvM for $\eta = (\theta, \ld)$]
\label{thm:nonparam-joint-bvm}
Let $\Pi$ be the independent product of the priors in {\rm \ref{prior:T}} and {\rm \ref{prior:W}}. Define the centering $T_n = (T_n^\theta, T_n^\ld)$. 
Let $\mathcal{M}_0$ be the separable multiscale subspace for some sequence $w_l \to \infty$ with $w_l \geq l$. 
Suppose {\rm{\ref{cond:P}}} is satisfied with $\epsilon_n$, $\zeta_n$, and a cut-off $L_n$ satisfy (\ref{rate-condition}) and 
$$
\rn \epsilon_n 2^{-L_n} = o\left(\min_{l\leq L_n} \{l^{-1/4} 2^{-l/2} w_l\}\right),
$$
and suppose {\rm \ref{cond:C1}} holds for $a = z$, $z \in \mathbb{R}^p$ is fixed, and
$b \in \mathcal{V}_{L} = \text{Vect}\{\psi_{lk}, \ l \leq L, 0 \leq k < 2^l\}$ with a fixed $L \geq 0$ and {\rm \ref{cond:C2}} holds uniformly for $a = z$, $z \in \mathbb{R}^p$ is fixed, and $b = \psi_{LK}$ with $L \leq L_n$ and $0 \leq K \leq 2^L$.
Then, for $\mathbb{Z}_{Q_0}$, as in (\ref{WN-cov}), is independent of the random variable $\mathbb{V} \sim N(0, \tIinv)$,
\begin{align}
\label{thm:nonparam-1}
\mathcal{B}_{\mathbb{R}^p \times \mathcal{M}_0} \left( \Pi((\theta, \lambda) \in \cdot \given X) \circ \tau_{T_n}^{-1}, \mathcal{L}(\mathbb{V} , \mathbb{Z}_{Q_0} - \mathbb{V}' \gamm \ld_0 ) \right) \stackrel{P_{\eta_0}}{\to} 0,
\end{align}
where
$\mathcal{B}_{\mathbb{R}^p \times \mathcal{M}_0}$ is the bounded-Lipschitz metric on $\mathbb{R}^p \times \mathcal{M}_0$.
\end{theorem}

By applying the delta method on both probability measures in (\ref{thm:nonparam-1}), Theorem \ref{thm:nonparam-joint-bvm} immediately implies the nonparametric BvM theorem for {the hazard function conditional on $z$}, which is given in the following corollary.  

\begin{corollary}[Nonparametric BvM for {the hazard function conditional on $z$}]
\label{cor:nonparam-bvm-chzd}
Under the same conditions as in Theorem \ref{thm:nonparam-joint-bvm}, define $\tilde \tau_{S_n}$ as the map $\tilde \tau_{S_n}: \eta \to \rn (\ld e^{\theta'\tz} - T_n^\ld e^{{T_n^\theta}'\tz})$, for a fixed $\tz \in \mathbb{R}^p$,
then
$$
\mathcal{B}_{\mathcal{M}_0} 
\left( \Pi(\lambda e^{\theta'\tz} \in \cdot \given X) \circ \tilde \tau_{T_n^\ld e^{{T_n^\theta}'\tz}}^{-1}, 
\mathcal{L} \left( e^{{T_n^\theta}'\tz} (\mathbb{W} + \mathbb{Z}_{Q_0}) \right) \right) \stackrel{P_{\eta_0}}{\to} 0,
$$
where $\mathbb{W}$ and $\mathbb{Z}_{Q_0}$ are independent, and 
$\mathbb{W} \sim N(0, \Delta)$ with
$
\Delta = 
(T_n^\ld)^2 \tz' \tilde I_{\eta_0}^{-1} \tz
- 2T_n^\ld \tz'\tilde I_{\eta_0}^{-1}  \gamma_{M} \ld_0 
+ \lambda_0^2 \gamma_{M}' \tilde I_{\eta_0}^{-1} \gamma_{M}
$.
\end{corollary}

\begin{remark}
\rm 
When $\tz = 0$, Corollary \ref{cor:nonparam-bvm-chzd} implies the nonparametric BvM for $\lambda$ in the survival model; i.e., under the same conditions as in Corollary \ref{cor:nonparam-bvm-chzd}, we have
$$
\mathcal{B}_{\mathcal{M}_0} \left( \Pi(\ld \in \cdot \given X) \circ \tilde\tau_{T_n^\ld}^{-1}, \
\mathcal{L}(\tilde{\mathbb{W}} + \mathbb{Z}_{Q_0}) \right) \stackrel{P_{\eta_0}}{\to} 0,
$$
where $\tilde{\mathbb{W}}$ and $\mathbb{Z}_{Q_0}$ are independent, and $\tilde{\mathbb{W}} \sim N(0, \tilde \Delta)$ with $\tilde \Delta = \lambda_0^2 \gamma_{M}' \tilde I_{\eta_0}^{-1} \gamma_{M}$.
\end{remark}

\subsection{Proof of nonparametric BvM results}
\label{pf-nonparametric-BvM}

In this section, we prove Theorem \ref{thm:nonparam-joint-bvm}.
We apply the general framework proposed by \citet{cast14a} to prove the theorem. A key is to establish the tightness criterion in space of $\mathbb{R}^p \times \mathcal{M}_0(w)$, which is given in Proposition \ref{prop-tightness}. In the proposition, we have to modify Proposition 6 in \citet{cast14a} in the space of $\mathcal{M}_0(w)$ to the product space $\mathbb{R}^p \times \mathcal{M}_0(w)$.

We need to verify the two conditions in Proposition \ref{prop-tightness}: 1) the BvM theorem for finite-dimensional distributions in (\ref{eqn:cond1}) and 2) tightness of $\ld$ at the rate $1/\rn$ in (\ref{eqn:cond2}). We first present a proposition in the next subsection, which is used for the proof of the tightness condition. 
The verification of the first condition is given in Section \ref{pf-npbvm}.
 
\subsubsection{Controlling Laplace transforms of linear functionals}
\label{sec:bvm-finite-distribution}

The following proposition will be used to verify the tightness criterion. 

\begin{prop}
\label{prop-bvm}
Suppose $b \in L^2(\Ld)$ possibly depends on $n$, for some positive constants $d_1$ and $d_2$ such that $\|b\|_\Ltwo \leq d_1$ and $\|b\|_\Linfty \leq d_2 2^{L_n/2}$. Assume {\rm{\ref{cond:P}}} and {\rm \ref{cond:C2}} hold, then for $\epsilon_n$, $\zeta_n$, and $A_n$ given in {\rm \ref{cond:P}}, a fixed $\tz \in \mathbb{R}^p$, and any $t \vee s \leq \log n$, 
\begin{align}
\label{prop:eqn-1}
& \log \mathbb{E}
\left(
e^{\rn\left( t(\theta'\tz - \theta_0'\tz) + s(\Ld\{b\} - \Ld_0\{b\})\right)} \given X, A_n
\right) \lesssim J_n(t, s),
\end{align}
where 
for
\begin{align*}
W_n^{(1)}(\tz) &= W_n(\tIinv \tz, \ - \gamm' \tIinv \tz),\\
W_n^{(2)}(b) &= W_n(- \tIinv \Ld_0\{b \gamm\},\  \gamb + \gamm' \tIinv \Ld_0\{b\gamm\}),
\end{align*}
and some constant $C > 0$, 
\begin{align*}
J_n(t,s) & = C\left(1+t^2 + s^2 + |s| (\rn \epsilon_n + O_{P_{\eta_0}}(1)) \|\gamb - \gambln\|_\Linfty
+ (|t| + |s|) \rn \epsilon_n^2\right)  \\
& \quad 
+ |s| p^2 \rn \epsilon_n 2^{-L_n} + (t^2 + s^2) \zeta_n
+ (|t| + |s|) O_{P_{\eta_0}}(\zeta_n) + t W_n^{(1)}(\tz) + sW_n^{(2)}(b).
\end{align*}
\end{prop}

\begin{proof}
By following the proof of Theorem \ref{thm:bvm-linear-fns}, one could bound the expectation on the left hand side of (\ref{prop:eqn-1}) by 
\begin{align}
& \exp
\Bigg(
\sup_{\eta \in A_n} \Big|
D_n + s\rn B_2(\eta, \eta_0) \Big| + \rn W_n(\theta - \theta_h, r - r_h)
\label{pf-cond1-1} \\
& \quad + \sup_{\eta \in A_n} \Big | R_{n}(\eta, \eta_0) - R_{n}(\eta_h, \eta_0)
- s \rn B_3(\eta, \eta_0) \Big|
\Bigg) 
\times
 \frac{
e^{h'\Sigma_{\tz,b}h/2}
\int_{A_n} e^{\ell_n(\eta_h) - \ell_n(\eta_0)} d\Pi(\eta)
}{
\Pi(A_n \given X) \int e^{\ell_n(\eta) - \ell_n(\eta_0)} d\Pi(\eta)
},
\label{pf-cond1-2} 
\end{align}
where 
$A_n = \{\|\theta - \theta_0\| \leq \epsilon_n, \|\ld - \ld_0\|_\Lone \leq \epsilon_n, \|r - r_0\|_\Linfty \leq \zeta_n\}$,
$h = (t,s)$, $\Sigma_{\tz, b}$ is $\Sigma_{a,b}$ in (\ref{Sigma-ab}) with $a = \tz$, and the expressions of
$B_2(\eta, \eta_0)$, $B_3(\eta,\eta_0)$, and $D_n$ are given in (\ref{B2}), (\ref{B3}), and (\ref{Dn}) respectively.

First, we bound (\ref{pf-cond1-1}). 
By Lemma \ref{lemma:bounding-B2} and $\|b\|_1 \leq \|b\|_2 \leq d_1$, 
\begin{align}
s\rn \sup_{\eta \in A_n} |B_2(\eta, \eta_0)|
& \lesssim 
\rn \epsilon_n \|\gamb - \gambln\|_\Linfty + p^2 \epsilon_n 2^{-L_n} \|b\|_\Lone
\nonumber \\
& \leq \rn \epsilon_n \|\gamb - \gambln\|_\Linfty + p^2 \epsilon_n 2^{-L_n} d_1.
\label{pf-cond1-3}
\end{align}

By Lemma \ref{lemma:bounding-Dn}, using $\|b\|_{\Ltwo} \leq d_1$ and $\|b\|_\Linfty \leq d_2 2^{L_n/2}$, 
\begin{align}
\sup_{\eta \in A_n} |D_n| 
& \lesssim 
s^2 \|\gamb^2 - \gambln^2\|_\Lone + (t^2 + s^2) p^2 2^{-L_n} (\|\gambln\|_\Lone + p^2 2^{-L_n})
\nonumber \\
& \lesssim s^2 \|\gamb - \gambln\|_\Linfty \|\gamb + \gambln\|_\Lone 
+ (t^2 + s^2)p^2 2^{-L_n} \left(C  +p^2 2^{-L_n}\right)
\nonumber \\
& \lesssim s^2 \|\gamb - \gambln\|_\Linfty + (t^2 + s^2) p^2 2^{-L_n},
\end{align}
where the second inequality in the last display is obtained by using $\|ab\|_\Lone \leq \|a\|_\Lone \|b\|_\Linfty$ for any two functions $a,b \in L^2$ and $\|\gambln\|_\Lone \leq C$ for some constant $C$. The last line is obtained by using the inequality $\|\gamb + \gambln\|_\Lone \leq \|\gamb\|_\Lone + \|\gambln\|_\Lone$ and both $\|\gamb\|_\Lone$ and $\|\gambln\|_\Lone$ are bounded by some constants. 

To bound the third term in (\ref{pf-cond1-1}),
by \ref{asp:v} and the second point in Lemma \ref{lemma:psib-2} (replacing $b$ with $M_{1j}$), due to the linearity of $W_n(\cdot)$, 
we obtain
\begin{align}
& \rn W_n(\theta - \theta_h, r - r_h) - tW_n^{(1)}(\tz) - sW_n^{(2)}(b) \nonumber \\
& \quad 
\lesssim O_{P_{\eta_0}}\left( (|t| + |s|) \max_j\|\gamm^j - \gammln^j\|_\Linfty + |s| \|\gamb - \gambln\|_\Linfty\right) 
\nonumber \\
& \quad \leq O_{P_{\eta_0}}\left( (|t| + |s|) 2^{-L_n} + |s| \|\gamb - \gambln\|_\Linfty\right).
\label{pf-cond1-4}
\end{align}
Therefore, 
$\rn W_n(\theta - \theta_h, r - r_h) \leq O_{P_{\eta_0}}\left( (|t| + |s|) 2^{-L_n} + |s| \|\gamb - \gambln\|_\Linfty\right) + tW_n^{(1)} (\tz) + sW_n^{(2)} (b)$. 

Next, we bound (\ref{pf-cond1-2}). We first bound $\sup_{\eta \in A_n} |R_{n}(\eta, \eta_0) - R_n(\eta_h, \eta_0) - s\rn B_3(\eta, \eta_0)|$.
Recall that $R_{n}(\eta, \eta_0) = R_{n,1}(\eta, \eta_0) + R_{n,2}(\eta, \eta_0)$. By (\ref{rn1-ub}), one can bound 
$$
\sup_{\eta \in A_n}|R_{n,1}(\eta, \eta_0) - R_{n,1}(\eta_h, \eta_0)|
\leq \|\mathbb{G}_n\|_{\mathcal{F}_{n,1}} + \|\mathbb{G}_n\|_{\mathcal{F}_{n,2}} , 
$$
where $\mathcal{F}_{n,1}$ and 
$\mathcal{F}_{n,2}$ are given in (\ref{F-n1}) and (\ref{F-n2}) respectively and $\mathcal{L}_n$ is replaced with $A_n$. 
We use $\|\cdot\|_\infty$-consistency for $\ld - \ld_0$ and Lemma \ref{R1-version2} to bound the last display. 
We first check the conditions in Lemma \ref{R1-version2}: for $\Dl_1$ and $\Dl_{2,L_n}$ given in (\ref{K1}) and (\ref{K2}) respectively, where 
\begin{align*}
\Delta_1 = - t \tIinv \tz + s \tIinv \Ld_0\{b\gamm\}, \quad 
\Delta_{2,L_n} = t \gammln' \tIinv \tz - s\gambln + s\gammln' \tIinv \Ld_0\{b\gamm\},
\end{align*}
by \ref{asp:i}, \ref{asp:iii}, \ref{asp:v}, and $\|b\|_\Ltwo \leq d_1$, 
$$
\|\Delta_1\|_\infty \leq p |t| \|\tIinv\|_{(\infty, \infty)} \|\tz\|_\infty + |s| p \|\tIinv\|_{(\infty, \infty)}
\|\Ld_0\{b\gamm\}\|_\infty
\leq (|t| + |s|) O(1).
$$
Since $t,s \leq \log n$, $|\Delta_1'\tz| /\rn = o(1)$.
With the same assumptions as bounding $\|\Dl_1\|_\infty$, by the first and the third points in Lemma \ref{lemma:psib-2}, we have 
\begin{align*}
\|\Delta_{2,L_n}\|_\Linfty
& \leq 
|t| p^2 \|\tz\|_\infty \|\tIinv\|_{(\infty, \infty)} \max_j \|\gammln^j\|_\Linfty 
+ |s| \|\gambln\|_\Linfty \\
& \quad + |s| \max_j\|\gammln^j\|_\Linfty \|\tIinv\|_{(\infty, \infty)} \|\Ld_0\{b\gamm\}\|_\Linfty\\
& \lesssim (|t| + |s|)L_n +  |s|L_n 2^{L_n/2},
\end{align*}
where we applied the inequality $ \|\Ld_0\{b\gamm\}\|_\Linfty \leq \|\ld_0\|_\infty \|b\|_1 \max_j\|\gamm^j\|_\infty \leq C$ for some constant $C > 0$.
Since $L_n (|t| + |s|) /\rn  = o(1)$ and $|s| L_n2^{L_n/2} /\rn  = o(1)$ for $t, s < \log n$ with the choice of $L_n$ in (\ref{choicel}), we have verified $\|\Delta_{2,L_n}\|_\infty/\rn = o(1)$.

By using the same assumptions as above and $\|b\|_2\leq d_1$, we have
\begin{align*}
\|\Delta_{2,L_n}\|_\Ltwo 
& \leq 
|t| p^2 \|\tz\|_\infty \|\tIinv\|_{(\infty, \infty)} \max_j \|\gammln^j\|_\Ltwo 
+ |s| \|\gambln\|_\Ltwo \\
& \quad + |s| \max_j\|\gammln^j\|_\Ltwo \|\tIinv\|_{(\infty, \infty)} \|\Ld_0\{b\gamm\}\|_\Linfty\\
& \lesssim (|t| + |s|)O(1).
\end{align*}
Thus, by Lemma \ref{R1-version2}, we obtain
\begin{align*}
\mathbb{E}^\star_{\eta_0} 
\left[\|\mathbb{G}_n\|_{\mathcal{F}_{n,1}} \right]
\lesssim (t^2 + s^2) /\rn,
\quad
\mathbb{E}^\star_{\eta_0} 
\left[\|\mathbb{G}_n\|_{\mathcal{F}_{n,2}} \right]
\lesssim \zeta_n (|t| + |s|),
\end{align*}
and then obtain
$$
\sup_{\eta \in A_n} |R_{n,1}(\eta, \eta_0) - R_{n,2}(\eta_h, \eta_0)| = O_{P_{\eta_0}} \left(\frac{t^2 + s^2}{\rn} + (|t| + |s|) \zeta_n\right),
$$
as $(t^2 + s^2)/\rn = o(1)$ as $n \to \infty$ for $t,s\leq \log n$.

Next, we bound $\sup_{\eta \in E_n} |R_{n,2}(\eta, \eta_0) - R_{n,2}(\eta_h, \eta_0) - s\rn B_3(\eta, \eta_0)|$. We apply Lemma \ref{lemma:bounding-R2-using-supnorm}. 
Define $\tilde K_{\tz,b,t,s} = p^2(|t| \|\tz\|_\infty + |s| \|b\|_\Lone) + s\|b\|_\Ltwo$,
then by assumption, $\|b\|_\Ltwo \leq d_1$, 
$K_{\tz,b,t,s} = O(|t| + |s|)$. 
By Lemma \ref{lemma:bounding-R2-using-supnorm},  
\begin{align*}
& \sup_{\eta \in A_n} |R_{n,2}(\eta_h, \eta_0) - R_{n,2} (\eta, \eta_0) - s\rn B_3(\eta, \eta_0)|\\
& \quad \lesssim 
\frac{|t|^3 + |s|^3}{\rn} + (t^2 + s^2) \zeta_n + |s|p^2 \rn \epsilon_n 2^{-L_n} + \rn \epsilon_n^2 (|t| + |s|). 
\end{align*}
Since $(|t|^3 + |s|^3)/\rn = o(1)$ as $t, s \leq \log n$, $n\to \infty$, the upper bound in the last display can be simplified to $(t^2 + s^2) \zeta_n + \rn \epsilon_n^2 (|t| + |s|) + |s| p^2 \rn \epsilon_n 2^{-L_n} + o(1)$. 

What left is to bound the last term in the product in (\ref{pf-cond1-2}), by the change of variable condition \ref{cond:C2},
$$
\frac{\int_{A_n} e^{\ell_n(\eta) - \ell_n(\eta_0)} d\Pi(\eta)}
{\int e^{\ell_n(\eta) - \ell_n(\eta_0)} d\Pi(\eta)}
\lesssim e^{C_1(1+t^2 + s^2)}
$$
for some constant $C_1$.
Also, by plugging-in the expression of $\Sigma_{\tz, b}$, we obtain
$h'\Sigma_{\tz,b}h \lesssim t^2 \|\tz\|_\infty^2 + s^2 \|b\|_\Ltwo^2 \lesssim t^2 + s^2$. 
By collecting all the relevant upper bounds derived above, we then complete the proof.
\end{proof}

\subsubsection{Tightness at rate $1/\rn$ for the hazard rate}
\label{sec:sqrt-n-rate}

We verify (\ref{eqn:cond2}) in the tightness criterion.
Consider the function $f = \lambda$ and the centering $T_n^f = T_n^\ld$,
we need to show there exists a divergence sequence $\bar w = (\bar w_l) \to \infty$ and $\bar w_l \geq \sqrt{l}$ such that 
$$
\mathbb{E}\left[
\|\lambda - T_n^\lambda\|_{\mathcal{M}_0(\bar w)} \given X\right] = O_{P_{\eta_0}}(1/\rn).
$$
We choose $\bar w_l = w_l/l^{1/4}$ such that $\rn \epsilon_n 2^{-L_n} \lesssim \bar w_l 2^{-l/2}$
and denote
$$
T_n^\ld = \lambda_{0, L_n} + \frac{1}{\rn} \sum_{L \leq L_n} \sum_{0 \leq K \leq 2^L}
W_n^{(2)} (\psi_{LK}) \psi_{LK},
$$
where 
$W_n^{(2)}(\psi_{LK}) = W_n(-\tIinv \Ld_0\{\psi_{LK}\gamm\}, \psi_{LK}/M_0 + \gamm' \tIinv \Ld_0\{\psi_{LK}\gamm\})$. 
Applying the inequality $\mathbb{E}(x) \leq M + \int_M^\infty P(x \geq \varkappa) d\varkappa$ 
for a constant $M > 0$ and any real-valued variable $\varkappa$ and by the definition of $\mathcal{M}_0(w)$-norm, 
we arrive at
\begin{align*}
\mathbb{E}\big[ 
\rn \| & \ld - T_n^\ld \|_{\mathcal{M}_0(\bar w)} \given X
\big] 
\leq 
M + \int_M^\infty P \left(\rn \| \ld - T_n^\ld \|_{\mathcal{M}_0(\bar w)} \geq \varkappa \given X \right) d\varkappa\\
&  \leq 
M + \int_M^\infty P \left(\rn \max_{l \leq L_n} \bar w_l^{-1} \max_{0 \leq k < 2^{l}} |\langle \ld - T_n^\ld, \psi_{lk} \rangle| \geq \varkappa \given X \right) d\varkappa\\
&  \leq M + 
\sum_{l \leq L_n} \sum_{k =0}^{2^l - 1} \int_M^\infty 
P \left( z_l^{-1} \rn | \langle \ld - T_n^\ld, \psi_{lk} \rangle | > \sqrt{l} \varkappa \given X \right) d\varkappa\\
&  \leq 
M + \sum_{l \leq L_n} \sum_{k=0}^{2^l - 1} \int_M^\infty e^{-\sqrt{l} \varkappa \sqrt{l}} \mathbb{E} \left[
e^{\sqrt{l} z_l^{-1} \rn |\langle \ld - T_n^\ld, \psi_{lk} \rangle|} \given X
\right] d\varkappa,
\end{align*}
where $z_l = \bar w_l/\sqrt{l}$.
The last inequality in the last display is obtained by simply applying Markov's inequality. 
Let $s = \sqrt{l}/z_l$ and $b = \psi_{lk}$ and applying Proposition \ref{prop-bvm}, the logarithm of the expectation in the last line of the last display can be further bounded by
$$
C_1 (1 + s^2 + |s| (\rn \epsilon_n + O_{P_{\eta_0}}(1))\|\gamb - \gambln\|_\Linfty + |s| \rn \epsilon_n^2) + |s| o_{P_{\eta_0}}(\zeta_n) + s^2 \zeta_n + |s| p^2 \rn \epsilon_n 2^{-L_n},
$$
for some positive constant $C_1$ with $s = \sqrt{l}/z_l$ and 
$b = \psi_{lk}$. 
To bound the last display,
by the first point of Lemma \ref{lemma:psib-2}, $\|\gamb - \gambln\|_\Linfty \leq 2^{l/2 - L_n}$. 
Thus $|s| \rn \epsilon_n\|\gamb - \gambln\|_\Linfty \leq |s| \rn \epsilon_n 2^{l/2 - L_n}
\leq \sqrt{l} \rn \epsilon_n 2^{-L_n} 2^{l/2} / z_l \leq l$, as $s = \sqrt{l}/z_l$ and $\bar w_l \geq \sqrt{l}$. 
Also, $s^2 = l/z_l^2 = l^2/ \bar w_l^2 \leq l$ and $|s| \leq l$.
Then by assumptions $\rn \epsilon_n^2 |s|= o(1)$ for $|s| < l < L_n$, $L_n \rn \epsilon_n 2^{-L_n} = o(1)$, and $\zeta L_n^2 = o(1)$, the last display is bounded by $C_2 l$.
Thus, the last line in the penultimate display is bounded by $M+ C_3\sum_{l \leq L_n} \sum_k \int_M^\infty e^{-l \varkappa + C_2l} d \varkappa$ for some constants $C_2$ and $C_3$. The second term in the summation is a constant if choosing $M > C_3$ for a large enough but fixed $M$. This leads to $\mathbb{E}[\rn \|\ld - T_n^\ld\|_{\mathcal{M}_0(\bar w)} \given X] \leq M + O(1)$. Thus we verified (\ref{eqn:cond2}).

%Next, we also need to verify (\ref{eqn:cond2-theta}). 
%By using the inequality $\mathbb{E}(x) \leq M + \int_M^\infty P(x \geq \varkappa) d\varkappa$ again, and let $M= 0$, 
%we obtain
%\begin{align*}
%\mathbb{E} \left(
% \rn |\theta'a - {T_n^\theta}'a| \given X
%\right)
%& \leq \int_0^\infty P(\rn |\theta'a - {T_n^\theta}'a| \geq \varkappa \given X ) d\varkappa\\
%& \leq \int_0^\infty e^{-\varkappa} \mathbb{E}
%\left(
%e^{\rn |\theta'a - {T_n^\theta}'a|} \given X
%\right) d\varkappa.
%\end{align*}
%Applying Proposition \ref{prop-bvm} and choosing $s = 0$ and $t = 1$ (if $\theta'a \geq {T_n^\theta}'a$) or $t = -1$ (if $\theta'a < {T_n^\theta}'a$), one obtains that 
%$\mathbb{E}
%\left(
%e^{\rn |\theta'a - {T_n^\theta}'a} \given X
%\right) \leq C'$ for some constant $C'$.
%Therefore, the last display is bounded by $C' \int_0^\infty e^{-\nu} d\nu = C'$.
%We thus verified (\ref{eqn:cond2-theta}).

\subsubsection{Proof of the main theorem}
\label{pf-npbvm}

With the tightness criterion established in Section \ref{sec:sqrt-n-rate}, what left is to check (\ref{eqn:cond1}). It is sufficient to check Theorem \ref{thm:bvm-linear-fns} holds by letting
$a = \tz$ for any $\tz \in \mathbb{R}^p$ and $b = \psi_T$ with $\psi_T = \sum_{(l, k) \in T} t_{lk} \psi_{lk}$ for any finite set of indices $T$ and $t_{lk} \in \mathbb{R}$.
By following the proofs in Section \ref{cov-laplace}, \ref{cond:C1} holds for any $b = \psi_T$, as $\psi_T \in \mathcal{V}_\mathcal{L}$. What remains is to verify \ref{cond:B}. 
By the first point in Lemma \ref{lemma:psib-2}, $\|\gamb - \gambln\|_\infty \lesssim 2^{-L_n}$ and hence \ref{cond:B} holds as we assume $\rn \epsilon_n 2^{-L_n} = o(1)$. Therefore, (\ref{eqn:cond1}) is verified.

%Last, we check (\ref{eqn:cond2-theta}). Using the fact that $\|a\|\leq \|a\|_1$ for any $a \in \mathbb{R}^p$, $\mathbb{E}\left[ \|\theta - T_n^\theta \| \given X \right] \leq \mathbb{E} \left[ \|\theta - T_n^\theta\|_1 \given X \right] = \sum_{j=1}^p \mathbb{E} \left[ |\theta_j - T_{n,j}^\theta | \given X \right]$.
%Apply Jensen's inequality, for each $j$, 
%$$
%\mathbb{E} \left[ |\theta_j - T_{n,j}^\theta | \given X \right] \leq \log \mathbb{E} \left[ e^{|\theta_j - T_{n,j}^\theta|} \given X \right].
%$$
%The last display is bounded by $O_{P_{\eta_0}}(1/\rn)$ by applying Theorem \ref{thm:bvm-linear-fns} letting $b = 0$ and $a = a_j = (0, \dots, 1, \dots, 0)$ (i.e., a vector at the $j$-th location is 1, the rest are $0$). Note that $p$ is a fixed number, we thus verified (\ref{eqn:cond2-theta}) and hence complete the proof. 

%%%%%%%%%%%%%%%%%%%%%%%%%%%%%%%%%
%%%%%%%%%%%%%%%%%%%%%%%%%%%%%%%%%
%%%%%%%%%%%%%%%%%%%%%%%%%%%%%%%%%
\section{Proof of the Bayesian Donsker theorem}
\label{pf-donsker}

\subsection{Proof of Theorem \ref{thm:donsker-joint}}

First, consider the Haar wavelet prior in \ref{prior:W}. Define the primitive of $T_n^\ld(\cdot)$ as $\mathbb{T}_n^\ld(\cdot) = \int_0^\cdot T_n^\ld(u) du$. 
Then, the `integration' map
\begin{align}
\label{integration-map}
L: \{h_{lk}\} \to L_t(\{h_{lk}\}) = \sum_{l,k} h_{lk} \langle \psi_{lk}, \mathbbm{1}_{[0,t]} \rangle
= \langle h, \mathbbm{1}_{[0,t]} \rangle = \int_0^t h(u) du,
\end{align}
for $t \in [0,1]$,
is linear and continuous from $\mathcal{M}_0(w)$ to $\mathcal{C}([0, 1])$ and $\|\cdot\|_\infty$ for $h \in L^2([0,1])$ with wavelet coefficients $\{h_{lk}\}$ (see Page 1955 of \citet{cast14a}).
It suffices to apply the continuous mapping theorem to $L$ and $\|\cdot\|_\infty \circ L$ for the nonparametric part in the joint posterior distribution. 
{We now invoke Theorem \ref{thm:nonparam-joint-bvm}, which its proof uses Proposition \ref{prop-tightness}, by} checking that the limiting distribution under the map $L$, i.e., $\{\mathbb{Z}_{Q_0} - \mathbb{V}' (\gamm \ld_0)\} \circ L^{-1}$, coincides with $[0, 1] \ni t \rightarrow \mathbb{B}(U_0(t)) - \mathbb{V}'\Ld_0\{\gamm\}(t)$, which follows from Lemma \ref{gp}, then the two claimed processes converges in distribution.

{{Now we focus on}} the random histogram prior in \ref{prior:H}. The proof is similar, we also need to check the condition $\rn \epsilon_n 2^{-L_n} = o\left(\min_{l \leq L_n}\{l^{-1/4} 2^{-l/2} w_l\}\right)$, which holds by choosing $w_l = 2^{l/2}/(1+l^2)$, then $\rn \epsilon_n 2^{-L_n} = o(L_n^{-9/4})$. 

\begin{lemma}
\label{gp}
The Gaussian process $[0, 1] \ni t \rightarrow \{\mathbb{B}(U_0(t)) - \mathbb{V}'\Ld_0\{\gamm\}(t)\}$ and 
$[0, 1] \ni t \rightarrow \{\mathbb{Z}_{Q_0} - \mathbb{V}'\gamm \ld_0\} \circ L^{-1}_t$ coincide, where $L_t$ is the integration map defined in (\ref{integration-map}). 
\end{lemma}

\begin{proof}
We check the respective reproducing kernel Hilbert space (RKHS) attached to the two Gaussian processes coincide. 
This is straightforward by noting that $\mathbb{V}$ and $\mathbb{B}(\cdot)$ (and $\mathbb{Z}_{Q_0}$) are independent and $\mathbb{V}$ is a mean-zero multivariate normal density. 
\end{proof}

\subsection{Proof of Corollary \ref{cor-donsker}}

Followed by Theorem \ref{thm:donsker-joint}, it is sufficient to show that $\rn \| T_n^\theta - \hat \theta \|_\infty = o_{P_{\eta_0}}(1)$ and $\rn \|\mathbb{T}_n^\ld (\cdot) - \hat \Ld (\cdot)\|_\Linfty = o_{P_{\eta_0}}(1)$. 
Note that both $\hat \theta$ and the Breslow estimator $\hat \Ld$ are efficient estimators. In other words, they are both asymptotically linear in their efficient influence function respectively. One can quickly check from Section VIII.4.3 of \citet{andersen93} that $\rn \|\hat \theta - \theta_0 - W_n^{(1)}(1)\|_\infty = o_{P_{\eta_0}}(1)$, thus we have $\rn \| T_n^\theta - \hat \theta \|_\infty = o(1)$.
Also, from Section VIII.4.3, 
$$
\sup_{t \in [0,1]} \left|
\rn \left(\hat \Ld(t) - \Ld_0(t)\right) - W_n^{(2)}(\mathbbm{1}_{[0,t]})
\right| = o_{P_{\eta_0}}(1).
$$
Lemma \ref{lemma-centering-Lambda} shows that $\rn \|\mathbb{T}_n^\ld(\cdot) - \Ld^\star(\cdot)\|_\Linfty = o_{P_{\eta_0}}(1)$ for $\Ld^\star(t) = \Ld_0(t) + W_n^{(2)}(\mathbbm{1}_{\cdot \leq t})/\rn$.

%%%%%%%%%%%%%%%%%%%%%%%%%%%%%%%%%
%%%%%%%%%%%%%%%%%%%%%%%%%%%%%%%%%
%%%%%%%%%%%%%%%%%%%%%%%%%%%%%%%%%
\section{Proof of the supremum-norm rate}
\label{pf-supnorm}

\subsection{Proof of Theorem \ref{thm:supnorm-rate}}

By applying the triangular inequality for $\ell_\infty$-norm, one obtains
$$
\|\ld e^{\theta' \tz} - \ld_0 e^{\theta_0' \tz}\|_\Linfty
\leq \|\ld\|_\Linfty |e^{\theta'\tz} - e^{\theta_0'\tz}| + \|\ld - \ld_0\|_\Linfty e^{\theta_0'\tz}.
$$
The first term can be bounded by
$(\|\ld - \ld_0\|_\Linfty + \|\ld_0\|_\Linfty) |e^{\theta'\tz} - e^{\theta_0'\tz}|$. From Lemma \ref{lemma-4}, we have $\|\theta - \theta_0\| \leq \epsilon_n$, where $\epsilon_n$ is the Hellinger rate, hence $
|e^{\theta'\tz} - e^{\theta_0'\tz}| \lesssim \|\theta - \theta_0\|\|\tz\| \lesssim \epsilon_n$. 
By \ref{asp:iii}, $\|\ld_0\|_\infty \leq c_6$. 
From Lemma \ref{lemma-4} and since $\beta > 1/2$, $\zeta_n = o(1)$, we apply Taylor's theorem to obtain
$\|e^{r - r_0} - 1\|_\infty \lesssim  \|r - r_0\|_\infty \leq \zeta_n$.
Therefore, the first term in the last display is bounded by a constant times $\epsilon_n + \epsilon_n \zeta_n \leq (1 + o(1))\epsilon_n$. 

Since $\tz$ is fixed and hence bounded, by \ref{asp:ii}, $\|\theta_0\| \leq c_2$, by Lemma \ref{lemma-4}, the second term is bounded by a constant times $\zeta_n$. 
Since $\zeta_n \geq \epsilon_n$, the previous display is bounded by some constant times $\zeta_n$. 
This rate is not optimal, as $2^{L_n/2} \epsilon_n$ can be large for a divergent sequence $L_n \to \infty$, e.g., $L_n$ in (\ref{choicel}). 
In the following lemma, we obtain a sharper rate via invoking our nonparametric BvM result.

The following quantity will be used in the next lemma: Define 
\begin{align}
\label{ld-star-psi-lk}
\langle \ld^\star, \psi_{LK} \rangle
= \begin{cases}
\langle \ld_0, \psi_{LK} \rangle + W_{n, LK}^{(2)}(\psi_{LK}), &\quad \text{if}\  L \leq L_n\\
0, &\quad \text{if} \ L > L_n,
\end{cases}
\end{align}
where $0 \leq K < 2^{L}$ and 
$W_{n, LK}^{(2)}(\psi_{LK}) = \langle W_{n}^{(2)} (\psi_{LK}), \psi_{LK} \rangle$ with
$W_{n}^{(2)} (\cdot)$ defined in (\ref{Wn-2}).
We denote $\ld_{L_n}$ as the orthogonal projection of $\lambda$ onto $\mathcal{V}_{L_n} = \text{Vect}\{\psi_{lk}, \ l \leq L_n, \ 0 \leq k \leq 2^{l}\}$, that is the element of $\mathcal{V}_{L_n}$ 
of coordinates $\{\ld_{lk}\}$ in the basis $\{\psi_{lk}\}$. Similar notations are used for $\ld^\star_{L_n}$, $\ld_{0,L_n}$, and $P_{L_n} W_n^{(2)}(\psi_{lk})$. 

\begin{lemma}
\label{lemma:supnorm-lambda}
Under the same conditions as in Theorem \ref{thm:supnorm-rate}, 
for 
$
\xi_n = \sqrt{L_n2^{L_n}/n} + 2^{-\beta L_n} + \epsilon_n,
$ 
then
$
\Pi(\|\ld - \ld_0\|_\Linfty > \xi_n \given X) = o_{P_{\eta_0}}(1).
$
\end{lemma}

\begin{proof}
Applying the triangle inequality for $\ell_\Linfty$-norm, we have
\begin{align*}
\|\ld - \ld_0\|_\Linfty 
& \leq 
\|\ld_{L_n^c}\|_\Linfty + \|\ld_{L_n} - \ld_0\|_\Linfty \\
& \leq 
\underbrace{
\|\ld_{L_n^c}\|_\Linfty}_{(I)} 
+
\underbrace{
\|\ld_{0, L_n^c}\|_\Linfty}_{(II)} 
+
\underbrace{
\|\lambda_{0,L_n} - \lambda^\star_{L_n} \|_\Linfty}_{(III)}
+ 
\underbrace{
\|\lambda^\star_{L_n} - \ld_{L_n} \|_\Linfty}_{(IV)}.
\end{align*}

Term $(I)$ is almost surely zero as for any draw of $\lambda$, the prior is truncated at the level of $L_n$. Note that since $r_H = \Psi r_S$ (see (\ref{haar-transf})), 
the inner product $\langle \lambda, \psi_{lk} \rangle$ is zero when $l \geq L_n$ implies that 
$\langle \lambda, \psi_{lk}^H \rangle$ is also zero,
where $\psi_{lk}^H$ is the $l,k$-th bases function in $\Psi r_S$.

Term $(II)$ can be bounded by $\sum_{L > L_n} 2^{l/2}\max_k|\langle \ld_0, \psi_{lk} \rangle|
\leq \sum_{L > L_n} 2^{-l/2} 2^{-l(1/2+ \beta)} \lesssim 2^{-\beta L_n}$ as $\ld_0$ is $\beta$-H\"older. 

By invoking Lemma \ref{lemma:bounding-iii} and assumptions \ref{asp:i}-\ref{asp:v},
we obtain $(III) \lesssim \sqrt{L_n2^{L_n}/n}$. 

Last, to bound $(IV)$, we introduce the set
$$
E_n = A_n \cap \{\|\ld - \ld_0\|_\Linfty \leq \zeta_n\}.
$$
Recall that $\zeta_n = 2^{L_n/2} \epsilon_n + 2^{-\beta L_n}$. 
Denote $\mathbb{E}^{\Pi_n}_{\eta_0}$ as the expectation under the posterior of $\eta_0$ conditional on $E_n$, then conditioning on the $E_n$, we have
\begin{align*}
\mathbb{E}^{\Pi_n}_{\eta_0} \max_{0 \leq k < 2^{l}} \rn |\langle \ld - \ld^\star, \psi_{lk} \rangle|
\leq 
\frac{1}{s} \log \sum_{k=1}^{2^l - 1}
\mathbb{E}^{\Pi_n} \left(e^{s\rn \langle \ld - \ld^\star, \psi_{lk} \rangle} 
+ e^{-s\rn \langle \ld - \ld^\star, \psi_{lk} \rangle} \right).
\end{align*}
To bound the last display, we first bound $\mathbb{E}^{\Pi_n} e^{s\rn \langle \ld - \ld^\star, \psi_{lk} \rangle}$ (the other part,
$\mathbb{E}^{\Pi_n} e^{-s\rn \langle \ld - \ld^\star, \psi_{lk} \rangle}$, can be bounded using a similar strategy).
{{The proof starts as that of Theorem \ref{thm:bvm-linear-fns} (with $t= 0$). We take $s$ equal to $\sqrt{l}$ to control uniformly all Laplace transforms for $l \le L_n$.}}
One can write 
\begin{align}
& \mathbb{E}^{\Pi_n} \big(e^{s\rn \langle \ld - \ld^\star, \psi_{lk} \rangle} \given X, E_n \big)
= 
\frac{
\int_{E_n} e^{\ell_n(\eta) - \ell_n(\eta_h) + s\rn \langle \ld - \ld^\star, \psi_{lk} \rangle} d\Pi(\eta)
}{
\int e^{\ell_n(\eta) - \ell_n(\eta_0)} d\Pi(\eta)
}  \nonumber \\
& \leq 
\exp\left(
\sup_{\eta \in E_n} 
\left|
D_n + \rn s B_2(\eta, \eta_0) + R_n(\eta, \eta_0) - R_n(\eta_h, \eta_0) - \rn s B_3(\eta, \eta_0)
\right|
\right)  
\label{pf-supnorm-1}
\\
& \quad \times 
\frac{
\int_{E_n} e^{\ell_n(\eta_h) - \ell_n(\eta_0)} d\Pi(\eta)
}{
\int e^{\ell_n(\eta) - \ell_n(\eta_0)} d\Pi(\eta)
},
\label{pf-supnorm-2}
\end{align}
where $D_n$, $B_2(\eta, \eta_0)$, $B_3(\eta, \eta_0)$, and $R_n$ are defined in the proof of Theorem \ref{thm:bvm-linear-fns} in Section \ref{sec:support-lemma-thm1} with the choice $b = \psi_{lk}$. 
By \ref{cond:C2}, (\ref{pf-supnorm-2}) is bounded by $\exp(C_6(1 + s^2))$ for some constant $C_6$. 

To bound (\ref{pf-supnorm-1}), first, we bound $|D_n|$.
From Lemma \ref{lemma:bounding-Dn}, using the inequality $\|fg\|_1 \leq \|f\|_1 \|g\|_\infty$ and the triangular inequality for $\ell_1$-norm, we obtain
$$
|D_n| \lesssim s^2 \|\gamb-\gambln\|_\infty (\|\gamb\|_\Lone +\|\gambln\|_\Lone) + s^2 p^2 2^{-L_n}\|\gambln\|_\Lone + s^2 p^4 2^{-2L_n},
$$
where $b = \psi_{lk}$. 
By Lemma \ref{lemma:psib-2}, 
{as $s^2 = l \leq L_n$,
$$
|D_n| \lesssim s^2 2^{l/2 - L_n}+ s^2 2^{-L_n} \leq L_n 2^{-L_n/2} + L_n 2^{-L_n}.
$$
The last display is $o(1)$.
}
Next, using Lemma \ref{lemma:bounding-B2} and replacing $A_n$ with $E_n$, we obtain 
$$
\sup_{\eta \in E_n} |B_2(\eta, \eta_0)| 
\lesssim \epsilon_n \|\gamb - \gambln\|_\Linfty + p^2 \epsilon_n 2^{-L_n} \|b\|_\Lone.
$$
From the first point of Lemma \ref{lemma:psib-2}, 
$\|\gamb - \gambln\|_\Linfty \lesssim 2^{l/2}2^{-L_n}$;
also, 
$\|\psi_{lk}\|_1 \leq 2^{-l/2}$. Then, for $b = \psi_{lk}$, the last display is bounded by 
$\epsilon_n 2^{l/2 - L_n} + \epsilon_n 2^{-L_n-l/2}$
Thus, we obtain
$$
s\rn \sup_{\eta \in E_n} |B_2(\eta, \eta_0)| \lesssim {\sqrt{ln}} \epsilon_{n} 2^{l/2 - L_n}
{\leq \sqrt{L_n n} \epsilon_n 2^{-L_n/2} = o(1)}.
$$
Last, we bound the two remainder terms. 
First, we use Lemma \ref{R1-version2} (also, replacing $A_n$ with $E_n$). With the assumptions \ref{asp:i}-\ref{asp:v}, 
{
for $\Delta_1$ and $\Delta_{2,L_n}$ given in \eqref{K1} and \eqref{K2} respectively, we have
\begin{align*}
& \|\Delta_1\|_\infty /\rn \leq |s| p^2\|\Ld_0\{b\gamm\}\|_\infty \|\tIinv\|_{(\infty, \infty)}/\rn
\lesssim \sqrt{l}/\rn \leq \sqrt{L_n}\rn = o(1),\\
& \|\Delta_{2,L_n}\|_\infty/\rn \lesssim (p^2|s|L_n + |s|L_n)/\rn \leq (p^2 +1)L_n^2/\rn = o(1),
\end{align*}
thus the conditions in Lemma \ref{R1-version2} are verified.}
Hence, we obtain 
$
\sup_{\eta \in E_n} |R_{n,1}(\eta, \eta_0) - R_{n,1}(\eta_h, \eta_0)|= O_{P_{\eta_0}}(s (1 + \zeta_n)).
$
Next, using Lemma \ref{lemma:bounding-R2-using-supnorm}, 
note that $\tilde K_{a,b,t,s} \lesssim s \|\psi_{lk}\|_\Ltwo = s$ as $\|\psi_{lk}\|_\Ltwo^2 = 1$, we have 
\begin{align*}
& \sup_{\eta \in E_n} |R_{n,2} (\eta, \eta_0) - R_{n,2}(\eta_h ,\eta_0) - \rn s B_3(\eta, \eta_0)|  \\
& \quad \lesssim s^3/\rn + s^2 \zeta_n + |s| \rn \epsilon_n 2^{-L_n} + |s| \rn \epsilon_n^2 \\
& \quad \leq s^2 (1+ \zeta_n) + o(1), 
\end{align*}
as $|s|/\rn \leq 1$ since {$|s| \leq \sqrt{L_n} \lesssim \sqrt{\log n}$} and $\rn \epsilon_{n}^2 L_n = o(1)$ by \ref{cond:P}. 

A simple calculation reveals that $s^2 \Ld_0\{\psi_{lk}^2 M_0\} + s^2 \Ld_0\{\psi_{lk} \gamm'\} \tIinv \Ld_0\{\psi_{lk} \gamm'\}
\lesssim s^2 2^{-l/2}$, then, by combining all the relevant bounds obtained above, 
let $s_l = s = \sqrt{l} \leq \sqrt{L_n}$ and $q_n = \rn \epsilon_n 2^{l/2 - L_n}$, 
then, for some constants $C_1$ and $C_2$, we have
\begin{align*}
& \int \|\ld - \hat \ld\|_\infty d\Pi_n(\eta) \\
& \quad \lesssim \frac{1}{\rn} \sum_{l \leq L_n} \frac{2^{l/2}}{s_l} \log 
\left\{
2C_1(1+o_{P_{\eta_0}}(1)) 
\sum_{k=1}^{2^l - 1} 
e^{C_2(s_l^2 (1+\zeta_n) + s_l  q_n + O_{P_{\eta_0}}(s_l (1 + \zeta_n))}
\right\}\\
& \quad \lesssim 
\frac{1}{\rn} \sum_{l\leq L_n} \frac{2^{l/2}}{s_l} \left[
l + s_l^2 ( 1+ \zeta_n) + s_l q_n + s_l O_{P_{\eta_0}}(1)
\right] \\
& \quad \lesssim 
\sqrt{\frac{L_n 2^{L_n}}{n}} (1+ \zeta_n) + \epsilon_n + \sqrt{\frac{2^{L_n}}{n}} O_{P_{\eta_0}}(1) \\
& \quad \leq\epsilon_n (1+ \zeta_n) +  \epsilon_n + \frac{\epsilon_n}{L_n}
O_{P_{\eta_0}}(1) \\
& \quad \lesssim\epsilon_n \left(1 + o(1) + o_{P_{\eta_0}}(1)\right). 
\end{align*}
Thus, the last display can be bounded by a constant times $\epsilon_n$.
By combining the upper bounds of $(I)$, $(II)$, $(III)$, and $(IV)$, 
we obtain $\|\ld - \ld_0\|_\infty \lesssim 2^{-\beta L_n} + \sqrt{L_n 2^{L_n}/n} + \epsilon_n = \xi_n$.
\end{proof} 

%%%%%%%%%%%%%%%%%%%%%%%%%%%%%%%%%%%%%%%%%%%
\begin{lemma}
\label{lemma:bounding-iii}
Let $\ld_{0,L_n}$ be the orthogonal projection of $\ld$ onto $\mathcal{V}_{L_n} = \text{Vect}\{\psi_{lk}, \ l\leq L_n, \ 0 \leq k < 2^{l}\}$ and $\ld^\star_{L_n}$ be the orthogonal projection of $\ld^\star$ onto 
$\mathcal{V}_{L_n}$, where the inner product $\langle \ld^\star, \psi_{lk} \rangle$ is defined in (\ref{ld-star-psi-lk}). Then,
$$
\mathbb{E}_{\eta_0} \|\ld_{0,L_n} - \ld_{L_n}^\star \|_\infty \lesssim  \sqrt{\frac{L_n2^{L_n}}{n}}.
$$

\end{lemma}
\begin{proof}
By the definition of $\lambda^\star_{L_n}$ with the inner product $\langle \ld^\star, \psi_{lk}\rangle$ is defined in (\ref{ld-star-psi-lk}), we have
\begin{align}
\mathbb{E}_{\eta_0} & \|\ld_{0,L_n} - \ld^\star_{L_n}\|_\Linfty
\leq 
\frac{1}{\rn} \sum_{l \leq L_n} \frac{2^{l/2}}{t_l} \sum_{k=1}^{2^l - 1}
\mathbb{E}_{\eta_0} | t_l \langle W_{n,lk}^{(2)}(\psi_{lk}), \psi_{lk} \rangle|
\nonumber \\
& \quad \quad \leq \frac{1}{\rn} \sum_{l\leq L_n} \frac{2^{l/2}}{t_l}
\log \sum_{k=1}^{2^l - 1} \mathbb{E}_{\eta_0}
\left(
e^{t_l \langle W_n^{(2)}(\psi_{lk}), \psi_{lk} \rangle} 
+ e^{- t_l \langle W_n^{(2)}(\psi_{lk}), \psi_{lk} \rangle}
\right),
\label{pf-supnorm-3}
\end{align}
where the second inequality is obtained by using the inequality $\mathbb{E}(x) \leq \log \mathbb{E}(e^x)$.
It suffices to bound $\mathbb{E}_{\eta_0} e^{t_l \langle W_n^{(2)}(\psi_{lk}), \psi_{lk} \rangle}$ for bound  the expectation term in (\ref{pf-supnorm-1}), as bounding the term written with a negative sign is similar. 

Denote
$$
H_{lk}(X_i) = \delta_i (g_{lk}(Y_i)'Z_i + h_{lk}(Y_i)) - e^{\theta_0'Z_i} \Ld_0 (g_{lk}(Y_i)'Z_i + h_{lk}(Y_i)),
$$
where $X_i$ is the triple $(\delta_i, Y_i, Z_i)$, 
$$
g_{lk}(\cdot) = \tIinv \Ld_0\{\psi_{lk} \gamm\} \int_0^\cdot \psi_{lk}(u) du,
$$
$$h_{lk}(\cdot) = \gamma_{n,lk}(\cdot) - \gamma_{M_1,lk}(\cdot)' \tIinv \Ld_0\{\psi_{lk}(\cdot) \gamm(\cdot)\},
$$
$\gamma_{n,lk} = \langle \gamma_{lk}, \psi_{lk} \rangle$ and 
$\gamma_{M_1, lk} = \langle \gamma_{lk}M_1, \psi_{lk} \rangle$.
Then $W_{n,lk}^{(2)}(\psi_{lk}) = \frac{1}{\rn} \sum_{i=1}^n H_{lk}(X_i)$.
By construction, $\mathbb{E}_{\eta_0}(H_{lk}(X_i)) = 0$. Then for each $X_i$, 
\begin{align*}
\mathbb{E}_{\eta_0} \left(\frac{e^{t_l H_{lk}(X_i)}}{\rn}\right)
& = \mathbb{E}_{\eta_0} \left(
\sum_{k\geq 0} \left( \frac{t_l H_{lk}(X_i)}{\rn} \right)^k \frac{1}{k!} 
\right)
= 1 + \mathbb{E}_{\eta_0} \left(
\sum_{k\geq 2} \left(\frac{ t_l H_{lk}(X_i) }{\rn}\right)^k \frac{1}{k!}
\right) \\
& \leq 1 + 
\sum_{k \geq 2} \left(
\frac{| t_l |\|H_{lk}(X_i)\|_\Linfty}{\rn}
\right)^{k-2}
\frac{ t_l^2 \mathbb{E}_{\eta_0} (H_{lk}^2(X_i))}{n k!}  \\
& \leq 1 + 
\frac{ t_l^2}{2n} \mathbb{E}_{\eta_0} (H_{lk}^2(X_i))
\exp \left( \frac{| t_l| \|H_{lk}(X_i)\|_\Linfty}{\rn} \right).
\end{align*}
First, using the following three results:
$\|0, \gambln\|_L \to \|0, \gamb\|_L$ and $\|0, \gammln\|_L \to \|0, \gamm\|_L$ by Lemma \ref{lemma:bound-W}, $\|0, \gamb\|_\Ltwo \leq c\|b\|_\Ltwo = O(1)$ by Lemma \ref{lemma:psib-1},
and $\|\psi_{lk}\|_\Lone \lesssim 2^{-l/2}$, as $\psi_{lk}$ is a Haar bases, we have
$$
|g_{lk}'z| \leq p^2\|z\|_\infty \|\tIinv\|_{(\infty,\infty)} \|\Ld_0\{\psi_{lk}\gamm\}\|_\infty \|\mathbbm{1}_{[0,t]}\|_\Linfty \|\psi_{lk}\|_\Lone \lesssim 2^{-l}.
$$
Thus, $\mathbb{E}_{\eta_0} (H^2_{lk}(X_i)) = \|g_{lk}'Z_i, h_{lk}\|_L^2$ is bounded.
Next, by the definition of $H_{lk}(\cdot)$, we have
$$
\|H_{lk}(X_i) \|_\infty \lesssim C_1(1+ \|\Ld_0\|_\Linfty) (C_2 + \|\gamma_{n,lk}\|_\Linfty) \lesssim l 2^{L_n/2} \leq L_n2^{L_n/2},
$$ 
for some constants $C_1$ and $C_2$.
Therefore, 
by plugging the two last two upper bounds in the last line of the previous display, we obtain
$$
\mathbb{E}_{\eta_0} ({e^{t_l H_{lk}(X_i)}}/{\rn}) \leq 1 + \frac{C_3 t_l^2}{2n} e^{|t_l| L_n2^{L_n}/\rn}.
$$
 
Now, let $t_l = \sqrt{l}$ (for the other part in (\ref{pf-supnorm-1}), let $ t_l = -\sqrt{l}$ instead),
then, $C_3t_l^2/(2n) \leq C_3 L_n/(2n) \leq C_5$ as $L_n/n \leq 1$.
By the assumption that $L_n2^{L_n}/\rn \leq C_4$ for some constant $C_4$, 
we arrive at
\begin{align*}
\mathbb{E}_{\eta_0} \|\ld_{0,L_n} - \ld_{L_n}^\star \|_\Linfty
& \lesssim
\frac{1}{\rn} \sum_{l \leq L_n} \frac{2^{l/2}}{\sqrt{l}} \log \left(2^{l+1} e^{C_4 t_l}\right)
= \frac{1}{\rn} \sum_{l\leq L_n} \frac{2^{l/2}}{\sqrt{l}} \left(C_4 \sqrt{l} + (l+1)\log 2\right)\\
&\lesssim \frac{1}{\rn} \sum_{l \leq L_n} \sqrt{l} 2^{l/2} 
\lesssim \sqrt{\frac{L_n2^{L_n}}{n}}.
\end{align*}
\end{proof}

\subsection{Lower bound for the supremum-norm rate of {the hazard function conditional on $z$}}
\label{sec:lower-bound}

\begin{lemma}[Lower bound for the sup-norm rate]
\label{lower-bound}
Let $\beta > 1/2$ and $L > 0$, if $\theta \in [-C, C]^p$ for some large but fixed $C$, then, for a given $\tz \in \mathbb{R}^p$, there exists a finite constant $M = M(\beta, L) > 0$ such that for large enough $n$, 
$$
\inf_{H} \sup_{\substack{\ld \in \mathcal{H}(\beta, L) \\ \theta \in [-C,C]^p}} \mathbb{E}_{\eta} \|H - \lambda e^{\theta'\tz}\|_\Linfty
\geq M \left(\frac{n}{\log n} \right)^{-\beta/(2\beta + 1)}.
$$
\end{lemma}
{The proof of above lemma is fairly similar to that for the hazard rate under right-censoring in \citet{cast20survival}; it is included for completeness.}
\begin{proof}
As $\theta$ is a parametric quantity and $\lambda$ is a nonparametric quantity, an estimator of $\lambda$ typically converges at a much slower rate than an estimator of $\theta$. In the proof, we fix $\theta$ in a compact set and should only consider $\ld$.

We follow the principle of lower bounds approach proposed by \citet{ibragimov77} to prove the result. Let $\ld_0, \dots, \ld_N$ with $N \geq 2$ be baseline hazard functions and recall that $K(P, Q)$ is the Kullback-Leibler divergence. Then, for some $\alpha \in (0, 1/8)$ and $C_\alpha > 0$, which is a constant depends on $\alpha$, the minimax risk is bounded by 
$$
\inf_H \sup_{\substack{\ld \in \mathcal{H}(\beta, L) \\ \theta \in [-C, C]^p}} \mathbb{E}_{\eta} \|H - \ld e^{\theta'\tz}\|_\Linfty \geq C_\alpha s,
$$ 
if the following two conditions are satisfied:
\begin{enumerate}[label=\textbf{(I)}] 
\item \label{lb:1} $\|\ld_i e^{\theta' \tz} - \ld_j e^{\theta' \tz}\|_\Linfty \geq 2 s > 0$, for $0 \leq i < j \leq N$;
\end{enumerate}
\begin{enumerate}[label=\textbf{(II)}]
\item \label{lb:2} $\sum_{j=1}^N K(P_{{(\theta, \lambda_j)}}^{\otimes n}, P_{{(\theta, \lambda_0)}}^{\otimes n}) \leq \alpha N \log N$.
\end{enumerate}
To verify \ref{lb:1}, using \ref{asp:i} and $\theta \in [-C, C]^p$, $e^{\theta'\tz}\|\ld_i - \ld_j \|_\Linfty  \gtrsim \|\ld_i - \ld_j\|_\Linfty$. 
Following the proof of Theorem S-4 of \citetalias{cast20survival}, set $\ld_0 = 1$, i.e., a constant baseline hazard function, and define $\ld_k$ such that 
$$
\ld_k = \ld_0 + Lh^\beta \psi\left(\frac{x- x_k}{h}\right), \ 1\leq k \leq N,
$$
where $x_k = (k - 1/2)/{N}$, $h = 1/N$, and $\psi(\cdot) \in \mathcal{H}(\beta, 1)$ such that $\psi$ has a compact support and $\psi(0) > c$, where $c$ is a small positive constant.
Then, $\ld_k \in \mathcal{H}(\beta, L)$ and $\|\ld_k - \ld_j\|_\Linfty = cLh^\beta$ by the construction in the last display. 
By choosing $h = (\delta \log n/n)^{1/(2\beta+1)}$ for some constant $\delta > 0$, we verified \ref{lb:1} by letting 
{$s= cLh^\beta/2$}. 

To verify \ref{lb:2}, we first obtain the Kullback-Leibler divergence between $P_{{(\theta, \lambda_j)}}$, $1\leq j\leq N$, and $P_{{(\theta, \lambda_j)}}$. 
From (\ref{eqn:density-function}), 
denote $S_{z}^j:= S_{z}^j(\cdot) = \exp(- e^{\theta'z} \Ld_j(\cdot))$, we have 
\begin{align*}
& K(P_{{(\theta, \lambda_j)}}, P_{{(\theta, \lambda_0)}}) \\
& \quad = 
\int_0^1 g_{z} S_{z}^0 \log \left(\frac{S_{z}^j}{S_{z}^0} \right) 
+ \int_0^1 \bar G_{z} \ld_0 e^{\theta_0'\tz} S_{z}^0 \log \left(\frac{\ld_j S_{z}^j}{\ld_0 S_{z}^0} \right) 
+ \bar G_{z}(1) S_{z}^0(1) \log \left(\frac{S_{z}^j(1)}{S_{z}^0(1)} \right).
\end{align*}
The second term in the previous display can be split into two parts by writing 
$\log \left({\ld_j S_{z}^j}/(\ld_0 S_{z}^0) \right) = \log \left({\ld_j}/{\ld_0}\right) + \log \left({S_{z}^j}/{ S_{z}^0}\right)$.
Using the integral by parts for the integral with the second part, the third term in the previous display cancels as $S_{z}^j(0) = S_{z}^0(0) = e^{-\theta'z}$. By rearranging other terms, we obtain 
\begin{align*}
K(P_{{(\theta, \lambda_j)}}, P_{{(\theta, \lambda_0)}})
= \int \int_0^1 \left[\log \left(\frac{\ld_0(u)}{\ld_j(u) }\right) - \frac{\ld_0(u) }{\ld_j(u) } + 1 \right] \ld_0(u) e^{\theta_0'\tz} S_{z}^0(u) \bar G_{z}(u) f_Z(\tz) du d\tz
\end{align*}
From the definition of $\ld_k$ and $h$ above, $\|\ld_j - \ld_0\|_\Ltwo^2 \lesssim L^2h^{2\beta+1} \lesssim \log n/n \leq 1/2$ for some sufficiently large $n$. Therefore, one can apply Taylor's theorem and assumptions \ref{asp:i} and \ref{asp:iv} to bound the last display. The bound is given by 
$
K(P_{{(\theta, \lambda_j)}}, P_{{(\theta, \lambda_0)}}) \lesssim \max_z \|\ld_j - \ld_0\|_\Ltwo^2
$
as $\|\theta\|\leq C$ by assumption and $\|z\|_\infty \leq c_1$ by \ref{asp:i}. 
Therefore, 
$$
\frac{1}{N} \sum_{j=1}^N K(P_{{{(\theta, \lambda_j)}}}^{\otimes n}, 
P_{{{(\theta, \lambda_0)}}}^{\otimes n})
= \frac{1}{N} \sum_{j=1}^N n K(P_{{(\theta, \lambda_j)}}, P_{{(\theta, \lambda_0)}})
\lesssim L^2 n h^{2\beta + 1}. 
$$
{We now claim that by choosing $\delta = CL^{-2}\alpha/(2\beta+1)$ for a sufficiently small but fixed constant $C$, the upper bound in the last display is bounded by $\alpha \log N$ for $\alpha < 1/8$. We thus verified \ref{lb:2}. }

{Here, we prove the claim. Recall that $h = (\delta \log n/n)^{1/(2\beta + 1)}$ and $N = 1/h$. By plugging-in the two expressions, we then need to verify  
$
L^2 \delta \log n \leq \frac{\alpha}{2\beta+1} (\log n - \log \delta - \log \log n),
$
which is the same as to verify
$$
\delta \leq \frac{\alpha}{L^2(2\beta+1)} (1 - \log \delta/\log n - \log\log n/\log n).
$$
Since $n \to \infty$, $\delta, \alpha, \beta, L$ are all fixed constant, $\log \delta/\log n = o(1)$ and 
$\log\log n / \log n = o(1)$. Thus, we obtain
$
\delta \leq C\frac{\alpha}{L^2(2\beta+1)}
$
by choosing $C$ to be sufficiently small (e.g., $C < 1/2$).
}

%The upper bound in the last display can be made smaller than $\alpha \log N$ for $\alpha < 1/8$ by choosing $\delta$ as small as desired. 
%We thus verified \ref{lb:2}. 

Thus, the minimax risk as in the statement of this Lemma is bounded from below by $M {{\nu}}_{n,\beta}$ for some constant $M$ which depends on $\alpha$, $L$, $p$, $C$, and $c_1$.
\end{proof}

%%%%%%%%%%%%%%%%%%%%%%%%%%%%%%%%%
%%%%%%%%%%%%%%%%%%%%%%%%%%%%%%%%%
\section{Proof of Theorem \ref{thm:specific-sup-rate}}
\label{verify-conditions}

In this section, we prove Theorem \ref{thm:specific-sup-rate}. 
One has to verify the conditions \ref{cond:P}, \ref{cond:B}, and the two change of variables conditions \ref{cond:C1} and \ref{cond:C2} for our specific choice of priors in Section \ref{sec:priors}. Since \ref{cond:P} has already been verified in Section \ref{sec:boundparameters}. 
Below we verify the rest three conditions.

%%%%%%%%%%%%%%%%%%%%%%%%%%%%%%%%%%%%%%%
%%%%%%%%%%%%%%%%%%%%%%%%%%%%%%%%%%%%%%%
\subsection{Verifying {\rm \ref{cond:B}}}
\label{verifyingB}

From the second point in Lemma \ref{lemma:psib-2}, 
$\|\gamb - \gambln\|_\Linfty \lesssim 2^{-\mu' L_n}$ with $\mu' = \mu \wedge 1$. 
Then, 
$\rn \epsilon_n \|\gamb - \gambln\|_\Linfty \lesssim \rn \epsilon_n 2^{-\mu'L_n}$. 
By plugging-in the rate $\epsilon_n = {{\nu}}_n$ in Theorem \ref{thm:specific-sup-rate}
and $L_n$ in \eqref{choicel},  
$\rn \epsilon_n 2^{-\mu'L_n} = \rn (n/\log n)^{- (\beta + \mu')/(2\beta+1)} = o(1)$ as $\beta > 1/2$ and $b \in \mathcal{H}(\mu, D)$, $\mu > 1/2$.

%%%%%%%%%%%%%%%%%%%%%%%%%%%%%%%%%%%%%%%
%%%%%%%%%%%%%%%%%%%%%%%%%%%%%%%%%%%%%%%
\subsection{Verifying the two change of variables conditions {\rm \ref{cond:C1}} and {\rm \ref{cond:C2}}}
\label{verify-cov}

This section has two subsections.
In Section \ref{cov-laplace}, we verify the two conditions for the Haar wavelet priors in \ref{prior:W}. The prior for each $Z_{lk}$ is chosen as an independent 1) Laplace density (i.e., $Z_{lk} \sim \text{Laplace}(0, 1)$) and 2) Gaussian density (i.e., $Z_{lk} \sim N(0, 1)$).
Both densities are non-conjugate. 
In Section \ref{cov-histogram}, we verify the two conditions for the random histograms priors \ref{prior:H}, 
including the independent gamma prior, which is conjugate, and the dependent gamma prior, which is non-conjugate. 

\subsubsection{Verifying the two conditions for {\rm \ref{prior:T}} and {\rm \ref{prior:W}}}
\label{cov-laplace}

We first verify \ref{cond:C2} by proving the following:
\begin{align}
\label{cv-1}
    \frac{\int_{A_n} e^{\ell_n(\eta_h) - \ell_n(\eta_0)} d\Pi(\eta)}
    {\int e^{\ell_n(\eta) - \ell_n(\eta_0)} d\Pi(\eta)} 
    \leq e^{C(1 + t^2 + s^2)},
\end{align}
for $A_n$ given in \ref{cond:P} and some constant $C$. 
The verification of \ref{cond:C1} is similar and is given after the proof of \ref{cond:C2}. 

Recall that $\eta_h = (\theta_h, r_h)$, where 
\begin{align*}
& \theta_h = \theta - \frac{t \tIinv a}{\rn} + \frac{s \tIinv \Ld_0\{b \gamm\}}{\rn},\\ 
r_h = r & + \frac{t \gammln' \tIinv a}{\rn} - \frac{s \gambln}{\rn} - \frac{s \gammln' \tIinv \Ld_0\{b \gamm\}}{\rn}.
\end{align*}
For simplicity, let's denote $\Delta_{1} = t \tIinv a - s\tIinv \Ld_0\{b \gamm\}$ and
$\Delta_{2, L_n} = - t \gammln' \tIinv a + s\gambln + s\gammln'\tIinv \Ld_0\{b \gamm\}$
and write $\theta_h = \theta - \Delta_{1}/\rn$ and $r_h = r - \Delta_{2,L_n}/\rn$ accordingly.
We further denote the projection
$\Delta_{2,lk} = \langle \Delta_{2, L_n}, \psi_{lk} \rangle$.
%
%$\Delta_{2,lk} = t \gamma_{M_1, lk}' \tIinv a - s \gamma_{b, lk} - s \gamma_{M_1, lk}' \tIinv \Ld_0\{b \gamm\}$,
%where $\gamma_{M_{1}, lk} = (\langle \gamma_{M_{11}}, \psi_{lk} \rangle, \dots, 
%\langle \gamma_{M_{1p}}, \psi_{lk} \rangle)'$ and $\gamma_{b, lk} = \langle \gamma_b, \psi_{lk} \rangle$.

{\it 1. The Laplace prior on $Z_{lk}$}

Define $\Theta_n = \{\theta: \|\theta - \theta_0\| \leq \epsilon_n\}$ and 
$\mathcal{H}_n = \{\ld: \|\ld - \ld_0\|_\Lone \lesssim \epsilon_n\}$,
then, $A_n = \{(\theta, \ld), \theta \in \Theta_n, \ld \in \mathcal{H}_n\}$.
Using the fact that $d\Pi(r) = \prod_{l \leq L_n;k} d\Pi(r_{lk})$,
the numerator in (\ref{cv-1}) can be written as
\begin{align*}
N_n & = \int_{\Theta_n}\int_{\mathcal{H}_n} 
e^{\ell_n(\eta_h) - \ell_n(\eta_0)} 
\prod_{l\leq L_n;k} \pi(r_{lk}) dr_{lk}
\prod_{j=1}^p \pi(\theta_j) d \theta_j,
\end{align*}
where $r_{lk} = \langle r, \psi_{lk} \rangle$. In fact, we can write $\pi(r_{lk}) = \phi(r_{lk}/\sigma_l)/\sigma_l$, where $\phi(\cdot)$ is denoted as the standard Laplace density.

Let $\vartheta_j = \theta_j - \Delta_{1j}/\rn$ ($\Delta_{1j}$ is the $j$-th coordinate of $\Delta_1$) and $\rho_{lk} = r_{lk} - \Delta_{2,lk}/\rn$ (hence, $\rho_{lk} = \langle \rho, \psi_{lk} \rangle$). By  
applying the change of variables from 
$\theta$ to $\vartheta$ and $r_{lk}$ to $\rho_{lk}$.
Due to the invariance of the Lebesgue measure, $dr_{lk} = d\rho_{lk}$ and $d\vartheta = d\theta$, 
%one can write 
%$\pi(\theta) = \pi(\theta - \Delta_1/\rn + \Delta_1/\rn)$ and $\phi(\sigma_l^{-1} r_{lk}) = \phi(\sigma_l^{-1} (r_{lk} - \Delta_{2,lk}/\rn + \Delta_{2,lk}/\rn) )$. 
%Define $\Theta_h = \{\theta: \|\theta - \theta_0\| \leq \epsilon_n\}$ and 
%$\mathcal{H} = \{\ld: \|\ld - \ld_0\|_\Lone \lesssim \epsilon_n\}$. 
then
\begin{align*}
N_n = & \int_{\Theta_n - \frac{\Delta_1}{\rn}} \int_{\mathcal{H}_n - \frac{\Delta_{2,L_n}}{\rn}} e^{\ell_n((\vartheta, \rho)) - \ell_n((\theta_0, r_0))} 
\prod_{l\leq L_{n}; k} \frac{1}{\sigma_l} \phi\left(\frac{\rho_{lk} + {\Delta_{2,lk}}/{\rn}}{\sigma_l}\right) d\rho_{lk} \\
& \times \prod_{j=1}^p \pi \left(\vartheta_j + \frac{\Delta_{1j}}{\rn}\right)  d\vartheta_j.
\end{align*}

We also apply the change of variable to the denominator in (\ref{cv-1})
and obtain that
$$
D_n = \int\int e^{\ell_n((\vartheta, \rho)) - \ell_n((\theta_0, r_0))} 
\prod_{l \leq L_{n};k}
\frac{1}{\sigma_l} \phi\left(\frac{\rho_{lk}}{\sigma_l}\right)d\rho_{lk} 
\prod_{j=1}^p \pi(\vartheta_j) d\vartheta_j.
$$
To bound the ratio of $N_n$ and $D_n$, using the fact that $\phi(\rho_{lk}/\sigma_l)$ is the Laplace density and hence is Lipschitz, then for some positive constant $C_1$, 
\begin{align*}
\phi\left(\frac{\rho_{lk} + \Delta_{2,lk}/\rn}{\sigma_l}\right) 
& = \phi\left(\frac{\rho_{lk}}{\sigma_l}\right) 
\exp\left(
\log \phi\left(\frac{\rho_{lk}+\Delta_{2,lk}/\rn}{\sigma_l}\right)
- \log \phi\left(\frac{\rho_{lk}}{\sigma_l}\right)
\right)\\
& \leq 
\phi\left(\frac{\rho_{lk}}{\sigma_l}\right) 
\exp\left(\frac{C_1 |\Delta_{2,lk}|}{\sigma_l \rn} \right). 
\end{align*}
For the prior of $\theta_j$, we consider either the uniform prior in $[-C, C]$ and the truncated Subbotin density given by 
$
f(\theta_j) = \frac{\tau \kappa}{2 \Gamma(1/\tau)} e^{- |\kappa \theta_j|^\tau}
$
for any $\tau \in [1, 2]$. For both priors,
\begin{align*}
\pi(\vartheta_j + \Delta_{1j}/\rn) 
& = \pi(\vartheta_j + \Delta_{1j}/\rn)
\leq 
\pi(\vartheta_j) 
\exp\left(
\frac{C_2 |\Delta_{1j}|}{\rn}
\right)
\end{align*}
for some positive constant $C_2$.
Therefore, 
\begin{align}
\label{cv-2}
\frac{N_n}{D_n} 
\leq \Pi\left( \left(\Theta_n - \frac{\Delta_1}{\rn}, 
\mathcal{H}_n - \frac{\Delta_{2,L_n}}{\rn} \right) \given X^n \right)
\exp\left( \sum_{j=1}^p \frac{C_1 |\Delta_{1j}| }{\rn} + \sum_{l \leq L_n;k} \frac{C_2|\Delta_{2, lk}|}{\sigma_l \rn} \right).
\end{align}
We will use the the following results to bound the last display. First, 
\begin{align*}
\|\Delta_1\|_\infty  
&\leq |t|\| \tIinv a\|_\infty + |s| \|\tIinv \Ld_0\{b\gamm\}\|_\infty\\
& \leq |t| \|\tIinv\|_{(\infty, \infty)} \|a\|_\infty + |s| \|\tIinv\|_{(\infty, \infty)} \|\Ld_0\{b \gamm\}\|_\infty \\
& \leq C_1' (|t| + |s|) 
\leq C_1' (1 + (t+s)^2) \leq C_1'' (1 + t^2 + s^2),
\end{align*}
for some constant $C_1'' \geq 2C_1'$.
Second, by Lemma \ref{lemma:cov-2}, $\sum_{0 \leq k \leq 2^{l}} \max_j |\langle \gamma_{M_{1j}, lk}, \psi_{lk} \rangle| \lesssim 2^{l/2}$ and $\sum_{0 \leq k \leq 2^{l}} |\langle \gamma_{b, lk}, \psi_{lk}\rangle| \lesssim 2^{-(1/2 - \mu')l}$, where $\mu' = \mu \wedge 1$, for any $b \in \mathcal{H}(\mu, D)$. 
Therefore, 
\begin{align*}
\sum_{l \leq L_n; k} |\Delta_{2, lk}| 
& \leq 
\sum_{l \leq L_n; k} \left(|t \gamma_{M_1,lk}' \tIinv a| + |s\gamma_{b,lk}| + |s\gamma_{M_{1},lk}' \tIinv \Ld_0\{b \gamm\}|\right)\\
& \leq 
\sum_{l \leq L_n; k} \Big(p^2|t| \max_j |\gamma_{M_{1j}, lk}| \|\tIinv\|_{(\infty, \infty)} \|a\|_\infty
+ |s| |\gamma_{b, lk}| \\
& \qquad + p^2 |s| \max_j |\gamma_{M_{1j}, lk}| \|\tIinv\|_{(\infty, \infty)} \|\Ld_0\{b \gamm\}\|_\infty \Big)\\
& \lesssim 
(|t| + |s|) p^2 L_n 2^{L_n/2} + |s| L_n 2^{-(1/2 - \mu')L_n}.
\end{align*}

Then, for $t, s \leq \log n$ and a fixed $p$, with the assumption that $(|t| + |s|) p^2 L_n 2^{L_n/2} \leq \rn$ and $|s| L_n 2^{(1/2 - \mu')L_n} \leq \rn$ as $\mu' > 0$, 
if choosing a value for $\sigma_l$ that it does not decrease to 0 too fast with $l$,
$\sum_{l \leq L_n;k} C_1|\Delta_{2,lk}|/(\sigma_l \rn) \lesssim (|t| + |s|) \leq C_1'(1 + t^2 + s^2)$, 
for a sufficient large $C_1'$. 
Also, 
$\sum_j C_2 |\Delta_{1j}|/\rn \leq C_2 p \|\Delta_1\|_\infty/\rn \leq C_2 p (1+t^2 + s^2)/\rn$.
Therefore, let $C = \max(C_1', C_2)$, we have
$$
\exp\left(\sum_{j=1}^p \frac{C_1|\Delta_{1j}|}{\rn} + \sum_{l\leq L_n;k} \frac{C_2|\Delta_{2,lk}|}{\rn} \right)
\leq \exp(C(1+t^2 + s^2)). 
$$

To bound (\ref{cv-2}), what remains to show is that $\Pi((\Theta_n - \Delta_1/\rn, \mathcal{H}_n - \Delta_{2,L_n}/\rn)\given X^n) = 1 + o_{P_0}(1)$.
Since a posterior probability is at most $1$, it is sufficient to show that the posterior distribution is bounded from below by $1 + o_{P_0}(1)$. From Lemma \ref{lemma:psib-2}, one can easily check that $\|\Delta_{2,L_n}\|_\infty/\rn \lesssim (|t| + |s|)p^2L_n/\rn = o(\epsilon_n)$. Therefore, one can define $\Theta_n'$ and $\mathcal{H}_n'$ such that $\epsilon_n$ is replaced by  $\epsilon_n/2$, then $\Theta_n' \subset \Theta_n - \Delta_1/\rn$ and $\mathcal{H}_n' \subset \mathcal{H}_n - \Delta_{2,L_n}/\rn$.  
One can choose $\epsilon_n$ to twice of its original value $\epsilon_n/2$ and still has $\Pi((\Theta_n', \mathcal{H}_n') \given X^n) = 1 + o_{p_0}(1)$. 
Therefore, we have showed that (\ref{cv-2}) is bounded by $\exp(C(1+t^2 + s^2))$.

To verify \ref{cond:C1}, one can use the same argument as above except for letting $t$ and $s$ be constants. 
Then, $C_1\|\Delta_1\|_\infty/\rn = o(1)$ and $C_2 \sum_{l\leq L_n;k} |\Delta_{2,lk}|/\rn = o(1)$. 
By following the proofs of the \ref{cond:C2} case, we have $\exp(C_1\sum_{j=1}^p |\Delta_{1j}|/\rn + C_2 \sum_{l\leq L_n;k}|\Delta_{2,lk}|/\rn) = 1 + o(1)$.
Since $\Pi((\Theta_n - \Delta_1/\rn, \mathcal{H}_n - \Delta_2/\rn)\given X^n) = 1 + o_{P_0}(1)$, we obtain \ref{cond:C1}.

{\it 2. The standard normal prior on $Z_{lk}$}

Consider the standard normal prior on $Z_{lk}$. 
Let $\mathbb{H}$ be the reproducing kernel Hilbert space (RKHS) associate to the Gaussian prior and let $\|\cdot\|_\mathbb{H}$ be the associated norm. For independent Gaussian wavelet prior on $r$, from Page 336 of \citet{ghosal17}, Lemma 11.43, $\|f\|_\mathbb{H}^2 = \sum_{l,k} \sigma_l^{-1} f_{lk}$ for any $f \in L^2[0,1]$. 

Let
$
\Delta_{2,n} = -t \gamma_{M_1,n}' \tIinv a + s \gamma_{b,n} + s \gamma_{M_1, n}' \tIinv \Ld_0\{b \gamm\},
$
where $\gamma_{M_1,n} = (\gamma_{n,lk}(M_1))$, $\gamma_{b,n} = (\gamma_{n,lk}(b))$, 
and $(\gamma_{n,lk}(b)) = \langle b/M_0, \psi_{lk} \rangle$.
Define $\varsigma_n = \Delta_{2,n}/\rn$, then we first verify that 
$$
\|\varsigma_n \|_\mathbb{H}^2 = O(t^2 + s^2).
$$
Using the bound for $\gamma_{n,lk}(b)$ for $b = \psi_{LK}$, $0\leq L\leq L_n$, $k < 2^L$, in Lemma \ref{lemma:cov-1}, note that there are $2^{l-L}$ indices $l$ such that $S_{lk} \subset S_{LK}$, we have
\begin{align*}
\|\gamma_{b,n}\|_\mathbb{H}^2
& \lesssim 
\sum_{l\leq L_n} \sigma_l^{-2} 2^{(l-L)/2} + \sum_{l = L+1}^{L_n}\sigma_l^{-2} 2^{l-L} 2^{L/2-3l/2}\\
& \leq L \sigma_l^{-2} + \sum_{l = L+1}^{L_n} \sigma_l^{-2} 2^{-l/2 - L/2} 
\end{align*}
As we choose $\sigma_l = 2^{-l/2}$, the last display is bounded by $L_n 2^{L_n}$. 
Therefore, $s^2\|\gamma_{b,n}\|_\mathbb{H}^2/n = L_n 2^{L_n}/n \leq s^2$. 
Similarly, for each $j$-th coordinate in $M_1$, we have $t^2\|\gamma_{M_{1j},n}\|_\mathbb{H}^2/n = O(t^2)$. By assumptions \ref{asp:i}-\ref{asp:v}, we can bound the squared RKHS-norm of the first term in $\Delta_{2,n}$ divided by $n$ by $O(t^2)$ and the squared RKHS-norm of the third term divided by $n$ by $O(s^2)$. Therefore, we obtain $\|\varsigma_n \|_\mathbb{H}^2 = O(t^2 + s^2)$. 

Next, we change of variables by letting {$\rho = r - \varsigma_n$}. Define the set $B_n$ such that 
$$
B_n = \{r: |\langle \varsigma_n, {\rho} \rangle| \leq M \rn \epsilon_n \|\varsigma_n\|_\mathbb{H} \},
$$
for a suitably large constant $M > 0$. 
Using the fact that $\langle r, \varsigma_n \rangle_\mathbb{H} \sim N(0, \|\varsigma_n\|_\mathbb{H}^2)$, then $\Pi(B_n^c) \leq e^{-Cn\epsilon_n^2}$ for some constant $C > 0$. This implies that $\Pi(B_n^c \given X) = o_{P_{\eta_0}}(1)$. 

We then modify the conditions \ref{cond:C2} (as well as \ref{cond:C1}) by replacing $A_n$ by $\tilde A_n  = A_n \cap B_n$. Then $\Pi({\tilde A_n}^C \given X) = o_{P_{\eta_0}}(1)$. 
The remaining proof proceeds similar as the Laplace prior case, the numerator in (\ref{cv-1}) after $A_n$ is replaced by $\tilde A_n$ (and let $\tilde{\mathcal{H}}_n = \mathcal{H}_n \cap B_n$) can be written as
\begin{align*}
N_n & = \int_{\Theta_n}\int_{\tilde{\mathcal{H}}_n} 
e^{\ell_n(\eta_h) - \ell_n(\eta_0)} 
\pi(r) dr
\prod_{j=1}^p \pi(\theta_j) d \theta_j \\
& = 
\int_{\Theta_n}\int_{\mathcal{H}^\tau_n} 
e^{\ell_n(\eta_h) - \ell_n(\eta_0)} 
\exp\{- \|\varsigma_n\|_\mathbb{H}^2/2 - \langle \varsigma_n, \rho \rangle_\mathbb{H} \} \pi (\rho) d\rho
\prod_{j=1}^p \pi(\theta_j) d \theta_j,
\end{align*}
where $\mathcal{H}^\tau_n = \tau(\tilde{\mathcal{H}}_n)$ with $\tau$ the translation map $\tau: f \to f - \varsigma_n$. 
On the other hand, the denominator in (\ref{cv-1}) is given by 
$$
D_n = \int \int e^{\ell_n(\eta_h) - \ell_n(\eta_0)} \pi(r) dr \pi(\theta) d\theta.
$$
Then to bound the ratio $N_n/D_n$, we need to control $\|\varsigma_n\|_\mathbb{H}^2$ and $\langle \varsigma_n, \rho\rangle$. 
On the set $B_n$, we immediately obtain $\langle \varsigma_n, \rho\rangle \lesssim \rn \epsilon_n \|\varsigma_n\|_\mathbb{H}$ and using a similar derivation as above, we have $\rn\epsilon_n \|\varsigma_n\|_{\mathbb{H}} = O((|t| + |s|) \epsilon_n L_n \sigma_l^{-1}) = O((|t| + |s|) \epsilon_n L_n 2^{L_n/2}) = O(|t| + |s|)$ as $\beta > 1/2$. In fact, one can choose $\sigma_l$ to be a fixed constant, then the condition $\beta > 1/2$ can be dropped. 
Also, we already showed that $\|\varsigma_n\|_\mathbb{H}^2 =O(t^2 + s^2)$. 
Therefore, the expression in the second exponential in $N_n$ is bounded by $\exp(1+ t^2 + s^2)$. 
Last, define $\Theta_n'$ and $\tilde{\mathcal{H}}_n'$ such that $\epsilon_n$ in $\Theta_n$ and $\tilde{\mathcal{H}}_n$ is replaced by $\epsilon_n/2$ and then apply the same argument as in the Laplace prior case, $N_n/D_n \lesssim \exp(1 + t^2 + s^2)$ and thus \ref{cond:C2} is verified. 
The verification of \ref{cond:C1} is similar except that one should use the fact that $t, s$ are fixed constants.

\subsubsection{Verifying the two conditions for {\rm \ref{prior:T}} and {\rm \ref{prior:H}}}
\label{cov-histogram}

Since we have verified \ref{cond:C1} and \ref{cond:C2} for the Haar wavelet prior, we can prove the results for the histograms prior by using the relation between the coefficients of the Haar wavelets and the histogram heights.
Denote the Haar coefficients $r_S$ and the histograms heights $r_H$, recall that $r_S = \Psi r_H$ for the matrix $\Psi$ given in Section \ref{sec:bg}. 

For the independent gamma prior on each $\lambda_k$, i.e., $\lambda_k \sim \text{Gamma}(\alpha_0, \beta_0)$. 
The density function for $r_k^H = \log \ld$ is 
$$
f(r_k^H \given \alpha_0, \beta_0) = \frac{\beta_0^{\alpha_0}}{\Gamma(\alpha_0)} \exp(\alpha_0 r_k^H - \beta_0 e^{r_k}),
$$ 
for $k = 0, \dots, 2^{L_n}$. 

For the dependent gamma prior on each $\lambda_k$, we have 
$\lambda_0 \sim \text{Gamma}(\alpha_0, \beta_0)$ and 
$\lambda_k \given \lambda_{k-1} \sim \text{Gamma}(\alpha, \alpha/\lambda_{k-1})$,
for $k = 1, \dots, 2^{L_n}$.
Therefore, apply the change of variables from $\lambda_k$ to $r_k^H = \log \ld_k$, $k \geq 1$,
we obtain 
$$
f(r_k^H \given \alpha, r_{k-1}^H) 
= \frac{\alpha^\alpha}{\Gamma(\alpha)} 
\exp\left(\alpha(r_k^H - r_{k-1}^H) - \alpha e^{r_k^H - r_{k-1}^H}\right).
$$ 
Then the numerator of the posterior with the dependent prior can be written as
$$
N_n' = \int_{\Theta_n} \int_{\mathcal{H}_n} 
e^{\ell_n(\eta_h)-\ell_n(\eta_0)} 
\prod_{k =1}^{2^{L_n}} f(r_k^H \given r_{k-1}^H) dr_k^H \prod_{j=1}^p \pi(\theta_j) d\theta_j,
$$
also, the denominator can be written as
$$
D_n' = \int\int e^{\ell_n(\eta_h)-\ell_n(\eta_0)} 
\prod_{k =1}^{2^{L_n}} f(r_k^H \given r_{k-1}^H) dr_k^H 
\prod_{j=1}^p \pi(\theta_j) d\theta_j.
$$
Let's denote $\vartheta_j = \theta_j - \Delta_{1j}/\rn$ and $\rho_k^H = r_k^H - \Delta_{2,k}^H/\rn$, 
where $\Delta_2^H = \Psi^{-1} \Delta_2$ and $\Delta_{2,k}^H$ is the $k$-the coordinate of $\Delta_2$. 
The remaining proof is similar to the proof for the Haar wavelet prior. 
Apply the change of variables from $\theta$ to $\vartheta$ and from $r^H$ to $\rho^H$. Due to the invariance of the Lebesgue measure, the denominator becomes
$$
N_n' = \int_{\Theta_n - \Delta_1/\rn} \int_{\mathcal{H}_n - \Delta_{2,L_n}^H/\rn} e^{\ell_n((\vartheta, \rho^H)) - \ell_n(\eta_0)} 
|\det(\Psi)|^{-1} F(\Psi^{-1} \rho^H + z) d\rho^H \cdot \pi(\vartheta + \Delta_1/\rn) d\vartheta,
$$
where 
$\Theta_n$ and $\mathcal{H}_n$ are the same as they defined in Section \ref{cov-laplace}, 
$z = \Psi^{-1} \Delta_{2,L_n}^H/\rn$, and $F(r^H) = f(r_0^H \given \alpha_0, \beta_0) \prod_{k=1}^{2^{L_n} } f(r_k^H \given \alpha, r_{k-1}^H)$.
The denominator can be written as
$$
D_n' = \int \int e^{\ell_n((\vartheta, \rho^H)) - \ell_n(\eta_0)} 
|\det(\Psi)|^{-1} F(\Psi^{-1} \rho^H) d\rho^H \cdot d\Pi(\vartheta). 
$$
We need to bound the ratio $N_n'/D_n'$, which requires to control the ratio 
$F(\Psi^{-1}\rho^H + z)/F(\Psi^{-1}\rho^H)$.
By plugging-in the expression for $F(\cdot)$, this ratio can be written as a product of $2^{L_n+1}$ individual terms, where the first term is
\begin{align*}
& \log( f((\Psi^{-1}\rho^H)_0 + z_1 \given \alpha_0, \beta_0)) - 
\log( f((\Psi^{-1}\rho^H)_0 \given \alpha_0, \beta_0))  =
\alpha_0 z_1 + \beta_0 \left(1 - e^{z_1}\right)e^{(\Psi^{-1}\rho^H)_0},
\end{align*}
and each of the remaining terms is
\begin{align*}
&\log (f( (\Psi^{-1}\rho^H)_k + z_k \given \alpha, (\Psi^{-1}\rho^H)_{k-1} + z_{k-1})
- \log (f((\Psi^{-1}\rho^H)_k \given \alpha, (\Psi^{-1}\rho^H)_{k-1}))\\
& \quad 
= \alpha (z_k  - z_{k-1}) + \alpha (1 - e^{z_k - z_{k-1}}) e^{(\Psi^{-1}\rho^H)_k - (\Psi^{-1}\rho^H)_{k-1}}
\end{align*}
for each $k = 1, \dots, 2^{L_n}$.
Denote $h_k = (\Psi^{-1}\rho^H)_k$, the $k$-th element in the vector $\Psi^{-1}\rho$,
then, 
\begin{align*}
& \log F(\Psi^{-1}\rho^H+z) - \log F(\Psi^{-1}\rho^H) \\
& \quad = 
\alpha_0 z_0 + \beta_0(1-e^{z_0})e^{h_0}
+ \alpha \sum_{k=1}^{2^{L_n}} (z_k - z_{k-1}) 
+ \alpha \sum_{k=1}^{2^{L_n}} (1 - e^{z_k - z_{k-1}}) e^{h_k - h_{k-1}}.
\end{align*}
By Lemmas \ref{lemma:cov-3} and \ref{lemma:cov-4}, we thus obtain
\begin{align*}
\sum_{k=1}^{2^{L_n}} |z_k - z_{k-1}| 
& \lesssim 
\frac{1}{\rn}\sum_{k=1}^{2^{L_n}}\left(p^2 (|t| + |s|)\max_j |\tilde H_k - \tilde H_{k-1}| + |s| |H_k - H_{k-1}|\right)\\
& \lesssim p^2 (|t|+|s|) 2^{L_n/2}/\rn + |s| 2^{L_n/2}/\rn,
\end{align*}
which the last line is bounded by $C(1+t^2+s^2)$ for some constant $C$, as $2^{L_n/2}/\rn = o(1)$.
In addition, $|h_k - h_{k-1}| = O(1)$ for each $k = 1, \dots, 2^{L_n}$ as one can invoke the supremum-norm consistency of $r$, which implies that all the histogram heights are bounded by a universal constant. 
Therefore, the penultimate display is bounded by $C_1 (1 + t^2 + s^2)$ for some constant $C_1 > C$.
This implies that
$$
\frac{F(\Psi^{-1}\rho^H + z)}{F(\Psi^{-1}\rho^H)}
\lesssim e^{C_1(1 + t^2 + s^2)}.  
$$ 
From Section \ref{cov-laplace}, we have 
$\pi(\vartheta + \Delta_1/\rn) / \pi(\vartheta) \leq \exp(C_2 (1 + t^2 + s^2))$.
By using a similar argument as we prove the case for the Haar wavelet prior, we thus \ref{cond:C2}.
We also verified \ref{cond:C1} by using that $t$ and $s$ are fixed constants.

For the independent gamma prior, the proof is similar to and simpler than the dependent gamma prior, 
By plugging-in its density function, one can check that 
$$
\frac{F(\Psi^{-1}\rho + z)}{F(\Psi^{-1}\rho) }\leq \exp\left(\alpha_0 \sum_{k=1}^{2^{L_n+1}-1}{z_k} \right)
\lesssim (|t| + |s|) 2^{L_n}/\rn.
$$ 
Due to $2^{L_n}/\rn = o(1)$ for the choice of $L_n$ in (\ref{choicel}) and $\beta > 1/2$, the previous display is bounded by $1 + t^2 + s^2$ if $t, s \leq \log n$ and by $o(1)$ if $t$ and $s$ are fixed constants. The remaining proof is essentially the same as the proof for the dependent gamma prior. We thus verified \ref{cond:C1} and \ref{cond:C2} for the use of the independent gamma prior. 
%%%%%%%%%%%%%%%%%%%%%%%%%%%%%%%%%%%%%%%%%%%%%%
%%%%%%%%%%%%%%%%%%%%%%%%%%%%%%%%%%%%%%%%%%%%%%
\section{Auxiliary lemmata} 
\label{auxiliary}

\subsection{Approximation lemmata for wavelets and histograms}
\label{wavelets-lemma}

Let $\tz \in \{z_1, \dots, z_p\}$ and $\max |z_j| \leq c_1$ for some positive constant $c_1$, denote $z = (z_1, \dots, z_p)'$,
we define
\begin{align}
\label{tilde-m-u}
\tilde M(u) 
= \mathbb{E}(\tz e^{\theta_0'z} \mathbbm{1}_{u \leq T}) 
= \int \tz \bar G_z(u) e^{\theta_0'z - \Ld_0(u) e^{\theta_0'z}} f(z) dz.
\end{align}
The lemmas listed in below give bounds for Haar wavelet basis and their projections. 

\begin{lemma}
\label{lemma:bound-W}
Let $\psi_{lk}$ be a Haar wavelet bases and $L_n$ be a cut-off and let $l < L_n$ and $0 \leq k < 2^L$. For $b \in L^\infty([0,1])$, recall that $\gamb = b/M_0$ and $\gambln = P_{L_n}(b/M_0)$; similarly, define 
$\gamma_{\tilde M} = \tilde M/M_0$ and $\gamma_{\tilde M, L_n} = P_{L_n}(\tilde M/M_0)$, $\tilde M$ given in (\ref{tilde-m-u}). Recall the LAN-norm $\|\cdot ,\cdot \|_L$ in Section \ref{sec:bg}, then
for any fixed and bounded function $b$, as $n\to \infty$ and $L_n \to \infty$, 
$$
\|0, \gambln\|_L \to \|0,\gamb\|_L, \quad 
\|0, \gamma_{\tilde M, L_n}\|_L \to \|0, \gamma_{\tilde M} \|_L.
$$
\end{lemma}
\begin{proof}
The proof of the first statement can be found in Lemma 12 of \citetalias{cast20survival}. For the second statement, since $\tilde M / M_0 \leq |\tz| \leq c_1$ by assumption, $\|\gamma_{\tilde M, L_n} - \gamma_{\tilde M}\|_{L^2} \to 0$ by definition. 
\end{proof}
%%%%%%%%%%%%%%%%%%%%%%%%%%%%%%%%%%%%%%%%%%%%%%%%

\begin{lemma}[Lemma 12 of \citetalias{cast20survival}]
\label{lemma:psib-1}
Under the same condition as in Lemma \ref{lemma:bound-W}, the following results hold:
$$
\|\gambln\|_\Linfty \leq C L_n \|b\|_\Linfty, \qquad 
\|\gambln\|_\Ltwo \leq C \|b\|_\Ltwo.
$$
\end{lemma}

%%%%%%%%%%%%%%%%%%%%%%%%%%%%%%%%%%%%%%%%%%%%%%%%
\begin{lemma}
\label{lemma:M0M1}
Under the same assumptions in Section \ref{sec:assumption}, $M_0(\cdot)$, $M_0^{-1}(\cdot)$, and $\tilde M(\cdot)$ are all Lipschitz functions on $[0, 1]$.
\end{lemma}
\begin{proof}
By the definition of $M_0(u)$ in Section \ref{sec:bg}, 
$M_0(u) = \int \bar G_z(u) e^{\theta_0'z} e^{-\Lambda_0(u) e^{\theta_0'z}} f(z) dz$ for $u \in [0, 1]$.
Since $\bar G_z(u) = 1 - \int_0^u g_z(v) dv$ and $g_z$ is bounded by assumption, $\bar G_z(u)$ is Lipschitz on $[0, 1]$. 
Due to $\ld_0$ is continuous, $e^{-\Ld_0(u)} = e^{-\int_0^u \lambda_0(v)}$ is $\mathcal{C}^1$. Thus, $e^{\theta_0'z} e^{-\Ld_0(u)} e^{\theta'z} \bar G_z(u)$ is a product of Lipschitz maps. 
If $f(z)$ is continuous and also bounded by assumption, $e^{\theta_0'z} e^{-\Ld_0(u)} e^{\theta'z} \bar G_z(u)f(z)$ is Lipschitz, so does $M_0(u)$. 
If $f(z)$ is discrete, the integral is the summation of multiple Lipschitz functions, each of which contains a different value of $z$, which $e^{\theta_0'z} e^{-\Ld_0(u)} e^{\theta'z} \bar G_z(u)f(z)$ is again Lipschitz. Therefore, $M_0(\cdot)$ is Lipschitz on $[0, 1]$. 

Next, by assumptions (3), (4), and (8) in Section \ref{sec:assumption}, $M_0(u)^{-1}$ is bounded. Hence,
we have $|M_0(u_1)^{-1} - M_0(u_2)^{-1}| \lesssim |M_0(u_1) - M_0(u_2)|$ and $M_0(\cdot)^{-1}$ is Lipschitz. 

Last, note that 
$|\tilde M(u_1) - \tilde M(u_2)| = |(\tilde M/M_0)(u_1) \cdot M_0(u_1) - (\tilde M/M_0)(u_2) \cdot M_0(u_2)| 
\leq |\tz| |M_0(u_1)  - M_0(u_2)| \leq c_1 |M_0(u_1)  - M_0(u_2)|$, hence $\tilde M(u_1)$ is a Lipschitz function.
 \end{proof}
 
 %%%%%%%%%%%%%%%%%%%%%%%%%%%%%%%%%%%%%%%%%%%%%%%%

\begin{lemma}
\label{lemma:psib-2}
Under the same condition as in Lemma \ref{lemma:psib-1}, 
\begin{enumerate}
\item If $b = \psi_{LK}$, uniformly over $L \leq L_n$ and $K$, for some $C > 0$,
\begin{align*}
	& \|\gambln\|_\Ltwo \leq C, \quad \|\gambln\|_\Linfty \leq C 2^{L/2}L_n, \quad 
	\|\gamb - \gambln\|_\Linfty \lesssim 2^{L/2} 2^{-L_n}.
\end{align*}
\item If $b \in \mathcal{H}(\mu, L)$ for some $\mu, L > 0$, with $\mu' = \mu \wedge 1$, 
$$
\|\gamb - \gambln\|_\Linfty \lesssim 2^{-\mu' L_n}.
$$
\item If $b = \tilde M$, with assumption {\rm \ref{asp:i}}, for some $C' > 0$ and $c_1$ in \ref{asp:i},
\begin{align*}
	& \|\gamma_{\tilde M, L_n}\|_\Ltwo \leq C', 
	\quad 
	\|\gamma_{\tilde M, L_n}\|_\Linfty \leq c_1L_n, 
	\quad
	\|\gamma_{\tilde M} - \gamma_{\tilde M, L_n}\|_\Linfty \lesssim 2^{- L_n}.
\end{align*} 
\item For any fixed bounded function $b$, suppose {\rm \ref{cond:B}} holds, then, as $n\to \infty$ and $L_n \to \infty$,
\begin{align*}
& W_n(0, \gamb - \gambln) = o_{P_{\eta_0}}(1),
\quad
W_n(0, \gamma_{\tilde M} - \gamma_{\tilde M, L_n}) = O_{P_{\eta_0}}(2^{-L_n}).
\end{align*}
\end{enumerate}
\end{lemma}

\begin{proof}
The first and second points are from Lemma 13 of \citetalias{cast20survival}. 
To prove the third point, first, by Lemma \ref{lemma:psib-1}, let $b = \tilde M$, we have
$\|\gamma_{\tilde M, L_n}\|_\Ltwo \leq C\|\tilde M\|_\Ltwo \leq C\|\tilde M/M_0\|_\Ltwo \|M_0\|_\Linfty \leq C'$ by \ref{asp:i}, where $C' > Cc_1$.
Next, by Lemma \ref{lemma:psib-1}, we also have
$\|\gamma_{\tilde M, L_n}\|_\Linfty \leq L_n \|\gamma_{\tilde M}\|_\Linfty \leq L_n \|z\|_\infty \leq c_1 L_n$.
Last, 
let $h = \tilde M/M_0$ and denote $\bar h$ as the mean of $h$ on the support of the wavelet $\gamma_{lk}$, then 
$$
\langle \tilde M/M_0, \psi_{lk} \rangle
= \langle h - \bar h, \psi_{lk} \rangle  + \bar h \langle 1, \psi_{lk} \rangle.
$$
The second term in the last display is 0. 
Since $h$ is a Lipschitz function by Lemma \ref{lemma:M0M1}, 
for all $x$ in the support of $S_{lk}$ of $\psi_{lk}$, there exist a $c$ in $S_{lk}$ such that 
$$
|h(x) - \bar h| = |h(x) - h(c)| \lesssim |x - c| \leq 2^{-l}.
$$
Therefore, 
\begin{align*}
\|\gamma_{\tilde M} - \gamma_{\tilde M, L_n}\|_\Linfty 
& \leq \sum_{l > L_n} 2^{l/2} \max \left|
\langle \tilde M/M_0, \psi_{lk} \rangle
\right|
\lesssim \sum_{l > L_n} 2^{l/2} 2^{-l} \|\psi_{lk}\|_\Lone \\
& \lesssim \sum_{l > L_n} 2^{l/2} 2^{-l} 2^{-l/2} \lesssim 2^{-L_n}.
\end{align*}
In fact, one could also obtain the result by using the second point. Since here $b = \tilde M/M_0$ is Lipschitz, thus is $\mathcal{1, L}$. Therefore, one can plug-in $\mu' = 1$ and obtain the upper bound $2^{-L_n}$.

To prove the fourth point,
since $W_n(0, \gamb - \gambln)$ is a collection of real variables that are centered at 0 with variance  equals to 
$\|0, \gamb - \gambln\|_L^2$ under $P_{\eta_0}$, by the definition of the LAN-norm, $\|0, \gamb - \gambln\|_L^2 = \Ld_0\{ (\gamb - \gambln)^2 M_0\}$ and since $\|M_0\lambda_0\|_\Linfty$ is bounded, we obtain
$\|0, \gamb - \gambln\|_L^2 \lesssim \|\gamb - \gambln\|_\Ltwo^2 = o(1)$ as we assume
$\rn \epsilon_n \|\gamb - \gambln\|_\Linfty = o(1)$ in \ref{cond:B}. 
Similarly, the variance of $W_n(0, \gamma_{\tilde M} - \gamma_{\tilde M,L_n})$ is $\|0, \gamma_{\tilde M} - \gamma_{\tilde M, L_n}\|_L^2$, thus it is bounded by $O_{\eta_0} (2^{-L_n})$. 
\end{proof}

%%%%%%%%%%%%%%%%%%%%%%%%%%%%%%%%%%%%%%%%%%%%%%%%

\begin{lemma}[Lemma 8 of \citetalias{cast20survival}]
\label{lemma:cov-1}
With $(\psi_{lk})$ a Haar wavelet bases, let $b = \psi_{LK}$ and set $\gambln = P_{L_n}(b/M_0)$ and $\gamma_{n,lk}: = \gamma_{n,lk}(b) = \langle b/M_0, \psi_{lk} \rangle$. Denote the support of $\psi_{lk}$ as $S_{lk}$. Suppose $L \leq L_n$, where $L_n$ is the cut-off such that $2^{L_n} = (n/\log n)^{1/(2\beta+1)}$, then
\begin{align*}
& |\gamma_{n,lk}| \lesssim 2^{-(L-l)/2}, \quad \text{if} \ l \leq L,\\
& |\gamma_{n,lk}| \lesssim 2^{L/2 - 3l/2}, \quad \text{if} \ l \geq L, S_{lk} \cap S_{LK} 
\neq \varnothing, (l,k) \neq (L,K),\\
& |\gamma_{n,lk}| = 0, \qquad \qquad \ \text{if} \ l \geq L, S_{lk} \cap S_{LK} = \varnothing.
\end{align*}
\end{lemma}
%%%%%%%%%%%%%%%%%%%%%%%%%%%%%%%%%%%%%%%%%%%%%%%%

\begin{lemma}
\label{lemma:cov-2}
For $b \in L^{\infty}[0,1]$, recall that $\gambln = P_{L_n}(b/M_0)$ and $\gamma_{n,lk}(b) = \langle b/M_0, \psi_{lk}\rangle$, for $l \leq L_n$ an integer, we have 
\begin{enumerate}
\item if $b = \psi_{LK}$, for any $L \leq L_n$ and $0 \leq K \leq 2^{L} - 1$, 
$$
\sum_{0\leq k\leq 2^l} |\gamma_{n,lk}(\psi_{LK})| \lesssim 2^{-(L - l)/2},
$$
\item if $b \in \mathcal{H}(\mu, D)$ for some positive constants $\mu$ and $D$, let 
$\mu' = 1 \wedge \mu$,
$$
\sum_{0 \leq k \leq 2^l} |\gamma_{n,lk}(b)| \lesssim 2^{-(2\mu' - 1)l/2},
$$
\item if $b = \tilde M(u)$ in (\ref{tilde-m-u}), 
$\displaystyle
\sum_{0 \leq k \leq 2^l} |\gamma_{n,lk}(\tilde M)| \lesssim 2^{-l/2}.
$
\end{enumerate}
\end{lemma}

\begin{proof}
The first two points are directly from Lemma 9 of \citetalias{cast20survival}. For the third point, note that $\tilde M(u) \in \mathcal{H}(1, D')$ for some constant $D' > 0$, we choose $\mu' = 1$ and obtain the upper bound using the second point.
\end{proof}

%%%%%%%%%%%%%%%%%%%%%%%%%%%%%%%%%%%%%%%%%%%%%%%%

\begin{lemma}[Lemma 10 \& 11 in \citetalias{cast20survival}]
\label{lemma:cov-3}
Let $\gambln = P_{L_n} (b/M_0)$ and $b = \psi_{LK}$ for some $L, K$, and set $H = \Psi^{-1} \gambln$, with $\Psi$ the matrix described in Section \ref{sec:bg}, then for $L \leq L_n$ and $1 \leq K \leq 2^{L+1}$, 
\begin{align*}
& |H_j| \leq C2^{L/2}, \quad \text{if} \ I_j^{L_n + 1} \cap S_{LK} \neq \varnothing,\\
& H_j = 0, \quad \quad \quad \ \ \text{if} \ I_{j}^{L_n + 1} \cap S_{LK} = \varnothing,\\
& \qquad \sum_{j = 1}^{2^{L_n+ 1}} | H_j - H_{j-1} | \leq C 2^{L/2}.
\end{align*}
\end{lemma}

%%%%%%%%%%%%%%%%%%%%%%%%%%%%%%%%%%%%%%%%%%%%%%%%

\begin{lemma}
\label{lemma:cov-4}
Let $\gamma_{\tilde M} = P_{L_n} (\tilde M/M_0)$ for $\tilde M$ given in (\ref{tilde-m-u}). Set $\tilde H = \Psi^{-1} \gamma_{\tilde M, L_n}$, with $\Psi$ the matrix given in Section \ref{sec:bg}, then for $l \leq L_n$ and $1 \leq j \leq 2^{L_n+1}$, then for some constant $D$,
$$
|\tilde H_j| \leq D,
\qquad 
\sum_{j = 1}^{2^{L_n+ 1}} | \tilde H_j - \tilde H_{j-1} | \leq 4D2^{L_n}.
$$
\end{lemma}
\begin{proof}
Recall that $\Psi$ is a $2^{L_n+1} \times 2^{L_n+1}$ matrix such that 
$$
\begin{pmatrix}
\Psi_{-1,1} & \Psi_{-1,2} & \cdots & \Psi_{-1, 2^{L_n+1}} \\
\Psi_{00,1} & \Psi_{00,2} & \cdots & \Psi_{00,2^{L_n+1}} \\
\Psi_{10, 1} & \Psi_{10,2} & \cdots & \Psi_{10,2^{L_n+1}} \\
\vdots & \vdots & \ddots & \vdots \\
\Psi_{L_n (2^{L_n}-1), 1} & \Psi_{L_n (2^{L_n}-1), 2} & \cdots & \Psi_{L (2^{L_n}-1), 2^{L_n+1}}
\end{pmatrix},
$$
where $\Psi_{-1,j} = 2^{-(L_n+1)}$ and 
$\Psi_{lk,j} = 2^{-(L_n+1) + l/2} \left[
\mathbbm{1}_{I_{j-1}^{L_n+1} \subset I_{2k}^{l+1}}
- \mathbbm{1}_{I_{j-1}^{L_n+1} \subset I_{2k+1}^{l+1}}
\right]$ for $1\leq l \leq L_n, 0 \leq k \leq 2^{l}-1$, and $j = 1, \dots, 2^{L_n+1}$.
Observe that $2^{(L_n+1)/2}\Psi$ is an orthogonal matrix, we denote $\tilde \Psi = 2^{(L_n+1)/2}\Psi$, then $\tilde \Psi^{-1} = \tilde \Psi'$. Thus, $\Psi^{-1} = 2^{(L_n+1)}\tilde \Psi'$.
Therefore, 
\begin{align*}
\tilde H_j &= 
\left( \Psi^{-1} \gamma_{\tilde M, L_n} \right)_j
= 2^{L_n+1} \sum_{l\leq L_n;k} \Psi_{lk,j} \left\langle {\tilde M}/{M_0}, \psi_{lk} \right\rangle
=  2^{L_n+1}\left \langle{\tilde M}/{M_0}, \sum_{l\leq L_n;k} \Psi_{lk,j} \psi_{lk} \right\rangle\\
& = 2^{L_n+1} \left\langle {\tilde M}/{M_0}, \psi_{r_{H,j}} \right\rangle
\leq 2^{L_n+1} \|\tilde M/M_0\|_\infty  \|r_{H,j}\|_\Lone
\leq c_12^{L_n+1} 2^{-L_n-1}
= c_1. 
\end{align*}
For the second inequality, let $D \geq c_1$, we have 
$
\sum_{j = 1}^{2^{L_n+ 1}} | \tilde H_j - \tilde H_{j-1} | 
\leq 2\sum_{j=1}^{2^{L_n + 1}} |\tilde H_j| \leq 4 D 2^{L_n}.
$
\end{proof}

%%%%%%%%%%%%%%%%%%%%%%%%%%%%%%%%%%%%%%%%%%%%%%%%

\subsection{Lemmata for bounding empirical processes}
\label{EPlemmas}
Let $N_{[]}$ be the usual bracketing number and $J_{[]}(\delta, \mathcal{F}, \|\cdot\|)$ be the bracketing integral of a class of function $\mathcal{F}$ equipped with a norm $\|\cdot\|$, from \citet{vdvEP}, 
\begin{align}
\label{bracketing-number}
J_{[]}(\delta, \mathcal{F}, \|\cdot\|) = \int_0^\delta \sqrt{1 + \log N_{[]}(\epsilon, \mathcal{F}, \|\cdot\|)} d\epsilon.
\end{align}

% Some empirical process result
\begin{lemma}[Example 19.7 of \citet{vdvAS}]
\label{lemma-ep-1}
Suppose $\mathcal{G} = \{g_{\theta}: \theta \in \Theta\}$ with $\Theta = \{\theta \in \mathbb{R}^p: \|\theta\|_\infty \leq M_1\}$ is a class of functions satisfies that 
$|g_{\theta_1}(\cdot) - g_{\theta_2}(\cdot)| \leq L \|\theta_1 - \theta_2\|$, 
If $L \leq \infty$, then there exists a constant $K_1 > 0$ such that for every $\epsilon$ such that $0 < \epsilon < M_1$, 
$$
N_{[]}(L \epsilon , \mathcal{G}, L^2(P)) \leq K_1 \left(\frac{M_1}{\epsilon}\right)^p.
$$
\end{lemma}

\begin{lemma}[Lemma 18 of \citetalias{cast20survival}] 
\label{lemma-ep-2}
Let $\mathcal{F}(M_2)$ be the set of all functions $f: \mathbb{R} \to \mathbb{R}$ with $f(0) = 0$ which have bounded total variation $M_2$ such that $0 < \nu < M_2$. Then, there exists a constant $K_2 > 0$ such that for every distribution $P$, 
$$
\log N_{[]}(\nu, \mathcal{F}(M_2), L^2(P)) \leq \frac{K_2M_2}{\nu}.
$$
\end{lemma}

\begin{lemma}[Lemma 3.4.2 in \citet{vdvEP}]
\label{lemma-ep-4}
Let $\mathcal{F}$ be a class of measurable functions such that for any $f \in \mathcal{F}$,
$\int f^2 dP \leq \delta$ and $\|f \|_\infty \leq M_3$, then for $j(\delta) = J_{[]}(\delta, \mathcal{F}, L^2(P_{\eta_0}))$, 
$$
\mathbb{E}_{P_{\eta_0}}^\star \|\mathbb{G}_n\|_\mathcal{F} \lesssim j(\delta) \left(1 + \frac{j(\delta)M_3}{\delta^2 \rn}\right).
$$
\end{lemma}

\begin{lemma}
\label{lemma-ep-3}
Let $\mathcal{F}$ be a class of measurable functions such that $f \in \mathcal{F}$ and $f = gh$, where 
$g \in \mathcal{G}$ satisfies the conditions in Lemma \ref{lemma-ep-1} and $h \in \mathcal{H}$ satisfies the conditions in Lemma \ref{lemma-ep-2}, 
if
for some positive $\mu_1, \mu_2$, and $D$,
$$
|g_{\theta_1}(\cdot) - g_{\theta_2}(\cdot)| \leq D \|\theta_1 - \theta_2\|, 
\quad 
\|g_\theta\|_\Linfty \leq \mu_1,
$$
and
$$
h(0) = 0, \quad \|h\|_{BV} := \int_0^1 |h'(u)|du \leq \mu_2,
$$
then, there exist some constants $M_1, K_1, K_2 >0$, for a constant $S$ such that $K_1p (\log D + 1) + K_2 \leq S^2 \leq \rn$, 
$$
\mathbb{E}^\star_{\eta_0} \|\mathbb{G}_n\|_\mathcal{F} \leq 4 S \mu_1\mu_2,
$$
where $\|\mathbb{G}_n\|_\mathcal{F} = \sup_{f \in \mathcal{F}}|\mathbb{G}_n f|$ and $\mathbb{E}^\star_{\eta_0}$ is the corresponding outer expectation under $P_{\eta_0}$.
\end{lemma}

\begin{proof}
Note that for any $f \in \mathcal{F}$, we have $\int f^2 dP_{\eta_0} \leq \|g\|_\Linfty^2 \|h\|_\Linfty^2$.
Since $\|g\|_\infty \leq \mu_1$ and $\|h\|_\Linfty \leq |h(0)| + \|h\|_{BV} \leq \mu_2$ by assumption,
$\int f^2 dP_{\eta_0} \leq \mu_1^2 \mu_2^2$. 

We remark that it is sufficient to prove the lemma when $\mu_1\mu_2 = 1$, as otherwise one can consider the set $\mathcal{F}' = \{f' = f/(\mu_1\mu_2), f\in \mathcal{F}\}$ and obtain $\|\mathbb{G}_n\|_{\mathcal{F}} = \mu_1\mu_2\|\mathbb{G}_n\|_{\mathcal{F}'}$. 
Let $\mu_1 = \mu_2 = 1$, the bracketing entropy number is bounded by
$$
\log N_{[]} (\varepsilon, \mathcal{F}, L^2(P_{\eta_0})) \leq 
\log N_{[]} (\varepsilon, \mathcal{G}, L^2(P_{\eta_0})) 
+ \log N_{[]} (\varepsilon, \mathcal{H}, L^2(P_{\eta_0})),
$$
From Lemma \ref{lemma-ep-1},
$\log N_{[]}(\varepsilon, \mathcal{G}, L^2(P_{\eta_0})) \leq K_1 p \log (D/\varepsilon)$,
as $\epsilon < 1$.
From Lemma \ref{lemma-ep-2}, we have
$\log N_{[]}(\varepsilon, \mathcal{H}(1), L^2(P_{\eta_0}))$ $\leq K_2 /\varepsilon$.
Thus, we obtain 
$\log N_{[]} (\varepsilon, \mathcal{F}, L^2(P_{\eta_0})) \leq K_2/\varepsilon + K_1p\log (D/\varepsilon)\leq S^2/\varepsilon$ as $\log (1/\varepsilon) \leq 1/\varepsilon$ and $S^2 \geq K_1p (\log D  + 1) + K_2$.
Using (\ref{bracketing-number}), the bracketing number is bounded by 
$$
J_{[]}(1, \mathcal{F}, L^2(P_{\eta_0})) 
\leq \int_0^1 S \sqrt{1/\varepsilon} d\varepsilon
\leq 2 S.
$$
Applying Lemma \ref{lemma-ep-4} with $\delta = \delta_n = 1$ and $M_3 = 1$ and using the assumption $S \leq \rn$, 
$$
\mathbb{E}^\star_{\eta_0} \|\mathbb{G}_n\|_\mathcal{F}
\leq J_{[]}(1, \mathcal{F}, L^2(P_{\eta_0})) 
\left( 1 + \frac{J_{[]}(1, \mathcal{F}, L^2(P_{\eta_0})) }{\rn} \right)
\leq 4 S.
$$
This concludes the proof of the case when $\mu_1 = \mu_2 = 1$. The proof for any $\mu_1 >0$ and $\mu_2 > 0$ is similar, as we argued above. Thus, the proof is completed. 
\end{proof}

%%%%%%%%%%%%%%%%%%%%%%%%%%%%%%%%%%%%%%%%%%%%%%%%%%%

\subsection{A proposition for establishing BvM in the product space $\mathbb{R}^d \times \mathcal{M}_0$}
\label{key-prop}

The proposition below extends Proposition 6 of \citet{cast14a} to the product space $\mathbb{R}^p \times \mathcal{M}_0(w)$. 
Define the operator,
$$
\tau_{V_L}^\times: (\theta, \ld) \to (\theta, \pi_{V_L} \ld),
$$
where $\pi_{V_L}\ld$ be the projection of $\ld$ onto $V_L$, $V_L$ is the subspace of $\mathcal{M}_0$ consisting of Haar wavelets functions up to level $L$.

\begin{prop}
\label{prop-tightness}
Let $(\theta, \ld) \sim \Pi(\cdot \given X)$
and $T_n^\theta$ and $T_n^\ld$ be the centerings of $\theta$ and $\ld$ respectively.
Define $\tilde \Pi(\cdot \given X)$ as the distribution of 
$\rn \left(\theta - T_n^\theta, \ \ld - T_n^\ld \right)'$ conditional on $X$.
Denote $\mathcal{N}$ the Gaussian probability measure on $\mathbb{R}^p \times \mathcal{M}_0$ with 
$N(0,1)^{\otimes p}$ the law on the first $p$ coordinates and, independently, the $\mathcal{M}_0$-part is the $P$-white noise $\mathbb{Z}_p$ from a bounded density $P$. 
Let us equip $\mathbb{R}^p \times \mathcal{M}_0(w)$ with the norm $\|\cdot\|_\times$ given by 
$\|(\theta, \ld)\|_\times = \|\theta\|+\|\ld\|_{\mathcal{M}_0}$, where $\|\cdot\|$ is the standard euclidean norm. 
Suppose, as $n \to \infty$,
\begin{enumerate}
\item the finite-dimensional distribution converges,
\begin{align}
\label{eqn:cond1}
\mathcal{B}_{\mathbb{R}^p \times \mathcal{M}_0(w)} 
\left(
\tilde \Pi_n \circ {\tau_{V_L}^\times}^{-1}, 
\mathcal{N} \circ {\tau_{V_L}^\times}^{-1}
\right) \to^{P_{\eta_0}} 0,
\end{align}
\item for some admissible sequence $\bar w_l = (\bar w_l) \to \infty$ and
$\bar w_l /\sqrt{l} \geq 1$, 
\begin{align}
\label{eqn:cond2}
& \mathbb{E}
\left[\|\ld - T_n^\ld \|_{\mathcal{M}_0(\bar w)} \given X \right] 
 = O_{P_{\eta_0}} (1/\rn)
\end{align}
\end{enumerate}
Then, for any sequence ($w_l$) such that $w_l/\bar w_l \to \infty$ as $n \to \infty$, 
\begin{align}
\label{eqn:cond2-result}
\displaystyle
\mathcal{B}_{\mathbb{R}^p \times \mathcal{M}_0(w)} \left( \tilde \Pi_n, \ \mathcal{N} \right) 
\to^{P_{\eta_0}} 0.
\end{align}
\end{prop}

\begin{proof}
For simplicity, we denote $\beta = \beta_{\mathbb{R}^p \times \mathcal{M}_0(w)}$.
By applying the triangle inequality,
\begin{align*}
\mathcal{B}(\tilde \Pi_n, \mathcal{N})
\leq 
& \underbrace{\mathcal{B}( \tilde \Pi_n, \tilde \Pi_n\circ {\tau_{V_L}^\times}^{-1})}_{(I)}
+ \underbrace{\mathcal{B}(\tilde \Pi_n\circ {\tau_{V_L}^\times}^{-1}, 
\mathcal{N} \circ {\tau_{V_L}^\times}^{-1})}_{(II)} 
+ \underbrace{\mathcal{B}(\mathcal{N}, \mathcal{N} \circ {\tau_{V_L}^\times}^{-1})}_{(III)}. 
\end{align*}
From (\ref{eqn:cond1}), we immediately obtain $(II) \to^{P_{\eta_0}} 0$.
To prove $(I)$ converges. By the definition of $\beta$  and denote $\mathcal{S} = \mathbb{R}^p \times \mathcal{M}_0(w)$, for any bounded function $F$ on $S$ such that $\|F\|_{BL} \leq 1$, 
\begin{align*}
(I) & = \sup_{F:\|F\|_{BL} \leq 1} 
\left|
	\int_\mathcal{S} F d\tilde \Pi_n - \int_\mathcal{S} F d\tilde \Pi_n \circ {\tau_{V_L}^\times}^{-1}
\right|
\leq \mathbb{E}
\left|
F(\tilde \theta_n, \tilde \ld_n) - F(\tilde \theta_n, \pi_{V_L} \tilde \ld_n)
\right| \\
& \leq  \mathbb{E} \left[
\| \tilde \theta_n - \tilde \theta_n, \tilde \ld_n - \pi_{V_L} \tilde \ld_n\|_\times  \given X
\right]
= \mathbb{E} \left[
\| \tilde \ld_n - \pi_{V_L} \tilde \ld_n \|_{\mathcal{M}_0(\tilde w)} \given X
\right],
\end{align*}
for $(\tilde \theta_n, \tilde \ld_n) \sim \tilde \Pi_n$ with $\tilde \theta_n := \rn (\theta - T_n^\theta)$ and 
$\tilde \ld_n := \rn (\ld - T_n^\ld)$.
By following the same argument of the proof of Proposition 6 on P. 1960 of \citet{cast14a}, the last display is bounded by $\sup_{l > L} (\bar w_l/w_l) \times O_{P_{\eta_0}}(1)$, which can be as small as desired by choosing a large but fixed $L$.

To show $(III)$ converges to 0 in probability, one can apply a similar argument as above but replacing $\tilde \Pi_n$ with $\mathcal{N}$,
then use the same argument as at the end of the proof of Theorem 1 on P. 1959 of \citet{cast14a}. 
\end{proof}

\subsection{Bounding the remainder using $\|\cdot\|_\Linfty$-consistency of $\lambda$}
\label{bound-use-linfty-norm}

The results in Lemmas \ref{lemma-R1} and \ref{thm1-lemma2} use $\|\cdot\|_\Lone$-consistency of $\lambda$. However, the upper bounds in those lemmas can be large if $t, s, b$ are divergent sequences as $n$ increases. 
This is especially problematic for obtaining the nonparametric BvM result in Section \ref{pf-nonparametric-BvM}.
The following two lemmas derive upper bounds using $\|\cdot\|_\Linfty$-consistency (instead of $\|\cdot\|_\Lone$ for $\lambda$.
The upper bounds can be smaller than those in Lemmas \ref{lemma-R1} and \ref{thm1-lemma2} when $t, s, b$ are divergent sequences. 

Let's consider two sequences of positive real numbers $(\epsilon_n)$ and $(\zeta_n)$, typically $v_n \geq \epsilon_n$, and define the set:
\begin{align}
\mathcal{L}_n' & = \{\eta = (\theta, \lambda): \theta \in \mathbb{R}^p, \lambda \in L^{\infty}[0,1], \
\|\theta - \theta_0\| \leq \epsilon_n,\
\|\ld - \ld_0\|_\Linfty \leq \zeta_n\}.
\label{Ln-infty}
\end{align}

\begin{lemma}
\label{R1-version2}
For $\mathcal{F}_{n,1}$ and $\mathcal{F}_{n,2}$ defined in (\ref{F-n1}) and (\ref{F-n2}) respectively with $\mathcal{L}_n$ is replaced with $\mathcal{L}_n'$ in (\ref{Ln-infty}) 
and 
$\Dl_1$ and $\Dl_{2,L_n}$ defined in (\ref{K1}) and (\ref{K2}) respectively, if $\|\Dl_1\|_\infty/\rn \leq d_1$ and $\|\Dl_{2,L_n}\|_\Linfty/\rn \leq d_2$ for some constants $d_1 + d_2 < 1$,
then 
\begin{align}
& \mathbb{E}_{\eta_0}^\star \left[ \|\mathbb{G}_n\|_{\mathcal{F}_{n,1}} \right] 
\lesssim (\|\Delta_1\|_\infty + \|\Delta_{2,L_n}\|_\Ltwo)^2/\rn,
\label{R1-v2-1}\\
& \mathbb{E}_{\eta_0}^\star \left[
\|\mathbb{G}_n\|_{\mathcal{F}_{n,2}} 
\right]
\lesssim \zeta_n (\|\Delta_1\|_\infty + \|\Delta_{2,L_n}\|_\Ltwo).
\label{R1-v2-2}
\end{align}
\end{lemma}

\begin{proof}
The proof is similar to that of Lemma \ref{R1-version1}, except that 
we bound $\|h_{n,11}\|_{BV}$ by
$$
\|h_{n,11}\|_{BV} 
\lesssim \|\ld - \ld_0\|_\infty 
\|\tilde \Delta_n\|_\Ltwo^2\rn
\leq \zeta_n (\|\Delta_1\|_\infty + \|\Delta_{2, L_n}\|_\Ltwo)^2/\rn.
$$ 
We have 
$$
\mathbb{E}_{\eta_0}^\star \left[ \|\mathbb{G}_n\|_{\mathcal{F}_{n,1}} \right]
\leq  
\mathbb{E}_{\eta_0}^\star \left[ \|\mathbb{G}_n\|_{\mathcal{F}_{n,11}} \right] + \mathbb{E}_{\eta_0}^\star \left[ \|\mathbb{G}_n\|_{\mathcal{F}_{n,12}} \right]$$ 
for each $\|\mathbb{G}_n\|_{\mathcal{F}_{n,1j}}$, $j = 1,2$, given in the proof of Lemma \ref{R1-version1}.
We immediately obtain
$\mathbb{E}_{\eta_0}^\star \|\mathbb{G}_n\|_{\mathcal{F}_{n,11}} 
\lesssim v_n (\|\Delta_1\|_\infty + \|\Delta_{2,L_n}\|_\Ltwo)^2/\rn.$
To bound the second term in the last display,
since $\|\Delta_{2, L_n}\|_\infty < d_2 < 1$ by assumption,
applying Taylor's theorem, we obtain $\|h_{n,12}\|_{BV} \lesssim (\|\Delta_1\|_\infty + \|\Delta_{2, L_n}\|_2)^2/\rn$. 
Thus, $ \mathbb{E}_{\eta_0}^\star \|\mathbb{G}_n\|_{\mathcal{F}_{n,12}} 
\lesssim (\|\Delta_1\|_\infty + \|\Delta_{2,L_n}\|_\Ltwo)^2/\rn$.
By combining the two bounds,
we obtain (\ref{R1-v2-1}).

The bound in (\ref{R1-v2-2}), we replace (\ref{h-n21}) with 
$$
\|h_{n,21}\|_{BV} \leq \|\ld - \ld_0\|_\Linfty \|\tilde \Delta_n\|_\Lone \leq \zeta_n \|\tilde \Delta_n\|_\Lone.
$$
Then, by following the same argument as in the proof of Lemma \ref{R1-version1}, we obtain (\ref{R1-v2-2}).
\end{proof}

\begin{lemma}
\label{lemma:bounding-R2-using-supnorm}
Suppose {\rm\ref{cond:P}} and assumptions 
{\rm \ref{asp:i}-\ref{asp:v}} hold, 
define $\tilde K_{a,b,t,s} = p^2(|t|\|a\|_\infty + |s| \|b\|_\Lone) + |s| \|b\|_\Ltwo$.
If $\tilde K_{a,b,t,s} /\rn = o(1)$, 
then 
\begin{align*}
& \sup_{\eta \in A_n} |R_{n,2}(\eta_h, \eta_0) - R_{n,2}(\eta, \eta_0) - s\rn B_3(\eta, \eta_0)| \\
& \quad \lesssim 
\tilde K^3_{a,b,n,p}/\rn + \tilde K^2_{a,b,n,p} \zeta_n + |s| p^2 \rn \epsilon_n 2^{-L_n} + \rn \epsilon_n^2 \tilde K_{a,b,t,s}.
\end{align*}
\end{lemma}

\begin{proof}
The proof is similar to that of Lemma \ref{thm1-lemma2}. 
We bound $(I)$ by a constant times 
$\|e^{(\theta - \theta_0)'z+ r-r_0} - 1\|_\Linfty \|(\theta_h - \theta)'z + r_h - r\|_\Ltwo^2$ instead. Using that $\|\ld - \ld_0\|_\Linfty \leq \zeta_n$ and $\|\gambln\|_\Ltwo \lesssim \|b\|_\Ltwo$ by Lemma \ref{lemma:psib-2}. Then by Lemma \ref{thm1-lemma7}, we have $(I) \lesssim \zeta_n \tilde K_{a,b,t,s}^2$. 
Also, one can bound (\ref{lem6-thm1-7}) by a constant times
$ \|(\theta_h - \theta)'z + r_h - r\|_\Lone 
\max_z\left(
\|(e^{r-r_0} - 1)(1 - e^{(\theta - \theta_0)'z})\|_\Linfty
+ |e^{(\theta- \theta_0)'z} -(\theta - \theta_0)'z - 1|
\right)$, which is bounded by a constant times $\tilde K_{a,b,t,s}\epsilon_n^2.$
The remaining proof is the same as that of Lemma \ref{thm1-lemma2}, and we thus obtain the result.
\end{proof}

\subsection{On centering and efficiency}
\label{centering-efficiency}

The two lemmas in this section enable us to center the posterior of $\Ld$ at an efficient estimator, e.g.,  the Breslow estimator.

For $\vartheta \in \mathbb{R}^p$ and $g(\cdot) \in L^2$, define the function 
$$
\Psi_\eta(\vartheta, g; \ X) = 
\delta (\vartheta'Z + g(Y)) - e^{\theta_0'Z} \int_0^Y (\vartheta'Z + g(u)) d\Ld_0(u).
$$
 From Section 12.3 of \citet{ghosal17}, the efficient influence function for estimating the linear function 
$\varphi_a(\theta)$ for $a \in \mathbb{R}^p$ is $\tilde \varphi_a = \tilde f_\eta(a' \tIinv, -a'\tIinv \gamm)$.
Similarly, the efficient influence function for estimating the linear function $\varphi_b(\ld) = \int_0^1 b\ld_0 = \Ld_0\{b\}$ for a function $b \in L^2(\Ld_0)$ is 
$\tilde \varphi_b = \tilde f_\eta (- \Ld_0\{\gamm'\} \tIinv, \gamb + \gamm' \tIinv \Ld_0\{b \gamm\})$.
Note that $W_n^{(1)}(a) = \tilde \varphi_a$,
an efficient estimator $\hat \varphi_a$ for $\varphi_a(\theta)$ should satisfy 
$$
\hat \varphi_a = \varphi_a + \frac{1}{\rn} W_n^{(1)}(a) + o_{P_{\eta_0}}(1).
$$
Also, $W_n^{(2)} = \tilde \varphi_b$ and an efficient estimator for $\varphi_b(\ld)$, $\hat \varphi_b$ should satisfy 
$$
\hat \varphi_b = \varphi_b + \frac{1}{\rn} W_n^{(2)}(b) + o_{P_{\eta_0}}(1).
$$

For $b = \psi_{lk}$, define
\begin{align}
\label{ld-star}
\ld_{L_n}^\star = \ld_{0, L_n} + \frac{1}{\rn} \sum_{L\leq L_n} \sum_{0 \leq K < 2^{L}} 
W_{n}^{(2)}(\psi_{LK}) \psi_{LK},
\end{align}
where $W_{n,LK}^{(2)}(\psi_{lk}) = \langle W_{n}^{(2)}(\psi_{lk}), \psi_{LK}\rangle$, $T_n^\ld$ is defined in (\ref{centering-Tn}). In the next lemma, we show that in the space of $\mathcal{M}_0$, the two quantities $T_n^\ld$ and $\ld_{L_n}^\star$ are close in probability. 

\begin{lemma}
\label{centering1}
Let $\ld^\star_{L_n}$ given by (\ref{ld-star}) and $T_n^\ld$ is given in (\ref{centering-Tn}), then for any admissible sequence $(w_l)$, 
$$
\mathbb{E}_{\eta_0} \|T_n^\ld - \ld^\star_{L_n}\|_{\mathcal{M}_0(w)}^2 = o(n^{-1}). 
$$
As a consequence, $\mathcal{B}_{\mathbb{R}^p \times \mathcal{M}_0} \left(\Pi(\cdot \given X) \circ \tau_{(T_n^\theta, T_n^\ld)}^{-1}, 
\Pi(\cdot \given X) \circ \tau_{(T_n^\theta, \ld_{L_n}^\star)}^{-1} \right) = o_{P_{\eta_0}}(1)$. 
\end{lemma}

\begin{proof}
By the definition of $W_n^{(2)}$, letting $\gamma_{n,lk} = P_{L_n}(\psi_{LK}/M_0)$
and
$\gamma_{M_1, nLK} = P_{L_n}(M_{1,LK}/M_0)$,
for $b = \psi_{LK}$, we have 
$$
\rn (T_n^\ld - \ld_{L_n}^\star) 
= \sum_{l \leq L_n, K} 
W_n\left(0, 
\gamb - \gamma_{n,LK} + 
(\gamm - \gamma_{M_1, nLK} )' \tIinv \Ld_0\{b\gamm\} \right) 
\psi_{LK}.
$$
By the definition of $\mathcal{M}(w)$-norm, we have $\|f\|_{\mathcal{M}(w)}^2 \leq \sum_{l,k} w_l^{-2} f_{lk}^2$ for any $f \in \mathcal{M}(w)$. Applying this inequality, we have
\begin{align*}
n \mathbb{E}_{\eta_0} \|T_n^\ld - \ld^\star_{L_n}\|_{\mathcal{M}_0(w)}^2 
\leq 
& \sum_{L\leq L_n, K} w_L^{-2} 
\|0, \psi_{LK}/{M_0} - P_{L_n}(\psi_{LK}/{M_0})\|_L^2 \\
& +
 \sum_{L\leq L_n, K} w_L^{-2} 
\left\|0, (\gamm - \gamma_{M_1, nLK})' \tIinv \Ld_0\{\psi_{LK}\gamm\}\right\|_L^2
\end{align*}

Using the $\|\cdot\|_\Linfty$-bound from Lemma \ref{lemma:psib-2}, 
the LAN-norm in the first line of the above display is bounded by 
$2^L 2^{-2L_n}$
and by the third point of Lemma \ref{lemma:psib-2}, the LAN-norm of the second line is bounded by 
$2^{-2L_n}$.
Therefore, we obtain that 
$$
n \mathbb{E}_{\eta_0} \|T_n^\ld - \ld^\star_{L_n}\|_{\mathcal{M}_0(w)}^2 
\lesssim 2^{-L_n} \sum_{L \leq L_n} w_L^{-2} 2^{L}2^{L- L_n}
\lesssim 1/L_n = o(1).
$$
Using the definition of the bounded Lipschitz metric and invoking Slutsky's theorem lead to 
$$
\mathcal{B}_{\mathbb{R}^p \times \mathcal{M}_0} \left(\Pi(\cdot \given X) \circ \tau_{(T_n^\theta, T_n^\ld)}^{-1}, 
\Pi(\cdot \given X) \circ \tau_{(T_n^\theta, \ld^\star)}^{-1} \right) 
\leq 
\rn \|T_n^\ld - \ld^\star\|_{\mathcal{M}_0(w)} = o_{P_{\eta_0}}(1).
$$
\end{proof}

\begin{lemma}
\label{lemma-centering-Lambda}
Let $T_n^\ld$ be defined as
\begin{align}
\langle T_n^\ld, \psi_{lk} \rangle =  
\begin{cases}
	\langle \ld_0, \psi_{lk} \rangle + W_n^{(2)}(\psi_{lk}) & \quad \text{if} \ l \leq L_n,\\
	0 & \quad \text{if} \ l > L_n,
\end{cases}
\end{align}
with the cut-off $L_n$ defined in (\ref{L_n}), where 
$$
W_n^{(2)}(\psi_{lk}) = W_n \left( - \tIinv \Ld_0\{\psi_{lk} \gamm\}, \ \psi_{lk}/M_0 + \gamm \tIinv \Ld_0\{\psi_{lk} \gamm\} \right).
$$
Let $\mathbb{T}_n^\ld(t) = \int_0^t T_n^\ld (u) du$, $t \in [0,1]$, and we set
$$
\Ld^\star(t) = \Ld_0(t) + \frac{1}{\rn} W_n^{(2)}(\mathbbm{1}_{\cdot \leq t}), \ t \in [0,1],
$$
Then, as $n \to \infty$, 
$
\rn \|\mathbb{T}_n^\ld(\cdot) - \Ld^\star(\cdot) \|_\Linfty = o_{P_{\eta_0}}(1). 
$
\end{lemma}

\begin{proof}
By the definition of $T_n^\ld$ and $\mathbb{T}_n$, one can write 
\begin{align}
\label{pf-lem13-1}
\mathbb{T}_n^\ld = \int_0^t \left(P_{L_n} \ld_0\right)(u) du 
+ \frac{1}{\rn} \sum_{L \leq L_n; K} W_n^{(2)}(\psi_{LK}) \int_0^t \psi_{LK}(u) du,
\end{align}
where $P_{L_n} \ld$ is the projection of $\ld_0$ onto the orthocomplement of $\mathcal{V}_n$. 
The first term in (\ref{pf-lem13-1}) can be written as 
$$
P_{L_n^c} \ld_0 = \ld_0 - P_{L_n} \ld_0. 
$$
Due to the linearity of $W_n$, the second term in (\ref{pf-lem13-1}) can be written as 
$$
W_n^{(2)}(\psi_{LK})
= \underbrace{W_n(- \tIinv \Ld_0\{\psi_{LK} \gamm\}, 0 ) }_{(I)}
+ \underbrace{W_n(0, \psi_{LK}/M_0)}_{(II)} 
+ \underbrace{W_n(0, \gamm \tIinv \Ld_0\{\psi_{LK} \gamm\})}_{(III)}. 
$$
Then, we have
$$
\sum_{L\leq L_n; K} (I) \int_0^t \psi_{LK}(u) du
= W_n \left( - \tIinv \sum_{L \leq L_n; K} \Ld_0\{\psi_{LK}\gamm\} \int_0^t \psi_{LK}(u) du, 0 \right).
$$
Since $\int_0^t \psi_{LK}(u) du = \langle \mathbbm{1}_{[0,t]} \psi_{LK} \rangle$, one can further write 
$
 \sum_{L \leq L_n; K} \Ld_0\{\psi_{LK}\gamm\} \int_0^t \psi_{LK}(u) du
 = 
\sum_{L \leq L_n; K}  \int_0^\cdot \psi_{LK}  \langle \mathbbm{1}_{[0,t]} \psi_{LK} \rangle \gamm(u) d\Ld_0(u)
 = \Ld_0\{P_{L_n} \mathbbm{1}_{[0,t]} \gamm\}.
$

Similarly, we can write 
$$
\sum_{L\leq L_n;K} (III) \int_0^t \psi_{LK}(u) du 
= W_n\left(0, 
\gamm' \tIinv \Ld_0\{ P_{L_n} \mathbbm{1}_{[0,1]} \gamm\}
\right).
$$

For the middle term, we have
\begin{align*}
\sum_{L \leq L_n;K} (II) \int_0^t \psi_{LK}(u) du
& = W_n \left(
0,  \sum_{L \leq L_n;K} \langle \mathbbm{1}_{[0,t]}, \psi_{LK} \rangle \psi_{LK} / M_0(\cdot)
\right) \\
& = W_n \left(0, P_{L_n} \mathbbm{1}_{[0,t]}(\cdot)/M_0(\cdot) \right).
\end{align*}
Therefore, using the results obtained above, (\ref{pf-lem13-1}) can be re-written as
\begin{align*}
\mathbb{T}_n^\ld(t) 
& = \frac{1}{\rn} W_n \left(\Ld_0\{P_{L_n} \mathbbm{1}_{[0,t]} \gamm\}, 
P_{L_n} \mathbbm{1}_{[0,t]}(\cdot)/M_0(\cdot) + 
\gamm' \tIinv \Ld_0\{ P_{L_n} \mathbbm{1}_{[0,t]} \gamm \}\right)\\
& \quad + \Ld_0(t) - \int_0^t (P_{L_n^c} \ld_0)(u) du.
\end{align*}
By comparing the expression of $\Ld^\star(t)$ with $\mathbb{T}_n^\ld(t)$ and note that $\Ld_0\{\gamm\} - \Ld_0\{P_{L_n} \mathbbm{1}_{[0,t]} \gamm\}
= \Ld_0\{P_{L_n^c} \mathbbm{1}_{[0,t]} \gamm\}$, we obtain 
\begin{align}
& \rn (\mathbb{T}_n^\ld (t) - \Ld^\star (t) ) \nonumber \\
& \quad= W_n \left(\Ld_0\{P_{L_n^c} \mathbbm{1}_{[0,t]} \gamm\}, 
P_{L_n^c} \mathbbm{1}_{[0,t]}/M_0 + 
\gamm' \tIinv \Ld_0\{ P_{L_n^c} \mathbbm{1}_{[0,t]} \gamm \}\right)
\label{pf-lem13-2}\\
& \quad \quad - \rn \int_0^t (P_{L_n^c} \ld_0)(u) du.
\label{pf-lem13-3}
\end{align}
The remaining proof is to bound $\rn \|\mathbb{T}_n^\ld (t) - \Ld^\star (t) \|_\Linfty$. 
First, we bound (\ref{pf-lem13-2}). By writing 
$
\int_0^t (P_{L_n^c} \ld_0)(u) du = \int_0^1 \left( P_{L_n^c} \mathbbm{1}_{[0,t]}(u) \right) P_{L_n^c}\ld_0(ud) du
$ and using the bound $\int f g \leq \|f\|_\Lone \|g\|_\Linfty$, we have
$$
\rn \int_0^t (P_{L_n^c} \ld_0)(u) du
\leq \| P_{L_n^c} \mathbbm{1}_{[0,t]}(u) \|_\Lone \|P_{L_n^c} \ld_0\|_\Linfty. 
$$
Using the fact that $\ld_0$ is $\beta$-H\"older, we have
$$
\|P_{L_n^c}\ld_0 \|_\Linfty \leq \sum_{l> L_n} 2^{l/2} \max_k|\langle \ld_0, \psi_{lk}\rangle| 
\leq \sum_{l \geq L_n} 2^{l/2} 2^{-\beta l} \|\psi_{lk}\|_\Lone \leq 2^{-\beta L_n}.
$$
On the other hand, we obtain 
$$\|P_{L_n^c} \mathbbm{1}_{[0,t]}\|_\Lone \leq \sum_{l \geq L_n;k} |\langle \mathbbm{1}_{[0,t]}, \psi_{lk} \rangle | \int_0^1 \psi_{lk}(u) du.
$$
Since $(\psi_{lk})$ is Haar basis, $|\langle \mathbbm{1}_{[0,1]}, \psi_{lk}| \rangle \leq 2^{-l/2}$. Therefore, the last display is bounded by $\sum_{l \geq L_n}2^{-l} \lesssim 2^{-L_n}$.

Thus, 
\begin{align}
\label{pf-lem13-4}
\rn  \int_0^t (P_{L_n^c} \ld_0)(u) du \lesssim 2^{-(1+\beta) L_n}.
\end{align}
What left is to bound (\ref{pf-lem13-2}), which we use the empirical process tools in Section \ref{EPlemmas}.
Define the function 
$
\Psi_t(\kappa_1, \kappa_2) = \delta (\kappa_1'z + \kappa_2) - e^{\theta_0'z} (\kappa_1'z \Ld_0(y) + (\Ld_0 \kappa_2)(y),
$
then, the empirical process (\ref{pf-lem13-2}) equals to $\mathbb{E}_{\eta_0}\Psi(\kappa_1, \kappa_2)$,
with 
$\kappa_1 = -\tIinv \Ld_0\{P_{L_n^c} \mathbbm{1}_{[0,t]} \gamm\}$ 
and $\kappa_2 = P_{L_n^c} \mathbbm{1}_{[0,t]}/M_0 + \gamm' \tIinv \Ld_0\{P_{L_n^c}\mathbbm{1}_{[0,t]} \gamm\}$.
Let $\mathcal{F}_n = \{f_t: \Psi_t(\kappa_1, \kappa_2)\}$. Then, for $f\in \mathcal{F}_n$, we obtain bounds for $\int f^2 dP_{\eta_0}$, $\|f\|_\Linfty$, and the bracketing intergral $J_{[]}(\delta, \mathcal{F}_n, L^2(P_{\eta_0}))$ and then invoke Lemma \ref{lemma-ep-4}.

Applying triangle inequality, 
\begin{align*}
& \int f_t^2 dP_{\eta_0} 
\lesssim \int_0^1 \left(\frac{P_{L_n^c} \mathbbm{1}_{[0,1]}}{M_0} \right)^2
+ \int_0^1 \left(\gamm' \tIinv \Ld_0\{P_{L_n^c} \mathbbm{1}_{[0,t]} \gamm\} \right)^2 \\
& \quad + \int_0^1 \left( \gamm' \tIinv \Ld_0\{P_{L_n^c} \mathbbm{1}_{[0,t]} \gamm\} \right)^2.
\end{align*}
By \ref{asp:i}-\ref{asp:v} and applying the inequality $\int f g \leq \|f\|_\Lone \|g\|_\Linfty$, the first term in the last display is bounded by a constant times $\|P_{L_n^c} \mathbbm{1}_{[0,t]}\|_\Ltwo^2$ and the second and third terms in the last display is bounded by a constant times $\|P_{L_n^c}\mathbbm{1}_{[0,t]}\|_\Lone^2 \leq \|P_{L_n^c}\mathbbm{1}_{[0,t]}\|_\Ltwo^2$, as $\|f\|_\Lone \leq \|f\|_\Ltwo$ for $f \in L^2$.
Since $(\psi_{lk})$ is Haar basis, 
$$
\|P_{L_n^c} \mathbbm{1}_{[0,t]}\|_\Ltwo^2
= 
\sum_{l \geq L_n,k} \langle \mathbbm{1}_{[0,t]}, \psi_{lk} \rangle^2
\leq \sum_{l \geq L_n} \sum_{0 \leq k < 2^{l}} 2^{-l/2} | \langle \mathbbm{1}_{[0,t]}, \psi_{lk} \rangle|
\lesssim \sum_{l \geq L_n} 2^{-l} \lesssim 2^{-L_n}.
$$
Thus, we obtain $\int f_t^2 dP_{\eta_0} \lesssim 2^{-L_n}$.

Bounding $\|f_t\|_\Linfty$ is similar, a simple calculation reveals that 
$$
\|f_t\|_\Linfty \lesssim \|P_{L_n^c} \mathbbm{1}_{[0,t]}\|_\Linfty + \|P_{L_n^c} \mathbbm{1}_{[0,t]}\|_\Lone \lesssim L_n + 2^{-L_n} \leq 2L_n,
$$ 
as $\mathbbm{1}_{[0,t]} = P_{L_n^c} \mathbbm{1}_{[0,t]} + P_{L_n} \mathbbm{1}_{[0,t]}$ for any $t \in [0,1]$. It is obvious that $\|\mathbbm{1}_{[0,t]}\|_\Linfty = 1$. Also,
$$
\|P_{L_n}\mathbbm{1}_{[0,t]}\|_\Linfty \leq \sum_{l \leq L_n} 2^{l/2} \max_k|\langle \mathbbm{1}_{[0,t]}, \psi_{lk} \rangle| \leq \sum_{l \leq L_n} 2^{l/2} \|\mathbbm{1}_{[0,t]}\|_\Linfty \max_k \|\psi_{lk}\|_\Lone
\lesssim L_n,
$$
thus, $\|P_{L_n^c}\mathbbm{1}_{[0,t]}\|_\Linfty \lesssim L_n$. 
Last, what remains is bounding the entropy $\mathcal{J}_{[]}(\delta, \mathcal{F}_n, L^2(P_{\eta_0})$. For any two $f_s, f_t \in \mathcal{F}_n$, $0 \leq s \leq t \leq 1$, 
noting that $|\mathbbm{1}_{[s,t]}, \psi_{lk}| \lesssim |s-t|^{1/2}$ by Cauchy-Schwarz inquality, the proceeding is similar as bounding $\|P_{L_n^c} \mathbbm{1}_{[0,t]}\|_2$, we have
$\|f_s - f_t\|_{L^2(P_{\eta_0})} \lesssim \|P_{L_n^c}\mathbbm{1}_{[0,t]}\|_\Ltwo$, and thus, by Lemma \ref{lemma-ep-2} and note that $N_{[]}(\epsilon, \mathcal{F}, L^2(P_{\eta_0}) \lesssim 2^{L_n}/\epsilon^4$, we obtain $J_{[]} (\delta, \mathcal{F}_n, L^2(P_{\eta_0})) \lesssim \sqrt{L_n} \delta + \delta \log (1/\delta)$ for a $\delta = o(1)$. 
By invoking Lemma \ref{lemma-ep-4} and choosing $\delta = 2^{-L_n}$, we have
$\|\mathbb{G}_n\|_{\mathcal{F}_n} \lesssim L_n 2^{-L_n} + L_n^3/\rn$. By the assumption $L_n^3 = o(\rn)$, this expression goes to 0 as $n \to \infty$. 
By combining the bounds of (\ref{pf-lem13-2}) and (\ref{pf-lem13-3}) together, we obtain 
$\rn \|\mathbb{T}_n^\ld(t) - \Ld^\star(t)\|_\Linfty \to^{P_{\eta_0}} 0$.
 
\end{proof}

%%%%%%%%%%%%%%%%%%%%%%%%%%%%%%%%%%%%%%%%
%%%%%%%%%%%%%%%% Acknowledgement %%%%%%%%%%%%%
%%%%%%%%%%%%%%%%%%%%%%%%%%%%%%%%%%%%%%%
%\section*{Acknowledgements}

%%%%%%%%%%%%%%%%%%%%%%%%%%%%%%%%%%%%%%%%
%%%%%%%%%%%%%%%% Bibliography %%%%%%%%%%%%%%%%%
%%%%%%%%%%%%%%%%%%%%%%%%%%%%%%%%%%%%%%%
\bibliographystyle{chicago}
\bibliography{citation.bib}

\end{document}